\documentclass[11pt,letterpaper, reqno]{amsart}
\usepackage[
top=3cm, bottom=3cm, inner=3cm, outer=3cm,
marginparwidth=2cm 
]{geometry}

\usepackage{amssymb,latexsym, amsmath, amsxtra, mathrsfs, bm}

\usepackage[dvips]{graphics}
\usepackage{xcolor}
\usepackage{subcaption}

\usepackage{verbatim}

\usepackage[abs]{overpic}
\usepackage{hyperref}
\usepackage{mathtools, tikz}
\usetikzlibrary{positioning, shapes.geometric}

\usepackage{graphicx}
\usepackage{tikz-cd}
\usepackage{tikz}
\usetikzlibrary{matrix}

\usepackage{amsmath,amsthm,amssymb,mathtools,latexsym,hyperref,xcolor,appendix,listings} 
\usepackage{tikz,overpic,graphicx,eepic,caption}
\usepackage{enumitem} 
\setlist[enumerate]{label=\rm{(\roman*)}}
\usepackage{algorithm,algorithmicx,multirow,algpseudocode,xargs}
\usepackage{verbatim,float,extarrows,bbm,cleveref}
\usepackage[matrix,arrow,curve,cmtip]{xy} 
\usepackage{ifthen} 

\usepackage{tikz-cd}
\usetikzlibrary{positioning}

\makeatletter
\@namedef{subjclassname@2020}{\textup{2020} Mathematics Subject Classification}
\makeatother

\DeclareMathAlphabet{\mathbbold}{U}{bbold}{m}{n}

\allowdisplaybreaks[4]

\makeatletter
\def\th@plain{%
	\thm@notefont{}
	\itshape 
}
\def\th@definition{%
	\thm@notefont{}
	\normalfont 
}
\makeatother

\theoremstyle{plain}
\newtheorem{theorem}{Theorem}[section]
\newtheorem*{theorem*}{Theorem}
\newtheorem*{conj*}{Conjecture}
\newtheorem{lemma}[theorem]{Lemma}
\newtheorem{prop}[theorem]{Proposition}

\newtheorem{thmx}{Theorem}

\theoremstyle{definition}
\newtheorem{definition}[theorem]{Definition}

\theoremstyle{remark}
\newtheorem*{remark}{Remark}

\newtheorem*{example}{Example}
\newtheorem*{claim}{Claim}

\newtheorem*{case1}{Case 1}
\newtheorem*{case2}{Case 2}

\numberwithin{equation}{section}
\numberwithin{theorem}{section}
\numberwithin{table}{section}
\numberwithin{figure}{section}

\renewcommand{\le}{\leqslant}
\renewcommand{\leq}{\leqslant}
\renewcommand{\ge}{\geqslant}
\renewcommand{\geq}{\geqslant}

\newcommand{\diam}  {\operatorname{diam}}

\newcommand{\inte}  {\operatorname{int}}
\newcommand{\inter}  {\operatorname{int}}

\newcommand{\card} {\operatorname{card}}

\newcommand{\R}{\mathbb{R}}

\newcommand{\C}{\mathbb{C}}

\newcommand{\Q}{\mathbb{Q}}
\newcommand{\N}{\mathbb{N}}
\newcommand{\Z}{\mathbb{Z}}

\newcommand{\CDach}{\widehat{\mathbb{C}}}
\newcommand{\RDach}{\widehat{\mathbb{R}}}
\newcommand{\D}{\mathbb{D}}

\providecommand{\abs}[1]{\lvert#1\rvert}

\providecommand{\norm}[1]{\|#1\|}

\renewcommand{\:}{\colon}

\renewcommand{\b}{\mathfrak{b}}
\newcommand{\w}{\mathfrak{w}}

\newcommand{\holderexp}{\beta} 
\newcommandx{\holderspace}[4][1=, 2=\holderexp, 3=X, 4=d, usedefault]{ C^{0, #2}\mathopen#1(#3, #4 \mathclose#1) }
\newcommandx{\holderspacesphere}[2][1=\holderexp, 2=d]{ C^{0, #1}(S^2, #2) }
\newcommandx{\holderseminorm}[6][1=, 3=\holderexp, 5=d, 6=]{ \mathopen #1| #2 \mathclose #1|_{#3, \, \mathopen #6(#4, #5 \mathclose #6)} }
\newcommandx{\holdernorm}[6][1=, 3=\holderexp, 5=d, 6=]{\mathopen #1\| #2 \mathclose #1\|_{C^{0, #3} \mathopen #6(#4, #5 \mathclose #6)} }

\newcommand{\I}{\mathbf{i}}

\newcommand{\crit}{\operatorname{crit}}
\newcommand{\post}{\operatorname{post}}

\newcommand{\mesh}{\operatorname{mesh}}

\newcommand{\X} {\mathbf{X}}

\newcommand{\E} {\mathbf{E}}

\newcommand{\V} {\mathbf{V}}

\newcommand{\describe}{:}  

\newcommand{\black}{\mathfrak{b}} 
\newcommand{\white}{\mathfrak{w}} 

\newcommand{\vertiii}[1]{{\left\vert\kern-0.25ex\left\vert\kern-0.25ex\left\vert #1
		\right\vert\kern-0.25ex\right\vert\kern-0.25ex\right\vert}}

\renewcommand{\=}{\coloneqq}

\newcommand{\cA}{\mathcal{A}}

\newcommand{\cC}{\mathcal{C}}

\newcommand{\cE}{\mathcal{E}}

\newcommand{\cK}{\mathcal{K}}

\newcommand{\sA}{\mathscr{A}}
\newcommand{\sB}{\mathscr{B}}

\newcommand{\oC}{\overline{C}}

\newcommand{\oS}{\overline{S}}

\newcommand{\oa}{\overline{a}}

\newif\ifFirstInitial
\FirstInitialfalse  

\newcommand{\prV}{\operatorname{pr}_{\V}}

\title{Zippers for expanding Thurston maps}
\makeatletter
\renewcommand{\@setauthors}{%
  \begingroup
  \trivlist
  \centering\normalsize \@topsep30\p@\relax
  \advance\@topsep by -\baselineskip
  \item\relax
  \authors
  \vspace{14pt}
  \endtrivlist
  \endgroup
}
\makeatother

\author[Yuke~An, Yan~Mary~He, Zhiqiang~Li, and Ruoyu~Zhong]{%
  \begin{tabular}{c}Yuke~An\\[-3pt]\Small Peking University
  \end{tabular}\;
  \begin{tabular}{c}Yan~Mary~He\\[-3pt]\Small University of Oklahoma
  \end{tabular}\;
  \begin{tabular}{c}Zhiqiang~Li\\[-3pt]\Small Peking University
  \end{tabular}\;
  \begin{tabular}{c}Ruoyu~Zhong\\[-3pt]\Small Peking University
  \end{tabular}
}

\begin{document}



\begin{abstract}
Motivated by the recent work of Calegari and Loukidou on zippers and
CaTherine wheels for closed hyperbolic $3$-manifolds, via Sullivan's dictionary, we study the
relation between zippers and invariant Peano curves for expanding
Thurston maps. For every
expanding Thurston map, we use Meyer's tile connections to construct
a P-zipper $(Z^+,Z^-)$ for a suitable iterate. This
construction allows perfect fits between the white and black
laminations. The two sides $Z^+$ and $Z^-$ admit ribbon tree models $Y^+$ and $Y^-$ respectively.
We show that at least one of $Y^+$ and $Y^-$ has an end circle, which we identify with the
circle $S^1$ parametrizing Meyer's approximations. When both end
circles exist, matching points of equal angle gives an
orientation-reversing identification between them and hence a
universal circle. The two tree models determine the fibers of
Meyer's invariant Peano curve on $S^1$, and the resulting quotient
recovers the Peano curve in its original parametrization.
\end{abstract}

\maketitle
\tableofcontents

\section{Introduction}
Universal circles arise in the study of several geometric and dynamical structures on hyperbolic $3$-manifolds, including taut foliations and quasigeodesic flows \cite{CD03,Calegari06}. These circles carry faithful actions of the fundamental group preserving pairs of laminations, connecting the dynamics of these structures with the geometry of the manifold. In several settings, this connection extends to the sphere at infinity through equivariant Peano curves, for example the Cannon--Thurston map for hyperbolic $3$-manifolds fibering over the circle \cite{CT07}, and constructions associated with pseudo-Anosov flows without perfect fits and quasigeodesic flows \cite{Fenley12,Frankel15}.

In \cite{CL24}, Calegari and Loukidou introduced \emph{zippers} to provide a unified approach to constructing universal circles. Their framework simplifies several existing constructions and yields new ones from uniform quasimorphisms and uniform actions. For a closed oriented hyperbolic $3$-manifold $M$, a zipper is a pair of disjoint, nonempty, path-connected, and $\pi_1(M)$-invariant subsets $Z^+,\, Z^-\subseteq S^2_\infty$. The two sides $Z^+,Z^-$ carry natural tree structures. For a minimal zipper, their end circles can be canonically identified to produce a universal circle; see \cite[Theorem~2.32]{CL24}. They also introduced \emph{P-zippers}. Each side of a P-zipper is the image of a path-connected space under an equivariant continuous surjection, and the images of paths from the two spaces do not cross; see \cite[Definition~2.34]{CL24}.

The theory of \emph{CaTherine wheels} developed by Calegari and Loukidou in \cite{CLWheels} further connects zippers with the geometry of sphere-filling curves. A CaTherine wheel is a continuous surjection from $S^1$ to $S^2$ for which the image of every nondegenerate proper closed interval is a closed topological disk, with the images of its endpoints lying on the boundary. For closed hyperbolic $3$-manifolds, Calegari and Loukidou establish a correspondence between minimal zippers, equivariant CaTherine wheels, and pseudo-Anosov flows without perfect fits; see \cite[Theorem~6.32]{CLWheels}. Their theory also includes constructions associated with expanding Thurston maps and probability theory.

Sullivan's dictionary \cite{Sullivan85} motivates dynamical analogues of this picture, with expanding Thurston maps providing dynamical analogs of Gromov hyperbolic groups with $S^2$ boundary \cite[Chapter~1]{BM}. Thurston maps arose in Thurston's topological characterization of postcritically finite rational maps \cite{DH93}. For an expanding Thurston map, the tiles cut out by successive pullbacks of a Jordan curve through the postcritical set are vertices of a Gromov hyperbolic graph whose boundary at infinity is naturally identified with $S^2$ \cite[Theorem~1.1 and Proposition~1.2]{Yin15}. Ha\"issinsky and Pilgrim obtained a related boundary construction for coarse expanding conformal systems \cite{HP09}. Meyer's invariant Peano curve theorem \cite[Theorem~1.1]{Meyer13} states that, for every expanding Thurston map $f_0\colon S^2\to S^2$ and every sufficiently large $N$, the iterate $f=f_0^N$ admits a continuous surjection $\gamma\colon S^1\to S^2$ satisfying
\begin{equation*}
f\circ\gamma=\gamma\circ\mu_d,
\quad \mu_d(z)=z^d,
\quad d=\deg f.
\end{equation*}
This semiconjugacy provides a dynamical analogue of the equivariant Peano curves arising in hyperbolic $3$-manifold geometry.

The purpose of this article is to investigate the zipper structures
underlying Meyer's invariant Peano curves for expanding Thurston
maps. For topological matings of dendritic Julia sets that admit a
CaTherine wheel, Calegari and Loukidou identify the associated
zippers with the unions of the iterated preimages of trimmed
Hubbard trees; see \cite[Proposition~7.9]{CLWheels}. The corresponding
laminar relations have no perfect fits; see
\cite[Theorem~2.13]{CLWheels}.  In this article, starting from the connections in Meyer's
construction, we construct
{\it P-zippers} for every expanding Thurston map after passing
to a suitable iterate, without assuming the absence of perfect
fits. Our construction thus also applies to examples with infinite
chains of perfect fits, such as those of Loukidou
\cite{Loukidou26}. We build the corresponding tree models, study their
ribbon structures and end circles, and recover Meyer's invariant
Peano curve from the associated laminar relations.

\subsection{P-zippers and their tree models}
Let $f\:S^2\to S^2$ be an expanding Thurston map.
Our first goal is to introduce \emph{P-zippers} for $f$ as in Definition~\ref{def_gen_zipper}. Roughly speaking, a P-zipper consists of two nonempty, path-connected, completely $f$-invariant subsets $Z^+,Z^-\subseteq S^2$ admitting continuous surjections $\pi^\sigma\:Y^\sigma\to Z^\sigma$ from path-connected spaces $Y^\sigma$. These maps satisfy $\pi^\sigma\circ f_\sigma=f\circ\pi^\sigma$ for continuous maps $f_\sigma\:Y^\sigma\to Y^\sigma$. The images of paths from $Y^+$ and $Y^-$ under $\pi^+$ and $\pi^-$ do not cross. We call $Z^+$ the \emph{white side} and $Z^-$ the \emph{black side}, and use the same color convention for $Y^+$ and $Y^-$. The sets $Z^+$ and $Z^-$ need not be disjoint, and the surjections from $Y^+$ and $Y^-$ need not be injective.

Our first theorem states that up to an iterate, every expanding Thurston map admits a P-zipper.

\begin{thmx}\label{A}
For every expanding Thurston map $f_0\:S^2\to S^2$, there exists an integer $N\geq1$ such that the iterate $f\=f_0^N$ admits a P-zipper $(Z^+,Z^-)$ with tree models $Y^+$ and $Y^-$. For each $\sigma\in\{+,\,-\}$, the set $Z^\sigma$ is dense in $S^2$, and the tree $Y^\sigma$ is the increasing union $\bigcup_{n=0}^{+\infty}Y_n^\sigma$ of finite subtrees, equipped with the weak topology. The projection $\pi^\sigma\:Y^\sigma\to Z^\sigma$ is continuous and surjective, and there is a continuous map $f_\sigma\:Y^\sigma\to Y^\sigma$ satisfying $\pi^\sigma\circ f_\sigma=f\circ\pi^\sigma$. For each pair of continuous paths $a^\pm\:[0,1]\to Y^\pm$, the paths $\pi^\pm\circ a^\pm$ do not cross in the sense of Definition \ref{def_gen_zipper}.
\end{thmx}

Theorem~\ref{A} gives a dynamical counterpart of the tree structures carried by geometric zippers in the sense of \cite{CL24}: the white and black connection data yield two tree models $Y^\pm$ whose dynamics project to $f$. These tree models provide the setting for passing from the sphere dynamics to end spaces and circle dynamics.

When $f_0$ has no periodic critical points, every sufficiently
high iterate $f=f_0^N$ is topologically conjugate to a mating
of two postcritically finite polynomials with dendritic Julia
sets $J^\pm$; see \cite[Theorem~1.2]{Meyer09}.
In this setting, if the white and black relations have no
perfect fits, the images of $\operatorname{Cut}(J^\pm)$ under
the mating projections are precisely the zippers of Calegari
and Loukidou, obtained from iterated preimages of trimmed
Hubbard trees; see \cite[Proposition~7.9]{CLWheels}.
Our sets $Z^\pm$ are obtained from these zippers by adjoining
$\bigcup_{n=0}^{+\infty}f^{-n}(\post f)$.
See the comparison at the end of Section~\ref{ConstructionofYpm}.

\subsection{End circles and their dynamics}
Our next goal is to construct the end circle of each hairy tree $Y^\sigma$. When both trees are hairy, we identify their end circles to obtain a universal circle. Our construction uses Meyer's connections of tiles.
For each $\sigma\in\{+,\,-\}$, we use the connection and the orientation of $S^2$ to equip $Y^\sigma$ with a canonical \emph{ribbon structure}, namely a cyclic order on the components of $Y^\sigma\smallsetminus\{y\}$ at each branch point $y$.

Let $\mathcal E(Y^\sigma)$ denote the \emph{interior end space} of $Y^\sigma$. Its elements are either \emph{hard ends}, corresponding to endpoints of $Y^\sigma$, or \emph{genuine ends}, represented by proper rays; see Proposition~\ref{prop:rooted-representatives}. When this space has at least three elements, the ribbon structure induces a \emph{circular order} on it, as described in Section~\ref{endcircles}.

We call a ribbon tree \emph{hairy} if, on each ribbon side of every nondegenerate oriented arc, there is a proper ray issuing from an interior point of the arc and meeting the arc only at that point; see Definition~\ref{def:ch5-hairy}. For a hairy $Y^\sigma$, the circular order completion of $\mathcal E(Y^\sigma)$, in the sense of \cite[Section~5]{Frankel13}, is a topological circle; see Theorem~\ref{thm:end-circle}. We call it the \emph{end circle} of $Y^\sigma$; see also \cite{CL24}.

Fix Meyer's parametrized approximations $\gamma^n\colon S^1\to S^2$ associated with the chosen connection, with the normalization as in Subsection~\ref{Pseudo_isotopies}. In particular, $\gamma^0$ follows the positive boundary orientation of the white $0$-tile. Their common parameter circle carries the map $\mu_d(z)=z^d$, where $d=\deg f$, and their uniform limit $\gamma$ satisfies $f\circ\gamma=\gamma\circ\mu_d$. All angle coordinates refer to these fixed parametrizations. The next theorem describes the end circles and the dynamics induced on them by the tree maps.

\begin{thmx}\label{B}
Let $f\colon S^2\to S^2$ be an expanding Thurston map of degree $d$, and let $\mathcal C\supseteq\post  f$ be an $f$-invariant Jordan curve. Fix a connection of $1$-tiles satisfying \textup{(C1)}--\textup{(C2)} of Definition~\ref{Propertiesofconnections}, together with the associated parametrized approximations $\gamma^n$ described above. Let $Y^\pm$ and $f_\pm$ be the tree models and maps constructed from this connection in the proof of Theorem~\ref{A}, and equip $Y^\pm$ with the ribbon structures. Then the following hold:
\begin{enumerate}
\smallskip
\item For each $\sigma\in\{+,\,-\}$, either $Y^\sigma=Y_0^\sigma$ is a compact interval or $Y^\sigma$ is hairy. At least one of $Y^+$ and $Y^-$ is hairy.

\smallskip
\item If $Y^\sigma$ is hairy, the circular order completion $S_\sigma\=S^1(Y^\sigma)=\overline{\mathcal E(Y^\sigma)}$ is a topological circle, oriented by its circular order, and the genuine ends of $Y^\sigma$ form a dense subset of $S_\sigma$. The map induced by $f_\sigma$ on genuine ends extends uniquely to a continuous map $\overline f_\sigma\colon S_\sigma\to S_\sigma$, which is a topologically exact covering of degree $d$. The fixed parametrizations of $\gamma^n$ determine a homeomorphism $\eta_\sigma\colon S_\sigma\to S^1$ to their common parameter circle satisfying
\begin{equation*}
\eta_\sigma\circ\overline f_\sigma
=\mu_d\circ\eta_\sigma \quad\text{ and } \quad
\mu_d(z)=z^d.
\end{equation*}
With respect to the standard orientation of $S^1$, the map $\eta_\sigma$ is orientation preserving for $\sigma\=+$ and orientation reversing for $\sigma\=-$.

\smallskip
\item If both trees are hairy, the homeomorphism
$
h\=\eta_-^{-1}\circ\eta_+\colon S_+\to S_-
$
pairs points with the same angle, reverses orientation, and satisfies
$
h\circ\overline f_+=\overline f_-\circ h.
$
\end{enumerate}
\end{thmx}

If $Y^\sigma$ is a compact interval, $\mathcal E(Y^\sigma)$ consists of two hard ends. If $Y^\sigma$ is hairy, $\mathcal E(Y^\sigma)$ is infinite, and its circular order completion is the end circle $S_\sigma$. Hairiness has a combinatorial characterization, namely, $Y^\sigma$ is hairy if and only if the induced $m$-connection graph of color $\sigma$ has a vertex of degree at least three for some $m\ge1$; see Proposition~\ref{lem:hairy}~\textup{(iii)}. Here the induced graph is obtained by deleting the degree-one block vertices from the connection graph; see Subsection~\ref{subsec:hairiness}. In particular, branching in both induced $1$-connection graphs suffices to obtain both end circles and the identification in Theorem~\ref{B}~\textup{(iii)}.

\subsection{Recovery of the invariant Peano curve}
We now relate the tree models to Meyer's invariant Peano curve. The fixed parametrizations of the approximations $\gamma^n$ determine compatible position maps $x_n^\sigma\colon S^1\to Y_n^\sigma$; see Section~\ref{The_end_circle_and_laminations}. For each color, two angles on $S^1$ are identified when their positions agree in every finite tree. The next theorem identifies the resulting relations with Meyer's closed white and black relations.

\begin{thmx}\label{C}
Let $f\colon S^2\to S^2$ be an expanding Thurston map of degree $d$, and let $\mathcal C\supseteq\post f$ be an $f$-invariant Jordan curve. Fix a connection of $1$-tiles satisfying \textup{(C1)}--\textup{(C2)} of Definition~\ref{Propertiesofconnections}, together with the associated parametrized approximations $\gamma^n$ and their Peano limit $\gamma$. For each $\sigma\in\{+,\,-\}$ and $n\in \N_0$, let $x_n^\sigma$ be the map defined in Lemma~\ref{lem:ch6-finite-positions-ref} and define
\begin{equation*}
    R_\sigma\=\bigcap_{n=0}^{+\infty}\ker x_n^\sigma.
\end{equation*}

Then $R_+$ and $R_-$ agree with Meyer's closed white and black relations, respectively. Consequently,
$
R\=\widehat{R_+\vee R_-}=\ker\gamma,
$
where $\widehat{R_+\vee R_-}$ is the smallest closed equivalence relation containing both $R_+$ and $R_-$.

If $q\colon S^1\to S^1/R$ is the quotient map, then the map
$H\colon S^1/R\to  S^2$ given by $H([t])\=\gamma(t)$
is a homeomorphism satisfying $H\circ q=\gamma$. Moreover, $\mu_d(z)\=z^d$ induces a continuous map $F\colon S^1/R\to S^1/R$ satisfying $F\circ q=q\circ\mu_d$ and $H\circ F=f\circ H$.
\end{thmx}

Theorem~\ref{C} recovers the invariant Peano curve associated with the chosen connection in its original angle parametrization. The two relations are determined by the finite tree positions, and their smallest common closed equivalence relation recovers the fibers of $\gamma$, together with the sphere dynamics on the quotient.

\subsection{Strategy of the proofs}
The proof of Theorem~\ref{A} constructs tree models from closed equivalence relations on disjoint unions of parameter intervals prescribed by the connection graphs. This retains distinctions between branches whose spherical images may overlap. The key step is to prove that the finite quotients are trees and that refinement induces embeddings; see Propositions~\ref{prop:compact-finite-trees} and~\ref{prop:compact-stage-stability}. Their weak union carries the lifted dynamics, while expansion controls the spherical realizations and yields the projections onto $Z^\pm$; see Propositions~\ref{prop:compact-projections} and~\ref{prop:pi-weak-tree}.

The proof of Theorem~\ref{B} combines compatible ribbon structures with a dynamical criterion for hairiness; see Propositions~\ref{Yribbon} and~\ref{lem:hairy}. The end circle of a hairy tree is realized as the inverse limit of the boundary circles of its finite subtrees; see Theorem~\ref{thm:end-circle}. The essential step is to identify this circle with Meyer's original parameter circle compatibly with the tree dynamics; see Theorem~\ref{prop_z^d}. This yields the conjugacy to $z\mapsto z^d$ and, when both trees are hairy, the orientation-reversing identification by equal angles.

The proof of Theorem~\ref{C} compares the equivalence relations defined by tree positions with Meyer's white and black relations. The key step is a finite level approximation argument relating angle classes to connection blocks and their images in the trees; see Proposition~\ref{relationsameprop}. Keeping Meyer's original angle coordinates makes these identifications compatible across the two colors. Meyer's characterization of the Peano fibers then identifies their smallest common closed equivalence relation with $\ker\gamma$; see Lemma~\ref{closureequi}. Finally, $\gamma$ induces a homeomorphism from the quotient onto $S^2$ that conjugates the induced dynamics to $f$.

\subsection{Organization of the article}
In Section~\ref{Preliminaries}, we review expanding Thurston maps, Meyer's tile connections, and invariant Peano curves.
In Section~\ref{sec:roads} we construct the road equivalence relations and the initial trees $Y_0^\pm$.
In Section~\ref{ConstructionofYpm} we construct the P-zippers $Z^\pm$ and its tree models $Y^\pm$, proving Theorem~\ref{A}.
In Section~\ref{endcircles} we develop the ribbon structures and studies the end circles and their dynamics, proving Theorem~\ref{B}.
In Section~\ref{The_end_circle_and_laminations} we recover Meyer's invariant Peano curve from the white and black tree relations, proving Theorem~\ref{C}.
When both induced connection graphs branch, Section~\ref{Connection_graphs_with_branching} describes the boundary of each $0$-arc image as the union of its white and black roads, in analogy with CaTherine wheels.

\subsection*{Acknowledgments} 
We would like to thank Hanyun Zheng for numerous comments on an early version of the manuscript.
Y.~A., Z.~L., and R.~Z.\ were partially supported by Beijing Natural Science Foundation~(JQ25001) and National Natural Science Foundation of China~(12471083).

\subsection*{AI Statement}
The authors confirm that no artificial intelligence (AI) or machine-learning tools were used in the research and writing of this manuscript.

\section{Preliminaries}
\label{Preliminaries}

\subsection{Notations}
We first fix some notation.
We follow the convention $\N \= \{1,\, 2,\, 3,\, \dots\}$ and $\N_0 \= \N\cup \{ 0\}$.
We write $\Z$ for the set of integers, and $\Q$, $\R$, and $\C$ for the set of rational, real, and complex numbers, respectively. We let $\CDach \= \C \cup \{\infty\}$ be the Riemann sphere and we denote $\RDach \= \R \cup \{\infty\}$. The symbol $\I$ stands for the imaginary unit in the complex plane $\C$. Let $\D\=\{z\in \C:\abs{z}<1\}$ be the unit disk in the complex plane.

Two nonnegative quantities $a$ and $b$ are said to be \emph{comparable} if there is a constant $C\geq 1$ (possibly depending on some ambient parameters) such that $\frac{1}{C}a\le b\le Ca.$
We then write $a\asymp b$. The constant $C$ is referred to as $C(\asymp)$. Similarly, we write $a\lesssim b$ or $b\gtrsim a$, if there is a constant $C>0$ such that $a\leq Cb$, and refer to the constant $C$ as $C(\lesssim)$ or $C(\gtrsim)$.

\subsection{Thurston maps}\label{sub:Thurston_maps}
In this subsection, we review some key concepts and results concerning Thurston maps, with particular emphasis on expanding Thurston maps. For a more comprehensive treatment of this subject, we refer the reader to \cite{BM}.

We denote by $S^2$ an oriented topological $2$-sphere. A continuous and surjective map $f \: S^2 \to S^2$ is called a \emph{branched covering map} on $S^2$ if, for each point $x\in S^2$, there exist a positive integer $d\in \N$, open neighborhoods $U$ of $x$ and $V$ of $y \= f(x)$, and orientation-preserving homeomorphisms $\varphi \: U \to \D$ and $\eta \: V \to \D$ with $\varphi(x) = 0$, $\eta(y) = 0$ such that $\bigl(\eta \circ f\circ \varphi^{-1}\bigr)(z) = z^d$ for all $z \in \D$.
The positive integer $d$ above is called the \emph{local degree} of $f$ at $x$ and is denoted by $\deg_f(x)$.

The \emph{degree} of $f$ is defined by $\deg f \= \sum_{x\in f^{-1}(y)} \deg_{f}(x)$ for every $y\in S^2$, and this quantity is independent of the choice of $y$. If $f \: S^2 \to S^2$ and $g \: S^2 \to S^2$ are two branched covering maps on $S^2$, then $f\circ g$ is also a branched covering map on $S^2$, and $\deg_{f\circ g}(x) = \deg_g(x) \deg_f(g(x))$ for each $x\in S^2$. Moreover, we have $\deg(f\circ g) = (\deg f)\cdot (\deg g)$.

A point $x\in S^2$ is called a \emph{critical point} of $f$ if $\deg_f(x) \geq 2$. The set of critical points of $f$ is denoted by $\crit f$. A point $y\in S^2$ is called a \emph{postcritical point} of $f$ if $y = f^n(x)$ for some $x \in \crit f$ and $n\in \N$. The set of postcritical points of $f$ is denoted by $\post f$. Note that $\post f = \post(f^n)$ for all $n\in \N$.

\begin{definition}[Thurston maps]
    A \emph{Thurston map} is a branched covering map $f \: S^2 \to S^2$ on $S^2$ with $\deg f \geq 2$ and a finite set of postcritical points.
\end{definition}

A basic approach to the investigation of Thurston maps, as well as a key tool in this article, is the cell decomposition. In what follows, we give a brief introduction to cell decompositions of $S^2$ and their relation to Thurston maps; see \cite[Section~5]{BM} for a detailed discussion. A \emph{cell of dimension $n$} in $S^2$, where $n \in \{1, \, 2\}$, is a subset $c \subseteq S^2$ that is homeomorphic to the closed unit ball $\overline{\mathbb{B}}^n$ in $\R^n$, where $\mathbb{B}^{n}$ denotes the open unit ball in $\R^{n}$. We define the \emph{cell boundary} (or simply \emph{boundary} when no ambiguity arises) of $c$, denoted by $\partial_{\circ} c$, to be the set of points corresponding to $\partial \mathbb{B}^n$ under such a homeomorphism between $c$ and $\overline {\mathbb B}^n$. This is independent of the homeomorphism chosen, and the set $\partial_{\circ} c$ is well-defined. The \emph{cell interior} (or simply \emph{interior} when no ambiguity arises) of $c$ is defined by $\inte_\circ (c) \= c \smallsetminus \partial_\circ c$. Note that $\partial_{\circ}$ and $\inte_\circ$ are different from the usual topological notions of interior and boundary. A \emph{cell of dimension} $0$ in $S^2$ is a set $c\subseteq S^2$ consisting of a single point. We set $\partial_\circ c \= \emptyset$ and $\inte_\circ (c) \= c$ in this case.

We recall the following definition of cell decompositions from \cite[Definition~5.1]{BM}.

\begin{definition}[Cell decompositions]    \label{def:cell decomposition}
    Let $\mathbf{D}$ be a collection of cells in $S^2$. We say that $\mathbf{D}$ is a \emph{cell decomposition of $S^2$} if the following hold:
    \begin{enumerate}
        \smallskip
        \item The union of all cells in $\mathbf{D}$ is equal to $S^2$.
        \smallskip
        \item If $c \in \mathbf{D}$, then $\partial_\circ c$ is a union of cells in $\mathbf{D}$.
        \smallskip
        \item For all $c_1,\,c_2 \in \mathbf{D}$ with $c_1 \neq c_2$, we have $\inte_\circ (c_1) \cap \inte_\circ(c_2) = \emptyset$.
        \smallskip
        \item Every point in $S^2$ has a neighborhood that meets only finitely many cells in $\mathbf{D}$.
    \end{enumerate}
\end{definition}

We now recall the notion of refinements of two cell decompositions.

\begin{definition}[Refinements]
    Let $\mathbf{D}'$ and $\mathbf{D}$ be two cell decompositions of $S^2$. We say that $\mathbf{D}'$ is a \emph{refinement} of $\mathbf{D}$ if the following hold:
    \begin{enumerate}
        \smallskip

        \item Every cell $c \in \mathbf{D}$ is the union of all cells $c' \in \mathbf{D}'$ with $c' \subseteq c$.

        \smallskip

        \item For every cell $c' \in \mathbf{D}'$ there exists a cell $c \in \mathbf{D}$ with $c' \subseteq c$.
    \end{enumerate}
\end{definition}

We also recall the notion of cellular maps on $S^2$, which are a class of maps that are ``compatible'' with pairs of cell decompositions.

\begin{definition}[Cellular maps and cellular Markov partitions]
    Let $\mathbf{D}'$ and $\mathbf{D}$ be two cell decompositions of $S^2$. We say that a continuous map $f \: S^2 \to S^2$ is \emph{cellular} for $(\mathbf{D}', \mathbf{D})$ if for every cell $c \in \mathbf{D}'$, the restriction $f|_c$ of $f$ to $c$ is a homeomorphism of $c$ onto a cell in $\mathbf{D}$. We say that $(\mathbf{D}',\mathbf{D})$ is a \emph{cellular Markov partition} for $f$ if $f$ is cellular for $(\mathbf{D}', \mathbf{D})$ and $\mathbf{D}'$ is a refinement of $\mathbf{D}$.
\end{definition}

Let $f \: S^2 \to S^2$ be a Thurston map, and $\mathcal{C}\subseteq S^2$ be a Jordan curve containing $\post f$. 
Then the pair $(f,\mathcal{C})$ induces a sequence of cell decompositions of $S^2$ in the following way.

By the Sch\"onflies theorem, the set $S^2 \smallsetminus \mathcal{C}$ has two connected components, whose closures are both homeomorphic to $\overline{\D}$ and have a common boundary $\cC$. We call the closure of one of them the \emph{white $0$-tile} for $(f,\mathcal{C})$, denoted by $X^0_+$, and the closure of the other one the \emph{black $0$-tile} for $(f,\mathcal{C})$, denoted by $X^0_-$. 
The set of \emph{$0$-tiles} is defined by $\X^0(f, \mathcal{C}) \= \bigl\{X^0_-, \, X^0_+ \bigr\}$. 
We call $\V^0(f, \mathcal{C}) \= \post f$ the set of \emph{$0$-vertices}. 
We set $\overline{\V}^0(f, \mathcal{C}) \= \bigl\{ \{x\} \describe x\in \V^0(f,\mathcal{C}) \bigr\}$. 
We call the set of the closures of the connected components of $\mathcal{C} \smallsetminus \post f$ the set of \emph{$0$-edges} and denote it by $\E^0(f,\mathcal{C})$. 
Then we get a cell decomposition $
    \mathbf{D}^0(f,\mathcal{C}) \= \X^0(f, \mathcal{C}) \cup \E^0(f,\mathcal{C}) \cup \overline{\V}^0(f,\mathcal{C})
$
of $S^2$. We call the cells in this decomposition \emph{cells of level} $0$, or simply $0$-\emph{cells}.

Starting with $\mathbf{D}^0(f,\mathcal{C})$, we can recursively define a sequence $\{\mathbf{D}^n(f,\mathcal{C})\}_{n\in\N_0}$ such that $f$ is cellular for all $\bigl(\mathbf{D}^{n+1}(f,\mathcal{C}), \mathbf{D}^n(f,\mathcal{C})\bigr)$. We call cells in $\mathbf{D}^n(f,\mathcal{C})$ the \emph{$n$-cells}. Note that such a sequence of cell decompositions is uniquely determined by $(f,\cC)$; cf. \cite[Corollary~5.13]{BM}.
We denote by $\X^n(f,\mathcal{C})$ the set of $n$-cells of dimension $2$, which are called \emph{$n$-tiles}; by $\E^n(f,\mathcal{C})$ the set of $n$-cells of dimension $1$, called \emph{$n$-edges}; by $\overline{\V}^n(f,\mathcal{C})$ the set of $n$-cells of dimension $0$; and by $\V^n(f)=\V^n(f,\mathcal{C})$ the set $\bigl\{x :\{x\} \in \overline{\V}^n(f,\mathcal{C})\bigr\}$, called the set of \emph{$n$-vertices}.

We recall some basic properties of the cell decompositions introduced above.

\begin{lemma}[{Consequences of \cite[Proposition~5.16]{BM}}]\label{lem:cell-pullback-counting}
For $n,k\in\N_0$, $f^k$ maps each $(n+k)$-tile or $(n+k)$-edge homeomorphically onto an $n$-tile or $n$-edge, respectively. For $U\in\X^n(f,\cC)$,
\begin{equation}\label{eq:tile-interior-pullback}
 f^{-k}(\inte_\circ U)=\bigsqcup_{\substack{T\in\X^{n+k}(f,\cC)\\f^k(T)=U}}\inte_\circ T.
\end{equation}
The vertex sets satisfy
\begin{equation}\label{eq:cell-vertices}
 \V^m=f^{-m}(\post f),\quad f^{-1}(\V^m)=\V^{m+1},\quad \crit f\subseteq\V^1
 \quad \text{ for all } m\in\N_0.
\end{equation}
\end{lemma}

\begin{lemma}[{Bonk \& Meyer \cite[Lemma~5.17]{BM}}]    \label{lem:cell mapping properties of Thurston map}
    Let $k,\,n \in \N_0$, $f \: S^2 \to S^2$ be a Thurston map, and $\mathcal{C} \subseteq S^2$ be a Jordan curve with $\post f \subseteq \mathcal{C}$.
    \begin{enumerate}
        \smallskip

        \item     \label{item:lem:cell mapping properties of Thurston map:i}
        If $c \subseteq S^2$ is a cell such that $f^{k}|_{c}$ is a homeomorphism onto its image and $f^{k}(c)$ is an $n$-cell, then $c$ is an $(n + k)$-cell.

        \smallskip

        \item \label{item:lem:cell mapping properties of Thurston map:ii}
        If $X$ is an $n$-tile and $p \in S^2$ is a point with $f^{k}(p) \in \inte_\circ(X)$, then there exists a unique $(n + k)$-tile $X'$ with $p \in X'$ and $f^{k}(X') = X$.
    \end{enumerate}
\end{lemma}

For each $n\in \N_0$, the \emph{set of black $n$-tiles} is defined as $
    \X^n_-(f,\mathcal{C}) \= \bigl\{ X \in \X^n (f,\mathcal{C}) \describe f^n(X) = X^0_- \bigr\}$,
and the \emph{set of white $n$-tiles} as $
    \X^n_+(f,\mathcal{C}) \= \bigl\{X\in \X^n(f,\mathcal{C}) \describe f^n(X) = X^0_+\bigr\}$.

From now on, if the map $f$ and the Jordan curve $\mathcal{C}$ are clear from the context, we will sometimes omit $(f,\mathcal{C})$ from the notation above.

We next consider the number and cyclic order of the $n$-tiles incident to a vertex.

\begin{lemma}[{Bonk \& Meyer \cite[Lemma~5.28~\textup{(i)}]{BM}}]\label{lem:incident-tile-count}
For $m\in \N_0$ and $v\in\V^m$, exactly~$2\deg_{f^m}(v)$ distinct $m$-tiles meet at $v$.
\end{lemma}

\begin{lemma}[{A consequence of \cite[Lemma~5.21~\textup{(iii)}]{BM}}]\label{lem:tile-color-alternation}
For $m\in \N_0$ and $v\in\V^m$, the $m$-tiles incident to $v$ alternate in color in positive cyclic order.
\end{lemma}

A Jordan curve $\mathcal{C} \subseteq S^2$ is \emph{$f$-invariant} if $f(\mathcal{C}) \subseteq \mathcal{C}$. If $\mathcal{C}$ is $f$-invariant with $\post f \subseteq \mathcal{C}$, then the cell decompositions $\mathbf{D}^{n}(f, \mathcal{C})$ have subdivision properties. More precisely, $\mathbf{D}^{n+k}(f, \mathcal{C})$ is a refinement of $\mathbf{D}^{n}(f, \mathcal{C})$, whenever $n,\,k \in \N_0$. Intuitively, this means that each cell decomposition $\mathbf{D}^{n}(f,\mathcal{C})$ is ``subdivided'' by cells in $\mathbf{D}^{n+k}(f, \mathcal{C})$. Each cell $c\in \mathbf{D}^{n}(f, \mathcal{C})$ is actually subdivided by the cells in $\mathbf{D}^{n+k}(f,\mathcal{C})$ in the same way that the cell $f^n(c) \in \mathbf{D}^{0}(f,\mathcal{C})$ is subdivided by the cells in $\mathbf{D}^{k}(f, \mathcal{C})$. 

These cell decompositions fit into the framework of finite subdivision rules developed by Cannon, Floyd, and Parry in \cite{CFP01}. For postcritically finite rational maps without periodic critical points, they proved that every sufficiently high iterate is the subdivision map of a finite subdivision rule \cite{CFP07}.

We recall \cite[Proposition~12.5~(ii)]{BM}, which is easy to check but useful.

\begin{prop} \label{prop:cell decomposition: invariant Jordan curve}
    Let $k,\,n \in \N_0$, $f \: S^2 \to S^2$ be a Thurston map, and $\mathcal{C} \subseteq S^2$ be an $f$-invariant Jordan curve with $\post f \subseteq \mathcal{C}$. Then every $(n+k)$-tile $X^{n+k}$ is contained in a unique $k$-tile $X^k$.
\end{prop}

In this article, we focus on the class of expanding Thurston maps defined as follows. We refer to \cite[Section~6.1]{BM} for more details.

\begin{definition}[Expansion]     \label{def:expanding_Thurston_maps}
    A Thurston map $f \: S^2 \to S^2$ is called \emph{expanding} if there exists a metric $d$ on $S^2$ that induces the standard topology on $S^2$ and a Jordan curve $\mathcal{C} \subseteq S^2$ containing $\post{f}$ such that
    \begin{equation*}    \label{eq:definition of expansion}
        \lim_{n \to +\infty} \max\{ \diam_{d}(X) : X \in \X^n(f,\mathcal{C}) \} = 0.
    \end{equation*}
\end{definition}
Here $\diam_{d}(X)$ denotes the diameter of $X$ with respect to the metric $d$.

\begin{remark}\label{rem:Expansion_is_independent}
It is clear that if $f$ is an expanding Thurston map, then so is $f^{n}$ for each $n \in \N$. We note that being expanding is a topological property of a Thurston map and is independent of the choice of the metric $d$ that induces the standard topology on $S^2$. By \cite[Lemma~6.2]{BM}, it is also independent of the choice of the Jordan curve $\mathcal{C}$ containing $\post{f}$. More precisely, if $f$ is an expanding Thurston map, then $\lim_{n \to +\infty} \max \bigl\{ \diam_{\widetilde{d}}{X} : X\in \X^n\bigl(f,\widetilde{\mathcal{C}}\bigr) \bigr\} = 0$, for each metric $\widetilde{d}$ that induces the standard topology on $S^2$ and each Jordan curve $\widetilde{\mathcal{C}} \subseteq S^2$ that contains $\post{f}$.

\end{remark}

For an expanding Thurston map $f$, $f$ induces a natural class of metrics called \emph{visual metrics for $f$}.
One major advantage of a visual metric $\rho$ is that in $\bigl(S^2,\rho\bigr)$ we have good quantitative control over the sizes of the cells in the cell decompositions discussed above. 
For a visual metric $\rho$ for $f$, there exists a unique constant $\Lambda > 1$, called the \emph{expansion factor} of $\rho$ (see \cite[Chapter~8]{BM} for more details), such that $\diam_{\rho}(\tau)\asymp \Lambda^{-n}$ for all $n$-edges and all $n$-tiles $\tau$ and the constant $C(\asymp)$ is independent of the cells and their level $n$.
For the existence and properties of such metrics, see \cite[Chapter~8]{BM} for a detailed discussion.

A key property of expanding Thurston maps is the existence of invariant Jordan curves for sufficiently high iterates.

\begin{theorem}[{Bonk \& Meyer \cite[Theorem~15.1]{BM}}]    \label{lem:invariant_Jordan_curve}
    Let $f \: S^2 \to S^2$ be an expanding Thurston map, and $\widetilde{\mathcal{C}} \subseteq S^2$ be a Jordan curve with $\post f \subseteq \widetilde{\mathcal{C}}$. Then there exists $N\bigl(f, \widetilde{\mathcal{C}}\bigr) \in \N$ such that for each integer $n \geq N\bigl(f,\widetilde{\mathcal{C}}\bigr)$ there exists an $f^n$-invariant Jordan curve $\mathcal{C}$ isotopic to $\widetilde{\mathcal{C}}$ rel.\ to $\post f$.
\end{theorem}

Related invariant graph constructions were developed by Hlushchanka \cite{Hlushchanka17}, who showed that every sufficiently high iterate of an expanding Thurston map admits an invariant planar embedded tree containing the postcritical set.

Here we record a basic example of expanding Thurston maps having an invariant Jordan curve.

\begin{figure}[htbp]
    \centering
    \begin{overpic}[
        width=0.9\textwidth,
        trim=-210bp -200bp -265bp -150bp
    ]{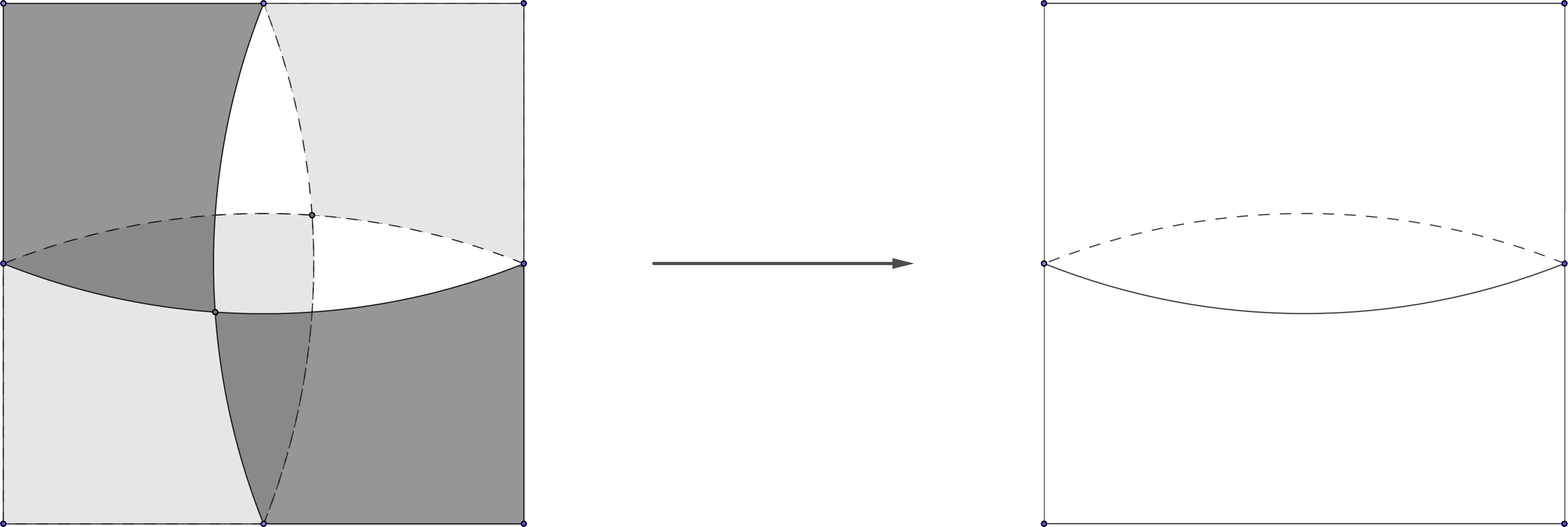}
    \setlength{\unitlength}{0.009\textwidth}
        \put(6.1,35.1){$-1$}
        \put(35.4,32.5){$\infty$}
        \put(35.4,3.3){$1$}
        \put(6.2,3.2){$0$}

        \put(64.6,35.3){$-1$}
        \put(94.0,32.5){$\infty$}
        \put(93.7,3.3){$1$}
        \put(64.6,3.2){$0$}

        \put(48.7,20.8){$g$}
    \end{overpic}
    \caption{The Latt\`es map $g$ in Example~\ref{ex:lattes-map}. The white $1$-tiles are mapped to the front side of the pillow on the right, and the black tiles are mapped to the back side. The vertices $-1,\,\infty,\,1,\,0$ are the postcritical points.}
    \label{2times2}
\end{figure}

\begin{example}[The $2\times2$ Latt\`es map]\label{ex:lattes-map}
    We consider the Latt\`es map $g\:\CDach\to \CDach$ defined in the following way (see \cite[Subsection~1.5]{Meyer13} for a more detailed construction). Let $L\=\Z+\mathbf{i}\Z$ and let $\wp\:\C\to\CDach$ be the normalized Weierstrass function with period lattice $L$. The map $\psi\:\C\to \C$, $\psi(z)\=2z$, induces $g$ by $g\circ\wp=\wp\circ\psi$. Explicitly, $g(z)\=\frac{4z(1-z^2)}{(1+z^2)^2}$.

    The map $g$ is an expanding Thurston map with degree $4$, and $\post g=\{0,\,1,\,\infty,\,-1\}$. The Jordan curve $\cC=\RDach$ is $g$-invariant and contains $\post g$.

    As shown in Figure~\ref{2times2}, the sphere $\CDach$ can be considered as a ``pillow'' obtained by gluing two squares together along their boundaries. Each square is subdivided into four smaller squares, giving eight $1$-tiles.
\end{example}

We conclude this subsection with a brief introduction to the invariant Peano curve.

\begin{theorem}[{Meyer \cite[Theorem~1.1]{Meyer13}}] \label{Peanocurves}
Let $f$ be an expanding Thurston map. Then for each sufficiently high iterate $F=f^n$ there is a Peano curve $\gamma\:S^1\to S^2$ (onto) such that $F(\gamma(z))=\gamma\bigl(z^d\bigr)$ (for all $z\in S^1$). Here $d=\deg F$. This means that the following diagram commutes:
\begin{equation*}
\begin{tikzcd}
S^1 \arrow[r,"z^d"] \arrow[d,"\gamma"'] &
S^1 \arrow[d,"\gamma"]\\
S^2 \arrow[r,"F"'] & S^2.
\end{tikzcd}
\end{equation*}

Furthermore, we can approximate the Peano curve $\gamma$ as follows. There is a homotopy $\Gamma \:S^2\times[0,1]\to S^2$, with $\Gamma(z,0)= z$, such that $\Gamma(z,1)=\gamma(z)$ for all $z\in S^1$. Here we view $S^1\subseteq S^2$ as a subset. In fact, $\Gamma$ may be chosen to be a pseudo-isotopy,
meaning that it is an isotopy on $[0,1)$.
\end{theorem}

\begin{figure}
    \centering
    \begin{overpic}[permil,width=0.9\linewidth]{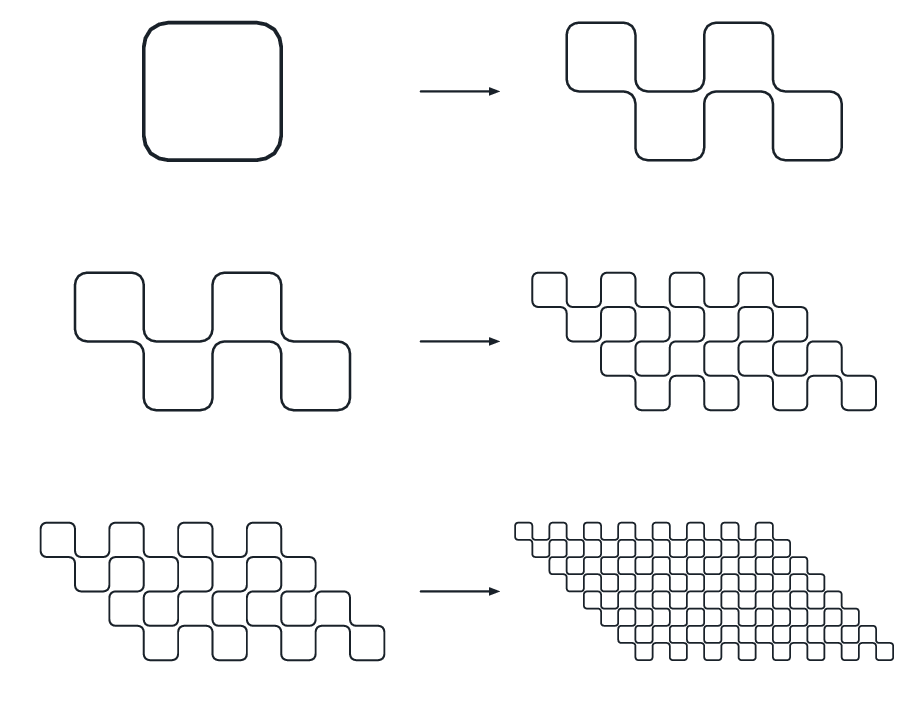}
    \small
    \put(71,605){\makebox(0,0){$\gamma^0$}}
    \put(499,702){\makebox(0,0){$H^0$}}
    \put(605,605){\makebox(0,0){$\gamma^1$}}
    \put(71,333){\makebox(0,0){$\gamma^1$}}
    \put(499,431){\makebox(0,0){$H^1$}}
    \put(605,333){\makebox(0,0){$\gamma^2$}}
    \put(71,62){\makebox(0,0){$\gamma^2$}}
    \put(499,160){\makebox(0,0){$H^2$}}
    \put(605,62){\makebox(0,0){$\gamma^3$}}
    \end{overpic}
    \caption{Successive approximations for Latt\`es map $g$ in Example~\ref{ex:lattes-map}, shown in the orbifold covering.}
    \label{fig:peano-approximations}
\end{figure}

Generally speaking, the construction in \cite{Meyer13} of the invariant Peano curve $\gamma$ proceeds as follows.
Let $F$ be a sufficiently high iterate such that there exists an $F$-invariant Jordan curve $\cC \supseteq \post F$.
The curve $\gamma$ is obtained as the uniform limit of a sequence $\{\gamma^n\}_{n\in\N_0}$ of curves, and $\gamma^n$ is called the \emph{$n$-th approximation} of $\gamma$. To construct the approximations, we start with an
$F$-invariant Jordan curve, letting $\gamma^0\bigl(S^1\bigr) = \cC$, and deform it via a suitable \emph{pseudo-isotopy} $H^0$ (rel.\ $\post F$) onto a curve $\gamma^1$. Then, recursively, we obtain $\gamma^{n+1}$ by deforming $\gamma^n$ using the lift $H^n$ of $H^0$ under $F^n$. Passing to the limit gives us the desired Peano curve $\gamma$. Figure~\ref{fig:peano-approximations} illustrates
the construction of approximation curves for Example~\ref{ex:lattes-map}.

A technical step is to obtain $\gamma^1$ and $H^0$, for which we give an intuitive description here (see also Figure~\ref{fig:peano-approximations} for an illustration). The idea is to join the white $1$-tiles at selected vertices
so that the resulting ``white region'' is a Jordan region, and we consider its boundary. By ``shrinking'' the ``white region” near the $1$-vertices, we turn the boundary into the desired $\gamma^1$ (with a suitable parameterization). The joins are also chosen so that we can deform $\gamma^0$ to $\gamma^1$ rel.\ $\post F$. To rigorously describe how tiles are joined, Meyer introduced the notions of \emph{connections} and \emph{connection graphs}.

From Subsection~\ref{subsec:connections} to \ref{wbequivalencerelations}, we recall from \cite{Meyer13} and \cite{Meyer09} the notions of the connection, the invariant Peano curve for expanding Thurston maps, and the corresponding equivalence relations for the convenience of the reader. We also list some results in Meyer's articles we may use later in the proof. Readers who are already familiar with this material may skip these subsections.

\subsection{Non-crossing connections at a vertex}
\label{subsec:connections}
We first recall the local combinatorial data, namely the connections, that describe how tiles are joined near a $1$-vertex. We refer to
\cite[Section~6]{Meyer13} for a detailed discussion. 

Following \cite[Section~6.2]{Meyer13}, let $v$ be a $1$-vertex at which $2n$ $1$-tiles meet. We label these tiles $X_0,\,X_1,\,\ldots,\,X_{2n-1}$ according to their positive cyclic order around $v$. Moreover, the labeling is chosen so that the white tiles have even indices and the black tiles have odd indices. Set  
\begin{equation}\label{evenoddlabel}
    \operatorname{even}_n\=\{0,\,2,\,\dots,\,2n-2\} 
    \quad \text{ and } \quad
    \operatorname{odd}_n\=\{1,\,3,\,\dots,\,2n-1\}.
\end{equation}

To describe the way that tiles $X_0,\,X_1,\,\ldots,\,X_{2n-1}$ join at $v$, we consider partitions of the set of labels.

For $A\subseteq \{0,\,1,\,\dots,\,2n-1\}$, a \emph{partition} $\pi(A)$ of $A$ is a set $\pi(A)\=\{b_1,\,b_2,\,\dots,\,b_N\}$, satisfying the following properties:
\begin{enumerate}
   \smallskip \item Each $b_i$ (called \emph{block}) is a nonempty subset of $A$.
   \smallskip \item $b_i\cap b_j=\emptyset$ for all distinct $i,\,j\in\{1,\,2,\,\dots,\,N\}$.
   \smallskip \item $\bigcup_{i=1}^{N} b_i=A$.
\end{enumerate}

Two distinct blocks $b,\,c\in\pi(A)$ are said to \emph{cross} if there exist $i,\,k\in b$ and $j,\,l\in c$ such that $i<j<k<l$ (considering the cyclic order on $\{0,\,1,\,\dots,\,2n-1\}$).

\begin{definition}[Non-crossing partition]
    A partition $\pi(A)$ of a finite set $A\subseteq \N_0$ is called \emph{non-crossing} if no two distinct blocks of $\pi(A)$ cross.
\end{definition}

Let $\pi_{\white}$ be a non-crossing partition of $\operatorname{even}_n$. By \cite[Lemma~6.1]{Meyer13}, there exists a unique maximal partition $\pi_{\black}$ of $\operatorname{odd}_n$ such that $\pi_{\white}\cup\pi_{\black}$ is a non-crossing partition of $\{0,\,1,\,\dots,\,2n-1\}$. Here, the terminology maximal partition means that whenever $\sigma$ is a partition of $\operatorname{odd}_n$ with $\pi_{\white}\cup\sigma$ non-crossing, every block of $\sigma$ is contained in a block of $\pi_{\black}$. We call this $\pi_{\black}$ the \emph{partition complementary to $\pi_{\white}$} and denote it by $\pi_{\black}(\pi_\white)$. A partition $\pi_{\white}\cup\pi_{\black}(\pi_{\white})$ as above is called a \emph{complementary non-crossing partition}, or \emph{cnc-partition}.

We consider a geometric realization of cnc-partitions.

\begin{definition}[Geometric cnc-partition]
Let $g_1,\,g_2,\,\dots,\,g_n\subseteq\overline{\mathbb{D}}$ be pairwise disjoint Jordan arcs such that
\begin{enumerate}
\smallskip
    \item the (distinct) endpoints of each $g_j$ lie in $S^1=\partial\mathbb{D}$ and the interior of $g_j$ lies in $\mathbb{D}$,
    \smallskip
    \item the arcs $\{g_j\}_{j=1}^n$ divide $S^1$ into $2n$ circular arcs $a_0,\,a_1,\,\dots,\,a_{2n-1}\subseteq S^1$ (labeled positively on $S^1$), where $a_j\=[e_j,e_{j+1}]\subseteq S^1$ and $e_j\=e^{\I 2j\pi/(2n)}$,
\end{enumerate}
and let $D_1,\,D_2,\,\dots,\,D_{n+1}$ be the (simply connected) domains obtained by dividing $\mathbb{D}$ with $g_1,\,g_2,\,\dots,\,g_n$. We call
\begin{equation*}
    \overline{\mathbb{D}}(g_1,g_2,\dots,g_n)\=\{\overline{D_1},\,\overline{D_2},\,\dots,\,\overline{D_{n+1}}\}
\end{equation*}
a \emph{geometric cnc-partition}.

Such a $\overline{\mathbb{D}}(g_1,g_2,\dots,g_n)$ induces a partition $\pi(\{g_m\}_{m=1}^n)$ of $\{0,\,1,\,\dots,\,2n-1\}$ in such a way that every $\overline{D}_l\in\overline{\mathbb{D}}(g_1,g_2,\dots,g_n)$ corresponds to exactly one block, which means $i,\,j\in\{0,\,1,\,\dots,\,2n-1\}$ are in the same block of $\pi(\{g_m\}_{m=1}^n)$ if and only if $a_i$ and $a_j$ are in the boundary of the same $\overline{D}_l\in\overline{\mathbb{D}}(g_1,g_2,\dots,g_n)$.
\end{definition}

\begin{lemma}[{Meyer \cite[Lemma~6.3]{Meyer13}}]\label{cncgeo}
The following hold:
\begin{enumerate}
\smallskip
    \item Let $\overline{\mathbb{D}}(g_1,g_2,\dots,g_n)$ be a geometric cnc-partition, then the induced partition $\pi(\{g_m\}_{m=1}^n)$ is a cnc-partition of $\{0,\,1,\,\dots,\,2n-1\}$.
\smallskip
    \item Conversely, for each cnc-partition $\pi$ of $\{0,\,1,\,\dots,\,2n-1\}$, there exists a geometric cnc-partition $\overline{\mathbb{D}}(g_1,g_2,\dots,g_n)$ such that $\pi$ is identical to $\pi(\{g_m\}_{m=1}^n)$.
\smallskip
    \item Let $\overline{D}_k,\,\overline{D}_l\in\overline{\mathbb{D}}(g_1,g_2,\dots,g_n)$ be distinct and their corresponding blocks be $b_k,\,b_l\in\pi(\{g_m\}_{m=1}^n)$. Then either $\overline{D}_k\cap\overline{D}_l=\emptyset$ or $\overline{D}_k\cap\overline{D}_l=g_m$ for some arc $g_m$. The latter case is equivalent to that $b_k,\,b_l$ are adjacent blocks in $\pi(\{g_m\}_{m=1}^n)$. Therefore, there is a one-to-one correspondence between pairs of adjacent blocks in $\pi(\{g_m\}_{m=1}^n)$ and the arcs $\{g_j\}_{j=1}^n$.
\end{enumerate}
\end{lemma}

\begin{definition}[Marking, marked arc]\label{markedarc}
A \emph{marking} of a cnc-partition $\pi_{\white}\cup\pi_{\black}$ is a pair $(b, c)$ of adjacent blocks $b,\,c\in\pi_{\white}\cup\pi_{\black}$. We also say that $\pi_{\white}\cup\pi_{\black}$ is \emph{marked by $(b, c)$} and $b,\,c$ are
the \emph{marked blocks}. Let $\overline{D}_l\in\overline{\mathbb{D}}(g_1,g_2,\dots,g_n)$ be a geometric cnc-partition with $\pi_{\white}\cup\pi_{\black}=\pi(\{g_m\}_{m=1}^n)$. When no ambiguity arises, we also call the arc $g_m$ corresponding to the marking $(b, c)$ (as in
Lemma~\ref{cncgeo}~\textup{(iii)}) the \emph{marked arc} of $\pi_{\white}\cup\pi_{\black}$.
\end{definition}

For each $\overline{D}_l\in\overline{\mathbb{D}}(g_1,g_2,\dots,g_n)$, color $\overline D_l$ white if the corresponding block $b_l\in\pi(\{g_m\}_{m=1}^n)$ belongs to $\pi_{\white}$, and black otherwise. Note that sets sharing a side $g_m$ (in which case the corresponding blocks are adjacent) have different colors.

\begin{figure}
    \centering
    \begin{overpic}[permil,width=0.9\linewidth]{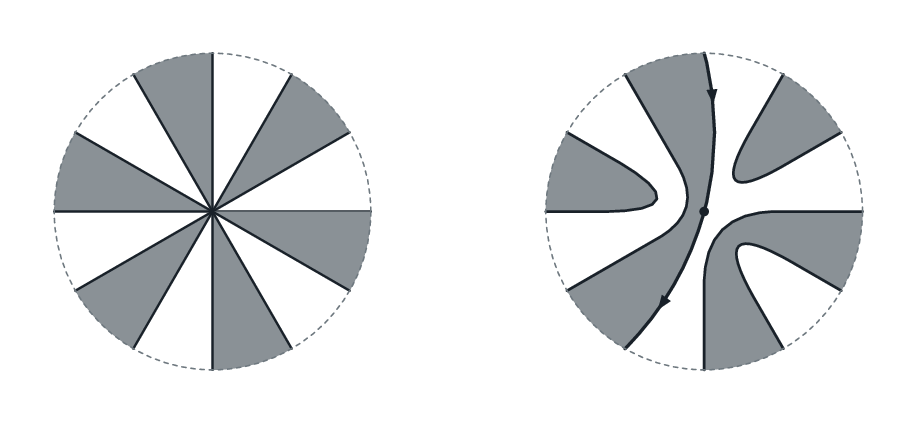}
        \footnotesize
        \put(424,282){\makebox(0,0){$X_0$}}
        \put(372,373){\makebox(0,0){$X_1$}}
        \put(282,425){\makebox(0,0){$X_2$}}
        \put(178,425){\makebox(0,0){$X_3$}}
        \put(87,373){\makebox(0,0){$X_4$}}
        \put(35,282){\makebox(0,0){$X_5$}}
        \put(35,178){\makebox(0,0){$X_6$}}
        \put(87,88){\makebox(0,0){$X_7$}}
        \put(178,36){\makebox(0,0){$X_8$}}
        \put(282,36){\makebox(0,0){$X_9$}}
        \put(372,88){\makebox(0,0){$X_{10}$}}
        \put(424,178){\makebox(0,0){$X_{11}$}}
        \put(958,282){\makebox(0,0){$X_0$}}
        \put(906,373){\makebox(0,0){$X_1$}}
        \put(816,425){\makebox(0,0){$X_2$}}
        \put(711,425){\makebox(0,0){$X_3$}}
        \put(621,373){\makebox(0,0){$X_4$}}
        \put(569,282){\makebox(0,0){$X_5$}}
        \put(569,178){\makebox(0,0){$X_6$}}
        \put(622,88){\makebox(0,0){$X_7$}}
        \put(711,36){\makebox(0,0){$X_8$}}
        \put(816,36){\makebox(0,0){$X_9$}}
        \put(906,88){\makebox(0,0){$X_{10}$}}
        \put(958,178){\makebox(0,0){$X_{11}$}}
        \put(764,432){\makebox(0,0){$E$}}
        \put(663,56){\makebox(0,0){$E'$}}
        \put(788,232){\makebox(0,0){$0$}}
    \end{overpic}
    \caption{A connection at a vertex and its geometric representation. White tiles have even indices and black tiles have odd indices. The arrows indicate the succeeding pair $E,\,E'$.}
    \label{fig:vertex-connection}
\end{figure}

\begin{definition}[Geometric representation of cnc-partition]\label{6.5}
Let $\pi_{\white}\cup\pi_{\black}$ be a cnc-partition. We call a geometric cnc-partition $\overline{\mathbb{D}}(g_1,g_2,\dots,g_n)$ satisfying $\pi_{\white}\cup\pi_{\black}=\pi(\{g_m\}_{m=1}^n)$ a geometric
representation of $\pi_{\white}\cup\pi_{\black}$, denoted by $\overline{\mathbb{D}}(\pi_{\white}\cup\pi_{\black})$. When $\pi_{\white}\cup\pi_{\black}$ is marked as in Definition~\ref{markedarc}, we further require that the marked arc contains the origin. The union of white (resp.\ black)
sets in $\overline{\mathbb{D}}(\pi_{\white}\cup\pi_{\black})$ is denoted by $\overline{\D}_{\white}(\pi_{\white}\cup\pi_{\black})$ (resp.\ $\overline{\D}_{\black}(\pi_{\white}\cup\pi_{\black})$). When no ambiguity arises,
we also denote the closed disk with $n$ disjoint Jordan arcs as above by $\overline{\mathbb{D}}(\pi_{\white}\cup\pi_{\black})$.
\end{definition}

Recall that the $2n$ $1$-tiles containing a $1$-vertex $v$ are labeled $X_0,\,X_1,\,\dots,\, X_{2n-1}$ in positive cyclic order, with even indices for white tiles and odd indices for black tiles.

\begin{definition}[Connection at a vertex]\label{def_Connection_at_a_vertex}
Let $f\:S^2\to S^2$ be an expanding Thurston map, $\cC$ be an $f$-invariant Jordan curve with $\cC\supseteq \post f$, and $v$ be an $1$-vertex where $2n$ $1$-tiles meet. A \emph{connection at $v$} consists of a labeling of $2n$ $1$-tiles containing $v$ (as in the beginning of this subsection), and a cnc-partition $\pi_{\white}(v)\cup\pi_{\black}(v)$, where $\pi_{\white}(v)$ is a partition of $\operatorname{even}_n$ representing the white $1$-tiles and $\pi_{\black}(v)$ is a partition of $\operatorname{odd}_n$ representing the black $1$-tiles. When $v\in\post f$, $\pi_{\white}(v)\cup\pi_{\black}(v)$ is marked in addition. 

Given a connection at $v$, two $1$-tiles $X_i$ and $X_j$ of the same color are said to be \emph{connected at $v$} (or $X_i$ is said to \emph{connect to $X_j$}) if $i$ and $j$ are contained in the same block of $\pi_{\white}(v)\cup\pi_{\black}(v)$. For each block $b\in\pi_{\white}(v)\cup\pi_{\black}(v)$ and each label
$i \in b$, we say that the $1$-tile $X_i$ is \emph{incident} (at $v$) to the block $b$.
\end{definition}

\begin{remark}
    Given a connection at $v$, the labeling of tiles implies that $\operatorname{even}_n$ (resp.\ $\operatorname{odd}_n$) represents the white (resp.\ black) $1$-tiles. It is clear from the definition that $1$-tiles of different colors are never connected.
\end{remark}

We use the geometric representations of cnc-partitions to describe the connections at vertices.
 
For $j\in\{0,\,1,\,\dots,\,2n-1\}$, a \emph{sector} in $\overline{\mathbb{D}}$ denoted by $S_j$ is defined as follows:
\begin{equation*}
    S_j\=\bigl\{re^{2\pi \mathbf{i}\theta}: j/(2n)\leq\theta\leq (j+1)/(2n)\text{ and }0\leq r\leq1\bigr\}.
\end{equation*}
Following \cite[Subsection~6.2]{Meyer13}, for every such $1$-vertex $v$ where $2n$ $1$-tiles meet, we can find a neighborhood $U(v)$, as small as we want, that is closed and homeomorphic to the closed disk $\overline{\mathbb{D}}$, such that there is a homeomorphism $h_v\: U(v)\to\overline{\mathbb{D}}$ that maps $v$ to $0$ and maps tiles to sectors, i.e., $h(X_j\cap U(v))=S_j$ for each $j\in\{0,\,1,\,\dots,\,2n-1\}$. Fix $U(v)$ and $h_v$ for each $v$, and we require that the neighborhoods $U(v)$ and $U(v')$ for distinct $1$-vertices $v$ and $v'$ are disjoint.

\begin{definition}[Geometric representation of a connection]\label{Geometric representation of a connection}
Let $f\:S^2\to S^2$ be an expanding Thurston map and $\cC$ be an $f$-invariant Jordan curve with $\cC\supseteq \post f$. Let $v$ be a $1$-vertex and $h_v$, $U(v)$ be as above. 

Let a connection at $v$ be given, with a cnc-partition $\pi_{\white}(v)\cup\pi_{\black}(v)$, geometrically represented by $\overline{\mathbb{D}}(\pi_{\white}(v)\cup\pi_{\black}(v))$ as in Definition~\ref{6.5}.  A \emph{geometric representation of the connection} at $v$ is a decomposition of $U(v)$ given by replacing $U(v)\cap(X_0\cup X_2\cup\dots\cup X_{2n-2})$ by $h_v^{-1}\bigl(\overline{\mathbb{D}}_{\white}(\pi_{\white}(v)\cup\pi_{\black}(v))\bigr)$ and replacing $U(v)\cap(X_1\cup X_3\cup\dots\cup X_{2n-1})$ by $h_v^{-1}\bigl(\overline{\mathbb{D}}_{\black}(\pi_{\white}(v)\cup\pi_{\black}(v))\bigr)$.

When $v$ is a postcritical point, the cnc-partition $\pi_{\white}(v)\cup\pi_{\black}(v)$ of the connection at $v$ is marked (see Definition~\ref{def_Connection_at_a_vertex}). Suppose that in the geometric representation $\overline{\mathbb{D}}(\pi_{\white}(v)\cup\pi_{\black}(v))=\overline{\mathbb{D}}(g_1,g_2,\dots,g_n)$, the marked arc is $g_m$.  
We define the \emph{geometric representation of the marked arc} to be the arc $h_v^{-1}(g_m)\subseteq U(v)$. 
\end{definition}

\begin{remark}
    Since the marked arc $g_m$ contains $0$ (see Definition~\ref{6.5}), and $h_v$ is chosen such that $h_v(v)=0$, it follows that $v\in h_v^{-1}(g_m)$, the geometric representation $h_v^{-1}(g_m)$ contains $v$.
\end{remark}

\begin{definition}[Connection of $1$-tiles]
Given an expanding Thurston map $f\:S^2\to S^2$ and an $f$-invariant Jordan curve $\cC$ with $\cC\supseteq \post f$, a \emph{connection of $1$-tiles} is an assignment of a connection at every $1$-vertex. Representing the connection at each $1$-vertex geometrically as in Definition~\ref{Geometric representation of a connection} gives a \emph{geometric representation} of this connection of $1$-tiles.
\end{definition}

\begin{lemma}[{Meyer \cite[Lemma~7.5]{Meyer13}}]\label{7.5}
Let $f_0$ be an expanding Thurston map. For each sufficiently high $n\in\N$ there is a Jordan curve $\cC$ with
$\post f_0\subseteq \cC$ satisfying the following condition:
\begin{enumerate}
\smallskip
    \item $\cC$ is invariant for the iterate $f\=f_0^n$. This means that $f(\cC)\subseteq \cC$.
\end{enumerate}  
    The $1$-tiles for $(f,\cC)$ satisfy the following conditions:
    \begin{enumerate}
    \smallskip
    \item\label{7.5(ii)} There is no $1$-tile $X$ that joins opposite sides of $\cC$. This means that
    no $1$-tile $X$ meets disjoint $0$-edges in the case $\card\post f\geq 4$, and no
    $1$-tile $X$ intersects all three $0$-edges in the case $\card\post f=3$.
    \smallskip
    \item The $1$-tiles do not form a link in the sense of \cite[Definition~7.12]{Meyer13}.
\end{enumerate}
\end{lemma}

For $(f,\cC)$ as in Lemma~\ref{7.5}, Meyer constructs a connection of $1$-tiles in \cite[Subsection~7.4 and Proposition~8.1]{Meyer13}, and extends it inductively to $n$-tiles in \cite[Subsection~8.1]{Meyer13}. By \cite[Proposition~8.1]{Meyer13}, these connections satisfy the following properties for each $n\in\N$:

\begin{definition}[{Meyer \cite[Definition~7.1]{Meyer13}}]\label{Propertiesofconnections}
\begin{enumerate}
    \item[(C1)] The associated white connection graph is a spanning tree $K$.
    \smallskip
    \item[(C2)] The Jordan curve $\partial K_\varepsilon$ is orientation-preserving
    isotopic to $\cC=\gamma^0\bigl(S^1\bigr)$ rel.\ $\post f$, where $K_\varepsilon$ is a geometric
    representation of $K$.
\end{enumerate}
\end{definition}

\subsection{Pseudo-isotopies and approximations of the Peano curve}
\label{Pseudo_isotopies}
With a suitable choice of connections of tiles, Meyer constructed the aforementioned pseudo-isotopies $H^n$, which further give rise to approximations $\gamma^n$, $n\in\N$. These approximations converge uniformly to the invariant Peano curve. We refer the reader to \cite[Subsections~3.3 and~4.2]{Meyer13} for a more comprehensive discussion.

A \emph{pseudo-isotopy} of $S^2$ is a homotopy $H\:S^2\times[0,1]\rightarrow S^2$
such that $H_t(x)\=H(x,t)$ is a homeomorphism for every $t\in[0,1)$ and, in particular, $H_0=\operatorname{id}_{S^2}$. Note that the map $H_1$ need not be a homeomorphism.

For a subset $A\subseteq S^2$, the pseudo-isotopy $H$ is said to be \emph{relative to $A$} if $H_t(a)=a$ for each $a \in A$ and each $t\in[0,1]$.
The connection between such pseudo-isotopies and polynomial mating is developed further in Meyer's study of pseudo-equators and unmating algorithms \cite[Sections~6--9]{Meyer14}.

Let $f\: S^2\to S^2$ be an expanding Thurston map and let $\cC\subseteq S^2$ be an $f$-invariant Jordan curve containing $\post f$. In \cite{Meyer13}, Meyer constructs a suitable connection of $1$-tiles (cf. \cite[Definition~7.1]{Meyer13}; see also Definition~\ref{Propertiesofconnections}), which resulted in a pseudo-isotopy $H^0$ that yields the first approximation $\gamma^1$. See \cite[Definition~3.2]{Meyer13} for a detailed description of the properties of $H^0$. Here, we merely record some facts of $H^0$ needed in this article.

First, $H^0$ deforms $\cC$, relative to $\post f$, to the union of $1$-edges, which means, $H^0_1(\cC)=\bigcup \E^1$ (cf.\ \cite[Definition~3.2, Lemma~7.2, and Theorem~7.3]{Meyer13}). It follows that $H_1^0\bigl(X_+^0\bigr)=\bigcup\X_+^1$ since $H^0$ preserves the orientation of $S^2$.

By \cite[Lemma~3.4]{Meyer13}, for each $n\in\N$, $H^0$ admits a unique lift $H^n$ under $f^n$, meaning $f^n\circ H^n_t=H^0_t\circ f^n$ for all $t\in [0,1]$. Note that $H^{n+1}$ is also the lift of $H^n$ under $f$, i.e., $H^n_t\circ f=f\circ H^{n+1}_t$ for all $t\in [0,1]$. By \cite[Lemmas~3.4 and~3.5]{Meyer13}, $H^n$ is a pseudo-isotopy relative to $\V^n$.

We next describe the action of $H_1^m$ on $m$-edges and the preimage of the interior of an $(m+1)$-tile.

\begin{lemma}[{Meyer \cite[Lemma~3.6]{Meyer13}}]\label{lem:lifted-edge-subdivision}
For $m\in \N_0$ and $E\in\E^m(f,\cC)$, the finite set $E\cap(H_1^m)^{-1}\bigl(\V^{m+1}\bigr)$ contains both endpoints of $E$. The points of $E\cap(H_1^m)^{-1}\bigl(\V^{m+1}\bigr)$ subdivide $E$ into closed subarcs. The restriction of $H_1^m$ to each subarc is a homeomorphism onto an $(m+1)$-edge. As $E$ ranges over all $m$-edges, each $(m+1)$-edge occurs exactly once.
\end{lemma}

\begin{lemma}[{Meyer \cite[proof of Lemma~5.11~\textup{(2)}]{Meyer09}}]
\label{lem:lifted-tile-preimages}
For $m\in\N_0$ and $U\in\X^{m+1}(f,\cC)$, set
\begin{equation}\label{eq:lifted-tile-interior}
 Q_U=(H_1^m)^{-1}(\inte_\circ U).
\end{equation}
The restriction $H_1^m|_{Q_U}\:Q_U\to\inte_\circ U$ is a homeomorphism. The set $Q_U$ lies in the cell interior of an $m$-tile of the same color as $U$.

Let $J$ be one of the closed subarcs obtained by subdividing an $m$-edge in Lemma~\ref{lem:lifted-edge-subdivision}. If $H_1^m(J)\subseteq\partial_\circ U$, then $J\subseteq\partial Q_U$, where $\partial Q_U$ denotes the topological boundary of $Q_U$ in $S^2$.
\end{lemma}

Fix a visual metric $\rho$ for $f$ with expansion factor $\Lambda$.
Meyer also shows in \cite[Lemmas~3.4]{Meyer13} that
\begin{equation*}
   \diam(H^n)\=\max_{x\in S^2} \diam\{H^n(x,t):t\in[0,1]\}\lesssim \Lambda^{-n},
\end{equation*}
where the constant $C_1(\lesssim)$ is independent of the level $n$.
Consequently, for integers $n>m\geq0$ and $x\in S^2$, we have
\begin{equation}\label{diamgmsmall}
    \rho\bigl(\bigl(H^{n}_1\circ H^{n-1}_1\circ \dots\circ H^m_1\bigr)(x),x\bigr)\leq\sum_{i=m}^{n} \diam\bigl(H^i\bigr)\lesssim \Lambda^{-m},
\end{equation}
where the constant $C_2(\lesssim)$ is independent of integers $n$, $m$ and $x\in S^2$.

Next we recall some notions and basic facts concerning the approximations $\gamma^n$.
The initial approximation $\gamma^0\: S^1\to S^2$ is chosen as an embedding with image $\cC$. Furthermore, $\gamma^0\:S^1\to \cC$ is assumed to be orientation preserving, with $\cC$ positively oriented as the boundary of the white $0$-tile $X_\w^0$.

For $n\in \N$, $\gamma^n$ is inductively defined by
$\gamma^n(t)\= H_1^{n-1}\bigl(\gamma^{n-1}(t)\bigr).$ In fact, $\gamma^n\:S^1\to S^2$ parametrizes an Eulerian circuit that traverses each $n$-edge exactly once. In particular, $\gamma^n\bigl(S^1\bigr)=\bigcup\E^n$. See \cite[Definition~3.8]{Meyer09} for more details.

Fix $n\in \N_0$. A point $\alpha^n\in S^1$ satisfying $\gamma^n(\alpha^n)\in \V^n$ is called an \emph{$n$-angle}; see \cite[Subsection~3.1]{Meyer09}. The set of all n-angles is denoted by $\mathbf A^n\=(\gamma^n)^{-1}(\V^n)$. The set $\mathbf A^n$ is finite. Moreover, the approximations stabilize on $n$-angles: for every integer $m\geq n$, $\gamma^n|_{\mathbf A^n}=\gamma^m|_{\mathbf A^n}=\gamma|_{\mathbf A^n}$. Since the vertices satisfy $\V^0\subseteq \V^1\subseteq\cdots$, we also have $\mathbf A^0\subseteq \mathbf A^1\subseteq\cdots$.

The $n$-angles divide the circle $S^1$ into finitely many closed arcs, which are called the \emph{$n$-arcs}. Each $n$-arc is mapped homeomorphically by $\gamma^n$ onto an $n$-edge. Conversely, for every $n$-edge $E^n$, there exists a unique $n$-arc $a^n\subseteq S^1$ such that $\gamma^n(a^n)=\E^n$.

Put $d\=\deg f$ and let $\mu\:S^1\to S^1$ denote the map $\mu(t)\=dt\pmod1$ with $S^1$ identified with $\R/\Z$. Every $n$-angle is preperiodic under $\mu$ and hence is a rational point. Each $n$-arc is mapped homeomorphically onto an $(n-1)$-arc by the map $\mu$. More precisely, the following diagram commutes:
\begin{equation}\label{commutediaav}
\begin{tikzcd}
	{\mathbf A^{n+1}\subseteq\mathbb{R}/\mathbb{Z}} && {\mathbf A^{n}\subseteq\mathbb{R}/\mathbb{Z}} \\
	\\
	{\V^{n+1}\subseteq S^2} && {\V^{n}\subseteq S^2.}
	\arrow["\mu", from=1-1, to=1-3]
	\arrow["{\gamma^{n+1}}"', from=1-1, to=3-1]
	\arrow["{\gamma^n}", from=1-3, to=3-3]
	\arrow["f"', from=3-1, to=3-3]
\end{tikzcd}   
\end{equation}

The action of $\mu$ on the arcs gives the following identities for the angle sets.

\begin{lemma}[{A consequence of \cite[Section~3.1]{Meyer09}}]\label{lem:angle-preimages}
The angle sets $\mathbf A^n=(\gamma^n)^{-1}(\V^n)$ satisfy
\begin{equation}\label{eq:finite-angle-preimages}
 \mathbf A^{n+1}=\mu^{-1}(\mathbf A^n),\quad \mathbf A^n=\mu^{-n}\bigl(\mathbf A^0\bigr)\quad \text{ for all } n\in\N_0.
\end{equation}
\end{lemma}

The size of the arcs decreases exponentially. More precisely, there exists a constant $C>0$ such that for every $n$-arc $a^n$,
\begin{equation}\label{diamadn}
\diam_{S^1}(a^n)\lesssim d^{-n},
\end{equation}
where $S^1$ is identified with $\R/\Z$ and is equipped with its standard metric.

Each $n$-edge $E^n$ is contained in the boundary of a unique white $n$-tile $X^n_+$. We orient $\partial X^n_+$ positively, inducing an orientation on $\E^n$. From another perspective, $\E^n$ is the image of a unique $n$-arc $a^n\subseteq S^1$ under $\gamma^n$, and hence inherits the orientation of $S^1$.

Two descriptions above give the same orientation on $\E^n$. In other words, the orientation induced from the cell structure agrees with that induced from the parametrization by $\gamma^n$.

Fix a visual metric $\rho$ of $f$ with expanding factor $\Lambda$. Then for all $\varepsilon>0$, there exists a Jordan curve $\gamma^n_\varepsilon\:S^1\rightarrow S^2$ such that 
\begin{equation*}
    \norm{\gamma^n-\gamma^n_\varepsilon}_\infty\=\max\bigl\{\rho(\gamma^n(t),\gamma^n_\varepsilon(t)):{t\in S^1}\bigr\}<\varepsilon.
\end{equation*}
This means although the approximation $\gamma^n$ may touch itself, it has no genuine self-intersections.

The approximations converge uniformly to the Peano curve (see \cite[Subsection~4.3]{Meyer13} for more details). More precisely,
$\norm{\gamma^n-\gamma}_\infty\lesssim \Lambda^{-n}$.
Here $C(\lesssim)$ does not depend on $n$ (see \cite[Lemma~4.5(3)]{Meyer13}).

Moreover, with a suitable parametrization (see \cite[Remark~3.1]{Meyer09}), 
$f\circ\gamma=\gamma\circ\mu$,
where $\mu(t)=dt\pmod 1$ and $d=\deg f$. This gives an invariant Peano curve for $f$.

\subsection{Succeeding edges and connection graphs}\label{prelimiconnectiongraph}

In this subsection, we briefly recall the notions of succeeding edges and connection graphs introduced by Meyer in \cite[Subsection~3.4]{Meyer09} and \cite[Subsection~6.3 and~6.4]{Meyer13}. These notions can describe the connection in Subsection~\ref{subsec:connections}.

Let $v$ be a $1$-vertex, and suppose the $1$-tiles containing $v$ are cyclically indexed by $X_0$, $X_1$, $\dots$, $X_{2n-1}$. All indices below are understood modulo $2n$.

Let the connection at a $1$-vertex $v$ be given by a cnc-partition $\pi_{\white}(v)\cup\pi_{\black}(v)$. Let $b$ be a block of $\pi_{\white}(v)$ or $\pi_{\black}(v)$. Two indices $i,\,j\in b$ are said to be \emph{succeeding in $b$} if none of the indices $i+1,\,i+2,\,\dots,\,j-1$, read in the positive cyclic order, belongs to $b$.

Let $E$ and $E'$ be positively oriented (as boundaries of white $1$-tiles) $1$-edges such that $v$ is the terminal endpoint of $E$ and the initial endpoint of $E'$. We say that $E$ and $E'$ are \emph{succeeding at $v$} with respect to the connection if there exist white $1$-tiles $X_i$ and $X_j$ and a block $b\in\pi_{\white}$ such that $E\subseteq X_i$, $E'\subseteq X_j$, and $i,\,j$ are succeeding in $b$. We call $E'$ the \emph{successor} to $E$ at $v$. Note that we can define these things for black $1$-tiles similarly.

For the connection of $1$-tiles used in the construction of $\gamma^1$, it follows from \cite[Definition~6.21 and Lemma~7.4]{Meyer13} that two positively oriented $1$-edges $E$ and $E'$ are succeeding along $\gamma^1$ if and only if they are succeeding with respect to this connection.

For each $n>1$, the \emph{connection of the $n$-tiles} is an assignment of a connection at every $n$-vertex and is constructed inductively from the connection of the $1$-tiles (see \cite[Lemma~7.2 and Section~8]{Meyer13}). We denote its cnc-partition at an $n$-vertex $v$ by $\pi_{\white}^n(v)\cup\pi_{\black}^n(v)$. At a $1$-vertex $v$, we write $\pi_{\white}^1(v)\=\pi_{\white}(v)$ and $\pi_{\black}^1(v)\=\pi_{\black}(v)$. The succeeding $n$-edges at $v$ are defined by the same rule as for $1$-edges. For $n\in \N$, if the succeeding $n$-edges $E$ and $E'$ correspond to the marked arc in the geometric representation of the connection at $v$, we say that they are \emph{marked} succeeding $n$-edges. We call $E'$ the \emph{marked $n$-edge starting at $v$}. We also say that the block $b$ \emph{contains the marked point} or the block $b$ is \emph{marked}. By \cite[Subsection~3.5]{Meyer09}, two positively oriented $n$-edges $E$ and $E'$ are \emph{succeeding} along $\gamma^n$ if and only if they are succeeding with respect to this connection.

We next recall how the marked edges are obtained from the subdivision arcs. We also describe how each block at an $m$-vertex determines a block at level $m+1$.

\begin{lemma}[{Meyer \cite[Section~8.1, paragraph following Figure~9]{Meyer13}}]\label{lem:initial-marked-edge}
For $p\in \post f$, let $E$ be the positively oriented $0$-edge starting at $p$ and let $J\subseteq E$ be its initial subdivision arc. Then $H_1^0(J)$ is the marked $1$-edge starting at $p$.
\end{lemma}

\begin{prop}[{Meyer \cite[Section~8.1, equations~(8.3)--(8.4)]{Meyer13}}]
\label{prop:recursive-connection-rules}
For $m\in\N$ and $v\in\V^m$, put $p\=f^m(v)$, $D\=\deg_{f^m}(v)$, and $k\=\deg_f(p)$. Let $E$ be a positively oriented $m$-edge starting at $v$ whose image under $f^m$ is the positively oriented $0$-edge starting at $p$. In the subdivision of $E$ from Lemma~\ref{lem:lifted-edge-subdivision}, let $J$ be the unique closed subarc containing $v$.

Use positive cyclic order for all tile labelings below, with even indices white. At $v$, label the incident $m$-tiles by $X_i^m(v)$, $0\le i<2D$, and the incident $(m+1)$-tiles by $X_i^{m+1}(v)$, $0\le i<2kD$, so that we have $E\subseteq X_0^m(v)$ and $H_1^m(J)\subseteq X_0^{m+1}(v)$. At $p$, label the incident $1$-tiles by $X_i^1(p)$, $0\le i<2k$, starting with the white tile containing the marked $1$-edge starting at $p$.

Let $b_*\in\pi_\white^1(p)$ and $c_*\in\pi_\black^1(p)$ be the blocks of level-$1$ containing $0$ and $2k-1$, respectively. For a set of integers $b$ and an integer $q$, write $b+q\=\{i+q:i\in b\}$. Then
\begin{equation}\label{eq:recursive-block-merges}
\begin{aligned}
&b_*^{m+1}(B)\=\bigcup_{2j\in B}(b_*+2kj)\in\pi_\white^{m+1}(v)
&&\quad \text{ for all } B\in\pi_\white^m(v),\\
&c_*^{m+1}(C)\=\bigcup_{2j+1\in C}(c_*+2kj)\in\pi_\black^{m+1}(v)
&&\quad \text{ for all } C\in\pi_\black^m(v).
\end{aligned}
\end{equation}
\end{prop}

\begin{lemma}[{Meyer \cite[Section~8.1, Case~3]{Meyer13}}]\label{lem:continued-marked-edges}
Let $m\in \N$ and $v\in \post f$. If $E$ is the marked $m$-edge starting at $v$ and $J\subseteq E$ is its initial subdivision arc, then $H_1^m(J)$ is the marked $(m+1)$-edge starting at $v$. In the labelings of Proposition~\ref{prop:recursive-connection-rules}, the new marking selects the adjacent white and black blocks containing $0$ and $2\deg_{f^{m+1}}(v)-1$, respectively.
\end{lemma}

\begin{figure}
    \centering
    \begin{overpic}[permil,width=0.7\linewidth]{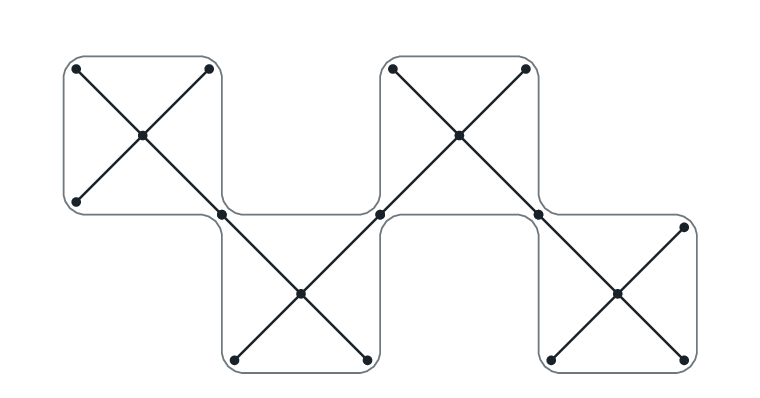}
    \small
        \put(99,483){\makebox(0,0){$c(v_1,b_1)$}}
        \put(274,483){\makebox(0,0){$c(v_2,b'_1)$}}
        \put(515,483){\makebox(0,0){$c(v_1,b_2)$}}
        \put(186,295){\makebox(0,0){$c(X_1)$}}
        \put(394,87){\makebox(0,0){$c(X_2)$}}
        \put(603,295){\makebox(0,0){$c(X_3)$}}
        \put(811,87){\makebox(0,0){$c(X_4)$}}
    \end{overpic}
    \caption{This is the white level-1 connection graph of the Latt\`es map $g$ in Example~\ref{ex:lattes-map}}
    \label{fig:connection-graph}
\end{figure}

We next recall the connection graph associated with a connection of $n$-tiles. Fix an integer $n\geq 1$. For each white $n$-tile $X$, introduce a vertex $c(X)$. For each $n$-vertex $v$ and block $b\in\pi_{\white}^n(v)$, introduce a vertex $c(v,b)$. The vertices $c(X)$ and $c(v,b)$ are joined by an edge if and only if $X$ is incident to $b$ at $v$. The resulting graph is called the \emph{$n$-th white connection graph} and is denoted by $\Gamma^+_n$. The black connection graph is defined analogously (see \cite[Subsection~6.3]{Meyer13} and \cite[Subsection~3.6]{Meyer09}).
Figure~\ref{fig:connection-graph} gives an example of the white level-$1$ connection graph of a $2\times2$ Latt\`es map.

We identify a $1$-tile $X$ with the corresponding vertex $c(X)$ of the level-$1$ white connection graph. Two white $1$-tiles $X$ and $Y$ are said to be \emph{connected} if the vertices $c(X)$ and $c(Y)$ belong to the same connected component of the white connection graph.

\begin{definition}[Cluster]
	A \emph{white (resp.\ black) cluster} $K$ is a component of the white (resp.\ black) connection graph.

	We say that a $1$-tile $X$ is \emph{contained in $K$}, denoted by $X\subseteq K$, if $c(X)\in K$. When no ambiguity arises, we also denote by $K$ the union of all $1$-tiles contained in it.
    
    A $1$-edge $E$ (resp.\ a $1$-vertex $v$) is \emph{contained in $K$} if there exists a $1$-tile $X$ with $E\subseteq X\subseteq K$ (resp.\ $v\in X\subseteq K$). 

	Every $1$-tile is contained in a unique cluster of the same color, and every $1$-edge is contained in exactly two clusters, one white and one black. A $1$-vertex $v$ can be contained in more than one cluster.

    A cluster $K$ is a \emph{tree} if the underlying component of the connection graph is a tree, i.e., contains no cycles. The white cluster $K$ is a \emph{spanning tree} if it is a tree and contains all white $1$-tiles.
\end{definition}

For the connection of $m$-tiles, with $m \in \N$, we take the geometric representation at each $m$-vertex as for the connection of $1$-tiles. We choose the vertex neighborhoods to be pairwise disjoint and leave the tiles unchanged outside them.

\begin{lemma}[{A consequence of \cite[Lemma~6.23 and its proof]{Meyer13}}]\label{lem:connection-tree-criterion}
At a fixed level, let $\mathcal G$ be the connection graph of either color and let $K_\varepsilon$ be the union of that color in a fixed geometric representation. Then $\mathcal G$ is a tree if and only if  $K_\varepsilon$ is a closed Jordan domain.
\end{lemma}

\begin{lemma}[{Meyer \cite[Proposition~8.1~\textup{(2)}]{Meyer13}}]\label{lem:white-connection-tree}
Fix the connection of $1$-tiles used to construct $H^0$. For every $m\in \N$, the white connection graph $\Gamma^+_m$ of the recursively defined connection of $m$-tiles is a tree.
\end{lemma}

\begin{lemma}
For each $m\in \N$, both $\Gamma_m^+$ and $\Gamma_m^-$ are trees.
\end{lemma}
\begin{proof}
    Fix $m\in\N$. In a geometric representation of the connection of $m$-tiles as in Subsection~\ref{subsec:connections}, let $K^+$ be the union of the closed white pieces, and let $K^-$ be the union of the closed black pieces. These two sets cover $S^2$, have disjoint interiors, and each is the closure of its interior. Hence, we have $K^-=S^2\smallsetminus\inter K^+$.

    By Lemma~\ref{lem:white-connection-tree}, $\Gamma_m^+$ is a tree. Lemma~\ref{lem:connection-tree-criterion} therefore implies that $K^+$ is a closed Jordan domain. The complementary closed Jordan domain bounded by $\partial K^+$ is $K^-$. Applying Lemma~\ref{lem:connection-tree-criterion} to $K^-$, we conclude that $\Gamma_m^-$ is a tree.
\end{proof}

\subsection{The white and black equivalence relations}\label{wbequivalencerelations}
We briefly recall the white and black equivalence relations introduced by Meyer in \cite[Section~5]{Meyer09}. We first describe the associated finite laminations and their gaps, and then recall how the white and black equivalence relations determine the equivalence relation $\sim$ induced by the invariant Peano curve $\gamma$.

Let $f$ be an expanding Thurston map, and let $\cC$ be an $f$-invariant Jordan curve with $\post f\subseteq \cC$. At each $n$-vertex $v\in \V^n$, let $\pi_\white^n(v)\cup\pi_\black^n(v)$ be the cnc-partition as above. Recall from Subsection~\ref{Pseudo_isotopies} that $\mathbf A^n=(\gamma^n)^{-1}(\V^n)$ is the set of $n$-angles.

Let $v\in \V^n$ and $b\in\pi_{\white}^n(v)\cup\pi_{\black}^n(v)$. We say that an $n$-edge $E$ with endpoint $v$ is \emph{incident to $b$ at $v$} if $E\subseteq X_i$ for an $n$-tile $X_i$ whose index $i$ at $v$ belongs to $b$. Let $\alpha\in \mathbf A^n$ with $\gamma^n(\alpha)=v$. Let $a_-^n(\alpha)$ and $a_+^n(\alpha)$ denote the two $n$-arcs adjacent to $\alpha$, where $a_-^n(\alpha)$ ends at $\alpha$ and $a_+^n(\alpha)$ starts at $\alpha$ with respect to the orientation of $S^1$. We say that $\alpha$ is \emph{incident to $b$ at $v$} if $\gamma^n(a_-^n(\alpha))$ and $\gamma^n(a_+^n(\alpha))$ are incident to $b$ at $v$.

Now we recall the definiton of equivalence relations $\overset{n,\white}{\sim}$, $\overset{n,\black}{\sim}$ on $\mathbf A^n$.

\begin{definition}[Equivalence relations $\overset{n,\white}{\sim},\,\overset{n,\black}{\sim}$]
For $n\in \N$ and $\alpha,\,\alpha'\in \mathbf A^n$, define
\begin{align*}
\alpha\overset{n,\white}\sim\alpha' &\Longleftrightarrow \text{ there exist }v\in \V^n \text{ and }b\in\pi_{\white}^n(v) \text{ such that both }\alpha,\,\alpha' \text{ are incident to }b\text{ at }v,\\
\alpha\overset{n,\black}\sim\alpha'&\Longleftrightarrow \text{ there exist }v\in \V^n \text{ and }b\in\pi_{\black}^n(v) \text{ such that both }\alpha,\,\alpha' \text{ are incident to }b\text{ at }v.
\end{align*}

If follows from the definition that there is a one-to-one correspondence between blocks $b\in\bigcup_{v\in \V^n} \pi_{\white}^n(v)$ and equivalence classes with respect to $\overset{n,\white}\sim$ (analogously, between blocks $b\in\pi_{\black}^n$ and equivalence classes with respect to $\overset{n,\black}\sim$). Note that $n$-angles $\alpha,\,\alpha'$ such that $\gamma^n(\alpha)\neq\gamma^n\bigl(\alpha'\bigr)$ are never equivalent with respect to $\overset{n,\white}\sim$ or $\overset{n,\black}\sim$.

In addition, the equivalence relations $\overset{0,\white}\sim$, $\overset{0,\black}\sim$ on $A^0$ defined to be trivial, meaning that each equivalence class is a singleton.
\end{definition}

\begin{lemma}[{Meyer \cite[Lemma~5.7~$(L^n2)$]{Meyer09}}]\label{lemma_wbnoncrossing}
    The equivalence relations $\overset{n,\white}{\sim}$ and $\overset{n,\black}{\sim}$ are \emph{non-crossing}, i.e., given two distinct equivalence classes $[s]_{n,\white},\,[t]_{n,\white}\subseteq \mathbb R/\mathbb Z$ of $\overset{n,\white}{\sim}$, there are no points $a,\,c\in[s]_{n,\white}$, $b,\,d\in[t]_{n,\white}$ such that $0\leq a<b<c<d<1.$
The same statement holds for $\overset{n,\black}{\sim}$.
\end{lemma}

We next relate succeeding angles in a class to the oriented boundaries of tiles and record how $\mu$ maps the classes.

\begin{lemma}[{Meyer \cite[Lemma~5.3, Remark~5.6]{Meyer09}}]\label{lem:finite-class-successors}
Let $n\in\N$ and distinct $\alpha,\,\beta\in\mathbf A^n$ satisfy $\gamma^n(\alpha)=\gamma^n(\beta)$. In the positive circle order, $\beta$ succeeds $\alpha$ in $[\alpha]_{n,\white}$ if and only if $\gamma^n(a_-^n(\beta))$ is followed by $\gamma^n(a_+^n(\alpha))$ on the positively oriented boundary of a white $n$-tile. In the reverse circle order, $\beta$ succeeds $\alpha$ in $[\alpha]_{n,\black}$ if and only if the reverse of $\gamma^n(a_+^n(\beta))$ is followed by the reverse of $\gamma^n(a_-^n(\alpha))$ on the positively oriented boundary of a black $n$-tile.
\end{lemma}

\begin{lemma}[{Meyer \cite[Lemma~5.7~\textup{($\mathcal L^n$4)}, Remark~5.6]{Meyer09}}]\label{lem:complete-angle-class-images}
For either color $\sigma\in\{\white,\,\black\}$, the complete angle classes satisfy
\begin{equation}\label{eq:finite-class-image}
 \mu([\alpha]_{n+1,\sigma})=[\mu(\alpha)]_{n,\sigma}\quad \text{ for all } n\in\N_0,\ \alpha\in\mathbf A^{n+1}.
\end{equation}
\end{lemma}

The equivalence relations $\overset{n,\white}{\sim}$ and $\overset{n,\black}{\sim}$ induce the white and black laminations, respectively, as in \cite[Subsection~5.3]{Meyer09}. View the parameter circle $S^1$ as the equator of a copy of $S^2$. Let $S_{\white}^2$ and $S_{\black}^2$ denote the two components of $S^2\smallsetminus S^1$, equipped with the hyperbolic metric, with ideal boundary $S^1$.

\begin{figure}
    \centering
    \begin{overpic}[permil,width=0.65\linewidth]{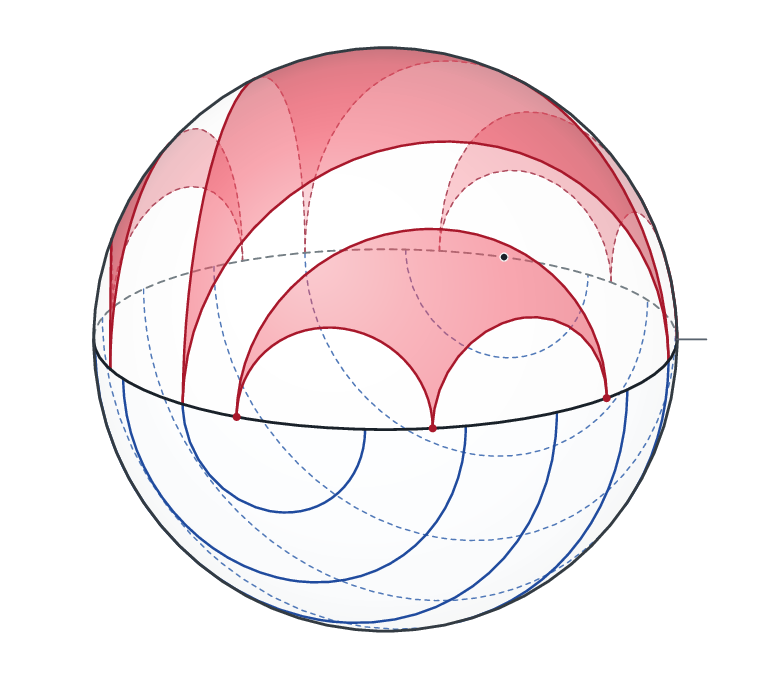}
        \normalsize
        \put(178,753){\makebox(0,0){$S_{\white}^2$}}
        \put(825,140){\makebox(0,0){$S_{\black}^2$}}
        \put(957,444){\makebox(0,0){$S^1$}}
        \put(656,522){\makebox(0,0){$0$}}
    \end{overpic}
    \caption{The parameter circle as the common boundary of two hemispheres, with convex hulls drawn on opposite sides. This picture shows a possible level-$2$ lamination of the map $R_4$ from \cite{Meyer09}.}
    \label{fig:hemispheres}
\end{figure}

The circle and the two hemispheres are illustrated in Figure~\ref{fig:hemispheres}.

\begin{definition}[Laminations $\mathcal{L}_{\white}^n,\,\mathcal{L}_{\black}^n$]
Let $\alpha\in\mathbf A^n$ be an $n$-angle and let $[\alpha]_{n,\white}\subseteq S^1$ be its equivalence class with respect to $\overset{n,\white}\sim$. We call the hyperbolically convex hull of $[\alpha]_{n,\white}$ in $S_{\white}^2$ a \emph{leaf}, denoted by $L=L([\alpha]_{n,\white})$. If $[\alpha]_{n,\white}$ corresponds to a block $b\in\pi_{\white}^n(v)$, we denote this leaf by $L(v,b)$. The set of all such leaves is the \emph{white lamination} $\mathcal{L}_{\white}^n$.

If $\card[\alpha]_{n,\white}=m$, $m\geq3$, the leaf $L=L([\alpha]_{n,\white})$ is an ideal $m$-gon. If $[\alpha]_{n,\white}$ consists of two points, the corresponding leaf $L$ is the hyperbolic geodesic (in $S_{\white}^2$) between these points. If $[\alpha]_{n,\white}$ consists of a single point the corresponding leaf is $L=\{\alpha\}$.

The \emph{black lamination} $\mathcal{L}_{\black}^n$ (containing $m$-gons in $S_{\black}^2$) is constructed similarly from the equivalence classes of $\overset{n,\black}\sim$.
\end{definition}

\begin{figure}
    \centering
    \begin{overpic}[permil,width=0.9\linewidth]{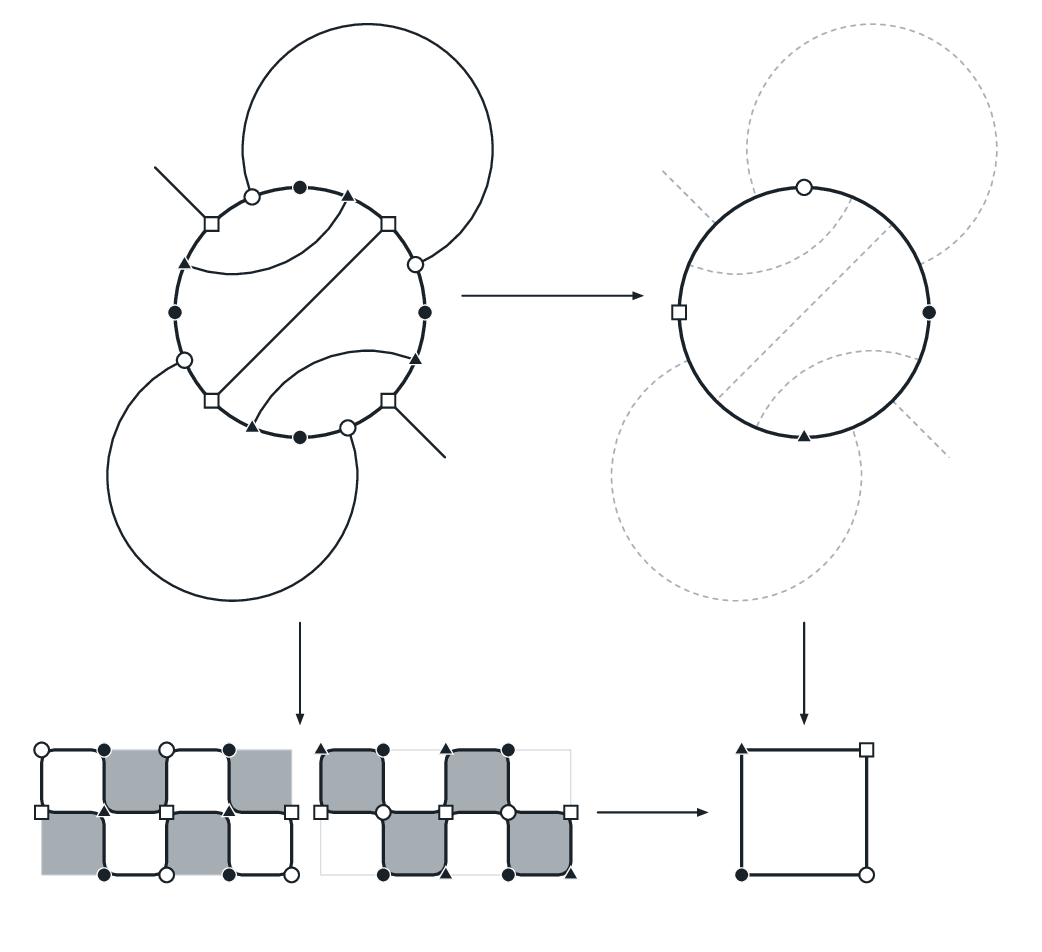}
        \scriptsize
        \put(430,614){\makebox(0,0){$0$}}
        \put(379,653){\makebox(0,0){$1$}}
        \put(388,715){\makebox(0,0){$2$}}
        \put(342,746){\makebox(0,0){$3$}}
        \put(287,757){\makebox(0,0){$4$}}
        \put(250,706){\makebox(0,0){$5$}}
        \put(218,685){\makebox(0,0){$6$}}
        \put(157,669){\makebox(0,0){$7$}}
        \put(145,614){\makebox(0,0){$8$}}
        \put(196,577){\makebox(0,0){$9$}}
        \put(188,514){\makebox(0,0){$10$}}
        \put(233,484){\makebox(0,0){$11$}}
        \put(288,473){\makebox(0,0){$12$}}
        \put(326,524){\makebox(0,0){$13$}}
        \put(358,545){\makebox(0,0){$14$}}
        \put(419,560){\makebox(0,0){$15$}}
        \put(914,614){\makebox(0,0){$0$}}
        \put(529,657){\makebox(0,0){\small $z^4$}}
        \put(327,271){\makebox(0,0){\small $\gamma^1$}}
        \put(811,271){\makebox(0,0){\small $\gamma^0$}}
        \put(625,161){\makebox(0,0){\small $g$}}
        \put(130,154){\makebox(0,0){$2,10$}}
        \put(70,116){\makebox(0,0){$11,15$}}
        \put(190,116){\makebox(0,0){$3,7$}}
        \put(398,154){\makebox(0,0){$6,14$}}
        \put(338,116){\makebox(0,0){$9,13$}}
        \put(458,116){\makebox(0,0){$1,5$}}
        \put(99,53){\makebox(0,0){$p_0=0$}}
        \put(487,217){\makebox(0,0){$p_0=0$}}
        \put(711,53){\makebox(0,0){$p_0=0$}}
        \put(160,21){\makebox(0,0){\footnotesize\sffamily White connections}}
        \put(428,21){\makebox(0,0){\footnotesize\sffamily Black connections}}
    \end{overpic}
    \caption{The first-level white and black laminations for the Latt\`es map $g$ in Example~\ref{ex:lattes-map}. In the upper-left panel, white leaves are drawn inside the parameter circle and black leaves outside it.}
    \label{fig:white-black-laminations}
\end{figure}

Figure~\ref{fig:white-black-laminations} illustrates the first-level white and black equivalence classes and their corresponding leaves for a Latt\`es map.

We next define the gaps of the finite laminations and describe their correspondence with the $n$-tiles.

\begin{definition}[{Meyer \cite[Definition~5.8]{Meyer09}}]\label{def:finite-gap-data}
For $n\in \N_0$ and $\sigma\in\{\white,\,\black\}$, an \emph{$n$-gap of color $\sigma$} is the closure in $S_\sigma^2\cup S^1$ of a component of $S_\sigma^2\smallsetminus\bigcup\mathcal L_\sigma^n$, where $\mathcal L_\sigma^n$ is the finite lamination of color $\sigma$.
\end{definition}

\begin{prop}[{Meyer \cite[Lemma~5.9~\textup{(G1)--(G3)}, Remark~5.6]{Meyer09}}]\label{prop:finite-gap-properties}
For $n\in \N$ and $\sigma\in\{\white,\,\black\}$, the $n$-tiles and $n$-gaps of the same color $\sigma$ correspond bijectively. Let $G_\sigma(X)$ denote the gap corresponding to $X\in\X_\sigma^n$. Then
\begin{equation}\label{eq:finite-gap-trace}
 G_\sigma(X)\cap S^1
 =\bigcup_{\substack{I\text{ an }n\text{-arc}\\\gamma^n(I)\subseteq\partial X}}I.
\end{equation}
In positive boundary order, the $0$-edges $f^n(\gamma^n(I))$ follow the positive order on $\cC$ for white and the reverse order for black. Each $t\in S^1\smallsetminus\mathbf A^n$ belongs to exactly one color-$\sigma$ $n$-gap. Set $G_\sigma\bigl(X_\sigma^0\bigr)\=S_\sigma^2\cup S^1$.
\end{prop}

\begin{lemma}[{Meyer \cite[proof of Lemmas~5.7 and~5.9,~\textup{($\mathcal L^n$3)}]{Meyer09}}]\label{lem:finite-gap-incidence}
Consider $n\in \N$, $\sigma\in\{\white,\,\black\}$, $X\in\X_\sigma^n$, $v\in\V^n$, and $b\in\pi_\sigma^n(v)$ with angle class $[\alpha]_{n,\sigma}$. Then $[\alpha]_{n,\sigma}\cap G_\sigma(X)\ne\emptyset$ if and only if $X$ is incident to $b$ at $v$.
\end{lemma}

\begin{lemma}[{Meyer \cite[Lemma~5.9~\textup{(G4)}, Remark 5.6]{Meyer09}}]\label{lem:finite-gap-dynamics}
For $n\in \N_0$, $\sigma\in\{\white,\,\black\}$ and $X\in\X_\sigma^{n+1}$,
\begin{equation}\label{eq:finite-gap-dynamics}
 \mu\bigl(G_\sigma(X)\cap S^1\bigr)=G_\sigma(f(X))\cap S^1.
\end{equation}
The restriction of $\mu$ to $G_\sigma(X)\cap\bigl(S^1\smallsetminus\mathbf A^{n+1}\bigr)$ is injective.
\end{lemma}

\begin{lemma}[{Meyer \cite[Lemma~5.9~\textup{(G5)}, Remark 5.6]{Meyer09}}]\label{lem:finite-gap-nesting}
For $n\in \N$ and either color $\sigma\in\{\white,\,\black\}$, each $(n+1)$-gap of color $\sigma$ is contained in a $n$-gap of color $\sigma$.
\end{lemma}

To describe the equivalence relation induced by $\gamma$, we first recall the partial order and the join of equivalence relations.

Let $\overset{a}{\sim}$ and $\overset{\black}{\sim}$ be equivalence relations on a set $S$. We write $\overset{a}{\sim}\leq\overset{b}{\sim}$ if every equivalence class of $\overset{a}{\sim}$ is contained in some equivalence class of $\overset{b}{\sim}$. Equivalently, viewing equivalence relations as subsets of $S\times S$, this means $\overset{a}{\sim}\ \subseteq\ \overset{\black}{\sim}$.

The set of all equivalence relations on $S$ forms a lattice with respect to this partial order. Recall that the join $\overset{a}{\sim}\vee\overset{\black}{\sim}$ is the smallest equivalence relation containing both $\overset{a}{\sim}$ and $\overset{\black}{\sim}$. More precisely, if $\overset{\vee}{\sim}\=\overset{a}{\sim}\vee\overset{\black}{\sim}$, then for $s,t\in S$,
\begin{align*}
    s\overset{\vee}{\sim} t
    \quad\Longleftrightarrow\quad
    &\text{ there exist }m\in\N_0
    \text{ and }s=x_0,\,x_1,\,\dots\,,x_m=t\in S\\
    &\text{ such that, for every }1\leq j\leq m,\quad
    x_{j-1}\overset{a}{\sim}x_j
    \text{ or }
    x_{j-1}\overset{\black}{\sim}x_j.
\end{align*}

An equivalence relation $\sim$ on a compact metric space $S$ is called \emph{closed} if $\{(s,t)\in S\times S:s\sim t\}$ is closed in $S\times S$. By \cite[Lemma~4.5]{Meyer09}, there exists a unique smallest closed equivalence relation containing $\sim$, which is called the \emph{closure} of $\sim$. See \cite[Section~2.2]{PM12} for more details. This closure of an equivalence relation need not coincide with the closure of its graph in $S\times S$, since the latter might not be transitive.

We now pass to the corresponding closed equivalence relations.
We regard $\overset{n,\white}{\sim}$ and $\overset{n,\black}{\sim}$ as equivalence relations on $S^1$ by declaring every point of $S^1\smallsetminus \mathbf A^n$ to be equivalent only to itself.
By \cite[Lemma~5.7~$(L^n6)$]{Meyer09}, for $n\in \N_0$ and $\alpha,\,\alpha'\in \mathbf A^n$,
\begin{align}\label{anglenn+1}
    \alpha\overset{n,\white}{\sim}\alpha'
    \quad\Longleftrightarrow\quad
    \alpha\overset{n+1,\white}{\sim}\alpha',
\end{align}
and the corresponding statement holds for $\overset{\black}{\sim}$.
Thus, $\overset{n,\white}{\sim}\subseteq\overset{n+1,\white}{\sim}$ and $\overset{n,\black}{\sim}\subseteq\overset{n+1,\black}{\sim}$ as equivalence relations on $S^1$.

\begin{definition}[Equivalence relations $\overset{\white}{\sim},\,\overset{\black}{\sim}$]\label{oversetwsim}
Define $\overset{\infty,\white}{\sim}\=\bigvee\overset{n,\white}{\sim}$, meaning for $s,\,t\in S^1$, $s\overset{\infty,\white}{\sim}t$ if and only if $s\overset{n,\white}{\sim}t$ for some $n\in\N_0$. Define equivalence relation $\overset{\white}{\sim}$ to be the closure of $\overset{\infty,\white}{\sim}$.

Similarly, define $\overset{\infty,\black}{\sim}\=\bigvee\overset{n,\black}{\sim}$ and let $\overset{\black}{\sim}$ be the closure of $\overset{\infty,\black}{\sim}$.
\end{definition}

The equivalence relation $\sim$ induced by the invariant Peano curve
$\gamma$ (i.e., $s\sim t$ if and only if $\gamma(s)=\gamma(t)$ for $s,\,t\in S^1$) may be recovered from $\overset{\white}{\sim},\overset{\black}{\sim}$. More precisely, Meyer proves the following lemma:

\begin{lemma}[{Meyer \cite[Lemma~5.14]{Meyer09}}] \label{closureequi}
Let $\overset{\white}{\sim},\,\overset{\black}{\sim}$ be defined as above, and $\sim$ be the equivalence relation induced by $\gamma$. Then 
\begin{enumerate}
    \smallskip
    \item $\sim$ is the closure of $\overset{\white}{\sim}\vee\overset{\black}{\sim}$;

    \smallskip
    \item if $f$ has no periodic critical points, then $\overset{\white}{\sim}\vee\overset{\black}{\sim} = \sim$.
\end{enumerate}
\end{lemma}

\begin{theorem}[{Meyer \cite[Theorem~6.1]{Meyer09}}] \label{finitefiber}
Let $\sim$ be the equivalence relation induced by $\gamma$.
\begin{enumerate}
    \smallskip
    \item If $f$ has no critical periodic cycles, there is a number $N<+\infty$ such that
    $\card[s]\leq N$ for all $s\in S^1$.

    \smallskip
    \item If $f$ has critical periodic cycles, there is (at least) one finite equivalence
    class $[s]$, where $\gamma(s)\notin \post f$.
\end{enumerate}
\end{theorem}

\subsection{Circular orders and order completion}\label{sec_order_completion}
In this section, we recall the circular order and order completion from \cite[Section~5]{Frankel13}. 

\begin{definition}\label{def_circularorder}
    A \emph{circularly ordered set} is defined to be a set $S$ with at least three elements together with a map $\langle \,\cdot\,, \,\cdot\,,\,\cdot\, \rangle$ from the set of triples of distinct elements of $S$ to $\pm 1$ such that for pairwise distinct $x,\,y,\,z \in S$, the following conditions are true:
\begin{enumerate}
\smallskip
\item $\langle x,y,z \rangle = {\rm sgn}(\tau)\langle \tau(x),\tau(y),\tau(z)\rangle$ for any permutation $\tau$ of $\{x,\,y,\,z\}$.
\smallskip
\item For each $w\in S\smallsetminus\{x,\,y,\,z\}$, if $\langle x,y,z \rangle = 1$ and $\langle x,z,w \rangle = 1$, then $\langle x,y,w \rangle = 1$ and $\langle y,z,w \rangle = 1$.
\end{enumerate}
\end{definition}

The \emph{order topology} on a circularly ordered set is defined to be the topology generated by the basis of \emph{open intervals} $(x,y)\=\{t\in S: \langle x,t,y \rangle=1\}$, where $x,\,y\in S$ are distinct. The \emph{closed interval} $[x,y]$ is defined as $[x,y]\=(x,y)\cup\{x,\,y\}$. For $x\in S$, we set $[x,x]\=\{x\}$.

A circularly ordered set is called \emph{order complete} if for any sequence $I_1\supseteq I_2\supseteq \cdots$ of nested closed intervals, we have $\bigcap_{n=1}^{+\infty} I_n \neq \emptyset$. A \emph{gap} in a circularly ordered set $S$ is defined to be an ordered pair $(x,y)$ of distinct elements $x,\,y\in S$ such that $(x,y)=\emptyset$.

\begin{prop}[{Frankel \cite[Proposition~5.4]{Frankel13}}]\label{completionofcircularorder}
Every circularly ordered set $S$ has a canonical \emph{order completion} $\overline{S}$ which is order complete. Moreover, $S$ is dense in $\overline{S}$ with respect to the order topology.
\end{prop} 

We give a brief overview of the construction of $\overline{S}$ following \cite[Construction~5.3]{Frankel13}, and refer the reader to \cite[Proposition~5.4]{Frankel13} for the rest of the proof.

We first define a set $S'$ of equivalence classes of admissible sequences in $S$ as follows. An \emph{admissible sequence} in $S$ is an infinite sequence of nested closed intervals $\{I_i\}$ such that $I_{i+1} \subseteq I_i$ for all $i\in\N$ and $\bigcap_{i=1}^{+\infty} I_i = \emptyset$. Two admissible sequences $\{I_i\}$ and $\{J_j\}$ are \emph{equivalent} if for all $n\in\N$, there exists $m\in\N$ such that $I_m \subseteq J_n$. Let $S'$ be the set of equivalence classes of admissible sequences in $S$. The order completion $\overline{S}$ is defined as $\overline{S} \= S \cup S'$. By~\cite[Construction~5.3]{Frankel13}, the circular order on $S$ extends canonically to $\overline{S}$. Thus, $\overline{S}$ is the order completion of $S$ in the proposition above.

\section{Roads and road equivalence relations}
\label{sec:roads}
In this section, we use the pseudo-isotopies from Subsections~\ref{Pseudo_isotopies} to show that connection graphs of higher levels can be obtained by gluing together connection graphs of lower levels. We then construct the road equivalence relations, which are useful for the construction of the P-zipper.

\subsection{Levels of connection graphs}

Let $f$ be a degree $d$ expanding Thurston map and  $\cC$ be an invariant Jordan curve satisfying Lemma~\ref{7.5}. Fix a visual metric $\rho$ for $f$ with expansion factor $\Lambda>1$.
Fix a connection of $1$-tiles satisfying \textup{(C1)}--\textup{(C2)} of Definition~\ref{Propertiesofconnections}.
Given the connection, let $H^m$ be its pseudo-isotopies, normalized approximations $\gamma^m, m \in \N_0$, and Peano limit $\gamma\:S^1\to S^2$ from Subsections~\ref{Pseudo_isotopies}. For each $m \in \N_0$ and each $t\in[0,1]$, write $H_t^m\=H^m(\,\cdot\,,t)$. Define $\mu\:\R/\Z\to\R/\Z$ by $\mu(t)\=d\,t\pmod1$. 

\begin{definition}\label{def:ch4-reference-edges}
For $v\in\bigcup_{j=0}^{+\infty}\V^j$, its \emph{birth level} is defined by $m_0\=\min\bigl\{j\in\N_0:v\in\V^j\bigr\}$. As in Subsection~\ref{Pseudo_isotopies}, each $m$-edge inherits the positive boundary orientation of the unique white $m$-tile containing it. Define the \emph{reference edge} $E_v^{m_0}$ as follows. If $m_0=0$, let $E_v^0$ be the unique $0$-edge with initial endpoint $v$. If $m_0=1$, choose a $1$-edge $E_v^1$ with initial endpoint $v$ and keep this choice fixed. If $m_0\ge2$, let $E_v^{m_0}$ be the unique $m_0$-edge with initial endpoint $v$ such that $f(E_v^{m_0})=E_{f(v)}^{m_0-1}$. Make these choices in increasing order of $m_0$.

We now define $E_v^m$ by recursion. For each integer $m\ge m_0$, let $J$ be the initial subdivision arc of $E_v^m$ at $v$ in Lemma~\ref{lem:lifted-edge-subdivision}, and define $E_v^{m+1}\= H_1^m(J)$. For each integer $m\ge\max\{1,\, m_0\}$, define the \emph{tile indices} $X_i^m(v)$ by listing the $m$-tiles incident to $v$ in positive cyclic order, starting from the white tile containing $E_v^m$, and choosing the $i$-th one to be $X_i^m(v)$.
\end{definition}

We first construct copies of lower-level connection graphs inside higher-level connection graphs. For $n,r\in \N$, Proposition~\ref{lem:ch4-sector-copies} shows that the level-$(n+r)$ graph is covered by copies of the level-$r$ graph, one for each level-$n$ tile of the same color. The copies for two distinct level-$n$ tiles $X$ and $X'$ of the same color meet exactly when the tiles share a block. In that case, they meet at a single vertex, which is the continuation of the common block to level $(n+r)$.

\begin{definition}\label{def:ch4-signed-connection-data}
For each $p\in\post f$ and $\sigma\in\{+,\,-\}$, we denote by $\widehat p_0^\sigma$ a copy of $p$. Set $a_p\=\bigl(\gamma^0\bigr)^{-1}(p)$. Define $\sB_0^\sigma\=\{\widehat p_0^\sigma:p\in\post f\}$ and $\cA_0^\sigma(\widehat p_0^\sigma)\=\{a_p\}$. 

Recall the definition of $\pi_\sigma^m$ and $m$-angles from Subsection~\ref{subsec:connections} and Subsection~\ref{Pseudo_isotopies}.
For each integer $m\ge1$ and $\sigma\in\{+,\,-\}$, define $\sB_m^\sigma\=\{c(v,b):v\in\V^m,\ b\in\pi_\sigma^m(v)\}$. For $u=c(v,b)\in\sB_m^\sigma$, let $\cA_m^\sigma(u)$ be the set of $m$-angles at $v$ whose two adjacent edges are contained in tiles with labels belonging to $b$.

For each integer $m\ge0$ and $\sigma\in\{+,\,-\}$, define $\prV\:\sB_m^\sigma\to\V^m$ by $\prV(c(v,b))\=v$ when $m\ge1$ and $\prV(\widehat p_0^\sigma)\=p$ when $m=0$. For each integer $m\ge1$, $X\in\X_\sigma^m$, and $u\in\sB_m^\sigma$, write $X\leftrightarrow u$ if $c(X)$ and $u$ are adjacent in $\Gamma_m^\sigma$.
\end{definition}

\begin{remark}
    Note that distinct blocks of a partition are disjoint. We observe that for each $m\in \N$, if $X\in\X_\sigma^m$, $u,u'\in\sB_m^\sigma$, and $0\le j\le m$, then \begin{equation}\label{eq:ch4-incident-block-uniqueness}
        \bigl\{X\leftrightarrow u,\,X\leftrightarrow u',\,f^j(\prV(u))=f^j(\prV(u'))\bigr\}\quad\Longrightarrow \quad u=u'.
\end{equation}
\end{remark}

\begin{definition}\label{def:ch4-block-continuation}
For each integer $m\ge0$ and $\sigma\in\{+,\,-\}$, define the \emph{continuation map} $\operatorname{Lift}_m^\sigma\:\sB_m^\sigma\to\sB_{m+1}^\sigma$ as follows.
For each $m\in\N$ and $c(v,b)\in\sB_m^\sigma$, define $\operatorname{Lift}_m^\sigma(c(v,b))\=c\bigl(v,b_\star^{m+1}(b)\bigr)$, where $b_\star^{m+1}(b)$ is given by~\eqref{eq:recursive-block-merges}. For each $p\in\post f$, define $\operatorname{Lift}_0^\sigma(\widehat p_0^\sigma)\=c\bigl(p,b_\star^1(p)\bigr)$, where $b_\star^1(p)$ is the marked white $1$-block at $p$.

For each $p\in\post f$, each $\sigma\in\{+,\,-\}$, and each integer $m\ge0$, define $\widehat p_{m+1}^\sigma\=\operatorname{Lift}_m^\sigma(\widehat p_m^\sigma)$ and write $\widehat p^\sigma\=\{\widehat p_m^\sigma\}_{m=0}^{+\infty}$.
\end{definition}

\begin{lemma}
\label{lem:ch4-connection-inputs}
Fix a color $\sigma\in\{+,\,-\}$. 
Use the tile indices of Definition~\ref{def:ch4-reference-edges}. 
For each $m\in \N_0$, the map $\operatorname{Lift}_m^\sigma$ is well-defined and injective, and
\begin{equation}\label{eq:ch4-block-persistence}
 \prV(\operatorname{Lift}_m^\sigma u)=\prV(u),\quad
 \cA_m^\sigma(u)\subseteq
 \cA_{m+1}^\sigma(\operatorname{Lift}_m^\sigma u),
 \quad u\in\sB_m^\sigma.
\end{equation}
Moreover, for $m\in \N$ and $c(v, b)\in\sB_m^\sigma$, put $p\=f^m(v)$ and $k\=\deg_f(p)$. Let $ b_\star^{1,+}\in\pi_{\white}^1(p)$ contain $0$ and $ b_\star^{1,-}\in\pi_{\black}^1(p)$ contain $2k-1$. Set $\varepsilon\=0$ if $\sigma=+$ and $\varepsilon\=1$ if $\sigma=-$. Then
\begin{equation}\label{eq:ch4-block-continuation}
 \operatorname{Lift}_m^\sigma(c(v,b))
 =c\Bigl(v,\bigcup_{2j+\varepsilon\in b}\bigl\{i+2kj:i\in b_\star^{1,\sigma}\bigr\}\Bigr).
\end{equation}
\end{lemma}
\begin{proof}
Fix $m\ge1$ and $v\in\V^m$. Put $p\=f^m(v)$, $k\=\deg_f(p)$, and $D\=\deg_{f^m}(v)$. By Lemma~\ref{lem:cell-pullback-counting}, $f^m$ maps $E_v^m$ homeomorphically onto a $0$-edge. The map $f^m$ preserves the white tile boundary orientations, so the image edge starts at $p$ and equals $E_p^0$. Since $f^m\circ H_1^m=H_1^0\circ f^m$ and $\V^{m+1}=f^{-m}\bigl(\V^1\bigr)$, the initial subdivision arcs of $E_v^m$ and $E_p^0$ correspond under $f^m$. Thus $f^m\bigl(E_v^{m+1}\bigr)=E_p^1$. With these reference edges, the local model at $v$ (see~\cite[Lemma~5.29]{BM}) gives $f^m\bigl(X_{2kj+i}^{m+1}(v)\bigr)=X_i^1(p)$ for integers $0\le j<D$ and $0\le i<2k$.

Put $\ell\=\max b_\star^{1,+}$. Let $J_j$ and $J'_j$ be the subdivision arcs containing $v$ on the outgoing and incoming edges of $X_{2j}^m(v)$, respectively. We have $H_1^m(J_j)\subseteq\partial X_{2kj}^{m+1}(v)$ and $H_1^m(J'_j)\subseteq\partial X_{2kj+\ell}^{m+1}(v)$ (see \cite[proof of Proposition~8.1~\textup{(3)}]{Meyer13}). The two subdivision arcs at $v$ on $\partial X_{2j+1}^m(v)$ are $J'_j$ and $J_{j+1}$, where $(j+1)$ is taken modulo $D$. Their images end and start at $v$, respectively. In positive cyclic order around $v$, the black tile containing the first image follows its white tile, and the black tile containing the second image precedes its white tile. The black indices are therefore $2kj+\ell+1$ and $2k(j+1)-1$ modulo $2kD$. Here $\ell+1=\min b_\star^{1,-}$ (see \cite[proof of Proposition~8.1, step~(1e)]{Meyer13}).

For $u=c(v,b)\in\sB_m^\sigma$, let $b'$ be the union in~\eqref{eq:ch4-block-continuation}. By~\eqref{eq:recursive-block-merges}, we have $b'\in\pi_\sigma^{m+1}(v)$. Fix $a\in\cA_m^\sigma(u)$. For either $m$-edge $E$ of $\gamma^m$ at $a$, write $E\subseteq\partial X_{2j+\varepsilon}^m(v)$ with $2j+\varepsilon\in b$. Let $J$ be the subdivision arc of $E$ containing $v$. The corresponding $(m+1)$-edge at $a$ is $H_1^m(J)$ (see \cite[Section~4.2]{Meyer13}), and $H_1^m(J)\subseteq\partial X_{2kj+i}^{m+1}(v)$ for some $i\in b_\star^{1,\sigma}$. Since $2kj+i\in b'$ for both edges, we obtain $\cA_m^\sigma(u)\subseteq\cA_{m+1}^\sigma(c(v,b'))$. For $\sigma=-$, we have $\operatorname{Lift}_m^-(u)=c(v,b')$ by definition. For $\sigma=+$, the formula follows from~\eqref{eq:recursive-block-merges}, since a shift by $2q$ in the $m$-tile indices shifts the $(m+1)$-tile indices by $2kq$.

For $p\in\post f$, the outgoing $1$-edge at $a_p$ is $E_p^1$ by Lemma~\ref{lem:initial-marked-edge}. The white (resp. black) tile containing $E_p^1$ has index $0$ (resp. $2\deg_f(p)-1$). Hence, $\widehat p_1^\sigma=c\bigl(p,b_\star^{1,\sigma}\bigr)$ and $a_p\in\cA_1^\sigma(\widehat p_1^\sigma)$.

By the restriction identity~\eqref{anglenn+1}, we have $\cA_{m+1}^\sigma(\operatorname{Lift}_m^\sigma(u))\cap\mathbf A^m=\cA_m^\sigma(u)$ for all $m\in\N_0$ and $u\in\sB_m^\sigma$. Thus $\operatorname{Lift}_m^\sigma(u)$ determines $u$. For $a\in\cA_m^\sigma(u)$, we have $\prV(\operatorname{Lift}_m^\sigma(u))=\gamma^{m+1}(a)=\gamma^m(a)=\prV(u)$. Iterating the inclusion in~\eqref{eq:ch4-block-persistence}, we obtain $a_p\in\cA_m^\sigma(\widehat p_m^\sigma)$ for every $m\in\N_0$.
\end{proof}

\begin{definition}\label{def:ch4-endpoint-composites}
Let $T$ be a connected graph with no cycles. A finite path is called \emph{simple} if it has no repeated vertices. For vertices $x,y$ of $T$, let $[x,y]_T$ denote the unique simple path joining $x$ to $y$, with $[x,x]_T=\{x\}$. The \emph{median} of three vertices is the unique vertex common to their three pairwise paths.

For a topological $\mathbb R$-tree $T$, let $[x,y]_T$ denote the unique arc joining distinct points $x,\,y\in T$, and set $[x,x]_T\=\{x\}$. The \emph{median} of three points is the unique point common to their three pairwise joining arcs.
\end{definition}

\begin{lemma}\label{lem:ch4-sector-ancestors}
For integers $0\le m\le n$, define $g_{n,m}\=H_1^n\circ\cdots\circ H_1^m\:S^2\to S^2$. For each integer $m\ge0$, define $g_{m-1,m}\=\mathrm{id}_{S^2}$. For integers $0\le i\le j$, the following hold:
\begin{enumerate}
\smallskip
\item For $V\in\X^j$, the set
$Q\=g_{j-1,i}^{-1}(\inte_\circ V)$ is nonempty and connected, and the restriction $g_{j-1,i}|_Q\:Q\to\inte_\circ V$ is a homeomorphism. There is a unique tile in $\X^i$ denoted by $A_{i,j}(V)$ of the same color as $V$ such that $Q\subseteq\inte_\circ A_{i,j}(V)$.
\smallskip
\item For all $0\le i\le j\le k$ and $V\in\X^k$, we have $A_{i,k}(V)=A_{i,j}(A_{j,k}(V))$ and $V\subseteq g_{k-1,i}(A_{i,k}(V))$.
\smallskip
\item They commute with $f$: for $V\in\X^{j+1}$,
$(A_{i,j}\circ f)(V)=(f\circ A_{i+1,j+1})(V)$. For $W'\in\X^{i+1}$, $T\mapsto f(T)$ is a bijection from $\bigl\{T\in\X^{j+1}:A_{i+1,j+1}(T)=W'\bigr\}$ onto $\bigl\{U\in\X^j:A_{i,j}(U)=f(W')\bigr\}$.
\end{enumerate}
\end{lemma}
\begin{proof}
\textup{(i)} For $m\in \N_0$ and $U\in\X^{m+1}$, Lemma~\ref{lem:lifted-tile-preimages} states that $(H_1^m)^{-1}(\inte_\circ U)$ maps homeomorphically onto $\inte_\circ U$ and lies in the interior of a single $m$-tile of the same color. Compose the inverse homeomorphisms for $m=\{j-1,\,j,\,\ldots,\,i\}$. Their image is the nonempty connected set $Q$ contained in one $i$-tile of the same color, and their composite is the inverse of $g_{j-1,i}|_Q$. Since distinct $i$-tiles have disjoint interiors, this tile is unique. For $i=j$, take $Q=\inte_\circ V$ and $A_{i,i}(V)=V$.

\textup{(ii)} By part~\textup{(i)} and the definition of $g_{n,m}$, we have $\emptyset\ne g_{k-1,i}^{-1}(\inte_\circ V)\subseteq\inte_\circ A_{i,k}(V)$ and $g_{k-1,i}^{-1}(\inte_\circ V)\subseteq\inte_\circ A_{i,j}(A_{j,k}(V))$. Since distinct $i$-tiles have disjoint interiors, $A_{i,k}(V)=A_{i,j}(A_{j,k}(V))$. The compact image $g_{k-1,i}(A_{i,k}(V))$ contains $\inte_\circ V$ and hence its closure $V$.

\textup{(iii)} Fix $W'\in\X^{i+1}$ and $U\in\X^j$, and set $Q_U\=g_{j-1,i}^{-1}(\inte_\circ U)$. For each tile $T$ with $f(T)=U$, set $Q_T\=g_{j,i+1}^{-1}(\inte_\circ T)$. Since $f\circ g_{j,i+1}=g_{j-1,i}\circ f$, equation~\eqref{eq:tile-interior-pullback} gives
\begin{equation*}
    f^{-1}(Q_U)=g_{j,i+1}^{-1}\bigl(f^{-1}(\inte_\circ U)\bigr)=\bigsqcup_{\substack{T\in\X^{j+1}\\f(T)=U}}Q_T.
\end{equation*}

For each such tile $T$, part~\textup{(i)} gives $Q_T\subseteq\inte_\circ A_{i+1,j+1}(T)$. Distinct $(i+1)$-tiles have disjoint interiors. Intersecting with $\inte_\circ W'$ therefore gives
\begin{equation*}
    \inte_\circ W'\cap f^{-1}(Q_U)
=\bigsqcup_{\substack{f(T)=U\\A_{i+1,j+1}(T)=W'}}Q_T.
\end{equation*}
Since $f|_{W'}$ is a homeomorphism, $\inte_\circ W'\cap f^{-1}(Q_U)$ is nonempty and connected when $A_{i,j}(U)=f(W')$, and empty otherwise. Each $Q_T$ is nonempty and open, so $\inte_\circ W'\cap f^{-1}(Q_U)=Q_T$ for a unique $T$ when $A_{i,j}(U)=f(W')$. Taking $W'=A_{i+1,j+1}(V)$ and $U=f(V)$, we have $\emptyset\ne Q_V\subseteq\inte_\circ W'\cap f^{-1}(Q_U)$, so $(A_{i,j}\circ f)(V)=(f\circ A_{i+1,j+1})(V)$.
\end{proof}

For integers $n,\,r\ge1$ and each tile $X\in\X_\sigma^n$, apply the bijection $T\mapsto f(T)$ established above with $i=k-1$, $j=k+r-1$, and $W'=f^{n-k}(X)$ for $k=n,\,n-1,\,\dots,\,1$. Since $f^n(X)=X_\sigma^0$, for every $U\in\X_\sigma^r$ there is a unique $V\in\X_\sigma^{n+r}$ satisfying $f^n(V)=U$ and $A_{n,n+r}(V)=X$.

\begin{figure}[htbp]
    \centering
    \begin{overpic}[
        percent,
        width=0.9\linewidth
    ]{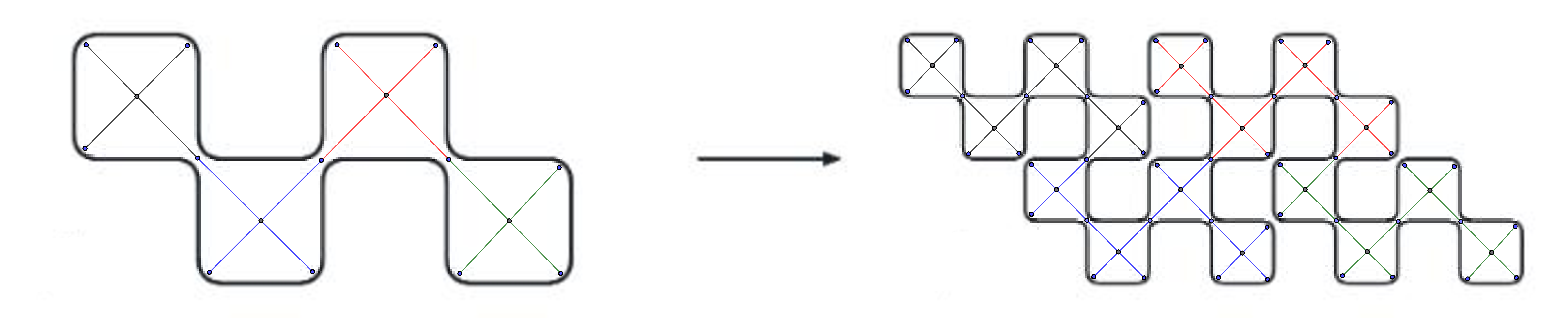}
        \setlength{\unitlength}{0.009\linewidth}
        \tiny

        \put(10.81,14.85){\makebox(0,0){\scalebox{1}{$c(X_1)$}}}
        \put(18.54,6.61){\makebox(0,0){\scalebox{1}{$c(X_2)$}}}
        \put(26.34,14.92){\makebox(0,0){\scalebox{1}{$c(X_3)$}}}
        \put(34.34,6.80){\makebox(0,0){\scalebox{1}{$c(X_4)$}}}

        \put(61.51,16.35){\makebox(0,0){\scalebox{1}{$c(U_4)$}}}
        \put(65.35,12.58){\makebox(0,0){\scalebox{1}{$c(U_3)$}}}
        \put(69.05,16.48){\makebox(0,0){\scalebox{1}{$c(U_2)$}}}
        \put(73.02,12.64){\makebox(0,0){\scalebox{1}{$c(U_1)$}}}

        \normalsize
        \put(49.17,13.85){\makebox(0,0){$H^1$}}
        \put(3.64,3.43){\makebox(0,0){$\gamma^1$}}
        \put(60.42,3.43){\makebox(0,0){$\gamma^2$}}

    \end{overpic}
    \caption{In this figure, we consider the expanding Thurston map in Example~\ref{ex:lattes-map} with connection of $1$-tiles shown in the left figure, and $n=r=1$. The $1$-tiles in the left figure are denoted by $X_1,\,X_2,\,X_3,\,X_4$. Then $L(X_1,1,X_i)=U_i$ as well as $L(X_1,1,c(X_i))=c(U_i)$ as shown in the right figure. Moreover, the edges in $\Gamma^\sigma_1[X_i]$ and $\Gamma^\sigma_2[X_i]$ have the same color.}
    \label{fig:copy}
\end{figure}

\begin{definition}
\label{def:ch4-sector-data}
\label{def:ch4-inverse-angle-map}
We call the tile $A_{i,j}(V)$ given by Lemma~\ref{lem:ch4-sector-ancestors} the \emph{sector ancestor} of a $j$-tile $V$ at level $i$. For integers $n,r\ge1$ and $X\in\X_\sigma^n$, let $\Gamma_{n+r}^\sigma[X]$ be the following subgraph of $\Gamma_{n+r}^\sigma$. Its vertices are all the tile vertices $c(V)$ satisfying $A_{n,n+r}(V)=X$ and all the block vertices adjacent to these tile vertices. Its edges are all the edges of $\Gamma_{n+r}^\sigma$ incident to these tile vertices. Define a map $b_X^\sigma\:\post f\to\sB_n^\sigma$ as follows. For $p\in \post f$, the homeomorphism $f^n|_X\:X\to X_\sigma^0$ determines a unique vertex $v=(f^n|_X)^{-1}(p)$. Write $X\=X_i^n(v)$ for an integer $0\leq i\leq 2\deg(f^n,v)$ and let $b$ be the unique block of the partition $\pi_\sigma^n(v)$ containing $i$. Set $b_X^\sigma(p)\=c(v,b)$. For $U\in\X_\sigma^r$, we denote by $L(X,r,U)$ the unique tile in $\X_\sigma^{n+r}$ (given by Lemma~\ref{lem:ch4-sector-ancestors}~\textup{(iii)}) satisfying 
\begin{equation*}
 A_{n,n+r}(L(X,r,U))=X,\quad f^n(L(X,r,U))=U,
\end{equation*}
and we call the map $L(X,r,\,\cdot\,)\:\X_\sigma^r\to\X_\sigma^{n+r}$ the \emph{tile map}. See Figure~\ref{fig:copy} for an example.

Write $G^+(T)\= G_\white(T)$ and $G^-(T)\= G_\black(T)$ for the gaps associated with tiles $T\in\X_\sigma^j$ of color $\sigma$ by the correspondence in Proposition~\ref{prop:finite-gap-properties}. The set $G^\sigma(T)\cap S^1$ is the union of the $j$-arcs $I$ satisfying $\gamma^j(I)\subseteq\partial T$.

For each integer $n\ge1$ and each tile $X\in\X_\sigma^n$, define the \emph{inverse angle map} $\ell_X^\sigma\:S^1\smallsetminus \mathbf A^0\to  G^\sigma(X)\cap\bigl(S^1\smallsetminus \mathbf A^n\bigr)$ as the inverse of the restriction of $\mu^n$ to $G^\sigma(X)\cap\bigl(S^1\smallsetminus \mathbf A^n\bigr)$, where $\mu\:S^1\to S^1$ satisfies $\mu(z)\=z^d$. 
\end{definition}

To extend $L(X,r,\,\cdot\,)$ from tiles to block vertices, we pull back the angle class of a level-$r$ block under $\mu^n$ and restrict the preimage to $G^\sigma(X)$. The next lemma identifies the equivalence class containing these preimages.

\begin{lemma}\label{lem:ch4-class-transport}
Fix a color $\sigma\in\{+,\,-\}$, integers $n,r\ge1$, a tile $X\in\X_\sigma^n$ and a block vertex $b\in\sB_r^\sigma$. For the gap $G^\sigma(X)$ and incidence map $b_X^\sigma$ of Definition~\ref{def:ch4-sector-data}, there exists a unique class $E\subseteq\mathbf A^{n+r}$ satisfying
\begin{equation}\label{eq:ch4-class-transport}
\mu^n(E)=\cA_r^{\sigma}(b)\quad\text{ and }\quad
E\cap G^{\sigma}(X)=G^{\sigma}(X)\cap\mu^{-n}(\cA_r^{\sigma}(b)).
\end{equation}
In addition, the following holds:
\begin{enumerate}
\smallskip
\item If $\cA_r^{\sigma}(b)\cap \mathbf A^0=\emptyset$, then $E\subseteq G^{\sigma}(X)$.
\smallskip
\item For $p\in \post f$ and $b\=\widehat p_r^{\sigma}$, we have $E=\cA_{n+r}^{\sigma}((\operatorname{Lift}_{n+r-1}^{\sigma}\circ \cdots\circ \operatorname{Lift}_n^{\sigma} \circ b_X^{\sigma})(p))$.
\end{enumerate}
\end{lemma}
\begin{proof}
Fix an integer $j\ge1$ and $T\in\X_\sigma^j$. For every $u\in\sB_j^\sigma$, Proposition~\ref{prop:finite-gap-properties} and Lemma~\ref{lem:finite-gap-incidence} give
\begin{equation}\label{eq:ch4-gap-incidence}
 \gamma^j\bigl(G^{\sigma}(T)\cap S^1\bigr)=\partial T \quad\text{ and }\quad
 T\leftrightarrow u\iff\cA_j^{\sigma}(u)\cap G^{\sigma}(T)\ne\emptyset.
\end{equation}
Let $T'\in\X_\sigma^{j+1}$ be a tile with ancestor $T$. It follows from Lemma~\ref{lem:lifted-tile-preimages} that each boundary edge of $T'$ pulls back to a subdivision arc $J\subseteq\partial T$. Let $I$ be the parameter arc for $H_1^j(J)$. By the construction of $\gamma^{j+1}$ in \cite[proof of Lemma~5.11~(1)]{Meyer09}, $I$ lies in the $j$-arc for the edge containing $J$, so $I\cap\bigl(G^\sigma(T)\smallsetminus\mathbf A^j\bigr)\ne\emptyset$ by~\eqref{eq:finite-gap-trace}. By Lemma~\ref{lem:finite-gap-nesting}, $G^\sigma(T')\subseteq G$ for some level-$j$ gap $G$ of color $\sigma$. Since $I\subseteq G^\sigma(T')$, we have $(G\cap G^\sigma(T))\smallsetminus\mathbf A^j\ne\emptyset$. Hence, $G=G^\sigma(T)$ by Proposition~\ref{prop:finite-gap-properties}, and $G^\sigma(T')\subseteq G^\sigma(T)$.

For $U\in\X_\sigma^r$, put $V_U^\sigma\=L(X,r,U)$. Then $A_{n,n+r}(V_U^\sigma)=X$ and $f^n(V_U^\sigma)=U$. Iterating the preceding inclusion $r$ times and~\eqref{eq:finite-gap-dynamics} $n$ times gives
\begin{equation}\label{eq:ch4-gap-transport}
 G^\sigma(V_U^\sigma)\subseteq G^\sigma(X) \quad\text{ and }\quad
 \mu^n\bigl(G^\sigma(V_U^\sigma)\cap S^1\bigr)=G^\sigma(U)\cap S^1.
\end{equation}
By Lemma~\ref{lem:cell-pullback-counting}, the homeomorphism $f^n|_X\:X\to X_\sigma^0$ maps the boundary edges of $X$ bijectively onto the $0$-edges. By Proposition~\ref{prop:finite-gap-properties} and Lemma~\ref{lem:finite-gap-dynamics}, $\mu^n$ restricts to a bijection from $G^\sigma(X)\cap\bigl(S^1\smallsetminus\mathbf A^n\bigr)$ onto $S^1\smallsetminus\mathbf A^0$. Its inverse is the map $\ell_X^\sigma$ prescribed in Definition~\ref{def:ch4-inverse-angle-map}. The continuous image $\mu^n\bigl(G^\sigma(X)\cap S^1\bigr)$ is compact and hence closed in $S^1$. Thus, since $\mu^n\bigl(G^\sigma(X)\cap S^1\bigr)$ contains $S^1\smallsetminus \mathbf A^0$, which is dense in $S^1$, we have $\mu^n\bigl(G^\sigma(X)\cap S^1\bigr)=S^1$.

For a finite angle class $A$ of color $\sigma$, let $L(A)$ be its closed hyperbolic convex hull in $S_{\white}^2\cup S^1$ if $\sigma=+$ and in $S_{\black}^2\cup S^1$ if $\sigma=-$ (see Subsection~\ref{wbequivalencerelations}). Put $C_b^\sigma\=\cA_r^\sigma(b)$. The restriction of $C_b^\sigma$ to $\mathbf A^0$ is $\{a_p\}$ if $b=\widehat p_r^\sigma$ and is empty otherwise. Since $\mu^n\bigl(G^\sigma(X)\cap S^1\bigr)=S^1$, the set $G^\sigma(X)\cap\mu^{-n}(C_b^\sigma)$ is nonempty. For each $a$ in this set, let $E_a^\sigma$ denote its complete level-$(n+r)$ class. Iterating~\eqref{eq:finite-class-image} in Lemma~\ref{lem:complete-angle-class-images} gives $\mu^n(E_a^\sigma)=C_b^\sigma$.

\begin{case1}
    Unmarked classes. Suppose $C_b^\sigma\cap \mathbf A^0=\emptyset$. Since $\mu^n(\mathbf A^n)=\mathbf A^0$, $E_a^\sigma\cap \mathbf A^n=\emptyset$. By~\eqref{eq:ch4-block-persistence}, each class $C\subseteq\mathbf A^n$ lies in a class $C'\subseteq\mathbf A^{n+r}$ with $C'\ne E_a^\sigma$. By non-crossing, $L(C')\cap L(E_a^\sigma)=\emptyset$. Since $L(C)\subseteq L(C')$, the hull $L(E_a^\sigma)$ is disjoint from every level-$n$ hull. If $\card E_a^\sigma\ge2$, the set $L(E_a^\sigma)\smallsetminus E_a^\sigma$ is connected by convexity. It is disjoint from $L(C)$ for every level-$n$ angle class $C$ of color $\sigma$, so it lies in a single component $D$ of the complement of these hulls. Each point of $E_a^\sigma$ is an ideal endpoint of a geodesic in $D$, so $E_a^\sigma\subseteq\overline D$. Both $\overline D$ and $G^\sigma(X)$ contain $a\notin\mathbf A^n$, so $\overline D=G^\sigma(X)$ by Proposition~\ref{prop:finite-gap-properties}. Thus, $E_a^\sigma\subseteq G^\sigma(X)$, which also holds if $E_a^\sigma=\{a\}$. Applying the inverse $\ell_X^\sigma$ to $\mu^n(E_a^\sigma)=C_b^\sigma$, we obtain $E_a^\sigma=\ell_X^\sigma(C_b^\sigma)$. The class $E\=\ell_X^\sigma(C_b^\sigma)$ is independent of $a$ and satisfies~\eqref{eq:ch4-class-transport} and part~\textup{(i)}.
\end{case1}

\begin{case2}
    Marked classes. Suppose $C_b^\sigma\cap\mathbf A^0=\{a_p\}$. Since $\mu^n(E_a^\sigma)=C_b^\sigma$ and $\mathbf A^n=\mu^{-n}\bigl(\mathbf A^0\bigr)$, we have $E_a^\sigma\cap\mathbf A^n\ne\emptyset$. By the restriction identity~\eqref{anglenn+1}, $E_a^\sigma\cap\mathbf A^n=\cA_n^\sigma\bigl(\widetilde b\bigr)$ for a unique $\widetilde b\in\sB_n^\sigma$. Hence,~\eqref{eq:ch4-block-persistence} gives $E_a^\sigma=\cA_{n+r}^\sigma\bigl((\operatorname{Lift}_{n+r-1}^\sigma\circ \cdots\circ \operatorname{Lift}_n^\sigma)\widetilde b\bigr)$.

Put $C\=\cA_n^\sigma\bigl(\widetilde b\bigr)$. By~\eqref{eq:ch4-gap-incidence}, to prove $X\leftrightarrow\widetilde b$, it suffices to show $C\cap G^\sigma(X)\ne\emptyset$. If $a\in\mathbf A^n$, then $a\in C\cap G^\sigma(X)$. Otherwise choose $c\in C$ and let $J\subseteq L(E_a^\sigma)$ be the closed hyperbolic geodesic from $a$ to $c$. Each class $B\subseteq\mathbf A^n$ with $B\ne C$ lies in a class $B'\subseteq\mathbf A^{n+r}$ by~\eqref{eq:ch4-block-persistence}. Since $E_a^\sigma\cap\mathbf A^n=C$, we have $B'\ne E_a^\sigma$. By non-crossing, $J\cap L(B)=\emptyset$. Let $z$ be the first point where $J$ meets $L(C)$, starting from $a$. The part of $J$ strictly between $a$ and $z$ lies in a single component $D$ of the complement of the level-$n$ hulls of color $\sigma$. Both $\overline D$ and $G^\sigma(X)$ contain $a\notin\mathbf A^n$, so, by Proposition~\ref{prop:finite-gap-properties}, $\overline D=G^\sigma(X)$. If $z\in S^1$, then $z\in C\cap G^\sigma(X)$. Otherwise, let $e$ be the geodesic side of $L(C)$ containing $z$. No other level-$n$ hull meets $e$, so $\partial D\cap\inte_\circ e$ is nonempty, open and closed in the connected arc $\inte_\circ e$. Hence, $e\subseteq\overline D=G^\sigma(X)$, and both endpoints of $e$ belong to $C\cap G^\sigma(X)$.

For $c\in\cA_n^\sigma\bigl(\widetilde b\bigr)$, we have $\mu^n(c)\in C_b^\sigma\cap\mathbf A^0=\{a_p\}$. Hence, $(f^n\circ\prV)\bigl(\widetilde b\bigr)=\bigl(\gamma^0\circ\mu^n\bigr)(c)=p$. Since $X\leftrightarrow\widetilde b$,~\eqref{eq:ch4-incident-block-uniqueness} gives $\widetilde b=b_X^\sigma(p)$. Thus, every $a\in G^\sigma(X)\cap\mu^{-n}(C_b^\sigma)$ has the same class $E_b^\sigma=\cA_{n+r}^\sigma((\operatorname{Lift}_{n+r-1}^\sigma\circ \cdots\circ \operatorname{Lift}_n^\sigma)(b_X^\sigma(p)))$. Consequently $G^\sigma(X)\cap\mu^{-n}(C_b^\sigma)\subseteq E_b^\sigma$, while $\mu^n(E_b^\sigma)=C_b^\sigma$ implies $E_b^\sigma\subseteq\mu^{-n}(C_b^\sigma)$. These inclusions prove~\eqref{eq:ch4-class-transport} and part~\textup{(ii)}.
\end{case2}

Since each complete equivalence class $E\subseteq\mathbf A^{n+r}$ satisfying~\eqref{eq:ch4-class-transport} contains $a$, we have $E=E_a^\sigma$.
\end{proof}

\begin{definition}\label{def:ch4-copy-maps}
For integers $n,r\ge1$ and $X\in\X_\sigma^n$, define $L(X,r,\,\cdot\,)$ on the vertices of $\Gamma_r^\sigma$ as follows. On tile vertices, set $L(X,r,c(U))\=c(L(X,r,U))$ for $U\in\X_\sigma^r$. On block vertices, for each $b\in\sB^\sigma_r$, let $L(X,r,b)$ be the element in $\sB_{n+r}^\sigma$ satisfying that $\cA_{n+r}^\sigma(L(X,r,b))$ is the unique complete class given by Lemma~\ref{lem:ch4-class-transport} for $X$ and $b$.

Set $L\bigl(X_\sigma^0,r,\,\cdot\,\bigr)\=\mathrm{id}_{\Gamma_r^\sigma}(\,\cdot\,)$. At a tile vertex $c(U)$, order the edges joining $c(U)$ to block vertices $c(v,b)$ by the positive cyclic order of the vertices $v$ along $\partial_\circ U$. At a block vertex $c(v,b)$, order the edges joining $c(v,b)$ to $c(X_i^r(v))$ by the cyclic order of the indices $i\in b$.
\end{definition}

We now check that the maps on tile and block vertices together define a graph isomorphism onto $\Gamma_{n+r}^\sigma[X]$.

\begin{prop}\label{lem:ch4-sector-copies}
Let $f\:S^2\to S^2$ be the expanding Thurston map with an $f$-invariant Jordan curve $\cC\supseteq\post f$ and a connection of $1$-tiles satisfying \textup{(C1)}--\textup{(C2)} of Definition~\ref{Propertiesofconnections}. Fix $\sigma\in\{+,\,-\}$, $n,r\ge1$, and $X\in\X_\sigma^n$. The map $L(X,r,\,\cdot\,)$ of Definition~\ref{def:ch4-copy-maps} have the following properties:
\begin{enumerate}
\smallskip
\item The vertex map extends to an isomorphism $L(X,r,\,\cdot\,)\:\Gamma_r^\sigma\to\Gamma_{n+r}^\sigma[X]$. For a subgraph $Q\subseteq\Gamma_r^\sigma$, write $L(X,r,Q)$ for its image under $L(X,r,\,\cdot\,)$. The union of the subgraphs $L(X,r,\Gamma_r^\sigma)$ over all $X\in\X_\sigma^n$ is $\Gamma_{n+r}^\sigma$. For distinct $X,X'\in\X_\sigma^n$,
\begin{equation*}
 L(X,r,\Gamma_r^{\sigma})\cap L(X',r,\Gamma_r^{\sigma})
 =\{(\operatorname{Lift}_{n+r-1}^{\sigma}\circ \cdots\circ\operatorname{Lift}_n^{\sigma})(b):
 b\in\sB_n^{\sigma},\ X\leftrightarrow b\leftrightarrow X'\}.
\end{equation*}

\item For each integer $s\ge1$ and $U\in\X_\sigma^r$,
\begin{equation*}
\begin{aligned}
L(X,r+1,\,\cdot\,)\circ\operatorname{Lift}_r^{\sigma}
&=\operatorname{Lift}_{n+r}^{\sigma}\circ L(X,r,\,\cdot\,)
&&\quad\text{ on }\sB_r^{\sigma},\\
L(X,r+s,\,\cdot\,)\circ L(U,s,\,\cdot\,)&=L(L(X,r,U),s,\,\cdot\,)
&&\quad \text{ on }\Gamma_s^{\sigma},\\
f(L(X,r,U))&=L(f(X),r,U).
\end{aligned}
\end{equation*}

\item The map $L(X,r,\,\cdot\,)$ preserves the cyclic order of incidences. For $p\in\post f$, put $v\=(f^n|_X)^{-1}(p)$ and $D_r\=\deg_{f^r}(p)$. Use the tile indices of Definition~\ref{def:ch4-reference-edges} and write $X\=X_{2j+\varepsilon}^n(v)$, where $\varepsilon=0$ for $\sigma=+$ and $\varepsilon=1$ for $\sigma=-$. For every index $i$ such that $X_i^r(p)\leftrightarrow\widehat p_r^\sigma$, we have
\begin{equation*}
L(X,r,X_i^r(p))=X_{2D_rj+i}^{n+r}(v).
\end{equation*}
\end{enumerate}
\end{prop}

\begin{proof}
\textup{(i)} For $b\in\sB_r^\sigma$ and $a\in\cA_{n+r}^\sigma(L(X,r,b))$, Lemma~\ref{lem:ch4-class-transport} and $f^n(\gamma^{n+r}(a))=\gamma^r( \mu^n(a))$ give the following identities:
\begin{equation}\label{eq:ch4-copy-values}
 \begin{aligned}
 f^n(\prV(L(X,r,b)))&=\prV(b)	&&\quad \text{ for all } b\in\sB_r^{\sigma},\\
 L(X,r,\widehat p_r^{\sigma})&=(\operatorname{Lift}_{n+r-1}^{\sigma}\circ \cdots\circ \operatorname{Lift}_n^{\sigma})(b_X^{\sigma}(p))&&\quad \text{ for all } p\in\post f.
 \end{aligned}
\end{equation}
Put $V_U^\sigma\=L(X,r,U)$, $C_b^\sigma\=\cA_r^\sigma(b)$ and $E_b^\sigma\=\cA_{n+r}^\sigma(L(X,r,b))$. By Lemma~\ref{lem:ch4-sector-ancestors}~\textup{(iii)}, the map $U\mapsto L(X,r,U)$ is a bijection from $\X_\sigma^r$ onto $\{V\in\X_\sigma^{n+r}:A_{n,n+r}(V)=X\}$, with inverse $V\mapsto f^n(V)$. Equations~\eqref{eq:ch4-gap-transport},~\eqref{eq:ch4-class-transport} and~\eqref{eq:ch4-gap-incidence} give
\begin{equation*}
 \begin{aligned}
 U\leftrightarrow b
 &\iff C_b^\sigma\cap G^\sigma(U)\ne\emptyset
 \iff\mu^{-n}(C_b^\sigma)\cap G^\sigma(V_U^\sigma)\ne\emptyset\\
 &\iff E_b^\sigma\cap G^\sigma(V_U^\sigma)\ne\emptyset
 \iff V_U^\sigma\leftrightarrow L(X,r,b).
 \end{aligned}
\end{equation*}
For a block vertex $\widetilde b\leftrightarrow V_U^\sigma$ in $\Gamma_{n+r}^\sigma[X]$, choose $a\in\cA_{n+r}^\sigma\bigl(\widetilde b\bigr)\cap G^\sigma(V_U^\sigma)$, and let $b$ denote the class $\mu^n(\cA_{n+r}^\sigma\bigl(\widetilde b\bigr))$. Since $a\in G^\sigma(X)$ and $\mu^n(a)\in C_b^\sigma$,~\eqref{eq:ch4-class-transport} implies $a\in E_b^\sigma$ and hence $\widetilde b=L(X,r,b)$. If $L(X,r,b)=L(X,r,b')$, then $C_b^\sigma=\mu^n(E_b^\sigma)=\mu^n(E_{b'}^\sigma)=C_{b'}^\sigma$, so $b=b'$. Thus, $L(X,r,\,\cdot\,)$ is a bijection from $\sB^\sigma_r$ to $\sB^\sigma_{n+r}\cap\Gamma_{n+r}^\sigma[X]$ and preserves incidence in both directions, so it is an isomorphism.

Each tile vertex $c(V)$ and all its edges belong to the copy indexed by $A_{n,n+r}(V)$. Every block vertex has a tile incidence, so the copies cover $\Gamma_{n+r}^\sigma$.

Each tile has a unique level-$n$ ancestor, so the copies with $X\ne X'$ can meet only at block vertices. Suppose $L(X,r,b)=L(X',r,b')$. By~\eqref{eq:ch4-class-transport}, $C_b^\sigma=\mu^n(E_b^\sigma)=C_{b'}^\sigma$, hence $b=b'$. If $b$ were unmarked, we would have $E_b^\sigma\cap\mathbf A^n=\emptyset$, whereas Lemma~\ref{lem:ch4-class-transport}~\textup{(i)} and Proposition~\ref{prop:finite-gap-properties} imply $E_b^\sigma\subseteq G^\sigma(X)\cap G^\sigma(X')\cap S^1\subseteq\mathbf A^n$. Since $E_b^\sigma\ne\emptyset$, this is impossible; thus $b=\widehat p_r^\sigma$ for some $p\in\post f$. Each $\operatorname{Lift}_j^\sigma$ with $n\le j<n+r$ is injective by Lemma~\ref{lem:ch4-connection-inputs}. Cancelling these maps in~\eqref{eq:ch4-copy-values} for $X$ and $X'$ gives $b_X^\sigma(p)=b_{X'}^\sigma(p)$.

Conversely, suppose $X\leftrightarrow\widetilde b\leftrightarrow X'$ and set $p\=f^n\bigl(\prV\bigl(\widetilde b\bigr)\bigr)$. Thus, by definition, $\widetilde b=b_X^\sigma(p)=b_{X'}^\sigma(p)$. By~\eqref{eq:ch4-copy-values}, $(\operatorname{Lift}_{n+r-1}^\sigma\circ\cdots\circ\operatorname{Lift}_n^\sigma)\bigl(\widetilde b\bigr)$ equals both $L(X,r,\widehat p_r^\sigma)$ and $L(X',r,\widehat p_r^\sigma)$, and hence belongs to $L(X,r,\Gamma_r^\sigma)\cap L(X',r,\Gamma_r^\sigma)$.

\smallskip
\textup{(ii)} Fix $b\in\sB_r^\sigma$ and put $B\=\operatorname{Lift}_{n+r}^\sigma(L(X,r,b))$. Since $\mu^n\bigl(G^\sigma(X)\cap S^1\bigr)=S^1$, there exists $a\in G^\sigma(X)\cap S^1$ with $\mu^n(a)\in C_b^\sigma$. By~\eqref{eq:ch4-class-transport} and~\eqref{eq:ch4-block-persistence}, we have $a\in E_b^\sigma\subseteq\cA_{n+r+1}^\sigma(B)$. Since $\mu^n(a)\in C_b^\sigma\subseteq\cA_{r+1}^\sigma(\operatorname{Lift}_r^\sigma b)$, applying~\eqref{eq:ch4-class-transport} at level $r+1$ gives $B=L(X,r+1,\operatorname{Lift}_r^\sigma b)$.

Fix $A\in\X_\sigma^s$ and put $T\=L(L(X,r,U),s,A)$. By Lemma~\ref{lem:ch4-sector-ancestors}~ \textup{(ii)--(iii)} and the definition of $L$, we have $A_{n,n+r+s}(T)=X$, $A_{r,r+s}(f^n(T))=U$, and $f^{n+r}(T)=A$. Thus, $f^n(T)=L(U,s,A)$ and $T=L(X,r+s,L(U,s,A))$. For each block $b'\leftrightarrow A$, put $B_1\=L(X,r+s,L(U,s,b'))$ and $B_2=L(L(X,r,U),s,b')$. For $i=1,\,2$, part~\textup{(i)} and~\eqref{eq:ch4-copy-values} imply $T\leftrightarrow B_i$ and $f^{n+r}(\prV(B_i))=\prV(b')$. Therefore, equation~\eqref{eq:ch4-incident-block-uniqueness} implies $B_1=B_2$. Since each edge is determined by its endpoints, these two compositions also send every edge to the same edge.

Put $V\=L(X,r,U)$. For $n\ge2$, Lemma~\ref{lem:ch4-sector-ancestors}~\textup{(iii)} gives $A_{n-1,n+r-1}(f(V))=f(X)$. Since $f^{n-1}(f(V))=U$, the definition of $L$ gives $f(V)=L(f(X),r,U)$. For $n=1$, we have $f(X)=X_\sigma^0$ and $f(V)=U=L\bigl(X_\sigma^0,r,U\bigr)=L(f(X),r,U)$.

\smallskip
\textup{(iii)} For $V=L(X,r,U)$, the homeomorphism $f^n|_V$ preserves the positive boundary order, hence the incidence order at the tile vertices. Read the $n$-arcs in $G^\sigma(X)\cap S^1$ in parameter-circle order. This agrees with the positive boundary order of the gap for white and reverses it for black. By Proposition~\ref{prop:finite-gap-properties}, $\mu^n$ maps these arcs to the $0$-arcs in order and preserves orientation on each arc. The terminal image of each arc is the initial image of the next. Thus $\mu^n$ preserves the cyclic order of points in $G^\sigma(X)\cap S^1$ with distinct images.

Fix $b\in\sB_r^\sigma$. For each $U\leftrightarrow b$, let $a_U$ be the angle over $\prV(b)$ on the parameter arc for the edge leaving $\prV(b)$ along the positively oriented boundary of $U$. Let $\widetilde a_U$ be the corresponding angle for $L(X,r,U)$ at $L(X,r,b)$. By~\eqref{eq:ch4-gap-transport}, $\widetilde a_U\in G^\sigma(X)$, and its arc maps to an arc of $G^\sigma(U)$. The correspondence $f^n\circ\gamma^{n+r}=\gamma^r\circ\mu^n$ on the endpoints identifies that arc with the one for $a_U$, so $\mu^n(\widetilde a_U)=a_U$. Distinct $U$ have distinct $a_U$. By Lemma~\ref{lem:finite-class-successors}, reading these angles positively for white and negatively for black gives the incidence order. The same rule applies to $\widetilde a_U$, so $L(X,r,\,\cdot\,)$ preserves the incidence order at $b$.

For $p\in\post f$, put $v\=(f^n|_X)^{-1}(p)$. Write $\widehat p_r^\sigma\=c\bigl(p,b_{p,r}^\sigma\bigr)$, $b_X^\sigma(p)\=c(v,b)$, and $X=X_{2j+\varepsilon}^n(v)$, with $\varepsilon=0$ for white and $\varepsilon=1$ for black. Put $D_r\=\deg_{f^r}(p)$ and $k\=\deg_f(p)$. For each of the two subdivision arc $J$ at $v$ on $\partial X$, let $T$ be the color-$\sigma$ tile incident to $H_1^n(J)$. Lemma~\ref{lem:lifted-tile-preimages} puts $J$ in $A_{n,n+1}(T)$. There is a unique color-$\sigma$ tile along the $n$-edge containing $J$, so $A_{n,n+1}(T)=X$. The boundary calculation in the proof of Lemma~\ref{lem:ch4-connection-inputs} identifies the two tiles $T$ as $X_{2kj+i}^{n+1}(v)$ for $i=\min\bigl(b_{p,1}^\sigma\bigr),\max\bigl(b_{p,1}^\sigma\bigr)$. Since $f^n(T)=X_i^1(p)$, these are $L\bigl(X,1,X_i^1(p)\bigr)$.

In what follows, we prove by induction on integers $r\ge1$ that
\begin{equation*}
 L(X,r,X_i^r(p))=X_{2D_rj+i}^{n+r}(v)	\quad \text{ for all } i\in b_{p,r}^\sigma.
\end{equation*}
For $r=1$, put $i_-\=\min\bigl(b_{p,1}^\sigma\bigr)$, $i_+\=\max\bigl(b_{p,1}^\sigma\bigr)$, and $T_i\=L\bigl(X,1,X_i^1(p)\bigr)$ for $i\in b_{p,1}^\sigma$. The preceding paragraph gives $T_{i_-}=X_{2kj+i_-}^{n+1}(v)$ and $T_{i_+}=X_{2kj+i_+}^{n+1}(v)$. When $\card b_{p,1}^\sigma\le2$, every index in $b_{p,1}^\sigma$ equals $i_-$ or $i_+$. Since $L(X,1,\,\cdot\,)$ preserves cyclic order, $T_i$ is between $T_{i_-}$ and $T_{i_+}$ in positive order around $v$ whenever $i_-<i<i_+$. Thus, every $T_i$ has the form $X_{2kj+s}^{n+1}(v)$ with $0\le s<2k$. Applying $f^n$ gives $X_i^1(p)=X_s^1(p)$, so $s=i$.

Make the induction hypothesis that $L(X,r,X_i^r(p))=X_{2D_rj+i}^{n+r}(v)$ for all $i\in b_{p,r}^\sigma$, and consider the case for $r+1$. For the induction step, put $q\=f^r(p)$ and $k'\=\deg_f(q)$, so $D_{r+1}=D_r\times k'$. Take $i=2\ell+\varepsilon\in b_{p,r}^\sigma$ and $h\in b_{q,1}^\sigma$. The case $r=1$ gives $L\bigl(X_i^r(p),1,X_h^1(q)\bigr)=X_{2k'\ell+h}^{r+1}(p)$. Using part~\textup{(ii)} and the induction hypothesis, we then apply the case $r=1$ at $X_{2D_rj+i}^{n+r}(v)$:
\begin{equation*}
L\bigl(X,r+1,X_{2k'\ell+h}^{r+1}(p)\bigr)
=L\bigl(X_{2D_rj+i}^{n+r}(v),1,X_h^1(q)\bigr)
=X_{2D_{r+1}j+2k'\ell+h}^{n+r+1}(v).
\end{equation*}
By~\eqref{eq:ch4-block-continuation}, every index in $b_{p,r+1}^\sigma$ has the form $2k'\ell+h$ with $2\ell+\varepsilon\in b_{p,r}^\sigma$ and $h\in b_{q,1}^\sigma$. Hence, the formula holds at level $r+1$.
\end{proof}

\subsection{Road equivalence relations}
We next use the copies of $\Gamma_1^\sigma$ to refine paths between marked blocks and thereby subdivide the road intervals. At each level, we identify subdivision nodes with the same block label. We take the closed equivalence relation generated by these identifications and show that the initial quotient $Y_0^\sigma$ is a finite tree.

\begin{definition}\label{def:ch4-initial-roads}
For $n\in \N_0$, define the \emph{road labels} by
\begin{equation*}
 \sA_0^{\sigma}\= \{\{\widehat p^{\sigma},\widehat q^{\sigma}\}:p,q\in \post f,\ p\ne q\},\quad
 \sA_n^{\sigma}\= \sA_0^{\sigma}\sqcup
 \bigsqcup_{k=1}^n\bigl\{(k,X,\alpha):X\in\X_\sigma^k,\ \alpha\in\sA_0^{\sigma}\bigr\}.
\end{equation*}
For each road label $\lambda\in\sA_n^\sigma$, define $I_\lambda^\sigma\=\{\lambda\}\times[0,1]$. When the road is fixed and no ambiguity arises, write $s$ for the point $(\lambda,s)$ and $[s,t]$ for the subinterval $\{\lambda\}\times[s,t]$. Define the \emph{birth level} of a road label by $\operatorname{birth}(\alpha)\=0$ for $\alpha\in\sA_0^\sigma$ and $\operatorname{birth}(k,X,\alpha)\=k$. For each $\alpha=\{\widehat p^\sigma,\,\widehat q^\sigma\}\in\sA_0^\sigma$, fix an order $(p,q)$ and label the endpoints $(\alpha,0)$ and $(\alpha,1)$ by $p$ and $q$, respectively.

For each integer $m\ge1$, set $P_m^\sigma(\alpha)\=[\widehat p_m^\sigma,\widehat q_m^\sigma]_{\Gamma_m^\sigma}$, with block vertices $B_{m,i}^{\alpha,\sigma}$ and tile vertices $c\bigl(X_{m,i}^{\alpha,\sigma}\bigr)$ in order. Write $k_m^{\alpha,\sigma}$ for the number of tile vertices of $P_m^\sigma(\alpha)$.

At level $1$, divide $I_\alpha^\sigma$ into $k_1^{\alpha,\sigma}$ equal subintervals and match them to the tile vertices of $P_1^\sigma(\alpha)$ in order. To subdivide an interval at level $m$, consider the block vertices immediately before and after its tile vertex on $P_m^\sigma(\alpha)$. Let $p_-$ and $p_+$ be their respective images under $f^m\circ\prV$. By~\eqref{eq:ch4-incident-block-uniqueness}, these images are distinct. Divide this interval into as many equal parts as there are tile vertices on $\bigl[\widehat{p_-}_1^\sigma,\widehat{p_+}_1^\sigma\bigr]_{\Gamma_1^\sigma}$. By Lemma~\ref{prop:ch4-road-data}~\textup{(i)}, the resulting intervals match the tile vertices of $P_{m+1}^\sigma(\alpha)$ in order. Write $t_{m,i}^{\alpha,\sigma}$ for the division points and $S_m^{\alpha,\sigma}$ for their set. Define the \emph{atoms}, \emph{node labels}, and \emph{tile labels} by
\begin{equation*}
 \begin{gathered}
 \mathcal P_m^{\alpha,\sigma}
 =\bigl\{\bigl[t_{m,i-1}^{\alpha,\sigma},t_{m,i}^{\alpha,\sigma}\bigr]:1\le i\le k_m^{\alpha,\sigma}\bigr\},\\
 b_m^{\alpha,\sigma}\bigl(t_{m,i}^{\alpha,\sigma}\bigr)=B_{m,i}^{\alpha,\sigma}\in\sB_m^\sigma,\quad\text{ and }\quad
 X_m^{\alpha,\sigma}\bigl(\bigl[t_{m,i-1}^{\alpha,\sigma},t_{m,i}^{\alpha,\sigma}\bigr]\bigr)=X_{m,i}^{\alpha,\sigma}.
 \end{gathered}
\end{equation*}

For a road label $\lambda=(k,X,\alpha)$, transfer the subdivision points and intervals from $I_\alpha^\sigma$ to $I_\lambda^\sigma$ via the map $(\alpha,s)\mapsto(\lambda,s)$. In these coordinates, for each integer $r\ge1$, set
\begin{equation}\label{eq:road-data}
 \begin{gathered}
 \mathcal P_{k+r}^{\lambda,\sigma}=\mathcal P_r^{\alpha,\sigma},\quad
 S_{k+r}^{\lambda,\sigma}=S_r^{\alpha,\sigma},\quad
 b_{k+r}^{\lambda,\sigma}=L(X,r,\,\cdot\,)\circ b_r^{\alpha,\sigma},\\
 X_{k+r}^{\lambda,\sigma}(J)=L(X,r,X_r^{\alpha,\sigma}(J)),\quad
 P_{k+r}^{\sigma}(\lambda)=L(X,r,P_r^{\sigma}(\alpha)).
 \end{gathered}
\end{equation}
For each $m>\operatorname{birth}(\lambda)$, we call the closed intervals in $\mathcal P_m^{\lambda,\sigma}$ \emph{level-$m$ atoms} and the points of $S_m^{\lambda,\sigma}$ \emph{level-$m$ nodes}. Their labels are $b_m^{\lambda,\sigma}\:S_m^{\lambda,\sigma}\to\sB_m^\sigma$ and $X_m^{\lambda,\sigma}\:\mathcal P_m^{\lambda,\sigma}\to\X_\sigma^m$.
\end{definition}

\begin{definition}[Atom coordinates and finite realizations]\label{def:ch4-finite-realizations}
For each initial road $I_\alpha^\sigma$ with fixed endpoint order $(p,q)$, fix a parametrization $a_\alpha^\sigma\:I_\alpha^\sigma\to X_\sigma^0$ of an arc whose interior lies in $\inte_\circ \bigl(X_\sigma^0\bigr)$, with $a_\alpha^\sigma(\alpha,0)=p$ and $a_\alpha^\sigma(\alpha,1)=q$.

For $m>\operatorname{birth}(\lambda)$ and an atom $J=[a_-,a_+]\in\mathcal P_m^{\lambda,\sigma}$, put $p_\pm\=f^m\bigl(\prV\bigl(b_m^{\lambda,\sigma}(a_\pm)\bigr)\bigr)$. By~\eqref{eq:ch4-incident-block-uniqueness}, we have $p_-\ne p_+$. Define $\alpha_J\=\bigl\{\widehat{p_-}^\sigma,\,\widehat{p_+}^\sigma\bigr\}\in\sA_0^\sigma$. Let $\kappa_J^\sigma\:J\to I_{\alpha_J}^\sigma$ be the affine map sending each $a_\pm$ to the endpoint labeled by $p_\pm$. On $J$, prescribe the finite realization by
\begin{equation}\label{eq:ch4-compatible-realization}
 F_m^{\lambda,\sigma}\big|_J
 =\bigl(f^m\big|_{X_m^{\lambda,\sigma}(J)}\bigr)^{-1}\circ a_{\alpha_J}^\sigma\circ\kappa_J^\sigma.
\end{equation}
\end{definition}

\begin{figure}[htbp]
    \centering
    \begin{overpic}[
        percent,
        width=0.9\linewidth,
        trim=0bp -16bp 0bp 0bp
    ]{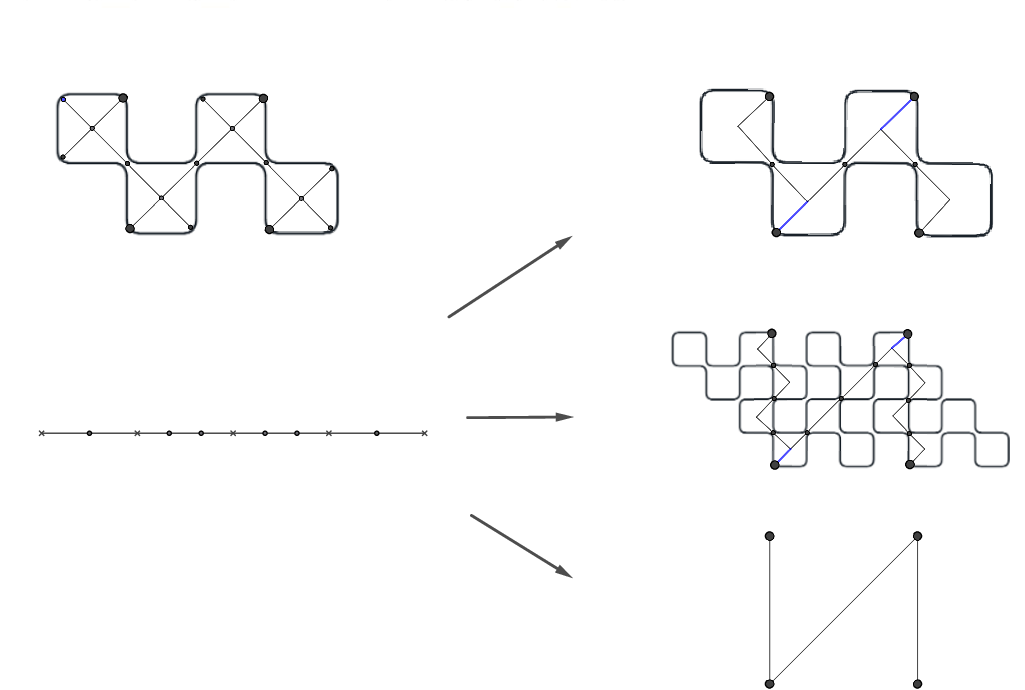}
        \setlength{\unitlength}{0.009\linewidth}
        \scriptsize

        \put(13.1,63.0){\makebox(0,0){$c(-1,b_1)$}}
        \put(13.4,56.4){\makebox(0,0){$c(v_1,b_1)$}}
        \put(19,52.8){\makebox(0,0){$c(v_2,b_2)$}}
        \put(27,56.4){\makebox(0,0){$c(v_3,b_3)$}}
        \put(27,45.9){\makebox(0,0){$c(1,b_4)$}}

        \put(77.1,60.7){\makebox(0,0){$-1$}}
        \put(91.1,59.8){\makebox(0,0){$\infty$}}
        \put(76.7,45.9){\makebox(0,0){$0$}}
        \put(90.4,45.9){\makebox(0,0){$1$}}

        \put(77.3,53.8){\makebox(0,0){$v_1$}}
        \put(84.0,53.5){\makebox(0,0){$v_2$}}
        \put(91.7,53.5){\makebox(0,0){$v_3$}}

        \put(76.8,37.7){\makebox(0,0){$-1$}}
        \put(90.2,37.4){\makebox(0,0){$\infty$}}
        \put(77.2,23.4){\makebox(0,0){$0$}}
        \put(89.9,23.4){\makebox(0,0){$1$}}

        \put(77.1,34.0){\makebox(0,0){$u_1$}}
        \put(77.1,27.5){\makebox(0,0){$u_3$}}
        \put(79.7,27.5){\makebox(0,0){$u_4$}}
        \put(86.5,34.0){\makebox(0,0){$u_6$}}
        \put(90.1,34.0){\makebox(0,0){$u_7$}}
        \put(90.8,27.5){\makebox(0,0){$u_9$}}

        \put(77.0,17.2){\makebox(0,0){$-1$}}
        \put(91.2,16.8){\makebox(0,0){$\infty$}}
        \put(76.9,2.3){\makebox(0,0){$0$}}
        \put(92.0,2.3){\makebox(0,0){$1$}}

        \put(49.4,45.0){\makebox(0,0){$F_1^{\alpha,+}$}}
        \put(50.9,31.7){\makebox(0,0){$F_2^{\alpha,+}$}}
        \put(52.8,19.0){\makebox(0,0){$F^{\alpha,+}$}}

        \put(22.8,21.7){\makebox(0,0){$I_\alpha^+$}}

        \tiny
        \put(4.1,30.2){\makebox(0,0){\scalebox{0.75}{$t_{1,0}^{\alpha,+}$}}}
        \put(13.4,30.2){\makebox(0,0){\scalebox{0.75}{$t_{1,1}^{\alpha,+}$}}}
        \put(22.8,30.2){\makebox(0,0){\scalebox{0.75}{$t_{1,2}^{\alpha,+}$}}}
        \put(32.0,30.2){\makebox(0,0){\scalebox{0.75}{$t_{1,3}^{\alpha,+}$}}}
        \put(41.4,30.2){\makebox(0,0){\scalebox{0.75}{$t_{1,4}^{\alpha,+}$}}}

        \put(4.1,26.3){\makebox(0,0){\scalebox{0.75}{$t_{2,0}^{\alpha,+}$}}}
        \put(8.6,26.3){\makebox(0,0){\scalebox{0.75}{$t_{2,1}^{\alpha,+}$}}}
        \put(13.4,26.3){\makebox(0,0){\scalebox{0.75}{$t_{2,2}^{\alpha,+}$}}}
        \put(16.5,26.3){\makebox(0,0){\scalebox{0.75}{$t_{2,3}^{\alpha,+}$}}}
        \put(19.6,26.3){\makebox(0,0){\scalebox{0.75}{$t_{2,4}^{\alpha,+}$}}}
        \put(22.8,26.3){\makebox(0,0){\scalebox{0.75}{$t_{2,5}^{\alpha,+}$}}}
        \put(25.8,26.3){\makebox(0,0){\scalebox{0.75}{$t_{2,6}^{\alpha,+}$}}}
        \put(29.0,26.3){\makebox(0,0){\scalebox{0.75}{$t_{2,7}^{\alpha,+}$}}}
        \put(32.0,26.3){\makebox(0,0){\scalebox{0.75}{$t_{2,8}^{\alpha,+}$}}}
        \put(36.8,26.3){\makebox(0,0){\scalebox{0.75}{$t_{2,9}^{\alpha,+}$}}}
        \put(41.4,26.3){\makebox(0,0){\scalebox{0.75}{$t_{2,10}^{\alpha,+}$}}}

    \end{overpic}
    \caption{In this figure, we consider the expanding map in Example~\ref{ex:lattes-map} with the connection of $1$-tiles in the upper-left figure. Let $\alpha \coloneqq \bigl\{\widehat{-1}^+,\,\widehat{1}^+\bigr\}\in\sA_0^+$. The block vertex at level $1$ is shown in the upper-left figure and $\widehat{-1}^+=c(-1,b_1)$. The interval $I_\alpha^+$ and the sets $S_1^{\alpha,+}$, $S_2^{\alpha,+}$ consisting of $t_{1,i}^{\alpha,+}$ and $t_{2,j}^{\alpha,+}$ for $i\in\{0,\,1,\,\dots,\,4\}$ and $j\in\{0,\,1,\,\dots,\,10\}$. The black curves in the right figures are images of $I_\alpha^+$ under $F_1^{\alpha,+}$, $F_2^{\alpha,+}$ and $F^{\alpha,+}$. Together with the blue curves, they compose $\bigcup_{\alpha\in\sA_0^+} F_1^{\alpha,+}(I_\alpha^+)$ and $\bigcup_{\alpha\in\sA_0^+} F_2^{\alpha,+}(I_\alpha^+)$.}
    \label{fig:Fm}
\end{figure}

In the proof of the next lemma, we verify that the restrictions in~\eqref{eq:ch4-compatible-realization} agree at shared endpoints. They therefore define a continuous map $F_m^{\lambda,\sigma}\:I_\lambda^\sigma\to S^2$. See Figure~\ref{fig:Fm} for an example.

\begin{lemma}\label{prop:ch4-road-data}\label{connection-path-compression}
Let $f\:S^2\to S^2$ be an expanding Thurston map, and let $\cC\supseteq\post f$ be an $f$-invariant Jordan curve. Fix a connection of $1$-tiles satisfying \textup{(C1)}--\textup{(C2)} of Definition~\ref{Propertiesofconnections}, a color $\sigma\in\{+,\,-\}$, and a road $I_\lambda^\sigma$ as in Definition~\ref{def:ch4-initial-roads}. In its level-$m$ subdivision $\mathcal P_m^{\lambda,\sigma}$, each node $s$ is labelled by the block vertex $b_m^{\lambda,\sigma}(s)$, and each atom $J$ is assigned the $m$-tile $X_m^{\lambda,\sigma}(J)$. Use the atom coordinates $\kappa_J^\sigma$ of Definition~\ref{def:ch4-finite-realizations}. The following hold:
\begin{enumerate}
\smallskip
\item For each integer $m>\operatorname{birth}(\lambda)$, $P_m^\sigma(\lambda)$ is a simple path. Its vertices $b_m^{\lambda,\sigma}(s)$ and $c\bigl(X_m^{\lambda,\sigma}(J)\bigr)$ occur in the order of the nodes $s$ and atoms $J$ along $I_\lambda^\sigma$. We have $S_m^{\lambda,\sigma}\subseteq S_{m+1}^{\lambda,\sigma}$ and $b_{m+1}^{\lambda,\sigma}(s)=\operatorname{Lift}_m^\sigma\bigl(b_m^{\lambda,\sigma}(s)\bigr)$ for every $s\in S_m^{\lambda,\sigma}$. For an integer $n>\operatorname{birth}(\lambda)$ and $J=[a_-,a_+]\in\mathcal P_n^{\lambda,\sigma}$, set $X\=X_n^{\lambda,\sigma}(J)$ and $p_\pm\=f^n\bigl(\prV\bigl(b_n^{\lambda,\sigma}(a_\pm)\bigr)\bigr)$. Then $p_-\ne p_+$ and, for every integer $r\ge1$,
\begin{equation}\label{eq:ch4-affine-refinement}
 \begin{aligned}
 \bigl\{K\in\mathcal P_{n+r}^{\lambda,\sigma}:K\subseteq J\bigr\}
 &=\bigl\{(\kappa_J^\sigma)^{-1}(D):D\in\mathcal P_r^{\alpha_J,\sigma}\bigr\},\\
 b_{n+r}^{\lambda,\sigma}(s)
 &=L(X,r,b_r^{\alpha_J,\sigma}(\kappa_J^\sigma(s)))
 \quad \text{ for all } s\in S_{n+r}^{\lambda,\sigma}\cap J,\\
 X_{n+r}^{\lambda,\sigma}\bigl((\kappa_J^\sigma)^{-1}(D)\bigr)
 &=L(X,r,X_r^{\alpha_J,\sigma}(D))
 \quad \text{ for all } D\in\mathcal P_r^{\alpha_J,\sigma}.
 \end{aligned}
\end{equation}
\smallskip
\item For each integer $m>\operatorname{birth}(\lambda)$, define $\mesh\mathcal P_m^{\lambda,\sigma}\=\max\bigl\{\abs{J}:J\in\mathcal P_m^{\lambda,\sigma}\bigr\}$, where $\abs{J}$ denotes the length of $J$ in the road coordinate. Then $\lim_{m\to+\infty}\mesh\mathcal P_m^{\lambda,\sigma}=0$.
\end{enumerate}
\end{lemma}

\begin{proof}
\textup{(i)} Fix an initial label $\alpha=\{\widehat p^\sigma,\,\widehat q^\sigma\}\in\sA_0^\sigma$ with endpoint order $(p,q)$ and let $n\in \N$. For each tile vertex $c(X)$ on $P_n^\sigma(\alpha)$, let $b_-$ and $b_+$ be the block vertices immediately before and after $c(X)$, respectively. Let $p_\pm=f^n(\prV(b_\pm))$. Since $b_-$ and $b_+$ are distinct and both adjacent to $c(X)$,~\eqref{eq:ch4-incident-block-uniqueness} gives $p_-\ne p_+$ and $b_X^\sigma(p_\pm)=b_\pm$. By Proposition~\ref{lem:ch4-sector-copies}~\textup{(i)} and~\eqref{eq:ch4-copy-values}, $L(X,1,\,\cdot\,)$ maps $[\widehat{p_-}_1^\sigma,\widehat{p_+}_1^\sigma]_{\Gamma_1^\sigma}$ to a simple path from $\operatorname{Lift}_n^\sigma(b_-)$ to $\operatorname{Lift}_n^\sigma(b_+)$. Concatenate these paths at their common endpoints in the order of the $n$-tiles on $P_n^\sigma(\alpha)$. At each join, $A_{n,n+1}$ maps the $(n+1)$-tiles immediately before and after the common endpoint to distinct $n$-tiles, so there is no backtracking. The concatenation is therefore the unique path $P_{n+1}^\sigma(\alpha)$ from $\widehat p_{n+1}^\sigma$ to $\widehat q_{n+1}^\sigma$ in the tree $\Gamma_{n+1}^\sigma$. At level $(n+1)$, the atoms of $\mathcal P_{n+1}^{\alpha,\sigma}$ correspond to the tile vertices of $P_{n+1}^\sigma(\alpha)$ in path order. The old nodes retain their parameters and their blocks continue by $\operatorname{Lift}_n^\sigma$.

The path of each label $\lambda=(k,U,\alpha)$ is simple by~\eqref{eq:road-data} and Proposition~\ref{lem:ch4-sector-copies}~\textup{(i)}. Thus, for each label $\lambda$, the endpoints of every atom $J$ on $I_\lambda^\sigma$ are labeled by distinct block vertices adjacent to the tile vertex for $J$. Equation~\eqref{eq:ch4-incident-block-uniqueness} implies $p_-\ne p_+$.

For $J\in\mathcal P_n^{\alpha,\sigma}$, put $X\=X_n^{\alpha,\sigma}(J)$. In $P_{n+1}^\sigma(\alpha)$, the block vertices labeling the level-$(n+1)$ nodes in $J$ and the tile vertices for its level-$(n+1)$ subintervals form the subpath $L(X,1,P_1^\sigma(\alpha_J))$. The affine map $\kappa_J^\sigma$ sends each endpoint $a_\pm$ of $J$ to the endpoint of $I_{\alpha_J}^\sigma$ labeled by $p_\pm$. Both subdivisions of $J$ and $I_{\alpha_J}^\sigma$ are uniform with the same number of intervals, so it sends the child intervals of $J$ onto $\mathcal P_1^{\alpha_J,\sigma}$. For $D\in\mathcal P_1^{\alpha_J,\sigma}$ and $K=(\kappa_J^\sigma)^{-1}(D)$, the tile and block labels are
\begin{equation}\label{eq:ch4-road-recursion}
 \begin{aligned}
 X_{n+1}^{\alpha,\sigma}(K)
 &=L\bigl(X,1,X_1^{\alpha_J,\sigma}(D)\bigr),\\
 b_{n+1}^{\alpha,\sigma}(s)
 &=L\bigl(X,1,b_1^{\alpha_J,\sigma}(\kappa_J^\sigma(s))\bigr)
 \quad \text{ for all } s\in(\kappa_J^\sigma)^{-1}\bigl(S_1^{\alpha_J,\sigma}\bigr).
 \end{aligned}
\end{equation}

Fix a road label $\lambda=(k,W,\alpha)$ and an integer $j\ge1$. Set $m\=k+j$. For $J=[a_-,a_+]\in\mathcal P_m^{\lambda,\sigma}$, equations~\eqref{eq:road-data} and~\eqref{eq:ch4-copy-values} give $p_\pm=f^j\bigl(\prV\bigl(b_j^{\alpha,\sigma}(a_\pm)\bigr)\bigr)$. By~\eqref{eq:road-data} and Proposition~\ref{lem:ch4-sector-copies}~\textup{(ii)}, we have
\begin{equation*}
 L(W,j+1,\,\cdot\,)\circ L\bigl(X_j^{\alpha,\sigma}(J),1,\,\cdot\,\bigr)
 =L\bigl(X_m^{\lambda,\sigma}(J),1,\,\cdot\,\bigr).
\end{equation*}
Applying $L(W,j+1,\,\cdot\,)$ to~\eqref{eq:ch4-road-recursion} at level $j$, we obtain~\eqref{eq:ch4-affine-refinement} for $\lambda$ with $n=m$ and $r=1$. For $K\in\mathcal P_{m+1}^{\lambda,\sigma}$ satisfying $K\subseteq J$, the definition of $L$ gives
\begin{equation}\label{eq:ch4-road-ancestry}
 A_{m,m+1}\bigl(X_{m+1}^{\lambda,\sigma}(K)\bigr)=X_m^{\lambda,\sigma}(J).
\end{equation}
For $s\in S_m^{\lambda,\sigma}$, we also have
\begin{equation*}
 \begin{aligned}
 b_{m+1}^{\lambda,\sigma}(s)
 &=L\bigl(W,j+1,b_{j+1}^{\alpha,\sigma}(s)\bigr)
 =L\bigl(W,j+1,\operatorname{Lift}_j^\sigma\bigl(b_j^{\alpha,\sigma}(s)\bigr)\bigr)\\
 &=\operatorname{Lift}_m^\sigma\bigl(L\bigl(W,j,b_j^{\alpha,\sigma}(s)\bigr)\bigr)
 =\operatorname{Lift}_m^\sigma\bigl(b_m^{\lambda,\sigma}(s)\bigr).
 \end{aligned}
\end{equation*}

For initial roads, the tile formula in~\eqref{eq:ch4-road-recursion} and the definition of $L$ give~\eqref{eq:ch4-road-ancestry}.

For an atom $D\in\mathcal P_1^{\alpha_J,\sigma}$, set $K\=(\kappa_J^\sigma)^{-1}(D)$. By~\eqref{eq:ch4-copy-values}, we have $\alpha_K=\alpha_D$. Since $\kappa_K^\sigma$ and $\kappa_D^\sigma\circ(\kappa_J^\sigma|_K)$ are affine and agree at both endpoints of $K$, we have $\kappa_K^\sigma=\kappa_D^\sigma\circ(\kappa_J^\sigma|_K)$. We iterate~\eqref{eq:ch4-road-recursion} and compose the maps $L$ by Proposition~\ref{lem:ch4-sector-copies}~\textup{(ii)}. This gives~\eqref{eq:ch4-affine-refinement} for every integer $r\ge1$.

Each restriction $F_m^{\lambda,\sigma}\big|_J$ is continuous by~\eqref{eq:ch4-compatible-realization}. At each endpoint $s$ of $J$, we have $\bigl(F_m^{\lambda,\sigma}\big|_J\bigr)(s)=\prV\bigl(b_m^{\lambda,\sigma}(s)\bigr)$. Hence, $F_m^{\lambda,\sigma}$ is continuous. By~\eqref{eq:ch4-block-persistence}, we have $\prV\circ\operatorname{Lift}_m^\sigma=\prV$. Thus, for every node $s\in S_m^{\lambda,\sigma}$ and integer $r\ge m$, we have
\begin{equation}\label{eq:node-realization}
b_{m+1}^{\lambda,\sigma}(s)=\operatorname{Lift}_m^\sigma\bigl(b_m^{\lambda,\sigma}(s)\bigr)\quad\text{ and }\quad
F_r^{\lambda,\sigma}(s)=\prV\bigl(b_m^{\lambda,\sigma}(s)\bigr).
\end{equation}

Fix a road label $\lambda=(k,X,\alpha)$ and an atom $J\in\mathcal P_r^{\alpha,\sigma}$. In the coordinates of~\eqref{eq:road-data}, equations~\eqref{eq:ch4-copy-values} and~\eqref{eq:ch4-compatible-realization} give $f^k\circ\bigl(F_{k+r}^{\lambda,\sigma}\big|_J\bigr)=F_r^{\alpha,\sigma}|_J$. Since $f^k$ maps $X_{k+r}^{\lambda,\sigma}(J)$ homeomorphically to $X_r^{\alpha,\sigma}(J)$, we obtain
\begin{equation}\label{eq:compact-born-road}
 F_{k+r}^{\lambda,\sigma}\big|_J
 =\bigl(f^k\big|_{X_{k+r}^{\lambda,\sigma}(J)}\bigr)^{-1}\circ(F_r^{\alpha,\sigma}|_J).
\end{equation}
For every initial road $\alpha$ and atom $J\in\mathcal P_m^{\alpha,\sigma}$,~\eqref{eq:ch4-compatible-realization} gives
\begin{equation}\label{eq:ch4-original-carriers}
 F_m^{\alpha,\sigma}(J)\subseteq X_m^{\alpha,\sigma}(J).
\end{equation}

\smallskip
\textup{(ii)} Fix a road label $\alpha\in\sA_0^\sigma$. Suppose an atom $J\=[a,b]\in\mathcal P_{m_0}^{\alpha,\sigma}$ satisfies $J\in\mathcal P_m^{\alpha,\sigma}$ for every integer $m\ge m_0$. Set $u_-\=F_{m_0}^{\alpha,\sigma}(a)$ and $u_+\=F_{m_0}^{\alpha,\sigma}(b)$. Part~\textup{(i)} gives $f^{m_0}(u_-)\ne f^{m_0}(u_+)$, so we have $u_-\ne u_+$. By~\eqref{eq:node-realization} and~\eqref{eq:ch4-original-carriers}, we have
\begin{equation*}
0<\rho(u_-,u_+)\le\diam_\rho X_m^{\alpha,\sigma}(J)\quad \text{ for all } m\ge m_0.
\end{equation*}
This contradicts $\lim_{m\to+\infty}\diam_\rho X_m^{\alpha,\sigma}(J)=0$.

Thus, every atom eventually splits. Each proper child $K$ of an atom $J$ satisfies $|K|\le |J|/2$. Since $\mathcal P_m^{\alpha,\sigma}$ is finite, there is a level $M>m$ by which all its atoms have split, and consequently $\mesh\mathcal P_M^{\alpha,\sigma}\le(1/2)\mesh\mathcal P_m^{\alpha,\sigma}$. The meshes are nonincreasing, so iteration gives $\lim_{m\to+\infty}\mesh\mathcal P_m^{\alpha,\sigma}=0$. For $\lambda=(k,X,\alpha)$ and $m>k$, we have $\mesh\mathcal P_m^{\lambda,\sigma}=\mesh\mathcal P_{m-k}^{\alpha,\sigma}$ by~\eqref{eq:road-data}. Thus, for every road label $\lambda$,
\begin{equation}\label{eq:ch4-road-mesh}
 \lim_{m\to+\infty}\mesh\mathcal P_m^{\lambda,\sigma}=0.
\end{equation}
Now the proof is complete.
\end{proof}

\begin{remark}
Fix $\sigma\in\{+,\,-\}$, a road label $\lambda$, an integer $m>\operatorname{birth}(\lambda)$, and an atom $J=[a_-,a_+]\in\mathcal P_m^{\lambda,\sigma}$. We denote the path $L\bigl(X_m^{\lambda,\sigma}(J),1,P_1^\sigma(\alpha_J)\bigr)$ by $R_{m,J}^{\lambda,\sigma}$ and call it the \emph{distinguished replacement packet} of the tile passage corresponding to $J$. By Lemma~\ref{prop:ch4-road-data}~\textup{(i)}, this is the subpath of $P_{m+1}^\sigma(\lambda)$ from $\operatorname{Lift}_m^\sigma\bigl(b_m^{\lambda,\sigma}(a_-)\bigr)$ to $\operatorname{Lift}_m^\sigma\bigl(b_m^{\lambda,\sigma}(a_+)\bigr)$. Orient each packet in this direction. Concatenating the packets in the order of the atoms of $\mathcal P_m^{\lambda,\sigma}$ gives $P_{m+1}^\sigma(\lambda)$.
\end{remark}

We now construct $Y_n^\sigma$ from the roads with labels in $\sA_n^\sigma$. Two nodes have the same image in $Y_n^\sigma$ whenever their block labels agree at a common level.

\begin{definition}\label{c:relation}
For each $n\in \N_0$, equip $K_n^\sigma\=\bigsqcup_{\lambda\in\sA_n^\sigma}I_\lambda^\sigma$ with the disjoint-union topology. For each integer $m\ge n+1$, set
\begin{equation}\label{eq:ps-generator}
 \mathcal R_{n,m}^{H,\sigma}
 \= \bigl\{((\lambda,s),(\tau,t)):\lambda,\tau\in\sA_n^{\sigma},
       \ s\in S_m^{\lambda,\sigma},\ t\in S_m^{\tau,\sigma},
       \ b_m^{\lambda,\sigma}(s)=b_m^{\tau,\sigma}(t)\bigr\}.
\end{equation}
For $E\subseteq(K_n^\sigma)^2$, let $\widehat E$ be the intersection of all closed equivalence relations on $K_n^\sigma$ containing $E$. Define a relation $\sim_n^\sigma$ on $K_n^\sigma$ by
\begin{equation*}
 \sim_n^\sigma\=\widehat{\bigcup_{m=n+1}^{+\infty}\mathcal R_{n,m}^{H,\sigma}}.
\end{equation*}
Denote by $Y_n^\sigma$ the quotient $K_n^\sigma/\sim_n^\sigma$, equipped with the quotient topology. Let $\eta_n^\sigma\:K_n^\sigma\to Y_n^\sigma$ be the quotient map given by $\eta_n^\sigma(\lambda,s)=[(\lambda,s)]$.
\end{definition}

\begin{definition}\label{def:ch4-path-matchings}
Fix integers $n\ge0$, $m\ge n+1$ and road labels $\lambda,\tau\in\sA_n^\sigma$. Let $Q$ be a common subpath of $P_m^\sigma(\lambda)$ and $P_m^\sigma(\tau)$ with endpoints $b_0,b_1\in\sB_m^\sigma$ (here $b_0$ and $b_1$ need not be distinct). For $i\in\{0,\,1\}$, let $(\lambda,s_i)\in S_m^{\lambda,\sigma}$ and $(\tau,t_i)\in S_m^{\tau,\sigma}$ be the unique nodes with block label $b_i$ on the two roads. Define
\begin{equation*}
 A\=\{\lambda\}\times[\min\{s_0,\,s_1\},\max\{s_0,\,s_1\}]\quad\text{ and }\quad
 B\=\{\tau\}\times[\min\{t_0,\,t_1\},\max\{t_0,\,t_1\}].
\end{equation*}

A homeomorphism $h\:A\to B$ is called a \emph{matching} if
\begin{equation}\label{eq:ch4-matching-definition}
h\bigl(S_r^{\lambda,\sigma}\cap A\bigr)=S_r^{\tau,\sigma}\cap B,\quad
b_r^{\tau,\sigma}\circ h=b_r^{\lambda,\sigma}
\quad\text{ on }S_r^{\lambda,\sigma}\cap A \text{ for all } r\ge m.
\end{equation}
Define the \emph{relation induced by $h$} by
\begin{equation*}
 \mathcal M(h)\=\{(x,h(x)):x\in A\}\subseteq K_n^\sigma\times K_n^\sigma.
\end{equation*}
For an atom $J=[a_-,a_+]\in\mathcal P_m^{\lambda,\sigma}$, put $X\=X_m^{\lambda,\sigma}(J)$. We have
\begin{equation*}
\bigl[b_{m+1}^{\lambda,\sigma}(a_-),b_{m+1}^{\lambda,\sigma}(a_+)\bigr]_{\Gamma_{m+1}^\sigma}=L(X,1,P_1^\sigma(\alpha_J)),\quad b_{m+1}^{\lambda,\sigma}(s) =\operatorname{Lift}_m^\sigma\bigl(b_m^{\lambda,\sigma}(s)\bigr)\quad \text{ for all } s\in S_m^{\lambda,\sigma}.
\end{equation*}
\end{definition}

\begin{lemma}\label{lem:ch4-path-matchings}
Let $f\:S^2\to S^2$ be an expanding Thurston map, and let $\cC\supseteq\post f$ be an $f$-invariant Jordan curve. Fix a connection of $1$-tiles satisfying \textup{(C1)}--\textup{(C2)} of Definition~\ref{Propertiesofconnections}. Fix a color $\sigma\in\{+,\,-\}$, integers $n\ge0$, $m\ge n+1$ and road labels $\lambda,\,\tau\in\sA_n^\sigma$. Let $Q$ be an arbitrary common subpath of $P_m^\sigma(\lambda)$ and $P_m^\sigma(\tau)$ whose endpoints lie in $\sB_m^\sigma$. Let $A\subseteq I_\lambda^\sigma$ and $B\subseteq I_\tau^\sigma$ be the intervals defined from $Q$ in Definition~\ref{def:ch4-path-matchings}. Then there is a unique matching $h\:A\to B$, and we have $\mathcal M(h)\subseteq\mathord{\sim_n^\sigma}$ for the relation $\sim_n^\sigma$ of Definition~\ref{c:relation}. For every integer $r\ge m$, there is a common subpath $Q_r\subseteq P_r^\sigma(\lambda)\cap P_r^\sigma(\tau)$ with $Q_m=Q$. The map $h\:A\to B$ is also the unique matching for every $Q_r$.
\end{lemma}
\begin{proof}
Let $s_0<\cdots<s_\ell$ be the nodes in $S_m^{\lambda,\sigma}\cap A$. For each index $i$, let $t_i\in S_m^{\tau,\sigma}\cap B$ be the unique node satisfying $b_m^{\tau,\sigma}(t_i)=b_m^{\lambda,\sigma}(s_i)$. Since $Q$ is a subpath of $P_m^\sigma(\tau)$, we have either $t_0<\cdots<t_\ell$ or $t_0>\cdots>t_\ell$ when $\ell\ge1$. Define $h(s_i)\=t_i$ and extend $h$ affinely on each atom $J_i\=[s_{i-1},s_i]$. This defines a homeomorphism $h\:A\to B$. For $\ell=0$, we have $A=\{s_0\}$, $B=\{t_0\}$, and $h(s_0)=t_0$. For each index $1\le i\le\ell$, the interval $K_i=h(J_i)$ is a level-$m$ atom of $B$. The identities $b_m^{\lambda,\sigma}(s_j)=b_m^{\tau,\sigma}(t_j)$ for $j\in \{i-1,\,i\}$ imply $\alpha_{J_i}=\alpha_{K_i}$. The affine map $(\kappa_{K_i}^\sigma)^{-1}\circ\kappa_{J_i}^\sigma$ sends $s_j$ to $t_j$, so we have
$ h|_{J_i}=(\kappa_{K_i}^\sigma)^{-1}\circ\kappa_{J_i}^\sigma.$
Equations~\eqref{eq:ch4-affine-refinement} and~\eqref{eq:node-realization} give, for every integer $r\ge m$,
\begin{equation}
\label{eq:ch4-atomic-matching}
\begin{aligned}
    h\bigl(S_r^{\lambda,\sigma}\cap A\bigr)
    &=S_r^{\tau,\sigma}\cap B,\\
    b_r^{\tau,\sigma}(h(s))
    &=b_r^{\lambda,\sigma}(s)
    \quad \text{ for all } s\in S_r^{\lambda,\sigma}\cap A.
\end{aligned}
\end{equation}

Write $A\=\{\lambda\}\times[a,b]$ and define $Q_r\=\bigl[b_r^{\lambda,\sigma}(a),b_r^{\lambda,\sigma}(b)\bigr]_{\Gamma_r^\sigma}$ for every integer $r\ge m$. By~\eqref{eq:ch4-atomic-matching}, we have $b_r^{\tau,\sigma}(h(a))=b_r^{\lambda,\sigma}(a)$ and $b_r^{\tau,\sigma}(h(b))=b_r^{\lambda,\sigma}(b)$. Since $\Gamma_r^\sigma$ is a tree, the path joining these two block vertices is unique. Hence, we have $Q_r\subseteq P_r^\sigma(\lambda)\cap P_r^\sigma(\tau)$ and $Q_m=Q$. By~\eqref{eq:ch4-atomic-matching}, the map $h\:A\to B$ is a matching for every $Q_r$.

Fix an integer $r\ge m$ and a point $s\in A$. By~\eqref{eq:ch4-road-mesh}, take nodes $u_j\in S_{r+j}^{\lambda,\sigma}\cap A$ satisfying $\lim_{j\to+\infty}u_j=s$. We have $(u_j,h(u_j))\in\mathcal R_{n,r+j}^{H,\sigma}$, and continuity of $h$ gives $\lim_{j\to+\infty}(u_j,h(u_j))=(s,h(s))$. Therefore,
\begin{equation}\label{eq:ch4-matching-graph}
 \mathcal M(h)\subseteq
 \overline{\bigcup_{\ell=1}^{+\infty}\mathcal R_{n,\ell}^{H,\sigma}}
 \subseteq\sim_n^\sigma.
\end{equation}
It remains to show that $h$ is the unique matching for $Q_r$. Let $g\:A\to B$ be another matching for $Q_r$. For every $i\in \N$, we have $b_{r+i}^{\tau,\sigma}(g(s_i))=b_{r+i}^{\lambda,\sigma}(s_i)=b_{r+i}^{\tau,\sigma}(h(s_i))$. Since $P_{r+i}^\sigma(\tau)$ is simple, we obtain $g(s_i)=h(s_i)$. Taking limits gives $g(s)=h(s)$ for every $s\in A$.
\end{proof}

For points $x,y$ in a topological $\mathbb R$-tree $T$, write $(x,y)_T\=[x,y]_T\smallsetminus\{x,\,y\}$. Let $A$ and $B$ be arbitrary sets, and let $u\:A\to B$ be a map. Write $\ker u\=\bigl\{(x,y)\in A^2:u(x)=u(y)\bigr\}$.

\begin{definition}\label{def:ch4-initial-quotient-data}
For $\alpha\in\sA_0^\sigma$, we call the restriction $q_\alpha^\sigma\=\eta_0^\sigma|_{I_\alpha^\sigma}\:I_\alpha^\sigma\to Y_0^\sigma$ the \emph{initial road map}. For $p\in\post f$, let $e_p$ be an endpoint of an initial road $I_\alpha^\sigma$ such that $\widehat p^\sigma\in\alpha$ and $b_1^{\alpha,\sigma}(e_p)=\widehat p_1^\sigma$. Set $\xi_p^\sigma\=\eta_0^\sigma(e_p)$.
\end{definition}
The following proposition will show that $\xi_p^\sigma$ is independent of the choice of $e_p$.

\begin{prop}\label{c:initial}
Let $f\:S^2\to S^2$ be an expanding Thurston map, and let $\cC\supseteq\post f$ be an $f$-invariant Jordan curve. Fix a color $\sigma\in\{+,\,-\}$ and a connection of $1$-tiles satisfying \textup{(C1)}--\textup{(C2)} of Definition~\ref{Propertiesofconnections}. Let $Y_0^\sigma\=K_0^\sigma/\sim_0^\sigma$ be the quotient space defined in Definition~\ref{c:relation}. For each road label $\alpha\in\sA_0^\sigma$, recall the map $q_\alpha^\sigma\:I_\alpha^\sigma\to Y_0^\sigma$ from Definition~\ref{def:ch4-initial-quotient-data}. For each point $p\in\post f$, recall the point $\xi_p^\sigma$ also from Definition~\ref{def:ch4-initial-quotient-data}.
\begin{enumerate}
\smallskip
\item[\textup{(i)}] Fix a road label $\alpha\in\sA_0^\sigma$ and a point $p\in\post f$. If $s$ is an endpoint of $I_\alpha^\sigma$ with $b_1^{\alpha,\sigma}(s)=\widehat p_1^\sigma$, then $q_\alpha^\sigma(s)=\xi_p^\sigma$.
\smallskip
\item[\textup{(ii)}] $Y_0^\sigma$ is a finite tree, and $\xi_p^\sigma\ne\xi_q^\sigma$ for distinct points $p,q\in\post f$.
\smallskip
\item[\textup{(iii)}] For each road label $\alpha=\{\widehat p^\sigma,\,\widehat q^\sigma\}\in\sA_0^\sigma$, the map $q_\alpha^\sigma$ embeds $I_\alpha^\sigma$ onto $[\xi_p^\sigma,\xi_q^\sigma]_{Y_0^\sigma}$. Moreover, we have
\begin{equation*}
 Y_0^\sigma
 =\bigcup_{\alpha\in\sA_0^\sigma}q_\alpha^\sigma(I_\alpha^\sigma).
\end{equation*}
\end{enumerate}
\end{prop}
\begin{proof}
\textup{(i)} Fix $p\in\post f$ and $\alpha\in\sA_0^\sigma$. Let $e_p$ be the endpoint chosen in Definition~\ref{def:ch4-initial-quotient-data}. Write $e_p\=(\tau,t)$, where $\tau$ is its road label and $t$ is its coordinate on $I_\tau^\sigma$. For any endpoint $s$ with $b_1^{\alpha,\sigma}(s)=\widehat p_1^\sigma=b_1^{\tau,\sigma}(t)$,~\eqref{eq:ps-generator} gives $((\alpha,s),e_p)\in\mathcal R_{0,1}^{H,\sigma}\subseteq\sim_0^\sigma$. Thus $q_\alpha^\sigma(s)=\xi_p^\sigma$.

\smallskip
\textup{(ii)} For an integer $m\ge1$, set $H_m^\sigma\=\bigcup_{\alpha\in\sA_0^\sigma}P_m^\sigma(\alpha)\subseteq\Gamma_m^\sigma$. The graph $H_m^\sigma$ is a tree, and every leaf equals $\widehat p_m^\sigma$ for some $p\in\post f$. Mark each $\widehat p_m^\sigma$ by $p$. In a copy of $H_m^\sigma$, for each unmarked vertex of degree two in $H_m^\sigma$, replace the vertex and its two incident edges by one edge joining its neighbors. Let $\ell$ be the number of leaves of the reduced tree. Then, $\sum_{\deg v\ge3}(\deg v-2)=\ell-2\le \card\post f-2$. Thus, at most $\card\post f-2$ unmarked vertices remain, and there are at most $2|\post f|-2$ vertices in total. Only finitely many graphs have this vertex bound up to isomorphism, and each has finitely many placements of the marks. We can therefore choose levels $1\le m_1<m_2<\cdots$ and graph isomorphisms $\varphi_j$ from the reduced tree at level $m_j$ to that at level $m_1$, with $\varphi_j\bigl(\widehat p_{m_j}^\sigma\bigr)=\widehat p_{m_1}^\sigma$ for every $p\in\post f$. Let $S$ be the reduced tree at level $m_1$, and denote $\widehat p_{m_1}^\sigma$ by $p$. For each vertex $v$ of $S$ and each $j\ge1$, denote $\varphi_j^{-1}(v)$ by $v_j$. If $v=p$ is a marked vertex of $S$, then $v_j=\widehat p_{m_j}^\sigma$. Note that an unmarked $v_j$ can be a block vertex or a tile vertex.

For each edge $e$ of $S$, set $\sA^\sigma(e)\=\{\{\widehat p^\sigma,\,\widehat q^\sigma\}\in\sA_0^\sigma:e\subseteq[p,q]_S\}$. If $x\notin e$, write $e\=[v,w]$ so that $x,w$ lie in different components of $S\smallsetminus\{v\}$. Take a mark $q$ in the component of $w$. Put $p\=v$ if $v$ is marked, and otherwise take a mark $p$ in a third component. Then, $e\subseteq[p,q]_S$ and $x\notin[p,q]_S$. Thus, for $e=[v,w]$, we have $\bigcap_{\alpha\in\sA^\sigma(e)}P_{m_j}^\sigma(\alpha)=[v_j,w_j]_{\Gamma_{m_j}^\sigma}$.

Let $b_0,\,b_1$ be the first and last block vertices on $[v_j,w_j]_{\Gamma_{m_j}^\sigma}$ from $v_j$ to $w_j$, and put $Q_{e,j}^\sigma\=[b_0,b_1]_{\Gamma_{m_j}^\sigma}$. For $\alpha\in\sA^\sigma(e)$, let $s_i^\alpha\in S_{m_j}^{\alpha,\sigma}$ satisfy $b_{m_j}^{\alpha,\sigma}(s_i^\alpha)=b_i$, $i=0,1$, and set $J_{e,j}^{\alpha,\sigma}\=\{\alpha\}\times[\min_i s_i^\alpha,\max_i s_i^\alpha]$. Lemma~\ref{lem:ch4-path-matchings} gives $h_{e,j}^{\alpha,\tau,\sigma}\:J_{e,j}^{\alpha,\sigma}\to J_{e,j}^{\tau,\sigma}$. For $i=0,1$ and $\alpha,\tau\in\sA^\sigma(e)$,~\eqref{eq:node-realization} gives $b_{m_{j+1}}^{\alpha,\sigma}(s_i^\alpha)=b_{m_{j+1}}^{\tau,\sigma}(s_i^\tau)\in Q_{e,j+1}^\sigma$. Hence, we have
\begin{equation}\label{eq:ch4-edge-core-refinement}
J_{e,j}^{\alpha,\sigma}\subseteq J_{e,j+1}^{\alpha,\sigma},\quad
h_{e,j+1}^{\alpha,\tau,\sigma}\big|_{J_{e,j}^{\alpha,\sigma}}
=h_{e,j}^{\alpha,\tau,\sigma}
\quad \text{ for all } \alpha,\tau\in\sA^\sigma(e).
\end{equation}
For each edge $e$ of $S$ and $\alpha\in\sA^\sigma(e)$, define $J_e^{\alpha,\sigma}  \=\overline{\bigcup_{j=1}^{+\infty}J_{e,j}^{\alpha,\sigma}}$.
Extend the maps $h_{e,j}^{\alpha,\tau,\sigma}$ and their inverses $\bigl(h_{e,j}^{\alpha,\tau,\sigma}\bigr)^{-1}$ to the closed intervals $J_e^{\alpha,\sigma}$ and $J_e^{\tau,\sigma}$, respectively, by endpoint limits. This defines homeomorphisms $h_e^{\alpha,\tau,\sigma}\:J_e^{\alpha,\sigma}\to J_e^{\tau,\sigma}$ satisfying
\begin{equation}\label{eq:ch4-edge-limit-matching}
 h_e^{\tau,\alpha,\sigma}=(h_e^{\alpha,\tau,\sigma})^{-1}\quad\text{ and }\quad
 \mathcal M(h_e^{\alpha,\tau,\sigma})
 =\overline{\bigcup_j\mathcal M\bigl(h_{e,j}^{\alpha,\tau,\sigma}\bigr)}
 \subseteq\sim_0^\sigma,
\end{equation}
where the inclusion uses the lemma and the closedness of $\sim_0^\sigma$.

Fix $\alpha=\{\widehat p^\sigma,\,\widehat q^\sigma\}\in\sA_0^\sigma$ and order $p,\,q$ so that $b_1^{\alpha,\sigma}(0)=\widehat p_1^\sigma$. List the vertices along $[p,q]_S$ as $p=v_0,\,v_2,\,\dots,\,v_k=q$, and put $e_i\=[v_{i-1},v_i]$. For $1\le i<k$, we have $0\le\min J_{e_{i+1},j}^{\alpha,\sigma}-\max J_{e_i,j}^{\alpha,\sigma}\le\mesh\mathcal P_{m_j}^{\alpha,\sigma}\to0$ by~\eqref{eq:ch4-road-mesh}. Therefore, we can define $s_\alpha$ by $J_{e_i}^{\alpha,\sigma}\=\{\alpha\}\times[s_\alpha(v_{i-1}),s_\alpha(v_i)]$ for $1\le i\le k$. These values satisfy $0=s_\alpha(v_0)\le\cdots\le s_\alpha(v_k)=1$.

For each edge $e=[v,w]$ of $S$ and $\alpha,\tau\in\sA^\sigma(e)$, the matching satisfies $h_e^{\alpha,\tau,\sigma}(s_\alpha(u))=s_\tau(u)$ for $u\in\{v,\,w\}$.

Let $\rho\:S\to T$ be the quotient map obtained by identifying all points of each edge $e$ satisfying $\card J_e^{\alpha,\sigma}=1$ for every $\alpha\in\sA^\sigma(e)$. Then $T$ is a finite tree. For each edge $e=[v,w]$, fix $\alpha_e\in\sA^\sigma(e)$ and a homeomorphism $\chi_e^\sigma\:J_e^{\alpha_e,\sigma}\to\rho(e)$ satisfying $\chi_e^\sigma(s_{\alpha_e}(u))=\rho(u)$ for $u\in\{v,\,w\}$. Define
\begin{equation}\label{eq:ch4-initial-model-map}
\eta(\alpha,s)\=\chi_e^\sigma(h_e^{\alpha,\alpha_e,\sigma}(s))
\quad \text{ for all } \alpha\in\sA^\sigma(e),\, (\alpha,s)\in J_e^{\alpha,\sigma}.
\end{equation}
For consecutive edges meeting at $v$, both formulas give $\eta(\alpha,s_\alpha(v))=\rho(v)$. Thus $\eta\:K_0^\sigma\to T$ is well-defined and continuous. The nondegenerate $J_e^{\alpha,\sigma}$ map homeomorphically onto consecutive edges of $\rho([p,q]_S)$, so $\eta|_{I_\alpha^\sigma}$ is injective. Since the paths $[p,q]_S$ cover $S$, we have $\eta(K_0^\sigma)=T$.

Fix $m\in \N$ and a pair $((\alpha,s),(\tau,t)) \in\mathcal R_{0,m}^{H,\sigma}$. Take an index $j$ satisfying $m_j\ge m$. By~\eqref{eq:node-realization}, the labels $b_{m_j}^{\alpha,\sigma}(s)$ and $b_{m_j}^{\tau,\sigma}(t)$ agree; denote their common value by $u$. If $u$ lies in $(v_j,w_j)_{\Gamma_{m_j}^\sigma}$ for an edge $e$ of $S$ with endpoints $v,w$, the matchings satisfy $h_{e,j}^{\alpha,\alpha_e,\sigma}(s)=h_{e,j}^{\tau,\alpha_e,\sigma}(t)$. Equations~\eqref{eq:ch4-edge-core-refinement} and~\eqref{eq:ch4-initial-model-map} then identify $\eta(\alpha,s)$ with $\eta(\tau,t)$. If $u$ equals $v_j$, write $\alpha$ as $\{\widehat p^\sigma,\,\widehat q^\sigma\}$. When $v$ lies in $(p,q)_S$, let $e_-,\,e_+$ be the two edges of $[p,q]_S$ incident to $v$. The node $(\alpha,s)$ satisfies
\begin{equation*}
(\alpha,s)\in J_{e_-,j}^{\alpha,\sigma}\cap J_{e_+,j}^{\alpha,\sigma}\subseteq J_{e_-}^{\alpha,\sigma}\cap J_{e_+}^{\alpha,\sigma}
 =\{(\alpha,s_\alpha(v))\}.
\end{equation*}
The equality $s=s_\alpha(v)$ holds at the endpoints of $[p,q]_S$ by the definition of $s_\alpha$, and the same argument on $I_\tau^\sigma$ gives $t=s_\tau(v)$. Substituting in~\eqref{eq:ch4-initial-model-map}, we obtain $\eta(\alpha,s)=\rho(v)=\eta(\tau,t)$. Hence $\ker\eta$ contains $\bigcup_{m=1}^{+\infty}\mathcal R_{0,m}^{H,\sigma}$. The kernel is closed because $\eta$ is continuous and $T$ is Hausdorff, so Definition~\ref{c:relation} gives $\sim_0^\sigma\subseteq\ker\eta$.

For the reverse inclusion, take a pair $((\alpha,s),(\tau,t))$ in $\ker\eta$ and let $z$ denote $\eta(\alpha,s)$. Take edges $e,e'$ of $S$ such that $(\alpha,s)\in J_e^{\alpha,\sigma}$ and $(\tau,t)\in J_{e'}^{\tau,\sigma}$. Since every leaf of $S$ is marked, there are marked vertices $p,\,q$ whose path $[p,q]_S$ contains $e\cup e'$. Let $\beta$ denote $\{\widehat p^\sigma,\,\widehat q^\sigma\}$; this label belongs to $\sA^\sigma(e)\cap\sA^\sigma(e')$. Define $r,\,r'$ by the formulas $r\=h_e^{\alpha,\beta,\sigma}(s)$ and $r'\=h_{e'}^{\tau,\beta,\sigma}(t)$. By~\eqref{eq:ch4-edge-limit-matching} and the inclusion $\sim_0^\sigma\subseteq\ker\eta$, we have $\eta(\beta,r)=z=\eta(\beta,r')$. Since $\eta$ is injective on $I_\beta^\sigma$, we obtain $r=r'$. Equation~\eqref{eq:ch4-edge-limit-matching} yields $ (\alpha,s)\sim_0^\sigma(\beta,r)=(\beta,r')\sim_0^\sigma(\tau,t)$. Therefore, $\ker\eta$ equals $\sim_0^\sigma$.

The induced continuous bijection $Y_0^\sigma\to T$ is a homeomorphism, since $K_0^\sigma$ is compact and $T$ is Hausdorff. Hence, $Y_0^\sigma$ is a finite tree. Each road has distinct endpoint images, so $\xi_p^\sigma\ne\xi_q^\sigma$ for $p\ne q$.

\smallskip

\textup{(iii)} For $\alpha=\{\widehat p^\sigma,\,\widehat q^\sigma\}\in\sA_0^\sigma$, the homeomorphism identifies $q_\alpha^\sigma$ with $\eta|_{I_\alpha^\sigma}$. Thus, $q_\alpha^\sigma$ embeds $I_\alpha^\sigma$ onto $[\xi_p^\sigma,\xi_q^\sigma]_{Y_0^\sigma}$, and $Y_0^\sigma=\bigcup_{\alpha\in\sA_0^\sigma}q_\alpha^\sigma(I_\alpha^\sigma)$.
\end{proof}

\begin{lemma}\label{lem:ch3-initial-road-dynamics}
Fix $\sigma\in\{+,\,-\}$, an initial label $\alpha\in\sA_0^\sigma$, and an atom $J\in\mathcal P_1^{\alpha,\sigma}$. For every integer $m\ge1$, we have
\begin{equation}\label{eq:ch3-initial-road-dynamics}
f\circ\bigl(F_{m+1}^{\alpha,\sigma}\big|_J\bigr)
=F_m^{\alpha_J,\sigma}\circ\kappa_J^\sigma.
\end{equation}
\end{lemma}
\begin{proof}
Set $X\=X_1^{\alpha,\sigma}(J)$ and $\beta\=\alpha_J$. Fix an integer $m\ge1$ and an atom $D\in\mathcal P_m^{\beta,\sigma}$. Set $K\=(\kappa_J^\sigma)^{-1}(D)$, $U\=X_m^{\beta,\sigma}(D)$, and $V\=X_{m+1}^{\alpha,\sigma}(K)$. Equation~\eqref{eq:ch4-affine-refinement} gives $V=L(X,m,U)$, so $f(V)=U$.

For each endpoint $s$ of $K$, the block formula in~\eqref{eq:ch4-affine-refinement} and equation~\eqref{eq:ch4-copy-values} give
\begin{equation*}
f^{m+1}\bigl(\prV\bigl(b_{m+1}^{\alpha,\sigma}(s)\bigr)\bigr)
=f^m\bigl(\prV\bigl(b_m^{\beta,\sigma}(\kappa_J^\sigma(s))\bigr)\bigr).
\end{equation*}
Thus $\alpha_K=\alpha_D$. By the definition of the atom coordinates, the affine maps $\kappa_K^\sigma$ and $\kappa_D^\sigma\circ(\kappa_J^\sigma|_K)$ agree at both endpoints of $K$, so they are equal.

Since $f\circ\bigl(f^{m+1}|_V\bigr)^{-1}=(f^m|_U)^{-1}$, substitution in~\eqref{eq:ch4-compatible-realization} gives
\begin{equation*}
f\circ\bigl(F_{m+1}^{\alpha,\sigma}\big|_K\bigr)
=\bigl(F_m^{\beta,\sigma}\big|_D\bigr)\circ(\kappa_J^\sigma|_K).
\end{equation*}
As $D$ ranges over $\mathcal P_m^{\beta,\sigma}$, the intervals $K$ cover $J$. This proves~\eqref{eq:ch3-initial-road-dynamics}.
\end{proof}

\begin{lemma}\label{lem:ch3-road-limits}
For each $\sigma\in\{+,\,-\}$ and road label $\lambda$, the continuous maps $F_m^{\lambda,\sigma}$, indexed by integers $m>\operatorname{birth}(\lambda)$, converge uniformly to a continuous map $F^{\lambda,\sigma}\:I_\lambda^\sigma\to S^2$.
\end{lemma}
\begin{proof}
Fix $\sigma$ and $\lambda$. At each node $s$, $F_m^{\lambda,\sigma}$ has value $\prV\bigl(b_m^{\lambda,\sigma}(s)\bigr)$ on either adjacent atom by~\eqref{eq:node-realization}. Hence, $F_m^{\lambda,\sigma}$ is continuous. Fix a visual metric $\rho$ for $f$ with expansion factor $\Lambda>1$. Recall from Subsections~\ref{sub:Thurston_maps} and~\ref{Pseudo_isotopies} that there exists a constant $C>0$ such that for every integer $m\ge0$,
\begin{equation}\label{eq:pi-visual-estimates}
\max_{U\in\X^m}\diam_\rho U\le C\Lambda^{-m}\quad\text{ and }\quad
\sup_{z\in S^2}\rho(H_1^m(z),z)\le C\Lambda^{-m}.
\end{equation}
For an integer $m>\operatorname{birth}(\lambda)$ and an atom $J\in\mathcal P_m^{\lambda,\sigma}$, equation~\eqref{eq:ch4-compatible-realization} gives
\begin{equation}\label{eq:ch4-road-containment}
F_m^{\lambda,\sigma}(J)\subseteq X_m^{\lambda,\sigma}(J)\quad\text{ and }\quad
F_m^{\lambda,\sigma}(\inte_\circ J)\subseteq\inte_\circ X_m^{\lambda,\sigma}(J).
\end{equation}

For a child interval $K\subseteq J$, equation~\eqref{eq:ch4-road-ancestry} and Lemma~\ref{lem:ch4-sector-ancestors}~\textup{(ii)} give $X_{m+1}^{\lambda,\sigma}(K)\subseteq H_1^m\bigl(X_m^{\lambda,\sigma}(J)\bigr)$. For $t\in K$, write $F_{m+1}^{\lambda,\sigma}(t)\=H_1^m(z)$ with $z\in X_m^{\lambda,\sigma}(J)$. Both $\rho(H_1^m(z),z)$ and $\rho\bigl(z,F_m^{\lambda,\sigma}(t)\bigr)$ are at most $C\Lambda^{-m}$ by~\eqref{eq:pi-visual-estimates}. For integers $s>m>\operatorname{birth}(\lambda)$, summing from $m$ to $s-1$ gives
\begin{equation}\label{eq:ch4-road-Cauchy}
\begin{split}
\rho\bigl(F_{m+1}^{\lambda,\sigma}(t),F_m^{\lambda,\sigma}(t)\bigr)
&\le2C\Lambda^{-m},\\
\sup_{t\in I_\lambda^\sigma}\rho\bigl(F_s^{\lambda,\sigma}(t),F_m^{\lambda,\sigma}(t)\bigr)
&\le C_0\Lambda^{-m},
\quad C_0=\frac{2C}{1-\Lambda^{-1}}.
\end{split}
\end{equation}
Thus, $F_m^{\lambda,\sigma}$ converges uniformly to a continuous $F^{\lambda,\sigma}$.
\end{proof}

\section{P-zipper and its tree models}\label{ConstructionofYpm}
In this section, we define and construct the P-zipper $Z^\pm$ and its tree models $Y^\pm$, proving Theorem~\ref{A}.

\subsection{Construction of the P-zipper}
\label{Construction_of_primary_zippers}
For distinct $p,q\in\post f$ and $\sigma\in\{+,\,-\}$, write $\alpha_\sigma\=\{\widehat p^\sigma,\,\widehat q^\sigma\}$. Let $\jmath_{p,q}^\sigma\:[0,1]\to I_{\alpha_\sigma}^\sigma$ be the affine map sending $0$ and $1$ to the endpoints labeled by $p$ and $q$, respectively. For each integer $n\ge1$, define
\begin{equation*}
F_n^{p,q}\=F_n^{\alpha_+,+}\circ\jmath_{p,q}^+\quad\text{ and }\quad
G_n^{p,q}\=F_n^{\alpha_-,-}\circ\jmath_{p,q}^-.
\end{equation*}
Define $S_n^{p,q}\=\bigl(\jmath_{p,q}^+\bigr)^{-1}\bigl(S_n^{\alpha_+,+}\bigr)$. Using the road limits from Lemma~\ref{lem:ch3-road-limits}, define
\begin{equation*}
F^{p,q}\=F^{\alpha_+,+}\circ\jmath_{p,q}^+\quad\text{ and }\quad
G^{p,q}\=F^{\alpha_-,-}\circ\jmath_{p,q}^-.
\end{equation*}
Write $I\=[0,1]$, and define $l_{p,q}\=F^{p,q}(I)$ and $g_{p,q}\=G^{p,q}(I)$.

For each color $\sigma\in\{+,\,-\}$, we call the set
\begin{equation}\label{defZ0+}
Z_0^\sigma\=\bigcup_{\alpha\in\sA_0^\sigma}F^{\alpha,\sigma}(I_\alpha^\sigma)
\end{equation}
the \emph{primary zipper} of color $\sigma$.
Thus $Z_0^+=\bigcup_{p\ne q}l_{p,q}$ and $Z_0^-=\bigcup_{p\ne q}g_{p,q}$, where $p,q$ range over $\post f$.

\begin{lemma}\label{Z_0invariant}
Let $Z_0^+$ and $Z_0^-$ be the primary zippers given above. Then we have $f(Z_0^+)\subseteq Z_0^+$ and $f(Z_0^-)\subseteq Z_0^-$.
\end{lemma}
\begin{proof}
We only consider $Z_0^+$, and the proof is the same for $Z_0^-$. Fix $\alpha\in\sA_0^+$. Since the first-level atoms cover $I_\alpha^+$, equation~\eqref{eq:ch3-initial-road-dynamics} gives, for every integer $n\ge2$,
\begin{equation*}
f(F_n^{\alpha,+}(I_\alpha^+))
\subseteq\bigcup_{\beta\in\sA_0^+}F_{n-1}^{\beta,+}\bigl(I_\beta^+\bigr).
\end{equation*}
By Lemma~\ref{lem:ch3-road-limits}, $F_n^{\alpha,+}$ converges uniformly to $F^{\alpha,+}$, and for each $\beta\in\sA_0^+$, $F_{n-1}^{\beta,+}$ converges uniformly to $F^{\beta,+}$. The road intervals are compact and $\sA_0^+$ is finite, so the continuity of $f$ allows passage to the limit in this inclusion. Hence,
\begin{equation*}
f(F^{\alpha,+}(I_\alpha^+))
\subseteq\bigcup_{\beta\in\sA_0^+}F^{\beta,+}\bigl(I_\beta^+\bigr)
=Z_0^+.
\end{equation*}
Taking the union over $\alpha\in\sA_0^+$, we have $f(Z_0^+)\subseteq Z_0^+$.
\end{proof}

\medskip

If a set $A\subseteq S^2$ satisfies $f^{-1}(A)=A=f(A)$, we say that $A$ is \emph{completely $f$-invariant}. By \cite[Lemma~6.1]{BM}, since $f$ is expanding, $\card\post f\geq 3$.
A \emph{topological $\R$-tree} $T$ is a Hausdorff space in which every two distinct points are joined by a unique embedded arc.

\begin{definition}\label{def_gen_zipper}
Let $f\:S^2\to S^2$ be an expanding Thurston map. A pair $(Z^+,Z^-)$ of subsets of $S^2$, equipped with their subspace topologies, is a \emph{P-zipper for $f$} if the following conditions hold.
\begin{enumerate}
\smallskip\item Each of $Z^+$ and $Z^-$ is nonempty and path-connected.
\smallskip\item Each of $Z^+$ and $Z^-$ is completely $f$-invariant.
\smallskip\item For each $\sigma\in\{+,\,-\}$, there are a path-connected space $Y^\sigma$, a continuous map $f_\sigma\:Y^\sigma\to Y^\sigma$, and a continuous surjection $\pi^\sigma\:Y^\sigma\to Z^\sigma$ satisfying
$\pi^\sigma\circ f_\sigma=f\circ\pi^\sigma$.
\smallskip\item For every pair of continuous paths $a\subseteq Y^+$ and $b\subseteq Y^-$ the images $\alpha\=\pi^+(a)$ and $\beta\=\pi^-(b)$ do not cross in $S^2$.
\end{enumerate}
The last condition means that there is no closed disk $D\subseteq S^2$ and paths $\alpha,\beta$ as above properly contained in $D$ so that $\alpha\cap\partial D$ and $\beta\cap\partial D$ are disjoint linked $S^0$s in $\partial D$ \cite[Definition~2.34]{CL24}.
\end{definition}

\begin{definition}
For each $\sigma\in\{+,\,-\}$, let $Z_0^\sigma$ be the primary zipper of color $\sigma$ given by~\eqref{defZ0+}. Define
\begin{equation}\label{defZpm}
Z^\sigma\=\bigcup_{m\in\N_0}f^{-m}(Z_0^\sigma).
\end{equation}
\end{definition}

\begin{lemma}\label{nonemptyandpathconnected}
Let $Z^+$ and $Z^-$ be given above. Then each of $Z^+$ and $Z^-$ is nonempty and path-connected.
\end{lemma}
\begin{proof}
First, since each of $Z_0^\pm$ is nonempty, each of $Z^\pm$ is nonempty. Next, for each point $v\in Z^+$, by~\eqref{defZpm}, we know that there exists $n_0\in\N_0$ such that $v\in f^{-n_0}(Z_0^+)$. By~\eqref{defZ0+}, we have that $\post f\subseteq Z_0^+\subseteq Z^+$ and $Z_0^+$ is the union of paths connecting different points in $\post f$. Then $f^{-n_0}(Z_0^+)$ contains all $n_0$-vertices and is a union of paths connecting different $n_0$-vertices. Since $\V^0\subseteq\V^{n_0}$, there is a path in $f^{-n_0}(Z_0^+)\subseteq Z^+$ connecting $v$ and a point in $\post f$, which means that $Z^+$ is path-connected. In a similar way, we can prove that $Z^-$ is also path-connected.
\end{proof}

\begin{lemma}\label{fcompletelyinvariant}
Let $Z^+$ and $Z^-$ be given above. Then each of $Z^+$ and $Z^-$ is completely $f$-invariant, i.e., $f^{-1}(Z^+)=Z^+=f(Z^+)$ and $f^{-1}(Z^-)=Z^-=f(Z^-)$.
\end{lemma}

\begin{proof}
By Lemma~\ref{Z_0invariant} we have 
\begin{equation*}
    f(Z^+)=f\Bigl(\bigcup_{m\in\N_0}f^{-m}(Z_0^+)\Bigr)
    =\bigcup_{m\in\N_0}f^{-m}(Z_0^+)\cup f(Z_0^{+})
    \subseteq \bigcup_{m\in\N_0}f^{-m}(Z_0^+).
\end{equation*}
Then 
$    Z^+
    =\bigcup_{m\in\N_0}f^{-m}(Z_0^+)
    \subseteq f(Z^+)
    \subseteq Z^+$.
It follows that $Z^+=f(Z^+)$.

On the other hand,
\begin{equation*}
    f^{-1}(Z^+)
    =f^{-1}\Bigl(\bigcup_{m\in\N_0}f^{-m}(Z_0^+)\Bigr)
    =\bigcup_{m\in\N}f^{-m}(Z_0^+)
    \subseteq \bigcup_{m\in\N_0}f^{-m}(Z_0^+).
\end{equation*}

Since $f(Z_0^+)\subseteq Z_0^+$, we have $Z_0^+\subseteq f^{-1}(Z_0^+)$. Then 
$    Z^+=\bigcup_{m\in\N_0}f^{-m}(Z_0^+)\subseteq f^{-1}(Z^+)\subseteq Z^+$.
It follows that $Z^+=f^{-1}(Z^+)$.

Thus, $Z^+$ is completely $f$-invariant. Similarly, $Z^-$ is also completely $f$-invariant.
\end{proof}

\begin{lemma}\label{dense}
Let $Z^+$ and $Z^-$ be given above. Then each of $Z^+$ and $Z^-$ is dense in $S^2$.
\end{lemma}
\begin{proof}
We prove the statement for $Z^+$ as the proof is the same for $Z^-$. Fix $p \in Z^+$. Since $Z^+$ is backward invariant, $\bigcup_{m\in\N}f^{-m}(p) \subseteq Z^+$ and therefore $\overline{\bigcup_{m\in\N}f^{-m}(p)} \subseteq \overline{Z^+}$. On the other hand, $\overline{\bigcup_{m\in\N}f^{-m}(p)} = S^2$ by \cite[Lemma~2.12]{Li17}. Therefore $S^2 = \overline{Z^+}$. The conclusion follows.
\end{proof}

\subsection{Construction of the corresponding trees}

To complete the proof of Theorem~\ref{A}, we construct the tree $Y^\sigma=\bigcup_{n=0}^{+\infty}Y_n^\sigma$ as an increasing union of finite trees and equip it with the weak topology. We then construct a continuous surjection $\pi^\sigma\:Y^\sigma\to Z^\sigma$ as required by Definition \ref{def_gen_zipper} and a map $f_\sigma\:Y^\sigma\to Y^\sigma$ satisfying $\pi^\sigma\circ f_\sigma=f\circ\pi^\sigma$.

For each integer $n\ge1$, we glue copies of $Y_0^\sigma$ along $\Gamma_n^\sigma$ to obtain $Y_n^\sigma$. We then construct the embeddings $\iota_n^\sigma$ and projections $\pi_n^\sigma$.

\begin{definition}\label{def:ch4-tile-copy-data}
For an integer $n\ge1$ and each tile $X\in\X_\sigma^n$, define $\theta_X^\sigma\:Y_0^\sigma\to Y_n^\sigma$ by $\theta_X^\sigma(q_\alpha^\sigma(s))\=\eta_n^\sigma((n,X,\alpha),s)$ for $\alpha\in\sA_0^\sigma$ and $s\in I_\alpha^\sigma$, where $\eta_n^\sigma$ is the quotient map of Definition~\ref{c:relation} and $q_\alpha^\sigma$ and $\xi_p^\sigma$ are as in Definition~\ref{def:ch4-initial-quotient-data}. We call $\theta_X^\sigma$ the \emph{tile-copy map}. Define the \emph{tile copy} associated with $X$ by $T_X^\sigma\=\theta_X^\sigma(Y_0^\sigma)$. Define the \emph{block attachment} $\beta_n^\sigma\:\sB_n^\sigma\to Y_n^\sigma$ as follows. For each $b\in\sB_n^\sigma$, let $\beta_n^\sigma(b)\=\theta_X^\sigma\bigl(\xi_{f^n(\prV(b))}^\sigma\bigr)$ for any $X\leftrightarrow b$.
\end{definition}

The following proposition establishes that $\theta_X^\sigma$ and $\beta_n^\sigma$ are well-defined and that $\theta_X^\sigma$ is an embedding.

\begin{prop}\label{prop:compact-finite-trees}
Fix a positive integer $n$ and a color $\sigma$ in $\{+,\,-\}$. Let $f\:S^2\to S^2$ be an expanding Thurston map, and let $\cC$ be an $f$-invariant Jordan curve containing $\post f$. Fix a connection of $1$-tiles satisfying \textup{(C1)}--\textup{(C2)} of Definition~\ref{Propertiesofconnections}. Let $\eta_n^\sigma\:K_n^\sigma\to Y_n^\sigma$ be the quotient map of Definition~\ref{c:relation}. Retain the initial road maps $q_\alpha^\sigma$ and points $\xi_p^\sigma$ of Definition~\ref{def:ch4-initial-quotient-data}. The following hold:
\begin{enumerate}
\smallskip\item The space $Y_n^\sigma$ is a finite tree. The prescriptions of Definition~\ref{def:ch4-tile-copy-data} are well-defined. For each tile $X\in\X_\sigma^n$, the map $\theta_X^\sigma\:Y_0^\sigma\to Y_n^\sigma$ is an embedding. For labels $\alpha\in\sA_0^\sigma$ and $\lambda\in\sA_{n-1}^\sigma$, these maps satisfy
\begin{equation}\label{eq:ch4-stage-model-map}
\begin{aligned}
(\theta_X^\sigma \circ q_\alpha^\sigma)(s)
&=\eta_n^\sigma((n,X,\alpha),s)
&&\quad \text{ for all } s\in I_\alpha^\sigma,\\
\eta_n^\sigma|_J
&=\theta_{X_n^{\lambda,\sigma}(J)}^\sigma
\circ q_{\alpha_J}^\sigma\circ\kappa_J^\sigma
&&\quad \text{ for all } J\in\mathcal P_n^{\lambda,\sigma}.
\end{aligned}
\end{equation}
The first identity determines $\theta_X^\sigma$ uniquely. For each block $b$ in $\sB_n^\sigma$, the point $\beta_n^\sigma(b)$ belongs to $T_X^\sigma$ exactly when $X$ is incident to $b$. The copies $T_X^\sigma$ cover $Y_n^\sigma$. For distinct tiles $X,Y$ in $\X_\sigma^n$, they satisfy
\begin{equation}\label{eq:ch4-stage-copy-intersection}
 T_X^\sigma\cap T_Y^\sigma
 =\{\beta_n^\sigma(b):b\in\sB_n^\sigma,\ X\leftrightarrow b\leftrightarrow Y\}.
\end{equation}
\smallskip\item For every label $\lambda$ in $\sA_n^\sigma$, the map $\eta_n^\sigma|_{I_\lambda^\sigma}$ is an embedding.
\end{enumerate}
\end{prop}
\begin{proof}
\textup{(i)} For each block $b$ in $\sB_n^\sigma$, let $p_b$ denote $f^n(\prV(b))$. Let $\equiv$ be the equivalence relation generated by identifying $\bigl(X,\xi_{p_b}^\sigma\bigr)$ with $\bigl(Y,\xi_{p_b}^\sigma\bigr)$ whenever both tiles are incident to $b$. Define
\begin{equation}\label{eq:compact-attachment}
T\=\bigsqcup_{X\in\X_\sigma^n}(\{X\}\times Y_0^\sigma)/{\equiv}
\quad \text{ and }\quad 
\widehat\theta_X(z)\=[(X,z)].
\end{equation}
For each tile $X$, call $\widehat\theta_X(Y_0^\sigma)\subseteq T$ the copy labeled by $X$. For each block $b$, choose a tile $X$ with $X\leftrightarrow b$ and define $\widehat\beta(b)\=\widehat\theta_X\bigl(\xi_{p_b}^\sigma\bigr)$. If $X\leftrightarrow b\leftrightarrow Y$, then $\bigl(X,\xi_{p_b}^\sigma\bigr)\equiv\bigl(Y,\xi_{p_b}^\sigma\bigr)$, so $\widehat\beta(b)$ is independent of the choice of $X$. By~\eqref{eq:ch4-incident-block-uniqueness} and Proposition~\ref{c:initial}, distinct blocks incident to $X$ determine distinct points $\xi_{p_b}^\sigma$ in $Y_0^\sigma$. Thus, the copies labeled by distinct tiles $X,\,Y$ intersect precisely at the points $\widehat\beta(b)$ with $X\leftrightarrow b\leftrightarrow Y$. Fix a tile $X_*$ and add the copies in order of increasing distance from $c(X_*)$ in $\Gamma_n^\sigma$. For $X\ne X_*$, the first block $b$ on the path from $c(X)$ to $c(X_*)$ is the only block that $X$ has in common with any of the preceding tiles. Hence its copy intersects the preceding union at the single point $\widehat\beta(b)$. It follows that $T$ is a finite tree and each $\widehat\theta_X$ is an embedding.

Define $\jmath_X\:K_0^\sigma\to K_n^\sigma$ by $\jmath_X(\alpha,s)\=((n,X,\alpha),s)$. Define $\eta\:K_n^\sigma\to T$ by
\begin{equation*}
\eta\circ \jmath_X\=\widehat\theta_X\circ \eta_0^\sigma,\quad\eta|_J\=\widehat\theta_{X_n^{\lambda,\sigma}(J)}\circ q_{\alpha_J}^\sigma\circ\kappa_J^\sigma\quad \text{ for all } \lambda\in\sA_{n-1}^\sigma,\ J\in\mathcal P_n^{\lambda,\sigma}.
\end{equation*}
For $\lambda\in\sA_{n-1}^\sigma$, the formulas on two adjacent level-$n$ atoms of $I_\lambda^\sigma$ both equal $\widehat\beta\bigl(b_n^{\lambda,\sigma}(s)\bigr)$ at their common node $s$. Thus, $\eta$ is continuous. The images of $\eta\circ \jmath_X$ cover $T$ since $\eta_0^\sigma$ is surjective. By~\eqref{eq:road-data}, the closed equivalence relation $(\jmath_X\times\jmath_X)^{-1}(\sim_n^\sigma)$ contains $\bigcup_{r=1}^{+\infty}\mathcal R_{0,r}^{H,\sigma}$. Therefore, Definition~\ref{c:relation} implies
$(\jmath_X\times\jmath_X)(\sim_0^\sigma)\subseteq\sim_n^\sigma.$
For $\lambda\in\sA_{n-1}^\sigma$ and $J\in\mathcal P_n^{\lambda,\sigma}$,~\eqref{eq:ch4-affine-refinement} and Lemma~\ref{lem:ch4-path-matchings} imply $x\sim_n^\sigma\bigl(\jmath_{X_n^{\lambda,\sigma}(J)}\circ\kappa_J^\sigma\bigr)(x)$ for all $x\in J$.

We first show that $\sim_n^\sigma\subseteq\ker\eta$. Take a positive integer $r$ and a pair $(x,y)$ in $\mathcal R_{n,n+r}^{H,\sigma}$. By~\eqref{eq:road-data} and~\eqref{eq:ch4-affine-refinement}, write $\eta(x)\=\widehat\theta_X
(q_\alpha^\sigma(s))$, $\eta(y)=\widehat\theta_Y(q_\tau^\sigma(t))$, and $L(X,r,b_r^{\alpha,\sigma}(s))=L(Y,r,b_r^{\tau,\sigma}(t))$,
where $X,\,Y$ lie in $\X_\sigma^n$ and $s,\,t$ are level-$r$ nodes on the initial roads $I_\alpha^\sigma,\, I_\tau^\sigma$. If the tiles coincide, injectivity of $L(X,r,\,\cdot\,)$ places $((\alpha,s),(\tau,t))$ in $\mathcal R_{0,r}^{H,\sigma}$. Otherwise, Proposition~\ref{lem:ch4-sector-copies}~\textup{(i)},~\eqref{eq:ch4-copy-values} and~\eqref{eq:ps-generator} identify $q_\alpha^\sigma(s)$ and $q_\tau^\sigma(t)$ with $\xi_{p_b}^\sigma$ for a block $b$ incident to both tiles. Thus, $\eta(x)$ and $\eta(y)$ agree. Since $\ker\eta$ is closed, Definition~\ref{c:relation} implies $\sim_n^\sigma\subseteq\ker\eta$.

Conversely, let $x,\,y\in K_n^\sigma$ satisfy $\eta(x)=\eta(y)$. The atom relations give representatives $\jmath_X(u),\, \jmath_Y(v)$, with $u,\,v$ in $K_0^\sigma$. If the tiles coincide, injectivity of $\widehat\theta_X$ identifies $\eta_0^\sigma(u)$ with $\eta_0^\sigma(v)$, so the lifted initial relation identifies $x$ with $y$. Otherwise, $\eta_0^\sigma(u)$ and $\eta_0^\sigma(v)$ both equal $\xi_{p_b}^\sigma$ for a block $b$ incident to both tiles. Equation~\eqref{eq:ch4-copy-values} places $(\jmath_X(e_{p_b}),\jmath_Y(e_{p_b}))$ in $\mathcal R_{n,n+1}^{H,\sigma}$, so $x\sim_n^\sigma\jmath_X(e_{p_b})\sim_n^\sigma\jmath_Y(e_{p_b})\sim_n^\sigma y$.
Thus, $\ker\eta=\sim_n^\sigma$.

Since $Y_n^\sigma$ is compact and $T$ is Hausdorff, the continuous bijection $\Phi\:Y_n^\sigma\to T$ induced by $\eta$ is a homeomorphism. The prescribed maps satisfy $\theta_X^\sigma=\Phi^{-1}\circ \widehat\theta_X$ and $\beta_n^\sigma=\Phi^{-1}\circ \widehat\beta$. Taking inverse images under $\Phi$ gives $Y_n^\sigma=\bigcup_{X\in\X_\sigma^n}T_X^\sigma$, the criterion $\beta_n^\sigma(b)\in T_X^\sigma\iff X\leftrightarrow b$, and~\eqref{eq:ch4-stage-copy-intersection}. Since $\eta_0^\sigma$ is onto, the identity $\theta_X^\sigma\circ \eta_0^\sigma=\eta_n^\sigma\circ \jmath_X$ determines $\theta_X^\sigma$ uniquely.

\smallskip
\textup{(ii)} Since $q_\alpha^\sigma$ and $\theta_X^\sigma$ are embeddings,~\eqref{eq:ch4-stage-model-map} shows that $\eta_n^\sigma$ embeds $I_{(n,X,\alpha)}^\sigma$. For $\lambda$ in $\sA_{n-1}^\sigma$, list the intervals of $\mathcal P_n^{\lambda,\sigma}$ in order as $J_1,\,J_2,\,\ldots,\, J_k$. By~\eqref{eq:ch4-stage-model-map}, each restriction $\eta_n^\sigma|_{J_i}$ is an embedding. Since $P_n^\sigma(\lambda)$ is simple,~\eqref{eq:ch4-stage-copy-intersection} gives $\eta_n^\sigma(J_i)\cap\eta_n^\sigma(J_j)
 =\eta_n^\sigma(J_i\cap J_j)\quad(1\le i<j\le k)$. Hence $\eta_n^\sigma|_{I_\lambda^\sigma}$ is continuous and injective, and therefore embeds the compact interval $I_\lambda^\sigma$.
\end{proof}

\begin{definition}\label{def:ch4-refined-copies}
Set $T_{X_\sigma^0}^{\sigma}\= Y_0^{\sigma}$, $\theta_{X_\sigma^0}^{\sigma}\= \mathrm{id}_{Y_0^{\sigma}}$, $\beta_0^{\sigma}(\widehat p_0^{\sigma})\= \xi_p^{\sigma}$ and $\Theta_{X_\sigma^0}^{\sigma}\= \mathrm{id}_{Y_1^{\sigma}}$. For each integer $n\ge1$ and $X\in\X_\sigma^n$, put
\begin{equation*}
 G_X^{\sigma}\= \bigcup_{\substack{V\in\X_\sigma^{n+1}\\A_{n,n+1}(V)=X}}T_V^{\sigma}\subseteq Y_{n+1}^{\sigma}.
\end{equation*}

For $n\in \N_0$, define the \emph{stage map} $\iota_n^\sigma\:Y_n^\sigma\to Y_{n+1}^\sigma$ by $\iota_n^\sigma\circ \eta_n^\sigma\=\eta_{n+1}^\sigma|_{K_n^\sigma}$ (see Proposition~\ref{prop:compact-stage-stability}~\textup{(i)} for well-definedness). For $n\in \N$, define $\Theta_X^\sigma\:Y_1^\sigma\to G_X^\sigma$ by $\Theta_X^\sigma\circ \theta_U^\sigma=\theta_{L(X,1,U)}^\sigma$ for $U\in\X_\sigma^1$.
\end{definition}

The following proposition shows that $\iota_n^\sigma$ is an embedding and $\Theta_X^\sigma$ is a homeomorphism.

\begin{prop}\label{prop:compact-stage-stability}
Let $f\:S^2\to S^2$ be an expanding Thurston map, and let $\cC\supseteq\post f$ be an $f$-invariant Jordan curve. Fix a connection of $1$-tiles satisfying \textup{(C1)}--\textup{(C2)} of Definition~\ref{Propertiesofconnections} and $\sigma\in\{+,\,-\}$. Use the finite trees $Y_n^\sigma\=K_n^\sigma/\sim_n^\sigma$ and copies $T_X^\sigma=\theta_X^\sigma(Y_0^\sigma)$ of Proposition~\ref{prop:compact-finite-trees}. Let $\iota_n^\sigma\:Y_n^\sigma\to Y_{n+1}^\sigma$, $n\in \N_0$, be prescribed by Definition~\ref{def:ch4-refined-copies}. For $n\in \N$ and $X\in\X_\sigma^n$, recall $\Theta_X^\sigma\:Y_1^\sigma\to G_X^\sigma$ from Definition~\ref{def:ch4-refined-copies}. Here, $G_X^\sigma$ is the union of $T_V^\sigma$ with $A_{n,n+1}(V)=X$. Then, the following hold:

\begin{enumerate}
\smallskip
\item The inclusion $K_n^\sigma\subseteq K_{n+1}^\sigma$ induces an embedding $\iota_n^\sigma\:Y_n^\sigma\hookrightarrow Y_{n+1}^\sigma$ for $n\in \N_0$, with
\begin{equation}\label{eq:compact-stage-stability}
 \iota_n^{\sigma}\circ \eta_n^{\sigma}=\eta_{n+1}^{\sigma}|_{K_n^{\sigma}}\quad\text{ and }\quad \sim_{n+1}^{\sigma}|_{K_n^{\sigma}\times K_n^{\sigma}}=\sim_n^{\sigma}.
\end{equation}
For $n\in \N$ and $X\in\X_\sigma^n$, the map $\Theta_X^\sigma\:Y_1^\sigma\to G_X^\sigma$ is a homeomorphism satisfying
\begin{equation}\label{eq:ch4-refined-copy}
(\Theta_X^{\sigma}\circ \theta_U^{\sigma})(z)=\theta_{L(X,1,U)}^{\sigma}(z) \quad \text{ for all } U\in\X_\sigma^1,\ z\in Y_0^{\sigma}.
\end{equation}
Moreover,
\begin{equation}\label{eq:ch4-copy-stability}
 \iota_n^{\sigma}\circ \theta_X^{\sigma}=\Theta_X^{\sigma}\circ \iota_0^{\sigma} \quad \text{ for all } n\ge0,\ X\in\X_\sigma^n.
\end{equation}
This properties gives the well-definedness of $\iota_n$ from Definition~\ref{def:ch4-refined-copies}.
\smallskip
\item For $m\in\N$ and $X\in\X_\sigma^m$, put $U\=f^{m-1}(X)\in\X_\sigma^1$. Then
\begin{equation}\label{eq:ch4-stage-sections}
\iota_{m-1}^\sigma(Y_{m-1}^\sigma)\cap T_X^\sigma=\bigl(\theta_X^\sigma\circ(\theta_U^\sigma)^{-1}\bigr)(\iota_0^\sigma(Y_0^\sigma)\cap T_U^\sigma).
\end{equation}
\end{enumerate}
\end{prop}

\begin{proof}
\textup{(i)} Since $Y_1^\sigma$ is Hausdorff, the continuous restriction $\eta_1^\sigma|_{K_0^\sigma}$ has a closed kernel. For every $m\in\N$,~\eqref{eq:node-realization} implies $\mathcal R_{0,m}^{H,\sigma}\subseteq\mathcal R_{1,m+1}^{H,\sigma}$, so this kernel contains $\bigcup_{m\ge1}\mathcal R_{0,m}^{H,\sigma}$. Thus, $\eta_1^\sigma|_{K_0^\sigma}$ descends to a continuous map $\iota_0^\sigma$. Any two points of $Y_0^\sigma$ lie in some $q_\alpha^\sigma(I_\alpha^\sigma)$ by Proposition~\ref{c:initial}. The restriction $\eta_1^\sigma|_{I_\alpha^\sigma}$ is injective by Proposition~\ref{prop:compact-finite-trees}~\textup{(ii)}. Hence, $\iota_0^\sigma$ is injective and embeds the compact tree $Y_0^\sigma$.

The map $L(X,1,\,\cdot\,)$ preserves incidences by Proposition~\ref{lem:ch4-sector-copies}~\textup{(i)} and attachment marks by~\eqref{eq:ch4-copy-values}. Thus,~\eqref{eq:ch4-refined-copy} defines a continuous bijection $\Theta_X^\sigma\:Y_1^\sigma\to G_X^\sigma$. Since $Y_1^\sigma$ is compact and $G_X^\sigma$ is Hausdorff, $\Theta_X^\sigma$ is a homeomorphism. For $X\ne Y$, Proposition~\ref{lem:ch4-sector-copies}~\textup{(i)} and~\eqref{eq:ch4-stage-copy-intersection} give
\begin{equation}\label{eq:ch4-refined-copy-intersection}
 G_X^\sigma\cap G_Y^\sigma =\{\beta_{n+1}^\sigma(\operatorname{Lift}_n^\sigma b): X\leftrightarrow b\leftrightarrow Y\}.
\end{equation}

For $X\leftrightarrow b$ and $p=f^n\prV(b)$,
\begin{equation}\label{eq:ch4-refined-mark}
 (\Theta_X^\sigma\circ\iota_0^\sigma)(\xi_p^\sigma)=(\Theta_X^\sigma\circ\beta_1^\sigma)(\widehat p_1^\sigma) =\beta_{n+1}^\sigma(L(X,1,\widehat p_1^\sigma)) =\beta_{n+1}^\sigma(\operatorname{Lift}_n^\sigma b).
\end{equation}
On $T_X^\sigma$, define $\iota_n^\sigma\=\Theta_X^\sigma\circ\iota_0^\sigma\circ (\theta_X^\sigma)^{-1}$ to satisfy~\eqref{eq:ch4-copy-stability}. For tiles $X,\,Y$ incident to a block $b$,~\eqref{eq:ch4-refined-mark} shows that the formulas on $T_X^\sigma$ and $T_Y^\sigma$ both send $\beta_n^\sigma(b)$ to $\beta_{n+1}^\sigma(\operatorname{Lift}_n^\sigma b)$. Thus, $\iota_n^\sigma$ is well-defined. By~\eqref{eq:ch4-refined-copy-intersection}, $\iota_n^\sigma(x)=\iota_n^\sigma(y)$ forces $x=y$. The map $\iota_n^\sigma$ is continuous on each $T_X^\sigma$, hence on $Y_n^\sigma$. Since $Y_n^\sigma$ is compact and $Y_{n+1}^\sigma$ is Hausdorff, $\iota_n^\sigma$ is an embedding.

Fix $\lambda\in\sA_{n-1}^\sigma$ and $s\in K\in\mathcal P_{n+1}^{\lambda,\sigma}$. By~\eqref{eq:ch4-affine-refinement} and~\eqref{eq:ch4-refined-copy}, both $\iota_n^\sigma\circ\eta_n^\sigma(\lambda,s)$ and $\eta_{n+1}^\sigma(\lambda,s)$ equal $\theta_{X_{n+1}^{\lambda,\sigma}(K)}^\sigma\circ q_{\alpha_K}^\sigma\circ\kappa_K^\sigma(s)$.

For $\lambda=(n,X,\alpha)$, let $D\in\mathcal P_1^{\alpha,\sigma}$ contain $(\alpha,s)$ and let $K$ be its copy in~\eqref{eq:road-data}. By~\eqref{eq:ch4-copy-values}, the endpoint marks agree, so $\alpha_K=\alpha_D$ and $\kappa_K^\sigma(\lambda,s)=\kappa_D^\sigma(\alpha,s)$. Equations~\eqref{eq:ch4-stage-model-map} and~\eqref{eq:ch4-refined-copy} give $\eta_{n+1}^\sigma(\lambda,s)=\Theta_X^\sigma(\eta_1^\sigma(\alpha,s))$. By~\eqref{eq:ch4-copy-stability} and $\eta_1^\sigma|_{K_0^\sigma}=\iota_0^\sigma\circ\eta_0^\sigma$, this equals $\iota_n^\sigma(\eta_n^\sigma(\lambda,s))$.

Thus $\iota_n^\sigma\circ\eta_n^\sigma=\eta_{n+1}^\sigma|_{K_n^\sigma}$. Taking kernels gives~\eqref{eq:compact-stage-stability}, since $\iota_n^\sigma$ is injective.

\smallskip

\textup{(ii)} For $m=1$, both sides of~\eqref{eq:ch4-stage-sections} equal $\iota_0^\sigma(Y_0^\sigma)\cap T_X^\sigma$. For $m\ge2$, put $W\=A_{m-1,m}(X)$ and $U\=f^{m-1}(X)$. Then $X=L(W,1,U)$ and $T_X^\sigma=\Theta_W^\sigma(T_U^\sigma)$. By~\eqref{eq:ch4-copy-stability}, $\iota_{m-1}^\sigma(Y_{m-1}^\sigma)=\bigcup_V(\Theta_V^\sigma\circ\iota_0^\sigma)(Y_0^\sigma)$, where $V\in\X_\sigma^{m-1}$. For $V\ne W$,~\eqref{eq:ch4-refined-copy-intersection} and~\eqref{eq:ch4-refined-mark} imply $(\Theta_V^\sigma\circ\iota_0^\sigma)(Y_0^\sigma)\cap G_W^\sigma\subseteq(\Theta_W^\sigma\circ\iota_0^\sigma)(Y_0^\sigma)$. Hence $\iota_{m-1}^\sigma(Y_{m-1}^\sigma)\cap T_X^\sigma=\Theta_W^\sigma(\iota_0^\sigma(Y_0^\sigma)\cap T_U^\sigma)$. Substituting $\Theta_W^\sigma|_{T_U^\sigma}=\theta_X^\sigma\circ(\theta_U^\sigma)^{-1}$ from~\eqref{eq:ch4-refined-copy} proves~\eqref{eq:ch4-stage-sections}.
\end{proof}

\begin{figure}[htbp]
    \centering
    \begin{overpic}[
        percent,
        width=0.9\linewidth,
        trim=-8bp -22bp -36bp -9bp
    ]{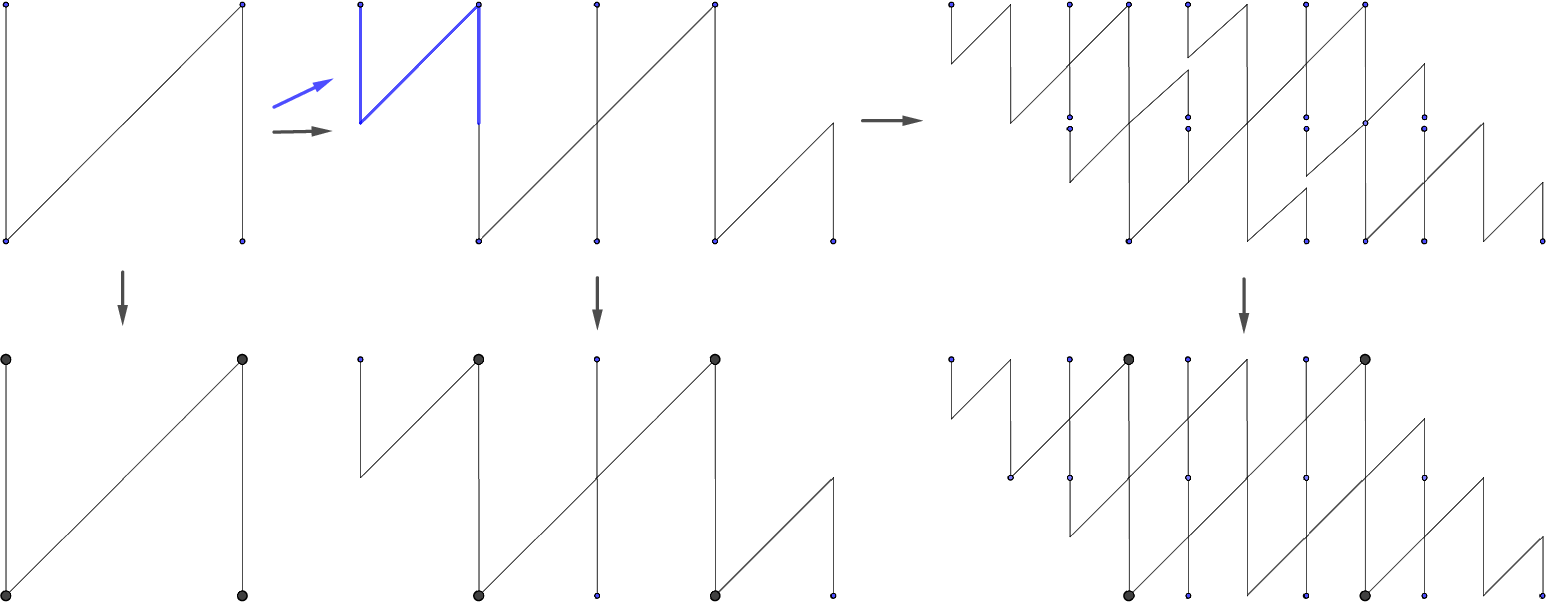}
        \setlength{\unitlength}{0.009\linewidth}
        \scriptsize

        \put(3.43,38.85){\makebox(0,0){$\xi_{-1}^+$}}
        \put(17.35,38.85){\makebox(0,0){$\xi_{\infty}^+$}}
        \put(2.54,23.09){\makebox(0,0){$\xi_{0}^+$}}
        \put(16.76,23.09){\makebox(0,0){$\xi_{1}^+$}}
        \put(18.38,26.11){\makebox(0,0){$Y_0^+$}}

        \put(53.81,25.30){\makebox(0,0){$Y_1^+$}}
        \put(97.50,25.15){\makebox(0,0){$Y_2^+$}}

        \put(19.44,35.37){\makebox(0,0){$\theta_X^+$}}
        \put(19.41,30.33){\makebox(0,0){$\iota_0^+$}}
        \put(55.51,30.92){\makebox(0,0){$\iota_1^+$}}

        \put(10.35,21.47){\makebox(0,0){$\pi_0^+$}}
        \put(39.60,20.95){\makebox(0,0){$\pi_1^+$}}
        \put(79.23,20.80){\makebox(0,0){$\pi_2^+$}}

        \put(2.54,16.46){\makebox(0,0){$-1$}}
        \put(16.91,16.46){\makebox(0,0){$\infty$}}
        \put(2.54,1.87){\makebox(0,0){$0$}}
        \put(16.32,1.65){\makebox(0,0){$1$}}
        \put(18.09,4.01){\makebox(0,0){$Z_0^+$}}

        \put(31.64,16.09){\makebox(0,0){$-1$}}
        \put(46.30,16.16){\makebox(0,0){$\infty$}}
        \put(31.64,2.17){\makebox(0,0){$0$}}
        \put(45.86,1.72){\makebox(0,0){$1$}}
        \put(53.81,3.79){\makebox(0,0){$Z_1^+$}}

        \put(71.13,16.68){\makebox(0,0){$-1$}}
        \put(85.78,16.46){\makebox(0,0){$\infty$}}
        \put(71.20,1.95){\makebox(0,0){$0$}}
        \put(85.78,1.95){\makebox(0,0){$1$}}
        \put(97.06,3.20){\makebox(0,0){$Z_2^+$}}

    \end{overpic}
    \caption{The finite trees ${Y_0}^+$, ${Y_1}^+$ and ${Y_2}^+$ for the $2\times2$ Latt\`es map $g$ in Example~\ref{ex:lattes-map}. The blue part of the upper middle figure is $\theta_X^+(Y_0^+)$.}
    \label{fig:finite-trees}
\end{figure}

We have identified $Y_n^\sigma$ with $\iota_n^\sigma(Y_n^\sigma)$. For the Latt\`es map $g$ in Example~\ref{ex:lattes-map}, the finite trees ${Y_0}^+$, ${Y_1}^+$, and ${Y_2}^+$ are shown in Figure~\ref{fig:finite-trees}. Recall the related N-shaped figure from Figure~\ref{fig:Fm}.

We next construct $\pi_n^\sigma$ from the road limits in $S^2$.

\begin{definition}\label{def:ch4-spherical-data}
Let $F^{\lambda,\sigma}$ denote the road limit from Lemma~\ref{lem:ch3-road-limits}, and let $Z_0^\sigma$ be the primary zipper from Subsection~\ref{Construction_of_primary_zippers}. For each integer $n\ge1$, define $Z_n^\sigma\=\bigcup_{k=0}^nf^{-k}(Z_0^\sigma)$. Equip $Z_n^\sigma$ with the subspace topology from $S^2$ for each $n\in\N_0$.

For each integer $n\ge0$, we call the map $\pi_n^\sigma\:Y_n^\sigma\to Z_n^\sigma$ given by $(\pi_n^\sigma\circ\eta_n^\sigma)(\lambda,s)\=F^{\lambda,\sigma}(s)$ the \emph{stage projection} (see~\ref{prop:compact-projections} for its well-definedness).
\end{definition}
The following proposition shows that this rule defines a continuous surjection $\pi_n^\sigma\:Y_n^\sigma\to Z_n^\sigma$.

\begin{prop}\label{prop:compact-projections}
Let $f\:S^2\to S^2$ be an expanding Thurston map, and let $\cC\supseteq\post f$ be an $f$-invariant Jordan curve. Fix a connection of $1$-tiles satisfying \textup{(C1)}--\textup{(C2)} of Definition~\ref{Propertiesofconnections} and $\sigma\in\{+,\,-\}$. Let $F_m^{\lambda,\sigma}\:I_\lambda^\sigma\to S^2$ be as in Definition~\ref{def:ch4-finite-realizations}. Use the quotients $Y_n^\sigma\=K_n^\sigma/\sim_n^\sigma$ of Definition~\ref{c:relation}, the embeddings $\iota_n^\sigma\:Y_n^\sigma\to Y_{n+1}^\sigma$ of Proposition~\ref{prop:compact-stage-stability}, and the sets $Z_n^\sigma=\bigcup_{k=0}^nf^{-k}(Z_0^\sigma)$ of Definition~\ref{def:ch4-spherical-data}. Then the following hold:
\begin{enumerate}
\smallskip\item For each road label $\lambda$, the maps $F_m^{\lambda,\sigma}$, $m>\operatorname{birth}(\lambda)$, are continuous and converge uniformly to a continuous map $F^{\lambda,\sigma}\:I_\lambda^\sigma\to S^2$.
\smallskip\item For $k\in \N$, $X\in\X_\sigma^k$ and $\alpha\in\sA_0^\sigma$,
\begin{equation}\label{eq:ch4-born-semiconjugacy}
 f^k\circ F^{(k,X,\alpha),\sigma}=F^{\alpha,\sigma}.
\end{equation}
\item For $n\in \N_0$, there is a continuous surjection $\pi_n^\sigma\:Y_n^\sigma\to Z_n^\sigma$ satisfying
\begin{equation}\label{eq:ch4-stage-projections}
\pi_n^{\sigma}\circ \eta_n^{\sigma}|_{I_\lambda^{\sigma}}=F^{\lambda,\sigma} \quad \text{ for all } \lambda\in\sA_n^{\sigma},\quad \pi_{n+1}^{\sigma}\circ \iota_n^{\sigma}=\pi_n^{\sigma}.
\end{equation}
\end{enumerate}
\end{prop}

\begin{proof}
\textup{(i)} The continuity and uniform convergence follow from Lemma~\ref{lem:ch3-road-limits}.

\textup{(ii)} For $k=1$,~\eqref{eq:compact-born-road} gives $f\circ F_{1+r}^{(1,X,\alpha),\sigma}=F_r^{\alpha,\sigma}$. For integer $k\ge2$, apply $f$ to~\eqref{eq:compact-born-road}. Use $f\circ\bigl(f^k|_V\bigr)^{-1}=\bigl(f^{k-1}|_{f(V)}\bigr)^{-1}$ and $f(L(X,r,U))=L(f(X),r,U)$. Since $f$ is continuous, taking limits gives
\begin{equation}\label{eq:ch4-road-one-step}
 f\circ F^{(k,X,\alpha),\sigma} =\begin{cases}
 F^{\alpha,\sigma} & \text{ if } k=1,\\
 F^{(k-1,f(X),\alpha),\sigma}&\text{ if } k\ge2.
 \end{cases}
\end{equation}
Iterating $k$ times proves~\eqref{eq:ch4-born-semiconjugacy}.

\textup{(iii)} By~\eqref{eq:node-realization}, the continuous rule $(\lambda,s)\mapsto F^{\lambda,\sigma}(s)$ on $K_n^\sigma$ identifies the pairs in~\eqref{eq:ps-generator}. The kernel is closed since $S^2$ is Hausdorff. It contains $\sim_n^\sigma$ by Definition~\ref{c:relation}. Hence, $F^{\lambda,\sigma}=\pi_n^\sigma\circ \eta_n^\sigma|_{I_\lambda^\sigma}$ for a continuous map $\pi_n^\sigma\:Y_n^\sigma\to S^2$.

Fix integers $k,\,r\ge1$ and an initial road $\alpha$. Fix $J\in\mathcal P_r^{\alpha,\sigma}$ and $t\in\inte_\circ J$, and put $U\=X_r^{\alpha,\sigma}(J)$. The tiles mapped onto $U$ by $f^k$ are $L(X,r,U)$, $X\in\X_\sigma^k$. Therefore, their inverse branches give
\begin{equation}\label{eq:mathematical-finite-cover}
 \bigl(f^k\bigr)^{-1}(F_r^{\alpha,\sigma}(t)) =\bigl\{F_{k+r}^{(k,X,\alpha),\sigma}(t):X\in\X_\sigma^k\bigr\}.
\end{equation}
To include the nodes, put $g\=f^k$ and let $y_j\to y$. Because $g$ is continuous, every subsequential limit of $x_j\in g^{-1}(y_j)$ lies in $g^{-1}(y)$. Since $S^2$ is compact, we have 
\begin{equation*}
\lim_{j\to+\infty}\sup_{x\in g^{-1}(y_j)}\rho\bigl(x,g^{-1}(y)\bigr)=0.
\end{equation*}
The map $g$ is open and $g^{-1}(y)$ is finite. Thus, $\bigcap_{x\in g^{-1}(y)}g(B_\rho(x,\varepsilon))$ is a neighborhood of $y$ for every $\varepsilon>0$. Hence, $g^{-1}(y_j)$ meets each $B_\rho(x,\varepsilon)$ for all large $j$. Consequently, $g^{-1}(y_j)$ converges to $g^{-1}(y)$ in Hausdorff distance, also over critical values.

The set $\X_\sigma^k$ is finite. By part~\textup{(i)} and the preceding argument for the fibers of $f^k$,~\eqref{eq:mathematical-finite-cover} extends by continuity from the dense set $\bigcup_{J\in\mathcal P_r^{\alpha,\sigma}}\inte_\circ J$ to every $t\in I_\alpha^\sigma$. Passing to the uniform limits of $F_r^{\alpha,\sigma}$ and $F_{k+r}^{(k,X,\alpha),\sigma}$ in~\eqref{eq:mathematical-finite-cover} and taking unions over $t$, we obtain $\bigcup_{X\in\X_\sigma^k}F^{(k,X,\alpha),\sigma}(I_{(k,X,\alpha)}^\sigma)=f^{-k}(F^{\alpha,\sigma}(I_\alpha^\sigma))$. By~\eqref{defZ0+}, the union of the initial road images is $Z_0^\sigma$. For each integer $1\le k\le n$, taking the union over $\alpha\in\sA_0^\sigma$ gives $f^{-k}(Z_0^\sigma)$ as the union of the images of the roads with birth level $k$. Taking the union over all road labels gives $\pi_n^\sigma(Y_n^\sigma)=\bigcup_{k=0}^nf^{-k}(Z_0^\sigma)=Z_n^\sigma$.

Take $x\in Y_n^\sigma$ be arbitrary, and choose $(\lambda,t)\in K_n^\sigma$ with $\eta_n^\sigma(\lambda,t)=x$. The definition gives $(\pi_{n+1}^\sigma\circ\iota_n^\sigma)(x)=(\pi_{n+1}^\sigma\circ\eta_{n+1}^\sigma)(\lambda,t)=F^{\lambda,\sigma}(t)=\pi_n^\sigma(x)$. This proves $\pi_{n+1}^\sigma\circ\iota_n^\sigma=\pi_n^\sigma$.
\end{proof}

We have shown that the maps $\pi_n^\sigma$ agree on the embedded trees $Y_n^\sigma$. We next extend them to $Y^\sigma=\bigcup_{n=0}^{+\infty} Y_n^\sigma$ and equip it with the weak topology.

\subsection{The weak tree and its dynamics}\label{sec:ch4-limit}

A point $p$ of a topological tree $T$ is called a \emph{cut point} if $T\smallsetminus\{p\}$ is disconnected.

Let $T$ be a topological $\mathbb R$-tree, and let $S\subseteq T$ be a nonempty finite subtree. For each $x\in T$, there is a unique point $y\in S$ satisfying $[x,y]_T\cap S=\{y\}.$
We call the map $r_{S,T}\:T\to S$ given by $r_{S,T}(x)\=y$ the \emph{first-point retraction} onto $S$. For each $x\in S$, we have $r_{S,T}(x)=x$.

Next, we extend $r_{S,Y_n^\sigma}$, $\pi_n^\sigma$ and $\phi_n^\sigma$ to $Y^\sigma$. To this end, we first show that every compact subset lies in some finite level $Y_n^\sigma$. We then verify the arcs in $Y^\sigma$ and prove that $f_\sigma$ is proper.

\begin{definition}[The weak tree and its maps]\label{def:ch4-limit-spaces}\label{def:ch4-cut-sets}\label{def:ch4-tree-maps}
Identify each point $x\in Y_n^\sigma$ with $\iota_n^\sigma(x)\in Y_{n+1}^\sigma$ using the embeddings in Proposition~\ref{prop:compact-stage-stability}. Set
\begin{equation*}
    Y^\sigma\=\bigcup_{n\in \N_0}Y_n^\sigma,
\end{equation*}
and equip $Y^\sigma$ with the weak topology
\begin{equation}\label{eq:pi-weak-union}
U\subseteq Y^\sigma\text{ is open}\iff U\cap Y_n^\sigma\text{ is open in }Y_n^\sigma\text{ for every }n\in \N_0.
\end{equation}
Equip $Z^\sigma=\bigcup_{n\in \N_0}Z_n^\sigma\subseteq S^2$ with the subspace topology.

For $x,\,y\in Y_n^\sigma$, let $[x,y]$ denote the unique arc joining them in $Y_n^\sigma\subseteq Y^\sigma$, with $[x,x]=\{x\}$. Since $\iota_n^\sigma$ preserve arcs, $[x,y]$ is independent of $n$. For distinct $x,z\in Y^\sigma$, set $D_z^\sigma(x)\=\{y\in Y^\sigma\smallsetminus\{z\}:z\notin[x,y]\}$.

For a nonempty finite subtree $S\subseteq Y^\sigma$, define $r_S^\sigma\:Y^\sigma\to S$ such that $r_S^\sigma(x)$ is the unique point $a\in S$ with $[x,a]\cap S=\{a\}$. Set $r_n^\sigma\=r_{Y_n^\sigma}^\sigma$.

Define $\pi^\sigma\:Y^\sigma\to Z^\sigma$ by $\pi^\sigma|_{Y_n^\sigma}=\pi_n^\sigma$. For $n\in \N$, define $\phi_n^\sigma\:Y_n^\sigma\to Y_{n-1}^\sigma$ and $f_\sigma\:Y^\sigma\to Y^\sigma$ by
\begin{equation}\label{eq:pi-local-dynamics}
\phi_n^\sigma|_{T_X^\sigma}\=\theta_{f(X)}^\sigma\circ(\theta_X^\sigma)^{-1}\quad \text{ for all } X\in\X_\sigma^n,\quad
f_\sigma|_{Y_n^\sigma}=\phi_n^\sigma.
\end{equation}
\end{definition}

Proposition~\ref{prop:pi-weak-tree} shows that $\pi^\sigma$ and $f_\sigma$ are well-defined.

\begin{prop}
\label{prop:pi-weak-tree}
Let $f\:S^2\to S^2$ be an expanding Thurston map, and let $\cC\supseteq\post f$ be an $f$-invariant Jordan curve. Fix a connection of $1$-tiles satisfying \textup{(C1)}--\textup{(C2)} of Definition~\ref{Propertiesofconnections} and $\sigma\in\{+,\,-\}$. Equip $Y^\sigma=\bigcup_{n\in \N_0}Y_n^\sigma$ with the weak topology defined in~\eqref{eq:pi-weak-union}. Recall the maps $\iota_n^\sigma$ and $\Theta_X^\sigma$ from Proposition~\ref{prop:compact-stage-stability}. Let $\theta_X^\sigma\:Y_0^\sigma\to T_X^\sigma$ and $\pi_n^\sigma\:Y_n^\sigma\to Z_n^\sigma$ be as in Definitions~\ref{def:ch4-tile-copy-data} and~\ref{def:ch4-spherical-data}. Recall the definition of the maps $\pi^\sigma$ and $f_\sigma$ from Definition~\ref{def:ch4-tree-maps}. The following hold:
\begin{enumerate}
\smallskip\item The space $Y^\sigma$ is Hausdorff and is a topological $\mathbb R$-tree. Each $Y_n^\sigma$ embeds as a subtree of $Y^\sigma$. Every compact subset $K\subseteq Y^\sigma$ is contained in some $Y_n^\sigma$. For $x,y\in Y_n^\sigma$, $[x,y]_{Y^\sigma}=[x,y]_{Y_n^\sigma}$.
\smallskip\item For every finite subtree $\emptyset\ne S\subseteq Y^\sigma$, $r_S^\sigma\:Y^\sigma\to S$ is continuous and fixes $S$. For finite subtrees $\emptyset\ne S\subseteq T\subseteq Y^\sigma$,
$
r_S^\sigma\circ r_T^\sigma=r_S^\sigma.
$
\smallskip\item The sets $D_z^\sigma(x)$ are the components of $Y^\sigma\smallsetminus\{z\}$ and are open and arcwise connected. Each connected subset $C\subseteq Y^\sigma$ contains $[x,y]$ whenever $x,y\in C$. For distinct $x,y\in Y^\sigma$, every $z\in(x,y)$ is a cut point.
\smallskip\item The map $\pi^\sigma\:Y^\sigma\to Z^\sigma$ is continuous and surjective.
\smallskip\item The map $f_\sigma\:Y^\sigma\to Y^\sigma$ is continuous, onto and proper, with
\begin{equation}\label{eq:pi-exact-pullback}
 \pi^\sigma \circ f_\sigma=f\circ \pi^\sigma,\quad f_\sigma^{-1}(Y_n^\sigma)=Y_{n+1}^\sigma\quad \text{ for all } n\in\N_0.
\end{equation}
\end{enumerate}
\end{prop}

\begin{proof}
\textup{(i)--(ii)} Let $\emptyset\ne R\subseteq S\subseteq T$ be finite trees. Add the vertices of $S$ to the vertex set of $T$. The map $r_{S,T}$ fixes each edge of $S$ and is constant on each remaining edge, so it is continuous. For $x\in T$, the arc $[x,r_{R,T}(x)]$ passes through $r_{S,T}(x)$. Hence $r_{R,S}\circ r_{S,T}=r_{R,T}$.

Fix an integer $N\ge0$. For integers $m\ge n\ge N$, the maps $r_{Y_N^\sigma,Y_m^\sigma}$ and $r_{Y_N^\sigma,Y_n^\sigma}$ agree on $Y_n^\sigma$. They define a map $r_N^\sigma\:Y^\sigma\to Y_N^\sigma$ by $r_N^\sigma|_{Y_n^\sigma}\=r_{Y_N^\sigma,Y_n^\sigma}$ for each integer $n\ge N$. This map fixes $Y_N^\sigma$. Each restriction to $Y_n^\sigma$ is continuous, so $r_N^\sigma$ is continuous by~\eqref{eq:pi-weak-union}. Let $x,\,y\in Y^\sigma$ be distinct and choose $N$ with $x,\,y\in Y_N^\sigma$. Choose disjoint open neighborhoods $U$ and $V$ of $x$ and $y$, respectively, in $Y_N^\sigma$. Then $(r_N^\sigma)^{-1}(U)$ and $(r_N^\sigma)^{-1}(V)$ separate $x,y$ in $Y^\sigma$. Hence $Y^\sigma$ is Hausdorff. The continuous inclusion $Y_N^\sigma\to Y^\sigma$ is therefore an embedding, since $Y_N^\sigma$ is compact.

Suppose a compact set $K\subseteq Y^\sigma$ lies in no $Y_n^\sigma$. Choose $x_j\in K$ with strictly increasing $\min\{n:x_j\in Y_n^\sigma\}$, and set $D\=\{x_j\}_{j\in\N}$. Each subset $A\subseteq D$ meets every $Y_n^\sigma$ in a finite set, so $A$ is closed in $Y^\sigma$ by the definition of the weak topology (see~\eqref{eq:pi-weak-union}). Thus $D$ is closed in $K$, while its singletons form an open cover with no finite subcover, which is a contradiction.

Since the union of two arcs with the same endpoints is compact, it is contained in some $Y_N^\sigma$. The arcs coincide because $Y_N^\sigma$ is a tree. Hence $Y^\sigma$ is uniquely arcwise connected. For a nonempty finite subtree $S\subseteq Y_N^\sigma$, the map $r_S^\sigma$ is the composition $r_{S,Y_N^\sigma}\circ r_N^\sigma$, so it is continuous and fixes $S$. For nonempty finite subtrees $S\subseteq T\subseteq Y_N^\sigma$, apply $r_{S,T}\circ r_{T,Y_N^\sigma}=r_{S,Y_N^\sigma}$ to $r_N^\sigma(x)$. This gives $(r_S^\sigma\circ r_T^\sigma)(x)=r_S^\sigma(x)$.

\smallskip
\textup{(iii)} Fix $x,\,y,\,z\in Y^\sigma$ with $x,\,y\ne z$. In every $Y_N^\sigma$ containing $x$, $y$, and $z$, the intersections $D_z^\sigma(x)\cap Y_N^\sigma$ and $D_z^\sigma(y)\cap Y_N^\sigma$ are connected components of $Y_N^\sigma\smallsetminus\{z\}$. They coincide precisely when $[x,y]$ avoids $z$. Taking unions over $N$ shows that $D_z^\sigma(x)$ and $D_z^\sigma(y)$ are equal or disjoint. Each $D_z^\sigma(x)$ contains $x$, so the distinct sets form a partition of $Y^\sigma\smallsetminus\{z\}$.

For $u,v\in D_z^\sigma(x)$, choose a tree $Y_N^\sigma$ containing $x,\,z,\,u,\,v$. The points $u,v$ lie in the component of $Y_N^\sigma\smallsetminus\{z\}$ containing $x$, so $[u,v]\subseteq D_z^\sigma(x)$. For every $n$, the intersection $D_z^\sigma(x)\cap Y_n^\sigma$ is either empty, equal to $Y_n^\sigma$, or a connected component of $Y_n^\sigma\smallsetminus\{z\}$. Each intersection is open in $Y_n^\sigma$, so $D_z^\sigma(x)$ is open by~\eqref{eq:pi-weak-union}. Its complement in $Y^\sigma\smallsetminus\{z\}$ is a union of the other $D_z^\sigma(y)$. Thus the connected set $D_z^\sigma(x)$ is a component.

Let $C\subseteq Y^\sigma$ be connected and contain $x,y$. If $C$ omits a point $z\in(x,y)$, then $C\cap D_z^\sigma(x)$ and $C\smallsetminus D_z^\sigma(x)$ separate $C$. Thus $[x,y]\subseteq C$. The sets $D_z^\sigma(x)$ and $D_z^\sigma(y)$ are distinct for $z\in(x,y)$, so every such $z$ cuts $Y^\sigma$.

\smallskip

\textup{(iv)} Equation~\eqref{eq:ch4-stage-projections} states that $\pi_{n+1}^\sigma\circ \iota_n^\sigma=\pi_n^\sigma$. Thus the maps $\pi_n^\sigma$ define $\pi^\sigma\:Y^\sigma\to Z^\sigma$. Each $\pi_n^\sigma$ is continuous, so $\pi^\sigma$ is continuous by~\eqref{eq:pi-weak-union}. Moreover, $\pi^\sigma(Y^\sigma)=\bigcup_{n=0}^{+\infty}\pi_n^\sigma(Y_n^\sigma)=\bigcup_{n=0}^{+\infty}Z_n^\sigma=Z^\sigma$.

\smallskip

\textup{(v)} For each integer $n\ge2$ and $b\in\sB_n^\sigma$, Lemma~\ref{lem:complete-angle-class-images} gives $b'\in\sB_{n-1}^\sigma$ with $\cA_{n-1}^\sigma(b')=\mu(\cA_n^\sigma(b))$. Here $\mu$ is the circle map $t\mapsto(\deg f)t\pmod1$. Evaluating $\gamma^{n-1}\circ \mu=f\circ \gamma^n$ on $\cA_n^\sigma(b)$ gives $\prV(b')=f(\prV(b))$. For every tile $X$ incident to $b$,~\eqref{eq:finite-gap-dynamics} gives $f(X)\leftrightarrow b'$. Hence $\bigl(\theta_{f(X)}^\sigma\circ (\theta_X^\sigma)^{-1}\bigr) (\beta_n^\sigma(b))=\beta_{n-1}^\sigma(b')$. For $b\in\sB_1^\sigma$ and a tile $X$ incident to $b$, the map $\theta_{f(X)}^\sigma\circ (\theta_X^\sigma)^{-1}$ takes $\beta_1^\sigma(b)$ to $\xi_{f(\prV(b))}^\sigma$. Thus, for every $n\in \N$, these maps agree on all intersections $T_X^\sigma\cap T_Y^\sigma$ by~\eqref{eq:ch4-stage-copy-intersection}. The finitely many $T_X^\sigma$ are closed in $Y_n^\sigma$, so $\phi_n^\sigma$ is well-defined and continuous. For each $U\in\X_\sigma^{n-1}$, choose $X\in\X_\sigma^n$ with $f(X)=U$. By~\eqref{eq:pi-local-dynamics}, $\phi_n^\sigma$ maps $T_X^\sigma$ homeomorphically onto $T_U^\sigma$. Hence $\phi_n^\sigma(Y_n^\sigma)=Y_{n-1}^\sigma$.

Fix an integer $n\ge1$, $X\in\X_\sigma^n$, and $U\in\X_\sigma^1$. By~\eqref{eq:pi-local-dynamics},~\eqref{eq:ch4-refined-copy}, and Proposition~\ref{lem:ch4-sector-copies}~\textup{(ii)},
\begin{equation*}
\phi_{n+1}^\sigma\circ \Theta_X^\sigma\circ \theta_U^\sigma=\theta_{f(L(X,1,U))}^\sigma=\theta_{L(f(X),1,U)}^\sigma=\Theta_{f(X)}^\sigma\circ \theta_U^\sigma.
\end{equation*}
Here $\Theta_{f(X)}^\sigma$ is the identity on $Y_1^\sigma$ when $n=1$. The sets $T_U^\sigma$ cover $Y_1^\sigma$. Thus~\eqref{eq:ch4-copy-stability} gives $\phi_{n+1}^\sigma\circ \iota_n^\sigma\circ \theta_X^\sigma=\Theta_{f(X)}^\sigma\circ \iota_0^\sigma=\iota_{n-1}^\sigma\circ \phi_n^\sigma\circ\theta_X^\sigma$. The sets $T_X^\sigma$ cover $Y_n^\sigma$, so
\begin{equation}\label{eq:finite-dynamics-compatibility}
\phi_{n+1}^\sigma\circ \iota_n^\sigma=\iota_{n-1}^\sigma\circ \phi_n^\sigma\quad \text{ for all } n\ge1.
\end{equation}

Fix an integer $n\ge1$ and a road $I_{(n,X,\alpha)}^\sigma$. Put $\lambda\=(n,X,\alpha)$. Equations~\eqref{eq:pi-local-dynamics} and~\eqref{eq:ch4-road-one-step} give
\begin{equation*}
\pi_{n-1}^\sigma\circ\phi_n^\sigma\circ\eta_n^\sigma(\lambda,s)=f\bigl(F^{\lambda,\sigma}(s)\bigr)=(f\circ\pi_n^\sigma\circ\eta_n^\sigma)(\lambda,s).
\end{equation*}
The sets $\eta_n^\sigma\bigl(I_{(n,X,\alpha)}^\sigma\bigr)$ cover $Y_n^\sigma$, so
\begin{equation}\label{eq:ch4-dynamics-projections}
\pi_{n-1}^\sigma\circ \phi_n^\sigma=f\circ \pi_n^\sigma\quad \text{ for all } n\ge1.
\end{equation}
By~\eqref{eq:finite-dynamics-compatibility}, the maps $\phi_n^\sigma$ agree on their common domains in $Y^\sigma$. Thus they induce a well-defined map $f_\sigma\:Y^\sigma\to Y^\sigma$ with $f_\sigma|_{Y_n^\sigma}\=\phi_n^\sigma$. Each restriction is continuous, so $f_\sigma$ is continuous by~\eqref{eq:pi-weak-union}. The map $f_\sigma$ sends each $Y_{n+1}^\sigma$ onto $Y_n^\sigma$, so it is onto $Y^\sigma$. Equation~\eqref{eq:ch4-dynamics-projections} gives $\pi^\sigma \circ f_\sigma=f\circ\pi^\sigma$.

It remains to show that $f_\sigma^{-1}(Y_n^\sigma)=Y_{n+1}^\sigma$ for each $n\in\N_0$ and that $f_\sigma$ is proper. Fix an integer $m\ge2$ and a tile $V\in\X_\sigma^m$. Put $U\=f^{m-1}(V)$. For every $z\in Y_0^\sigma$,~\eqref{eq:ch4-stage-sections} and~\eqref{eq:pi-local-dynamics} give
\begin{equation*}
\theta_V^\sigma(z)\in Y_{m-1}^\sigma
\iff\theta_U^\sigma(z)\in Y_0^\sigma
\iff\phi_m^\sigma(\theta_V^\sigma(z))\in Y_{m-2}^\sigma.
\end{equation*}
The copies $T_V^\sigma$ cover $Y_m^\sigma$, so $(\phi_m^\sigma)^{-1}(Y_{m-2}^\sigma)=Y_{m-1}^\sigma$. Consequently, $f_\sigma$ maps $Y_m^\sigma\smallsetminus Y_{m-1}^\sigma$ into $Y_{m-1}^\sigma\smallsetminus Y_{m-2}^\sigma$. It also maps $Y_1^\sigma$ onto $Y_0^\sigma$. Hence, $f_\sigma^{-1}(Y_n^\sigma)=Y_{n+1}^\sigma$ for each $n\in \N_0$.

Finally, let $K\subseteq Y^\sigma$ be compact. By~\textup{(i)}, choose $N$ with $K\subseteq Y_N^\sigma$. The set $K$ is closed in $Y_N^\sigma$. Therefore, $(\phi_{N+1}^\sigma)^{-1}(K)$ is closed in compact $Y_{N+1}^\sigma$ and hence compact. This set equals $f_\sigma^{-1}(K)$, so $f_\sigma$ is proper.
\end{proof}

\begin{prop}\label{prop:zipper-noncrossing}
Let $\pi^\pm\:Y^\pm\to Z^\pm$ be the projections given by Proposition \ref{prop:pi-weak-tree}. For continuous paths $a^\pm\:[0,1]\to Y^\pm$, the paths $\pi^\pm\circ a^\pm$ do not cross.
\end{prop}

\begin{proof}
We prove by contradiction. Put $E^\sigma\=a^\sigma([0,1])$. Suppose $\pi^\sigma(E^\sigma)\subseteq D$ for a closed disk $D\subseteq S^2$ and $\{p_0^\sigma,p_1^\sigma\}\=\pi^\sigma(E^\sigma)\cap\partial D$, with $\{p_0^+,p_1^+\}$ and $\{p_0^-,p_1^-\}$ disjoint and alternating. Choose pairwise disjoint closed arcs $J_i^\sigma\subseteq\partial D$ with $p_i^\sigma\in\inte_\circ J_i^\sigma$. Choose open disks $W_i^\sigma\ni p_i^\sigma$ with pairwise disjoint closures such that $\overline{W_i^\sigma}\cap\partial D\subseteq J_i^\sigma$ and $\overline{W_i^\sigma}\cap\pi^{-\sigma}(E^{-\sigma})=\emptyset$.

By Proposition~\ref{prop:pi-weak-tree}~\textup{(i)}, $E^\sigma\subseteq Y_m^\sigma$ for all large integers $m$. In the geometric representation of the $m$-tile connection from Subsections~\ref{subsec:connections} and~\ref{prelimiconnectiongraph}, choose disjoint closed disks $U(v)$, $v\in\V^m$, with $\diam_\rho U(v)\le\Lambda^{-m}$. Take the disks small enough that $U(v)$ meets only the $m$-tiles containing $v$. Write $Q_b$ for the disk representing $b$ in $U(\prV(b))$. For each $X\in\X_\sigma^m$, put
\begin{equation*}
O_X\=\inte\biggl[
\bigl(X\smallsetminus\bigcup_{v\in\V^m}\inte U(v)\bigr)
\cup\bigcup_{X\leftrightarrow b}Q_b
\biggr].
\end{equation*}
The arcs $X\cap\partial U(v)$ for $v\in X\cap\V^m$ cut off disjoint disks from $X$. Attaching $Q_b$ along $X\cap\partial U(\prV(b))$ gives a connected $O_X$. By \eqref{eq:ch4-stage-copy-intersection}, if $T_X^\sigma\cap T_Y^\sigma\ne\emptyset$ for $X\ne Y$, then $\inte Q_b\subseteq O_X\cap O_Y$ for a block $b$ with $X\leftrightarrow b\leftrightarrow Y$. Since $E^\sigma\cap T_X^\sigma$ form a finite closed cover of the connected set $E^\sigma$, the union $O_m^\sigma\=\bigcup_{E^\sigma\cap T_X^\sigma\ne\emptyset}O_X$ is connected. Each $O_m^\sigma$ lies in the interior of the union of the closed pieces of color $\sigma$, so $O_m^+\cap O_m^-=\emptyset$.

For $y\in T_X^\sigma$, Proposition~\ref{c:initial}~\textup{(iii)}, \eqref{eq:ch4-stage-model-map}, and \eqref{eq:ch4-stage-projections} give $\pi^\sigma(y)=F^{\lambda,\sigma}(s)$ with $\lambda\=(m,X,\alpha)$ and $s\in I_\lambda^\sigma$. For $J\in\mathcal P_{m+1}^{\lambda,\sigma}$ containing $s$, \eqref{eq:road-data} and the definition of $L$ give $A_{m,m+1}(X_{m+1}^{\lambda,\sigma}(J))=X$. Thus $F_{m+1}^{\lambda,\sigma}(s)\in H_1^m(X)$ by \eqref{eq:ch4-road-containment} and Lemma~\ref{lem:ch4-sector-ancestors}~\textup{(ii)}. Each $z\in O_X$ lies in $X$ or in some $U(v)$ with $v\in X$. Thus, for $z\in O_X$ and $y\in T_X^\sigma$, \eqref{eq:pi-visual-estimates} and \eqref{eq:ch4-road-Cauchy} give
\begin{equation*}
\rho(z,\pi^\sigma(y))\le(1+2C+C_0)\Lambda^{-m}.
\end{equation*}
For large $m$, the estimate and compactness give $O_m^\sigma\subseteq\inte D\cup W_0^\sigma\cup W_1^\sigma$ and $O_m^\sigma\cap\overline{W_i^{-\sigma}}=\emptyset$. Taking $y\in E^\sigma\cap T_X^\sigma$ with $\pi^\sigma(y)=p_i^\sigma$ also gives $O_m^\sigma\cap W_i^\sigma\ne\emptyset$.

The two sets $O_m^\sigma\cup W_0^\sigma\cup W_1^\sigma$ are disjoint, open, and connected. For each $\sigma\in\{+,\,-\}$, take an embedding $a_m^\sigma\:[0,1]\to O_m^\sigma\cup W_0^\sigma\cup W_1^\sigma$ with $a_m^\sigma(i)=p_i^\sigma$ for $i\in\{0,1\}$. Put $s_1^\sigma\=\min\bigl\{t\in[0,1]:a_m^\sigma(t)\in\overline{W_1^\sigma}\smallsetminus\inte D\bigr\}$ and $s_0^\sigma\=\max\bigl\{t\in[0,s_1^\sigma]:a_m^\sigma(t)\in\overline{W_0^\sigma}\smallsetminus\inte D\bigr\}$. Then $a_m^\sigma((s_0^\sigma,s_1^\sigma))\subseteq\inte D$. By continuity, $a_m^\sigma(s_i^\sigma)\in\partial D\cap\overline{W_i^\sigma}\subseteq J_i^\sigma$ for $i\in\{0,1\}$. The two disjoint subarcs have alternating endpoints, which contradicts the Jordan curve theorem.
\end{proof}

Now we give a proof of Theorem~\ref{A}.

\begin{proof}[\textbf{Proof of Theorem~\ref{A}}]
Let $f_0\:S^2\to S^2$ be an expanding Thurston map. Choose an integer $N\ge1$ and an $f_0^N$-invariant Jordan curve $\cC$ as in Lemma~\ref{7.5}, and write $f\=f_0^N$. The construction recalled before Definition~\ref{Propertiesofconnections} provides a connection of $1$-tiles satisfying \textup{(C1)}--\textup{(C2)}. Use this connection to construct the sets $Z^+$ and $Z^-$ in~\eqref{defZpm}.

By Lemmas~\ref{nonemptyandpathconnected},~\ref{fcompletelyinvariant}, and~\ref{dense}, both sets are nonempty, path-connected, dense in $S^2$, and completely $f$-invariant.

Fix $\sigma\in\{+,\,-\}$. The finite trees $Y_n^\sigma$ constructed above form an increasing union $Y^\sigma$, equipped with the weak topology in~\eqref{eq:pi-weak-union}. Proposition~\ref{prop:pi-weak-tree}~\textup{(i)} shows that $Y^\sigma$ is a topological $\R$-tree and that each $Y_n^\sigma$ is a subtree. By Proposition~\ref{prop:pi-weak-tree}~\textup{(iv)}, the projection $\pi^\sigma\:Y^\sigma\to Z^\sigma$ is continuous and surjective. Proposition~\ref{prop:pi-weak-tree}~\textup{(v)} provides a continuous map $f_\sigma\:Y^\sigma\to Y^\sigma$ satisfying $\pi^\sigma\circ f_\sigma=f\circ\pi^\sigma$. By Proposition~\ref{prop:zipper-noncrossing}, for each pair of continuous paths $a^\pm\:[0,1]\to Y^\pm$, the paths $\pi^\pm\circ a^\pm$ do not cross, as required by Definition \ref{def_gen_zipper}~\textup{(iv)}. Thus, $(Z^+,Z^-)$ satisfies Definition~\ref{def_gen_zipper}, and the tree models and their maps have all the properties asserted in Theorem~\ref{A}.
\end{proof}

\begin{remark}
We relate the P-zipper of Theorem~\ref{A} to Meyer's polynomial models and to the zippers of Calegari and Loukidou. When $f$ has no periodic critical points, Meyer realizes $f$ as a mating of two polynomials $P_\pm$ with dendritic Julia sets $J_\pm$ \cite[Theorem~1.2]{Meyer09}. For $\sigma\in\{+,\,-\}$, the tree model $Y^\sigma$ admits a continuous injection into $J_\sigma$. If, in addition, there are no perfect fits, $Z^\sigma$ is obtained from the corresponding branch of Calegari and Loukidou's zipper by adjoining the iterated preimages of $\post f$. We are aware of this similarity with Meyer's construction. However, we have chosen to work directly with the marked connection and the lifts $H^m$ of $H^0$ rather than with Meyer's polynomial models. This keeps the proof of Theorem~\ref{A} independent of the mating theorem and allows us to describe the finite trees $Y_n^\sigma$ and their dynamics directly in terms of the connections of $n$-tiles.

Suppose that $f$ has no periodic critical points. Take Meyer's polynomial $P_\sigma$ for the same invariant Jordan curve $\cC$ and connection of $1$-tiles, and let $q_\sigma\:S^1\to J_\sigma$ be the Carath\'eodory semi-conjugacy whose kernel is the corresponding equivalence relation of Definition~\ref{oversetwsim}. For each integer $m\ge1$ and $b\in\sB_m^\sigma$, define
\begin{equation*}
\Psi^\sigma(\beta_m^\sigma(b))\=q_\sigma(t)\quad\text{ for }t\in\cA_m^\sigma(b).
\end{equation*}
This map extends uniquely to a continuous injection $\Psi^\sigma\:Y^\sigma\to J_\sigma$, which maps each $Y_n^\sigma$ homeomorphically onto the subtree of $J_\sigma$ spanned by $q_\sigma(\mathbf A^n)$. Consequently,
\begin{equation*}
\Psi^\sigma(Y^\sigma)=\operatorname{Cut}(J_\sigma)\cup\bigcup_{n=0}^{+\infty}q_\sigma(\mathbf A^n),
\end{equation*}
where $\operatorname{Cut}(J_\sigma)$ denotes the set of cut points of $J_\sigma$. For the mating projection $M_\sigma\:J_\sigma\to S^2$ satisfying $M_\sigma\circ q_\sigma=\gamma$, we have $M_\sigma\circ\Psi^\sigma=\pi^\sigma$.

Suppose further that the relations of Definition~\ref{oversetwsim} have no perfect fits. Let $T(P_\sigma)$ be the Hubbard tree spanned in $J_\sigma$ by $(\post P_\sigma)\cap\mathbb C$, and let $H_{\mathrm o}^\sigma$ be the image under $M_\sigma$ of $T(P_\sigma)$ with its $1$-valent vertices removed \cite[Definition~7.8]{CLWheels}. By \cite[Subsection~7.1.2, Theorem~1.23, and Proposition~7.9]{CLWheels}, the color-$\sigma$ branch of Calegari and Loukidou's zipper is
\begin{equation*}
Z_{\mathrm{CL}}^\sigma\=\bigcup_{n=0}^{+\infty}f^{-n}(H_{\mathrm o}^\sigma)=M_\sigma(\operatorname{Cut}(J_\sigma)).
\end{equation*}
By the above description of $\Psi^\sigma(Y^\sigma)$ and \eqref{eq:cell-vertices}, we have
\begin{equation*}
Z^\sigma=Z_{\mathrm{CL}}^\sigma\cup\bigcup_{n=0}^{+\infty}\V^n.
\end{equation*}
Thus, if there are no perfect fits, our construction recovers the zipper example of Calegari and Loukidou up to the inclusion of additional vertices.
\end{remark}

\section{End circles}
\label{endcircles}
In this section, we discuss the end circles of the tree models $Y^\pm$ of the P-zipper. We first recall the notion of the end space and show that each end is represented by a rooted ray in Subsection~\ref{sub_Ends_and_arc_representations}. Next, in Subsection~\ref{sub_Spherical_ribbon_structures}, we show that $Y^\sigma$ defined above carries a ribbon structure, i.e., we can define a circular order for the germs at a branched point in $Y^\sigma$. Then in Subsection~\ref{subsec:hairiness}, we show that $Y^\sigma$ is hairy for the connection of $1$-tiles in Lemma~\ref{degree3}. Then we show that there is a circular order on the end space of $Y^\sigma$ so that it can be completed to a circle in Subsection~\ref{sub_Circular_order_on_the_end_spaces}. In Subsection~\ref{theendspace}, we describe more properties of the end circles of $Y^+$ and $Y^-$. We also show that there is an orientation-reversing homeomorphism between the two circles. Finally, in Subsection~\ref{sub_Dynamics_on_the_end_circle}, we show that $f$ induces a map on the end circle, which is conjugate to $z\mapsto z^d$.

Throughout this section, we use $\sigma\in\{+,\,-\}$ to represent the color as in the previous sections. A signed $n$-tile $X_\sigma^n$ is a white $n$-tile if $\sigma=+$ and a black $n$-tile if $\sigma=-$. Let $Y_n^\sigma$ be given by Definition~\ref{c:relation}, $\iota_n^\sigma$ be given by Definition~\ref{def:ch4-refined-copies}, and $Y^\sigma$, $\pi^\sigma$, $f_\sigma$ be given by Definition~\ref{def:ch4-limit-spaces}.

\subsection{Ends and arc representations}\label{sub_Ends_and_arc_representations} 
In this subsection, we recall the notion of the end space of a topological $\R$-tree. Then we focus on the end space of $Y^\sigma$ and show that each end of $Y^\sigma$ is represented by a unique ray; see Proposition~\ref{prop:rooted-representatives}.

Recall that every topological $\R$-tree is arc-wise connected. For two points $x,\,y\in T$, let $[x,y]_T\subseteq T$ denote the unique arc joining $x$ and $y$. We also set $(x,y]_T\=[x,y]_T\smallsetminus\{x\}$ to indicate that $x$ is not included.

By a \emph{ray} in a topological $\R$-tree $T$ we mean the image of a topological embedding $r\:[0,+\infty)\to T$; it is \emph{rooted at} $o\in T$ if $r(0)=o$. A continuous map is \emph{proper} if the preimage of every compact set is compact. A ray is \emph{proper} if one, or equivalently every, parametrization of it is proper. For each $N\in(0,+\infty)$, we call $r|_{[N,+\infty)}$ a \emph{tail} of the ray $r$.

\begin{definition}
    Let $T$ be a topological $\R$-tree. A point $p\in T$ is called a \emph{cut point} if $T\smallsetminus\{p\}$ is not connected. We write $\inter (T)\=\{x\in T:x\text{ is a cut point of }T\}$ for the \emph{interior} of $T$. A point in $T$ is called a \emph{hard endpoint} if it is not a cut point. An \emph{interior subtree} of $T$ is either a single cut point in $T$ or a finite subtree consisting entirely of cut points, and we denote $\mathcal K(T)\=\{K:K \text{ is an interior subtree of }T\}$.
\end{definition}

\begin{definition}\label{def_collectionend}
Let $T$ be a topological $\mathbb{R}$-tree. An \emph{end} of $T$ is a collection $\{C_K\}_{K\in \mathcal{K}(T)}$ satisfying the following conditions.
\begin{enumerate}
\smallskip
    \item For each $K\in\cK(T)$, $C_K$ is a connected component of $T\smallsetminus K$;
\smallskip
    \item For all $K_1,\,K_2\in\cK(T)$, if $K_1 \subseteq K_2$, then $C_{K_2} \subseteq C_{K_1}$.
\end{enumerate}
\end{definition}

\begin{definition}[{\cite[Definition~2.8]{CL24}}]
\label{def_endspace}
The {\it interior end space} of a topological $\mathbb R$-tree $T$ is the inverse limit $\mathcal{E}(T) \= \displaystyle \lim_{\longleftarrow} \pi_0(T\smallsetminus K)$
over $\mathcal K(T)$.
\end{definition}

The following proposition yields that the ends defined using all interior subtrees are canonically the same as the ends defined using only interior subtrees containing a point $p\in Y^+$.
We denote $C^+(K)\=\pi_0(Y^+\smallsetminus K)$ for each $K\in \mathcal K(Y^+)$.

\begin{prop} \label{p_interior_subtree}
    Let $\mathcal{K}(Y^+)$ denote the set of interior subtrees of $Y^+$, let $p$ be an arbitrary cut point of $Y^+$ and denote $\mathcal K_p(Y^+)\=\{K\in \mathcal K(Y^+):p\in K\}$. Then the restriction map 
    \begin{equation*}
\operatorname{res}_p\:\varprojlim_{K\in \mathcal K(Y^+)} C^+(K)\to \varprojlim_{K\in \mathcal K_p(Y^+)} C^+(K)
    \end{equation*}
    is a bijection. 
\end{prop}

\begin{proof}
Write $\mathcal K\=\mathcal K(Y^+)$ and $\mathcal K_p\=\mathcal K_p(Y^+)$. For $K\subseteq L$ in $\mathcal K$, let $\pi_{L,K}\:C^+(L)\to C^+(K)$ be the natural map sending each component of $T\smallsetminus L$ to the unique component of $T\smallsetminus K$ containing it.

We first prove that for every $K\in\mathcal K$, there exists $L\in\mathcal K_p$ such that $K\subseteq L$. Since $K$ is an interior subtree, it is a finite subtree that contains no hard ends. Thus, every point of $K$ is a cut point of $Y^+$. Let $x$ be a point in $K$, then every point in the arc $[p,x]$ is a cut point. Define $L\=K\cup [p,x]$. Then $L$ is an interior subtree with $K\subseteq L$ and $p\in L$. Hence, $L\in \mathcal K_p$.

Now we prove that $\operatorname{res}_p$ is injective. Fix an arbitrary $K\in\mathcal K$. Let $\{C_K\}_{K\in\mathcal K}$ and $\{D_K\}_{K\in\mathcal K}$ be two elements of $\varprojlim_{K\in\mathcal K} C^+(K)$ whose restrictions to $\mathcal K_p$ agree. Choose $L\in\mathcal K_p$ with $K\subseteq L$. Then $C_K=\pi_{L,K}(C_L)$ and $D_K=\pi_{L,K}(D_L)$, and we have $C_L=D_L$. Hence $C_K=\pi_{L,\,K}(C_L)=\pi_{L,\,K}(D_L)=D_K$. Since $K$ is arbitrary, $\{C_K\}_{K\in\mathcal K}=\{D_K\}_{K\in\mathcal K}$. Therefore, $\operatorname{res}_p$ is injective.

Next, we prove that $\operatorname{res}_p$ is surjective. Fix $\{D_L\}_{L\in\mathcal K_p}\in\varprojlim_{L\in\mathcal K_p} C^+(L)$. For each $K\in\mathcal K$, choose $L\in\mathcal K_p$ with $K\subseteq L$, and denote $C_K\=\pi_{L, K}(D_L)$. We claim that such a $C_K$ is independent of the choice of $L$. Suppose $L_1,\,L_2\in\mathcal K_p$ both contain $K$. Then there exists $M\in\mathcal K_p$ such that $L_1\subseteq M$ and $L_2\subseteq M$. Since $\{D_L\}_{L\in\mathcal K_p}$ is an inverse-limit element, we have $D_{L_i}=\pi_{M, L_i}(D_M)$ for $i\in\{1,\,2\}$. Therefore, $\pi_{L_i, K}(D_{L_i})=\pi_{L_i, K}(\pi_{M,L_i}(D_M))=\pi_{M, K}(D_M)$ for $i\in\{1,\,2\}$. Thus, the value of $C_K$ is independent of the choice of $L\in\mathcal K_p$ containing $K$. We now check that $\{C_K\}_{K\in\mathcal K}$ is an element of $\varprojlim_{K\in\mathcal K} C^+(K)$. Let $K\subseteq K'$ be two elements of $\mathcal K$. Choose $L\in\mathcal K_p$ with $K'\subseteq L$. We have $C_{K'}=\pi_{L, K'}(D_L)$ and $C_K=\pi_{L, K}(D_L)$. Hence $\pi_{K', K}(C_{K'})=\pi_{K', K}(\pi_{L, K'}(D_L))=\pi_{L, K}(D_L)=C_K$. The set $\{C_K\}_{K\in\mathcal K}$ satisfies Definition~\ref{def_collectionend}. Finally, we show that this set restricts to $\{D_L\}_{L\in\mathcal K_p}$. If $K\in\mathcal K_p$, then in the definition of $C_K$ we can choose $L=K$. Hence $C_K=\pi_{K, K}(D_K)=D_K$. Every element of $\varprojlim_{K\in\mathcal K_p} C^+(K)$ comes from an element of $\varprojlim_{K\in\mathcal K} C^+(K)$. This concludes that $\operatorname{res}_p$ is surjective.

Combining injectivity and surjectivity, the restriction map is a bijection.
\end{proof}

We now identify the ends of $Y^\sigma$ with rooted rays. Fix a cut point $o^\sigma \in Y_0^\sigma$ and let $\mathcal K_{o^\sigma}(Y^\sigma)$ be the family of interior subtrees containing $o^\sigma$. Note that if $K$ is an interior subtree and $x\in K$ is a point, then $K\cup [o^\sigma,x]$ is again an interior subtree. The following proposition shows that each end of $Y^\sigma$ is represented by a unique ray.

\begin{prop}\label{prop:rooted-representatives}
    For every $\omega=\{C_K^\omega\}_{K\in \mathcal K_{ o^\sigma}(Y^\sigma)} \in \mathcal E(Y^\sigma)$, there exists a unique ray $R_\omega\subseteq \inter (Y^\sigma)$ issuing from $o^\sigma$ and representing $\omega$, in the sense that for every $K\in\mathcal K_{o^\sigma}(Y^\sigma)$ a sufficiently far tail of $R_\omega$ lies in $C_K^\omega$, where $C_K^\omega$ is a connected component of $Y^\sigma\smallsetminus K$.

    Put $A_\omega\=\overline {R_\omega}^{Y^\sigma}$. If $A_\omega$ is compact, then $A_\omega=[o^\sigma, h_\omega]_{Y^\sigma}$ for a unique hard endpoint $h_\omega$. 
\end{prop}

\begin{proof}
Fix $K\in\mathcal K_{o^\sigma}(Y^\sigma)$ and write $C\=C_K^\omega$. 
We first associate to $C$ a unique attachment point $a_K(\omega)\in K$. Choose $x\in C$. Since $K$ is a finite subtree and $[o^\sigma,x]_{Y^\sigma}$ is an arc, both of them are compact. Hence, $[o^\sigma,x]_{Y^\sigma}\cap K$ is a compact initial subarc of $[o^\sigma,x]_{Y^\sigma}$. Let $a_K(\omega)$ be its terminal point. Then $(a_K(\omega),x]_{Y^\sigma}\cap K=\emptyset$, and $(a_K(\omega),x]_{Y^\sigma}\subseteq C$. It follows that $a_K(\omega)\in K\cap \oC$.

    We claim that $a_K(\omega)$ is independent of the choice of $x\in C$. We argue by contradiction and suppose $x,\,y\in C$ induce distinct attachment points $a$ and $b$. The unique arc $[x,y]_{Y^\sigma}$ would contain the nondegenerate arc $[a,b]_{Y^\sigma}\subseteq K$, or the attachment points would be the same. Each point of $(a,b)_{Y^\sigma}$ belongs to $K$ and separates $x$ from $y$ in the tree ${Y^\sigma}$, so $x$ and $y$ would lie in different components of $Y^\sigma\smallsetminus K$, which is a contradiction since $x,\,y\in C$. Thus, $a_K(\omega)$ is independent of the choice of $x\in C$.
    
    To prove $\overline{C}\cap K=\{a_K(\omega)\}$, it remains to rule out another point of $K$ in $\overline C$. If $K$ contains only one point, then we are done. Otherwise, suppose that $b\in K\smallsetminus\{a_K(\omega)\}$ and choose $z\in(a_K(\omega),b)_{Y^\sigma}$. Then $z$ separates $C$ from $b$. The component of $Y^\sigma\smallsetminus\{z\}$ containing $b$ is open in the weak topology (see Definition~\ref{def:ch4-limit-spaces}), since its intersection with each finite subtree $Y_n^\sigma$ for some $n\in\N$ is either empty, the whole subtree $Y_n^\sigma$, or a union of components of a finite tree with $z$ removed. Hence, it is an open neighborhood of $b$ disjoint from $C$, so $b\notin\overline C$. Therefore, $\overline{C}\cap K=\{a_K(\omega)\}$.

    Next, we define the ray $R_\omega$. If $K\subseteq L$, choose $x\in C_L^\omega\subseteq C_K^\omega$. Along the unique arc $[o^\sigma,x]_{Y^\sigma}$, one leaves $K$ through $a_K(\omega)$ and later leaves $L$ through $a_L(\omega)$. Hence, $a_K(\omega)\in [o^\sigma,a_L(\omega)]_{Y^\sigma}$. Consequently, the rooted arcs $[o^\sigma,a_K(\omega)]_{Y^\sigma}$ for each $K\in \mathcal K_{o^\sigma}(Y^\sigma)$ are nested. Define $R_\omega\=\bigcup_{K\in \mathcal K_{o^\sigma}(Y^\sigma)}[o^\sigma,a_K(\omega)]$.

    For every fixed $K$ and every $L\supseteq K$, we have $(a_K(\omega),a_L(\omega)]\subseteq C_K^\omega$. Take $x\in C_K^\omega\cap\inter (Y^\sigma)$ and set $L'=K\cup[a_K(\omega),x]_{Y^\sigma}$. We claim that $L'$ is an interior subtree and $a_{L'}(\omega)\ne a_K(\omega)$. We argue by contradiction and suppose $a_{L'}(\omega)= a_K(\omega)$. Choose  $y\in C_{L'}^\omega\subseteq C_K^\omega$. Then we have $[x,y]_{Y^\sigma}=[x,a_K(\omega)]_{Y^\sigma}\cup[y,a_{L'}(\omega)]_{Y^\sigma}$ and $[x,a_K(\omega)]_{Y^\sigma}\cap[y,a_{L'}(\omega)]_{Y^\sigma}=\{a_K(\omega)\}$. Since $(a_{L'}(\omega),y]_{Y^\sigma}\subseteq C_{L'}^\omega\subseteq Y^\sigma\smallsetminus L'$. This contradicts $x,\,y\in C_K^\omega$ and $a_K(\omega)\in K$. Hence, $R_\omega$ has a nonempty tail in $C_K^\omega$, and therefore represents $\omega$.

    Every point of $R_\omega$ is a cut point. Indeed, the attachment points are cut points, and every point of an arc joining two cut points is a cut point. The union $R_\omega$ has no terminal point in $\inter (Y^\sigma)$, i.e., there is no $p\in\inter (Y^\sigma) $ such that $R_\omega\subseteq [o^\sigma,p]_{Y^\sigma}$. If $p$ were such a terminal point, let $K\=[o^\sigma,p]_{Y^\sigma}$ be an interior subtree in $Y^\sigma$. However, since $R_\omega$ represents $\omega$ as shown above, we have a nonempty tail of $R_\omega$ outside $[o^\sigma,p]_{Y^\sigma}$, contradicting $R_\omega\subseteq[o^\sigma,p]_{Y^\sigma}$. 

    We next show that $R_\omega$ is a topological ray. Suppose first that $R_\omega$ is not contained in $Y_n^\sigma$ for every $n\in\N$. Choose inductively integers $N_1<N_2<\cdots$ and points $p_j\in R_\omega$ such that $[o^\sigma,p_j]_{Y^\sigma}\subseteq Y_{N_j}^\sigma$ and $p_{j+1}\notin Y_{N_j}^\sigma$. Then, $[o^\sigma,p_1]_{Y^\sigma}\subsetneq[o^\sigma,p_2]_{Y^\sigma}\subsetneq\cdots$ and the points $p_j$ are cofinal in the rooted order on $R_\omega$. That means, if $x\in R_\omega\cap Y_N^\sigma$ for some $N\in\N$, then for $N_j\geq N$ we have $p_{j+1}\notin[o^\sigma,x]_{Y^\sigma}$, since otherwise $p_{j+1}\in[o^\sigma,x]_{Y^\sigma}\subseteq Y_N^\sigma\subseteq Y_{N_j}^\sigma$, which is a contradiction.

    Choose homeomorphisms $r_0\:[0,1]\to [o^\sigma,p_1]$ and $r_j\:[j,j+1]\to [p_j,p_{j+1}]_{Y^\sigma}$ for $j\in\N$, such that they agree at their common endpoints. Since the intervals $[j,j+1]$ form a locally finite closed cover of $[0,+\infty)$, these maps paste to a continuous injective map $r\:[0,+\infty)\to  Y^\sigma$ with image $R_\omega$.

    Next, we prove that the map $r$ is proper. Let $Q\subseteq Y^\sigma$ be compact. By Proposition~\ref{prop:pi-weak-tree}, $Q\subseteq Y_N^\sigma$ for some $N$. Choose $j\in\N$ with $N_j\geq N$. Note that every point of $R_\omega\smallsetminus [o^\sigma,p_{j+1}]_{Y^\sigma}$ lies outside $Y_N^\sigma$ since $R_\omega\cap Y_N^\sigma\subseteq R_\omega\cap Y_{N_j}^\sigma \subseteq [o^\sigma,p_{j+1}]_{Y^\sigma}$. In addition, $Q\cap R_\omega=Q\cap (R_\omega\cap Y_N^\sigma)$ is closed and $r$ is continuous. Thus, $r^{-1}(Q)\subseteq r^{-1}(Y_N^\sigma)\subseteq r^{-1}([o^\sigma,p_{j+1}]_{Y^\sigma})$ is closed and bounded in $[0,+\infty)$, and hence compact.

    We use the properness to show that $r$ is a topological embedding. If $A\subseteq[0,+\infty)$ is closed, then for every $N$ the set $A\cap r^{-1}(Y_N^\sigma)$ is compact. Hence, $r(A)\cap Y_N^\sigma=r\bigl(A\cap r^{-1}(Y_N^\sigma)\bigr)$ is compact, and therefore closed, in the finite tree $Y_N^\sigma$. By the definition of the weak topology on $Y^\sigma$ (see Definition~\ref{def:ch4-limit-spaces}), $r(A)$ is closed in $Y^\sigma$. Thus, $r$ is a closed continuous injection. Hence, $r$ is an embedding onto its image.

    Suppose instead that $R_\omega$ is contained in one finite stage $Y_N^\sigma$ for some $N\in\N$. Its closure in $Y_N^\sigma$ is an arc $[o^\sigma,h]_{Y^\sigma}$. Since $R_\omega$ has no terminal point as above, we have $h\notin\inter (Y^\sigma) $ and $R_\omega=[o^\sigma,h)_{Y^\sigma}$, which is homeomorphic to $[0,1)$, and hence to $[0,+\infty)$. Thus, $R_\omega$ is again a topological ray.

    For uniqueness, let $R$ be another rooted ray in $\inter (Y^\sigma)$ representing $\omega$. Since $R\cap C_K^\omega\neq\emptyset$ for every $K\in\mathcal K_{o^\sigma}(Y^\sigma)$ and $R$ is a ray containing $o^\sigma$, so $R$ passes through every $a_K(\omega)$ and therefore $R_\omega\subseteq R$. Conversely, if $z\in R$, choose $K\in\mathcal K_{o^\sigma}(Y^\sigma)$ containing $z$. A sufficiently far tail of $R$ lies in $C_K^\omega$, so $R$ exits $K$ through $a_K(\omega)$. Hence $z\in R\cap K=[o^\sigma,a_K(\omega)]_{Y^\sigma}\subseteq R_\omega$. Thus $R=R_\omega$.

    Finally, suppose that $A_\omega\=\overline{R_\omega}^{\,Y^\sigma}$ is compact. By Proposition~\ref{prop:pi-weak-tree}, $A_\omega$ is in some finite tree. Since it is linearly ordered by the rooted order, $A_\omega=[o^\sigma,h_\omega]_{Y^\sigma}$ for a unique terminal point $h_\omega$. If $h_\omega$ were a cut point, then $[o^\sigma,h_\omega]_{Y^\sigma}$ would be an interior subtree. But the representation property proved above gives a nonempty tail of $R_\omega$ outside the interior subtree $[o^\sigma,h_\omega]_{Y^\sigma}$, contradicting $R_\omega\subseteq[o^\sigma,h_\omega]_{Y^\sigma}$. Hence, $h_\omega$ is a hard endpoint.
\end{proof}

\begin{remark}
    This proposition shows that if $A_\omega$ is compact, then $A_\omega=[o^\sigma, h_\omega]_{Y^\sigma}$ for a unique hard endpoint $h_\omega$. We call $\omega$ a \emph{hard end} and call $h_\omega$ a \emph{hard endpoint}. If $\omega$ is not a hard end, it is called a \emph{genuine end}.
\end{remark}

Now we introduce some properties of the ends of $Y^\sigma$. Let $\mathcal E_g(Y^\sigma)\subseteq\mathcal E(Y^\sigma)$ denote the set of genuine ends.

\begin{lemma}
    For each $\omega\in\mathcal E_g(Y^\sigma)$ and each compact set $Q\subseteq Y^\sigma$, every sufficiently far tail of $R_\omega$ is disjoint from $Q$.
\end{lemma}

\begin{proof}
    By Proposition~\ref{prop:pi-weak-tree}, there exists $N$ such that $Q\subseteq Y_N^\sigma$. Since $\omega$ is genuine, $A_\omega=\overline{R_\omega}^{\,Y^\sigma}$ is noncompact. Hence $R_\omega$ is not contained in $Y_N^\sigma$, so we can choose $p\in R_\omega\smallsetminus Y_N^\sigma$. If a point $q\in R_\omega\smallsetminus[o^\sigma,p]_{Y^\sigma}$ belonged to $Y_N^\sigma$, then $Y_N^\sigma$ would contain the whole arc $[o^\sigma,q]$, since $Y_N^\sigma$ is a subtree containing $o^\sigma$ and $q$. Hence, $p\in[o^\sigma,q]_{Y^\sigma}\subseteq Y_N^\sigma$, which is a contradiction. Thus, the tail $R_\omega\smallsetminus[o^\sigma,p]_{Y^\sigma}$ is disjoint from $Y_N^\sigma$, and therefore from $Q$.
\end{proof}

\begin{lemma}\label{lem:end-in-component}
    Let $K$ be an interior subtree, and let $C$ be a component of $Y^\sigma\smallsetminus K$. Then, there exists $\omega=\{C^\omega_K\}_{K\in\cK(Y^\sigma)}\in\mathcal E(Y^\sigma)$ with $C_K^\omega=C$.
\end{lemma}

\begin{proof}
Take $x_0\in C$ and let $a$ be the attachment point of $C$ to $K$, i.e., $\oC\cap K=\{a\}$. Define an order on $C\cup\{a\}$ by $x\leq_a y$ if and only if $x\in[a,y]$.

By Zorn's lemma, take a maximal chain $R\subseteq C\cup\{a\}$ containing $a$ and $x_0$. By maximality, $[a,x]\subseteq R$ for every $x\in R$. If $R$ has a terminal point $p$, then $R=[a,p]$. If $p$ were a cut point, a point $y$ in a component of $Y^\sigma\smallsetminus\{p\}$ distinct from the one containing $a$ would extend $[a,p]$ to a longer arc $[a,y]\subseteq C\cup\{a\}$. Hence $p$ is a hard endpoint. If $R$ has no terminal point and lies in some $Y_N^\sigma$, its closure there is an arc $[a,p]$. Since $[x_0,p]$ lies outside $K$, we have $p\in C$. Adding $p$ would enlarge $R$, contrary to maximality. Thus in this case $R$ leaves every finite stage. The parametrization by successive arcs in the proof of Proposition~\ref{prop:rooted-representatives}, with $a$ as the root, shows that $R$ is a proper ray.

In either case, for every interior subtree $L$ containing $a$, a sufficiently far tail of $R$ lies in a unique component $C_L$ of $Y^\sigma\smallsetminus L$. These components satisfy that $C_{L'}\subseteq C_L$ for $L\subseteq L'$ containing $a$.

These components induce an end on the family of interior subtrees containing $a$. They also determine a unique end $\omega$ of $Y^\sigma$: for an arbitrary interior subtree $K$, choose an interior subtree $L\supseteq K$ with $a\in L$, and let $C_K$ be the component of $Y^\sigma\smallsetminus K$ containing $C_L$. This $C_K$ is independent of the choice of $L$ since if $L\neq L'\supseteq K$ and $a\in L\cap L'$, we can consider $L\cup L'$ and $C_{L\cup L'}\subseteq C_L\cap C_{L'}$. Since the tail of $R$ lies in $C$, $C_K^\omega=C$.
\end{proof}

\subsection{Spherical ribbon structures}\label{sub_Spherical_ribbon_structures}
In this subsection we show that $Y^\sigma$ carries a ribbon structure. First we recall the notion of the ribbon structure from \cite[Section~3]{Mcmullen09}.

\begin{definition}
    A \emph{germ} at $x\in T$ in a topological $\R$-tree $T$ is a component of $T\smallsetminus\{x\}$. A \emph{ribbon structure} is a cyclic order (see Subsection~\ref{sec_order_completion}) on the germs at each point of degree at least three. For three germs $g_1,\,g_2,\,g_3$ at $x$, if they are in the positive order, we denote it by $\langle g_1,g_2,g_3\rangle=1$. If they are in the negative order, we denote it by $\langle g_1,g_2,g_3\rangle=-1$. A topological $\R$-tree with ribbon structure is called a \emph{ribbon tree}. An embedding is \emph{ribbon preserving} if it preserves every cyclically ordered triple of germs.
\end{definition}

\begin{definition}
    Let $T$ be a tree and let $x_1,\,x_2,\,x_3\in T$ be three distinct points. We call $\operatorname{med}_T(x_1,x_2,x_3)\=[x_1,x_2]_T\cap[x_2,x_3]_T\cap[x_3,x_1]_T$ the medium of the points in the tree $T$.
\end{definition}

Now we equip the topological $\R$-trees $Y_N^\sigma$ and $Y^\sigma$ with ribbon structures. 

\begin{prop}\label{Yribbon}
    For every $\sigma\in\{+,\,-\}$, the following hold:
    \begin{enumerate}
    \smallskip
        \item\label{finiteribbon(i)} Each $Y_n^\sigma$ carries a canonical ribbon structure $\langle \,\cdot\,,\,\cdot\,,\,\cdot\,\rangle_{n,\sigma}$, and every $\iota_n^\sigma$ is ribbon preserving.
    \smallskip
        \item\label{finiteribbon(ii)} There exists a unique ribbon structure $\langle \,\cdot\,,\,\cdot\,,\,\cdot\,\rangle_{\sigma}$ on $Y^\sigma$ whose restriction to $Y_n^\sigma$ is $\langle \,\cdot\,,\,\cdot\,,\,\cdot\,\rangle_{n,\sigma}$.
    \end{enumerate}
\end{prop}

\begin{proof}
     \textup{(i)} We first show that for $n\in\N$, the $n$-connection graph $\Gamma_n^\sigma$ is equipped with a ribbon structure. In fact, there is a canonical way to embed $\Gamma_n^\sigma$ into the sphere $S^2$. We consider the geometric representation of the connection of $n$-tiles. Let $\nu_n\: \Gamma_n^\sigma\to S^2$ satisfy that $\nu_n(c(X))$ is the center point of the $n$-tile $X$ and $\nu_n(c(v,b))$ is a point in the component of the geometric representation $\D(\pi_\w(v)\cup\pi_\b(v))$ of the connection at $v$ corresponding to the block $b$. For $X$ and $b$ with $X\leftrightarrow b$, we join $\nu_n(c(X))$ and $\nu_n(c(v,b))$ by an arc in $S^2$ that is disjoint from the Jordan curve $\partial K_\varepsilon$ in Definition~\ref{Propertiesofconnections}. Then $\nu_n$ maps $\Gamma_n^\sigma$ to a side of the Jordan curve $\partial K_\varepsilon$ determined by $\sigma$. Hence let $\Gamma_n^\sigma$ inherit the order of $\nu_n(\Gamma_n^\sigma)$ induced by the orientation of $S^2$. Moreover, we may assume $\nu_n(c(u,b))=\nu_{n+1}(c(u,\operatorname{Lift}_{n}^{\sigma}(b)))$ by Proposition~\ref{lem:ch4-sector-copies}~\textup{(i)}.  
    
    We then prove $Y_0^\sigma$ carries a canonical ribbon structure $\langle \,\cdot\,,\,\cdot\,,\,\cdot\,\rangle_{0,\sigma}$. By Proposition~\ref{c:initial}, the space $Y_0^\sigma$ is a finite tree and is the union of paths connecting (in the sense of standard topology) $\xi^\sigma_p$ and $\xi^\sigma_q$ for each pair of $p\neq q\in\post f$. Now we consider a vertex $v$ in $Y_0^\sigma$ with degree at least $3$. Fix $m\in\N$. Let $H^\sigma_{m}\=\bigcup_{\lambda\in\sA_0^{\sigma}}P^\sigma_m(\lambda)$ be a subgraph of the $m$-connection graph $\Gamma_m^\sigma$ connecting the vertices $\widehat p_m^\sigma$ representing marked block at postcritical points. Suppose that $C_1,\,C_2,\,C_3$ are three germs at $v$ and $\xi_{p_i}$ is a leaf contained in $C_i$ for $i\in\{1,\,2,\,3\}$. Since $H^\sigma_{m}$ is a tree, the paths $P_m^\sigma(\{\widehat p_i^\sigma,\,\widehat p_j^\sigma\})$ for $i\neq j\in\{1,\,2,\,3\}$ in $\Gamma_m^\sigma$ have a unique intersection vertex $c$, i.e., $c=\operatorname{med}_{\Gamma_m^\sigma}(\widehat p_1^\sigma,\widehat p_2^\sigma,\widehat p_3^\sigma)$. 

    There are two cases to consider. 
    
    \begin{case1}
        The first one is that $c$ is a block vertex $c(u,b)$ in $H^\sigma_{m}$. Here $u\in\V^m$ is an $m$-vertex and $b$ is a block in the connection at $u$. In the path from $c(u,b)$ to $\widehat p_{i,m}^\sigma$, we consider the vertex adjacent to $c(u,b)$. By the construction of the connection graph, this vertex is a tile vertex, and we denote it by $c(X_i)$ for some $X_i\in\X^m_\sigma$. All such tiles contain the vertex $u$, so we can define an order of these tiles which is identified with the standard orientation at $u$ on $S^2$. Letting the germs $C_1,\,C_2,\,C_3$ inherit this order, we obtain the ribbon structure at $v$.
    \end{case1}
    
    \begin{case2}
        The second case is that $c$ is a tile vertex $c(X^m)$ in $H^\sigma_{m}$, for an $m$-tile $X^m$. We still consider the vertex adjacent to $c(X^m)$ in the path from $c(X^m)$ to $\widehat p_{i,m}^\sigma$. By the construction of the connection graph, this vertex is a block vertex, and we denote it by $c(u_i,b_i)$ for some $u_i\in\V^m$. All such vertices are contained in the $m$-tile $X^m$, so we can define an order of these vertices on the boundary of $X^m$ which is identified with the standard orientation on $S^2$. Then let the germs $C_1,\,C_2,\,C_3$ inherit this order, which defines the ribbon structure at $v$.
    \end{case2}

    For the finite tree $Y_n^\sigma$, we consider a vertex $v$ with degree at least $3$. By Proposition~\ref{prop:compact-finite-trees}, the copies $T^\sigma_X$ cover $Y_n^\sigma$ for $X\in \X_\sigma^n$. Suppose $v$ is contained in $T^\sigma_{X_0}$ for $X_0\in \X_\sigma^n$. Suppose that $C_1,\,C_2,\,C_3$ are three germs at $v$ and $\xi_{i}$ is a leaf of $Y_n^\sigma$ contained in $C_i$ for $i\in\{1,\,2,\,3\}$. By Proposition~\ref{prop:compact-finite-trees}, there exists $b_i\in\sB_n^\sigma$ such that $\beta_n^\sigma(b_i)=\xi_i$. Then we consider the paths in $\Gamma_{m+n}^\sigma$ connecting $(\operatorname{Lift}_{n+m-1}^{\sigma}\circ \cdots\circ \operatorname{Lift}_n^{\sigma})(b_i)$ and $(\operatorname{Lift}_{n+m-1}^{\sigma}\circ \cdots\circ \operatorname{Lift}_n^{\sigma})(b_j)$ for $i\neq j\in\{1,\,2,\,3\}$. Similar to the discussion above, we consider the unique intersection vertex and define the order of $C_1,\,C_2,\,C_3$ by the order around the vertex with respect to the orientation on $S^2$.

    Now we show that $\iota_n^\sigma$ is ribbon preserving. We consider a vertex $v$ in $Y_n^\sigma$ with degree at least $3$. Suppose that $C_1,\,C_2,\,C_3$ are three germs of $Y_n^\sigma$ at $v$ and $\xi_{i}$ is a leaf of $Y_n^\sigma$ contained in $C_i$ for $i\in\{1,\,2,\,3\}$. As above, suppose $b_i\in\sB_n^\sigma$ satisfies $\beta_n^\sigma(b_i)=\xi_i$. Since $\iota_n^\sigma\:Y_n^\sigma\to Y_{n+1}^\sigma$ is an embedding, there are germs $C_1',\,C_2',\,C_3'$ at $\iota_n(v)$ with $\iota_n(C_i)\subseteq C_i'$. By Proposition~\ref{prop:compact-stage-stability}, we have $\iota_n(\xi_i)=\iota_n(\beta_n^\sigma(b_i))=\beta_{n+1}^\sigma(\operatorname{Lift}_n^{\sigma}(b_i))$. Let $B_i\=(\operatorname{Lift}_{n+m-1}^{\sigma}\circ \operatorname{Lift}_{n+m-2}^{\sigma}\circ\cdots\circ \operatorname{Lift}_n^{\sigma})(b_i)$ be a vertex in $\Gamma_{m+n}^\sigma$ and $B_i'\=(\operatorname{Lift}_{n+m}^{\sigma}\circ\operatorname{Lift}_{n+m-1}^{\sigma}\circ \cdots\circ \operatorname{Lift}_{n+1}^{\sigma})(\operatorname{Lift}_n^{\sigma}(b_i))$ be a vertex in $\Gamma_{m+n+1}^\sigma$. Put $c\=\operatorname{med}_{\Gamma_{n+m}^\sigma}(B_1,B_2,B_3)$ and we define the order of $C_1,\,C_2,\,C_3$ by the order around $c$. 
    
    First suppose that $c$ is a block vertex $c(u,b)$, then the intersection of the paths in $\Gamma_{m+n+1}^\sigma$ connecting $B_i'$ and $B_j'$ for $i\neq j\in\{1,\,2,\,3\}$ contains a unique point $c(u,\operatorname{Lift}_{n+m}^{\sigma}(b))$ and the order of $C_1',\,C_2',\,C_3'$ is defined by the order around $c(u,\operatorname{Lift}_{n+m}^{\sigma}(b))$. Since $\nu_{n+m}(c(u,b))=\nu_{n+m+1}(c(u,\operatorname{Lift}_{n+m}^{\sigma}(b)))$ and $\iota_n(C_i)\subseteq C_i'$,  we have $\langle C_1,C_2,C_3\rangle_{n,\sigma}=\langle C_1',C_2',C_3'\rangle_{n+1,\sigma}$. 
    
    Next suppose that $c$ is a tile vertex $c(X^{m+n})$, then $C_i$ corresponds to an $(n+m)$-vertex $v_i$ on the boundary of $X^{m+n}$ and $\langle C_1,C_2,C_3\rangle_{n,\sigma}$ is defined by the order of $v_i$ on the boundary of $X^{m+n}$. By Proposition~\ref{lem:ch4-sector-copies}~\textup{(i)}, $L(X^{m+n},1,\,\cdot\,)\:\Gamma_1^\sigma\to\Gamma_{m+n+1}^\sigma[X^{m+n}]$ is an isomorphism. It follows that $c'\=\operatorname{med}_{\Gamma_{n+m+1}^\sigma}(B_1',B_2',B_3')$ is $L(X^{m+n},1,c_0)$, where $c_0\=\operatorname{med}_{\Gamma_{1}^\sigma}(\widehat q_1^\sigma,\widehat q_2^\sigma,\widehat q_3^\sigma)$ $q_i\=f^{n+m}(v_i)\in\post f$. Then, by the discussion above, we know that $\langle C_1',C_2',C_3'\rangle_{n+1,\sigma}$ is defined by the order around $c'$, which is identified with the order of $[c_0,\widehat q_1^\sigma]_{\Gamma_1^\sigma},\,[c_0,\widehat q_2^\sigma]_{\Gamma_1^\sigma},\,[c_0,\widehat q_3^\sigma]_{\Gamma_1^\sigma}$ around $c_0$ in $\Gamma_1^\sigma$. Since the embedding $\nu_1(\Gamma_1^\sigma)$ is on one side of the Jordan curve $\partial K_\varepsilon$, which is isotopic to $\cC$, the order of $[c_0,\widehat q_1^\sigma]_{\Gamma_1^\sigma},\,[c_0,\widehat q_2^\sigma]_{\Gamma_1^\sigma},\,[c_0,\widehat q_3^\sigma]_{\Gamma_1^\sigma}$ around $c_0$ is in fact identified with the order of $q_1,\,q_2,\,q_3$ on $\cC$. Hence, we still have $\langle C_1,C_2,C_3\rangle_{n,\sigma}=\langle C_1',C_2',C_3'\rangle_{n+1,\sigma}$.

\smallskip
     \textup{(ii)} Every three germs in $Y^\sigma$ have a compact finite hull, hence lie in one finite stage $Y_n^\sigma$ for some $n\in\N$ by Proposition~\ref{prop:pi-weak-tree}; by the stage compatibility, we have the unique global ribbon structure $\operatorname{Rib}^\sigma$.
\end{proof}

\begin{remark}
   Although in the construction above of the ribbon structure $\langle \,\cdot\,,\,\cdot\,,\,\cdot\,\rangle_{0,\sigma}$ we fix $m\in\N$, $\langle \,\cdot\,,\,\cdot\,,\,\cdot\,\rangle_{0,\sigma}$ is independent of the choice of $m$ by arguments similar to those in the proof of ``every $\iota^\sigma_n$ is ribbon pre-serving''. 
\end{remark}

Now we can define the left and right sides of an interval in $Y^\sigma$ with respect to the ribbon structure.

\begin{definition}\label{sidedef}
    Let $I\=[a,b]_{Y^\sigma}\subseteq Y^\sigma$ be oriented from $a$ to $b$ and $x\in(a,b)_{Y^\sigma}$. If $g_x^-,\,g_x^+$ are the germs toward $a,\,b$, we say that a third germ $g$ lies on the \emph{left} of $I$ when $\langle g_x^+,g,g_x^-\rangle_{\sigma}=+1$, and on the \emph{right} when $\langle g_x^+,g,g_x^-\rangle_{\sigma}=-1$.
\end{definition}

\subsection{Hairiness}\label{subsec:hairiness}
First we show that for an expanding Thurston map $f$ with invariant Jordan curve $\cC\supseteq\post f$, we can construct a connection of $1$-tiles for $(f,\cC)$ (or $(f^n,\cC)$ for some $n\in\N$) with ``branching'' (see Lemma~\ref{degree3}). With this connection of $1$-tiles, we then show the associated $Y^\sigma$ is hairy.

Recall that the vertices in the associated connection graph are $c(X)$ representing a tile $X$, and $c(v,b)$ representing the block $b$ at a vertex $v$. By deleting the vertex $c(v,b)$ if the block $b$ contains only one label of a $1$-tile containing $v$, we obtain a tree, which is called the \emph{induced connection graph}. In what follows, we show that there exist connections of $1$-tiles that satisfy the conditions in Definition~\ref{Propertiesofconnections}, and, in addition, the induced connection graph has a vertex with degree at least $3$. Here, the \emph{degree} of a vertex in the graph is the number of edges incident to it.

\begin{lemma}\label{degree3}
    Let $f_0$ be an expanding Thurston map, with an invariant Jordan curve $\cC$ satisfying
$\post f_0\subseteq \cC$. Then there is a connection of $1$-tiles for $(f_0,\cC)$ (or $(f_0^m,\cC)$ for some $m\in\N$) satisfying the conditions in Definition~\ref{Propertiesofconnections}. Moreover, the induced white and black connection graphs determined by this connection of $1$-tiles both contain a vertex of degree at least $3$.
\end{lemma}
\begin{proof}
    Let $n\in\N$ be large enough for Lemma~\ref{7.5} and $\cC$ be the Jordan curve. \cite[Section~7]{Meyer13} gives a method to construct connections of $1$-tiles for $(f_1\=f_0^n,\cC)$ satisfying the conditions in Definition~\ref{Propertiesofconnections}.
    We make adjustments to this method, so that the connection additionally satisfies the conditions of this lemma. We may replace $f_1$ by an iterate when necessary.

    First, we briefly recall the construction of connections of $1$-tiles following \cite[Section~7]{Meyer13}. Consider all the white $1$-tiles for $(f_1,\cC)$ contained in the white $0$-tile $X_+^0$. We assume that these~$1$-tiles are connected at every~$1$-vertex $v$ in the interior of $X_+^0$ and disconnected at every~$1$-vertex on $\cC$. Then, by \cite[Lemma~7.6]{Meyer13}, there exists a unique (white) cluster $K$ that intersects all sides ($0$-edges), called the \emph{main cluster}. Here, we regard $K$ as the union of its constituent $1$-tiles. The other clusters are called the \emph{secondary clusters}. The next step is to choose an arbitrary spanning tree of $K$, called the \emph{main tree}. We also choose an arbitrary spanning tree of each secondary cluster, called the \emph{secondary tree}. The connections at all~$1$-vertices $v\in X_+^{0}\smallsetminus \cC$ are then fixed and will remain unchanged throughout the subsequent construction of a connection of $1$-tiles using Meyer's method.
    
    For each $l\in \N$, we consider $\bigl(f_1^l,\cC\bigr)$ and define the white main cluster $K_l$ and other white secondary clusters contained in $X_+^0$ as above. We then claim the following: 
    \begin{claim}
        For some $l_0\in\N$, there exists a white $1$-tile for $\bigl(f_1^{l_0},\cC\bigr)$ contained in the main cluster $K_{l_0}$ or a secondary cluster that connects to at least three other $1$-tiles. We denote this cluster by $K'$.
    \end{claim}
    
    \begin{proof}[Proof of the Claim]
    To prove this, we first show that, for some $l_1\in\N$, there exist at least three white $1$-tiles for $\bigl(f_1^{l_1},\cC\bigr)$ contained in $K_{l_1}$ that have nonempty intersection with $\cC$. We argue by contradiction and suppose that for all $l\in\N$, there are at most two white $1$-tiles for $\bigl(f_1^{l},\cC\bigr)$ contained in $K_{l}$ that have nonempty intersection with $\cC$. By~\ref{7.5(ii)} of Lemma~\ref{7.5}, a white $1$-tile can reach one $0$-edge or two adjacent $0$-edges. Then there are at most four~$0$-edges, so $\card\post f_1^{l}=\card\post f_1=\card\post f_0$ is $3$ or $4$, and there is a white $1$-tile contained in $K_{l}$ that reaches two adjacent $0$-edges. Let $\rho$ be a visual metric for $f_1$. Since the number of $0$-edges is finite, there exists a constant $\varepsilon_0$ such that $\rho(p,E)>\varepsilon_0$ for all $p\in\post f_1$ and~$0$-edge $E$ with $p\notin E$. 
    
    Now we consider a postcritical point $p\in\post f_1$ and $0$-edges $E_1\ni p$ and $E_2\ni p$. We define $z\:E_1\to \R$ by $z(u)\=\rho(u,E_2)$ for all $u\in E_1$. Then $z(u)>0$ whenever $u\neq p$. Since $E_1$ is closed and $z$ is continuous, for all real numbers $\eta>0$, we can define a positive real number $\delta_{p,E_1}(\eta)\=\sup\{\rho(u,p):u\in E_1\text{ and }z(u)\leq\eta\}$. Note that $\delta_{p,E_1}(\eta)\to0$ as $\eta\to 0$. Let $\varepsilon_1$ be a constant such that $\delta_{p,E_1}(\eta)<\varepsilon_0/2$ for all $\eta\leq \varepsilon_1$, all $p\in\post f_1$ and all $0$-edge $E_1\ni p$. 
    
    Since $f_1$ is expanding and $l_1$-tiles for $(f_1,\cC)$ are $1$-tiles for $\bigl(f_1^{l_1},\mathcal{C}\bigr)$ (see Lemma~\ref{lem:cell mapping properties of Thurston map}), when $l_1\in\N$ is large enough, we have that $\max\bigl\{ \diam_{\rho}(X) : X \in \X^1\bigl(f_1^{l_1},\mathcal{C}\bigr) \bigr\}<\min\{\varepsilon_0/2,\,\varepsilon_1\}$. Let $X$ be a white $1$-tile contained in $K_{l_1}$ that meets two adjacent $0$-edges $E$ and $E'$, and let $X'$ be the other contained in $K_{l_1}$ that intersects $\cC$. Suppose $E\cap E'=\{p\}\subseteq\post f_1$. Since $\card\post f_1\geq3$, $X'$ has nonempty intersection with a $0$-edge $E''$ different from $E$ and $E'$. It follows that $X'\cap E=\emptyset$ or $X'\cap E'=\emptyset$ by~\ref{7.5(ii)} of Lemma~\ref{7.5}. Without loss of generality, we suppose that $X'\cap E=\emptyset$. Suppose the other endpoint of $E$ is $p'\in\post f_1$. Consider the component $B$ (in the standard topological sense) of $X_+^0\smallsetminus K_{l_1}$ that contains $p'$. We call the set $a\=\partial B\cap K_{l_1}$ a boundary arc of $K_{l_1}$. Then $a$ is a closed arc joining $A\in E$ to $A'\in\cC\smallsetminus(E\cup E')$, since every $1$-tile in $K_{l_1}$ other than $X$ and $X'$ does not intersect $\cC$. It is shown in the proof of \cite[Lemma~7.6]{Meyer13} that every boundary arc $a$ as above is contained in a single black $1$-tile. Then $\diam_{\rho}(a)<\varepsilon_0/2$. Since $\rho(A,E')\leq\diam_{\rho}(X)\leq\varepsilon_1$, we have that $\rho(A,p)\leq \delta_{p,E}(\varepsilon_1)<\varepsilon_0/2$. Since $p\notin A'$, we have $\rho(p,A')\geq\varepsilon_0$. Thus, we have $\diam_{\rho}(a) \geq \rho (A,A')\geq \rho(p,A')-\rho(A,p)\geq\varepsilon_0-\varepsilon_0/2=\varepsilon_0/2$, which is a contradiction.

\begin{figure}[htbp]
    \centering
    \begin{overpic}[
        percent,
        width=\linewidth,
        trim=0bp -14bp 0bp -14bp
    ]{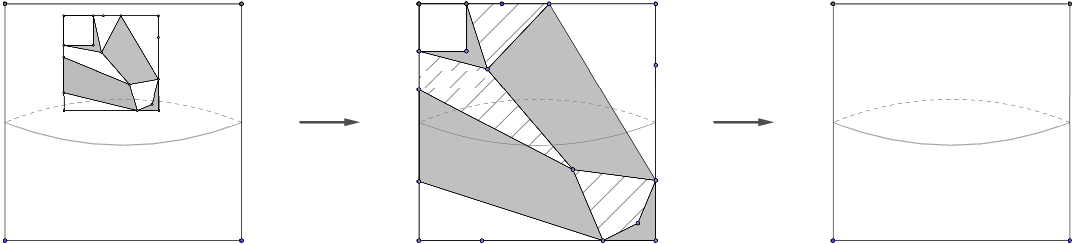}
        \setlength{\unitlength}{0.01\linewidth}

        \tiny
        \put(7.15,22.78){\makebox(0,0){\scalebox{0.75}{$U_1$}}}
        \put(9.85,22.78){\makebox(0,0){\scalebox{0.75}{$U_2$}}}
        \put(13.56,16.77){\makebox(0,0){\scalebox{0.75}{$U_3$}}}
        \put(4.65,21.32){\makebox(0,0){\scalebox{0.75}{$w_0$}}}
        \put(8.28,19.86){\makebox(0,0){\scalebox{0.75}{$W$}}}
        \put(15.86,16.04){\makebox(0,0){$Z$}}

        \scriptsize
        \put(48.79,26.94){\makebox(0,0){$K_{l_1}$}}
        \put(41.21,22.95){\makebox(0,0){$W'$}}
        \put(37.78,20.64){\makebox(0,0){$v_0$}}
        \put(43.40,17.44){\makebox(0,0){\scalebox{0.72}{$f_1^{kl_1}(W)$}}}

        \put(30.64,16.15){\makebox(0,0){$f_1^{kl_1}$}}
        \put(68.91,16.15){\makebox(0,0){$f_1^{l_1}$}}

        \put(75.43,4.74){\makebox(0,0){$X^0_{\mathfrak{w}}$}}

        \put(11.14,1.59){\makebox(0,0){$\mathcal{C}$}}
        \put(50.37,1.59){\makebox(0,0){$\mathcal{C}$}}
        \put(89.70,1.59){\makebox(0,0){$\mathcal{C}$}}

    \end{overpic}
    \caption{We identify $S^2$ with a ``pillow'' obtained by gluing two squares together along their boundary. Moreover, the two squares are the $0$-tiles: the square on the front side is the white $0$-tile $X_\w^0$, and the other one is the black $0$-tile $X_\b^0$. Their common boundary is the invariant Jordan curve $\cC$. The map $f_1^{l_1}$ maps the white $1$-tiles in the middle figure to $X_\w^0$ and the black $1$-tiles to $X_\b^0$. The union of the white tiles with diagonal hatching is the cluster $K_{l_1}$. The map $f_1^{kl_1}$ maps the white tile $Z$ in the left figure to $X_\w^0$. The map $f_1^{(k+1)l_1}$ maps the white tiles contained in $Z$ to $X_\w^0$.}
    \label{fig:branch}
\end{figure}

    If there is a white~$1$-tile for $\bigl(f_1^{l_1},\cC\bigr)$ contained in $K_{l_1}$ that connects to at least three other $1$-tiles, then let $l_0\=l_1$. If not, every white $1$-tile contained in $K_{l_1}$ is connected to at most two other $1$-tiles. Moreover, it is easy to see that at most two $1$-tiles contained in $K_{l_1}$ are connected to exactly one other $1$-tile. We consider the iterates $f'\=f_1^{kl_1}$ such that for $\bigl(f',\cC\bigr)$ there is a white $1$-tile $Z$ (a white $k$-tile for $\bigl(f_1^{l_1},\cC\bigr)$) with $Z\cap \cC=\emptyset$. Then we assume that the white $1$-tiles for $\bigl(f_1^{(k+1)l_1},\cC\bigr)$ (white $(k+1)$-tile for $\bigl(f_1^{l_1},\cC\bigr)$) in the white $0$-tile $X_+^{0}$ are all connected at all $1$-vertices $v$ in the interior of $X_+^{0}$, and disconnected at all $1$-vertices on $\cC$. Since $f_1^{kl_1}(Z)=X_+^{0}$, it follows that the tiles in $f_1^{-kl_1}(K_{l_1})\cap Z$ have the same connection pattern as $K_{l_1}$. Let $K'$ be the white cluster for $\bigl(f_1^{(k+1)l_1},\cC\bigr)$ that contains $f_1^{-kl_1}(K_{l_1})\cap Z$. Since there are at most two $1$-tiles for $\bigl(f_1^{(k+1)l_1},\cC\bigr)$ in $f_1^{-kl_1}(K_{l_1})\cap Z$ that connect to exactly one $1$-tile in $f_1^{-kl_1}(K_{l_1})\cap Z$, and there are at least three white $1$-tiles for $\bigl(f_1^{l_1},\cC\bigr)$ contained in $K_{l_1}$ that have nonempty intersection with $\cC$, we know that there is a white $1$-tile $W$ for $\bigl(f_1^{(k+1)l_1},\cC\bigr)$ in $f_1^{-kl_1}(K_{l_1})\cap Z$ that connect to two $1$-tile in $f_1^{-kl_1}(K_{l_1})\cap Z$ and has nonempty intersection with $f_1^{-kl_1}(\cC)\cap Z$. Note that the boundary point $v_0\in f_1^{kl_1}(W) \cap\cC$ is a critical point of $f_1^{l_1}$ since there is other white $1$-tile $W'$ (for $\bigl(f_1^{l_1},\cC\bigr)$) containing $v_0$, or $v_0$ is not the boundary point. By the construction of $K_{l_1}$, $f_1^{kl_1}(W)$, and $W'$ are not connected since $v_0\in\cC$. Note that $f_1^{-kl_1}(v_0)\cap Z$ contains only one point $w_0$, which is still a critical point for $f_1^{(k+1)l_1}$. Then this tile $W$ connects to a white $1$-tile for $\bigl(f_1^{(k+1)l_1},\cC\bigr)$ at $w_0$ (a preimage of $W'$ under $f_1^{-kl_1}$ that contains $w_0$) outside $f_1^{-kl_1}(K_{l_1})\cap Z$. Thus let $l_0\=(k+1)l_1$ we have that $W$ is a white $1$-tile for $\bigl(f_1^{l_0},\cC\bigr)$ contained in the cluster $K'$ that connects to at least three other~$1$-tiles, denoted by $U_1,\,U_2,\,U_3$. See Figure~\ref{fig:branch} for an example. Note that for all $n\in\N$, we also have that $f_1^{-n}(W)$ contains a white $1$-tile for $\bigl(f_1^{l_0+n},\cC\bigr)$ contained in $X_\w^0$ that connects to at least three other $1$-tiles in $f_1^{-n}(U_1)\cap X_\w^0,\,f_1^{-n}(U_2)\cap X_\w^0,\,f_1^{-n}(U_3)\cap X_\w^0$, respectively. Hence, we can take $l_0$ to be each integer larger than $(k+1)l_1$.

    \end{proof}
    
    We can choose one spanning tree of $K'$ whose induced connection graph contains a vertex of degree at least $3$. More precisely, suppose $X_1$, $X_2$, and $X_3$ are white $1$-tiles for $\bigl(f_1^{l_0},\cC\bigr)$ contained in $K'$ that are all connected to a white $1$-tile $X_0$. Besides, we suppose that $c(X_0)$ and $c(X_i)$ are both incident to $c(v_i,b_i)$, where $v_i\in X_0\cap X_i$ is a $1$-vertex and $b_i\in \pi^1_{\white}(v_i)$ for $i\in\{1,\,2,\,3\}$. If there is a loop in the induced connection graph corresponding to $K'$, suppose there are three adjacent vertices $c(Y_1)-c(v,b)-c(Y_2)$ in the loop such that $v\in Y_1\cap Y_2$ is a~$1$-vertex and $b\in \pi^1_{\white}(v)$. In addition, we suppose that $Y_1$ has a label $x$ in the connection at $v$. Then we define a new connection at $v$ with the same labeling as before. We break up $b$ into $\{x\}$ and $b\smallsetminus\{x\}$ and let $\pi^{1'}_{\white}(v)$ be the set of $\{x\}$ and $b\smallsetminus\{x\}$ together with blocks different from $b$ in $\pi^1_{\white}(v)$. Then let $\pi^{1'}_{\b}(v)$ be the partition complementary to $\pi^{1'}_{\white}(v)$. Then the cluster (still denoted by $K'$) remains the same as a union of $1$-tiles, and the corresponding induced connection graph has fewer loops than before. Repeating the procedure until there is no loop in the corresponding induced connection graph, we obtain a spanning tree $K'$. Since $X_1$, $X_2$, and $X_3$ are pairwise distinct, the loop described above cannot consist solely of paths of the form $c(X_i)-c(v_i,b_i)-c(X_0)$. Hence, we can retain all three subpaths $c(X_i)-c(v_i,b_i)-c(X_0)$, for $i\in\{1,\,2,\,3\}$, throughout the procedure. Consequently, these three subpaths are contained in the resulting spanning tree $K'$, and the corresponding induced connection graph has a vertex of degree at least $3$.
    
    Since the connections at all $1$-vertices $v\in X_+^{0}\smallsetminus \cC$ will remain unchanged in the subsequent construction of a connection of $1$-tiles using Meyer's method in \cite[Section~7]{Meyer13}, constructing the connections at the remaining $1$-vertices in this way ensures that the white induced connection graph satisfies the conditions stated in this lemma. 
    
    The connections at all $1$-vertices $v\in X_-^0\smallsetminus \cC$ are constructed analogously. We assume that the black $1$-tiles for $\bigl(f_1^{l_0},\cC\bigr)$ in the black $0$-tile $X_-^{0}$ for $\bigl(f_1^{l_0},\cC\bigr)$ are all connected at all $1$-vertices $v$ in the interior of $X_-^{0}$, and disconnected at all $1$-vertices on $\cC$. We also consider the black cluster $L$ that intersects all $0$-edges. As in the white case, we can take an iterate of $f_1^{l_0}$ if necessary such that there is a black $1$-tile for $\bigl(f_1^{l_0},\cC\bigr)$ contained in $L$ that connects to at least three other $1$-tiles. Then we can define a spanning tree of $L$ whose induced connection graph contains a vertex of degree at least $3$. As above, the connections at all $1$-vertices $v\in X_-^{0}\smallsetminus \cC$ will not be changed further during the subsequent construction. Thus, we can construct a connection of the $1$-tiles for $\bigl(f\=f_0^{nl_0},\cC\bigr)$ (or $(f^s,\cC)$ for some $s\in\N$) so that they satisfy all the conditions in this lemma simultaneously.
\end{proof}

\begin{definition}\label{def:ch5-hairy}
A ribbon tree $T$ is \emph{hairy} if, for each side of every nondegenerate oriented arc $I\subseteq T$, there exist $x\in\inte_\circ I$ and a proper ray $R$ from $x$ such that $R\smallsetminus\{x\}$ is contained in a germ on that side at $x$.
\end{definition}

Write $\mathbf A^\infty\=\bigcup_{m=0}^{+\infty}\mathbf A^m$, where $\mathbf A^m$ is the set of $m$-angles. For $\sigma\in\{+,\,-\}$ and $t\in\mathbf A^\infty$, choose an integer $m\ge1$ with $t\in\mathbf A^m$ and let $b\in\sB_m^\sigma$ be the unique block with $t\in\cA_m^\sigma(b)$. By~\eqref{eq:ch4-block-persistence}, $t\in\cA_{m+1}^\sigma(\operatorname{Lift}_m^\sigma b)$. Equations~\eqref{eq:ch4-copy-stability} and~\eqref{eq:ch4-refined-mark} give $\beta_{m+1}^\sigma(\operatorname{Lift}_m^\sigma b)=\iota_m^\sigma(\beta_m^\sigma(b))$. Thus $x^\sigma(t)\=\beta_m^\sigma(b)$ defines a point of $Y^\sigma$ independent of $m$.

\begin{lemma}\label{lem:ch5-marked-turns}
Fix $\sigma\in\{+,\,-\}$. Let $Y^\sigma$ be the topological $\R$-tree of Definition~\ref{def:ch4-limit-spaces} for the fixed connection of $1$-tiles. Use its canonical ribbon structure from Proposition~\ref{Yribbon}.

Let $I\subseteq Y^\sigma$ be a nondegenerate oriented arc. Suppose that no germ at any point of $\inte_\circ I$ lies on the right of $I$. Then there exist an integer $m\ge1$ and a nondegenerate closed arc $A\subseteq S^1$ with endpoints in $\mathbf A^m$ such that for every integer $n\ge m$ and every block $b\in\sB_n^\sigma$,
\begin{equation*}
A\cap\cA_n^\sigma(b)\ne\emptyset
\quad\to \quad
\beta_n^\sigma(b)\in I.
\end{equation*}
In particular, $x^\sigma(A\cap\mathbf A^\infty)\subseteq I$.
\end{lemma}

\begin{proof}
By Proposition~\ref{prop:pi-weak-tree}~\textup{(i)}, take an integer $N\ge1$ with $I\subseteq Y_N^\sigma$. Finitely many road images cover $Y_N^\sigma$, so some $\eta_N^\sigma(I_\lambda^\sigma)$ contains a nondegenerate subarc of $\inte_\circ I$. Since the mesh of the subdivisions of $I_\lambda^\sigma$ tends to zero by Lemma~\ref{prop:ch4-road-data}~\textup{(i)--(ii)}, we may take an integer $m>N$ so that this subarc contains four points $u,\,x,\,y,\,v\in\eta_N^\sigma\bigl(S_m^{\lambda,\sigma}\bigr)$ in this given order along $I$. For each such $z=\eta_N^\sigma(s)$, put $b_z\=b_m^{\lambda,\sigma}(s)$. Then $\beta_m^\sigma(b_z)=z$ by~\eqref{eq:ch4-stage-model-map} and~\eqref{eq:compact-stage-stability}.

For $z=x$ or $z=y$, let $c(X)$ and $c(X')$ be the tile vertices immediately before and after $b_z$ on the path from $b_u$ to $b_v$ in $\Gamma_m^\sigma$. Recall the notation $D_z^\sigma(\,\cdot\,)$ from Definition~\ref{def:ch4-cut-sets}. By~\eqref{eq:ch4-stage-model-map}, $T_X^\sigma\cap D_z^\sigma(u)\ne\emptyset$ and $T_{X'}^\sigma\cap D_z^\sigma(v)\ne\emptyset$. If $U\leftrightarrow b_z$ has index strictly between those of $X$ and $X'$ in positive cyclic order, take $w\in T_U^\sigma\smallsetminus\{z\}$. The copies meet pairwise only at $z$ by~\eqref{eq:ch4-stage-copy-intersection}, so $D_z^\sigma(w)$ is on the right of $I$ by Propositions~\ref{Yribbon} and~\ref{lem:ch4-sector-copies}~\textup{(iii)}. Thus the index of $X'$ succeeds that of $X$ at $b_z$.

Orient $S^1$ so that $\gamma^\ell$ traverses the edges of $\partial_\circ X$ positively for every integer $\ell\ge1$ and $X\in\X_\sigma^\ell$. Replace each $\ell$-arc mapped into $X$ by $b,\,c(X),\,b'$, where $b,b'\in\sB_\ell^\sigma$ are the blocks whose angle classes contain the initial and terminal angles of the arc, respectively. For $X\leftrightarrow b$, the two edges of $\partial_\circ X$ at $\prV(b)$ contribute $c(X)\to b$ and $b\to c(X)$ once each. The visits to $b$ occur at the angles in $\cA_\ell^\sigma(b)$, so $\card\cA_\ell^\sigma(b)=\deg_{\Gamma_\ell^\sigma}(b)$. The circuit follows $\partial_\circ X$ at $c(X)$ and takes the positive successor at $b$ (Subsection~\ref{prelimiconnectiongraph} and Lemma~\ref{lem:finite-class-successors}).

With $S^1$ oriented as above, for each $z\in\{x,\,y\}$, let $t_z\in\cA_m^\sigma(b_z)$ be the unique angle for which $\gamma^m$ maps the $m$-arcs ending and starting at $t_z$ into the tiles $X$ and $X'$ chosen above for $z$, respectively. Since $x^\sigma(t_x)=x$ and $x^\sigma(t_y)=y$, we have $t_x\ne t_y$. Let $A$ be the closed arc from $t_x$ to $t_y$ in this orientation.

Fix an integer $n\ge m$, $b\in\sB_n^\sigma$, and $t\in A\cap\cA_n^\sigma(b)$. For a term $r=m_j$ of the subsequence in the proof of Proposition~\ref{c:initial}~\textup{(ii)}, write $B'\=(\operatorname{Lift}_{n+r-1}^\sigma\circ\cdots\circ\operatorname{Lift}_k^\sigma)(B)$ for $B\in\sB_k^\sigma$ with $k\in\{m,\,n\}$. As in the proof of Proposition~\ref{Yribbon}~\textup{(i)}, take $j$ large enough that, for every triple of distinct leaves $\beta_n^\sigma(B_i)$ of $Y_n^\sigma$, the ribbon order of the three germs toward these leaves at their median agrees with the cyclic order of the edges toward $B_i'$ at $\operatorname{med}_{\Gamma_{n+r}^\sigma}(B_1',B_2',B_3')$.

By Lemma~\ref{prop:ch4-road-data}~\textup{(i)}, $b_x',\,b_y'$ occur in this order on $P\=[b_u',b_v']_{\Gamma_{n+r}^\sigma}$ from $b_u'$ to $b_v'$. Write $Q\=[b_x',b_y']_{\Gamma_{n+r}^\sigma}$. The circuit at level $n+r$ visits $b'$ at $t$ and $b_z'$ at $t_z$ by~\eqref{eq:ch4-block-persistence}. For each $z\in\{x,\,y\}$, use the tiles $X$ and $X'$ chosen above for $z$. By~\eqref{eq:ch4-affine-refinement}, $P$ enters $b_z'$ in $\Gamma_{n+r}^\sigma[X]$ and leaves in $\Gamma_{n+r}^\sigma[X']$. So does the circuit at $t_z$ by the refinement of the adjacent arcs in the proof of Lemma~\ref{lem:ch4-connection-inputs}. Since the incidences of each copy at $b_z'$ are consecutive by Proposition~\ref{lem:ch4-sector-copies}~\textup{(iii)}, the passage at $t_z$ occurs between the incoming and outgoing incidences of $P$ in positive cyclic order.

Since $\Gamma_{n+r}^\sigma$ is a tree, the circuit returns along each edge by which it leaves $P$. At an interior vertex it thus takes the other edge of $P$ after the intervening incidences in positive cyclic order. This places $t_x$ and $t_y$ on the traversal toward $b_v'$. The edges of $P$ traversed between $t_x$ and $t_y$ are therefore exactly those of $Q$, in order from $b_x'$ to $b_y'$.

By Proposition~\ref{lem:ch4-sector-copies}~\textup{(i)} and~\eqref{eq:ch4-copy-values}, the continued level-$n$ blocks span $G=\bigcup_{X\in\X_\sigma^n}L(X,r,H_r^\sigma)$. On each copy, apply $L(X,r,\,\cdot\,)^{-1}$, the reduction and contraction $H_r^\sigma\to Y_0^\sigma$ in the proof of Proposition~\ref{c:initial}~\textup{(ii)}, and then $\theta_X^\sigma$. These maps agree at shared blocks since each sends $B'$ to $\beta_n^\sigma(B)$ for $X\leftrightarrow B$. They therefore define $q\:G\to Y_n^\sigma$. By~\eqref{eq:ch4-stage-copy-intersection}, their connected fibers over the same point meet at a shared block, so $q$ sends paths onto arcs. For $z\in\{u,\,x,\,y,\,v\}$, let $B\in\sB_n^\sigma$ be the continuation of $b_z$ to level $n$, so $B'=b_z'$. By definition, $q(b_z')=\beta_n^\sigma(B)$. Iterating~\eqref{eq:ch4-copy-stability} and~\eqref{eq:ch4-refined-mark} gives $\beta_n^\sigma(B)=\beta_m^\sigma(b_z)=z$ in $Y^\sigma$. Hence, $q(Q)=[x,y]$.

If $q(b')\notin[x,y]$, the circuit leaves $P$ at $w\in Q$ on its way to $b'$. On the traversal toward $b_v'$, the edges at $w$ toward $b_u',\,b',\,b_v'$ occur in positive cyclic order. Since $q$ sends paths onto arcs, $z=q(w)\in[x,y]$ is the median of $u,\,q(b'),\,v$. Take leaves $\beta_n^\sigma(B_i)$ of $Y_n^\sigma$ in $D_z^\sigma(u),\,D_z^\sigma(q(b')),\,D_z^\sigma(v)$ for $i=1,\,2,\,3$, respectively. The paths from $B_i'$ to the corresponding blocks $b_u',\,b',\,b_v'$ avoid $w$ since their images avoid $z$. The order comparison at the common median $w$ for $r=m_j$ therefore yields $\langle D_z^\sigma(v),D_z^\sigma(q(b')),D_z^\sigma(u)\rangle_\sigma=-1,$
contrary to the hypothesis at $z\in\inte_\circ I$. Hence $\beta_n^\sigma(b)=q(b')\in[x,y]\subseteq I$.
\end{proof}

\begin{prop}[Hairiness]
\label{lem:hairy}
Let $f\:S^2\to S^2$ be an expanding Thurston map and $\cC\supseteq\post f$ an $f$-invariant Jordan curve. Fix a connection of $1$-tiles satisfying \textup{(C1)}--\textup{(C2)} of Definition~\ref{Propertiesofconnections}. For $\sigma\in\{+,\,-\}$, let $Y^\sigma\=\bigcup_{n\in\N_0}Y_n^\sigma$ be the topological $\R$-tree of Definition~\ref{def:ch4-limit-spaces} with the canonical ribbon structure of Proposition~\ref{Yribbon}. At each level $m\in\N$, use the induced connection graphs in Lemma~\ref{degree3}.

\begin{enumerate}
\smallskip\item For each $\sigma$, either $Y^\sigma=Y_0^\sigma$ with $Y_0^\sigma$ a compact interval, or $Y^\sigma$ is hairy in the sense of Definition~\ref{def:ch5-hairy}. The alternatives are mutually exclusive.
\smallskip\item For each $\sigma$, $Y^\sigma=Y_0^\sigma$ if and only if $Y_1^\sigma=Y_0^\sigma$, in which case $Y_n^\sigma=Y_0^\sigma$ for all $n\in\N_0$.
\smallskip\item For each $\sigma$, $Y^\sigma$ is hairy if and only if, for some $m\in\N$, the induced $m$-connection graph of color $\sigma$ has a vertex of degree at least three.
\smallskip\item At least one of $Y^+$ and $Y^-$ is hairy. Both are hairy if the induced $1$-connection graph of each color has a vertex of degree at least three.
\end{enumerate}
\end{prop}

\begin{proof}
Fix $\sigma\in\{+,\,-\}$. For $t\in\mathbf A^\infty$, take an integer $m\ge2$ and $b\in\sB_m^\sigma$ with $t\in\cA_m^\sigma(b)$. For the block $b'\in\sB_{m-1}^\sigma$ with class $\mu(\cA_m^\sigma(b))$, the proof of Proposition~\ref{prop:pi-weak-tree}~\textup{(v)} gives
\begin{equation}\label{eq:ch5-angle-point-dynamics}
f_\sigma(x^\sigma(t))=\beta_{m-1}^\sigma(b')=x^\sigma(\mu(t)).
\end{equation}

\textup{(i)} First, suppose that $Y^\sigma=Y_0^\sigma$ and let $L$ be its number of leaves. Recall that $Y_1^\sigma$ is constructed from $d$ copies of $Y_0^\sigma$ by attaching one copy at a time at a single point (see the proof of Proposition~\ref{prop:compact-finite-trees}~\textup{(i)}). Each of the $d-1$ attachments removes at most two leaves, so $L\ge dL-2(d-1)$. Since $d,\,L\ge2$, we have $L=2$, and $Y^\sigma$ is a compact interval with no proper ray.

Now suppose $Y^\sigma\ne Y_0^\sigma$. If a nondegenerate oriented arc $I$ has no germ on its right at any interior point, Lemma~\ref{lem:ch5-marked-turns} gives a nondegenerate arc $A\subseteq S^1$ with $x^\sigma(A\cap\mathbf A^\infty)\subseteq I$. For large $k$, we have $I\subseteq Y_k^\sigma$ and $\mu^k(A)=S^1$, so $\mu^k(A\cap\mathbf A^\infty)=\mathbf A^\infty$ by~\eqref{eq:finite-angle-preimages}. By~\eqref{eq:ch5-angle-point-dynamics} and~\eqref{eq:pi-exact-pullback},
\begin{equation*}
x^\sigma(\mathbf A^\infty)=f_\sigma^k(x^\sigma(A\cap\mathbf A^\infty))\subseteq Y_0^\sigma.
\end{equation*}
For an integer $m\in\N$ and $X\in\X_\sigma^m$, the marks $\theta_X^\sigma\bigl(\xi_p^\sigma\bigr)=\beta_m^\sigma(b_X^\sigma(p))$ for $p\in\post f$ belong to $x^\sigma(\mathbf A^\infty)$ and span $T_X^\sigma$ by Proposition~\ref{c:initial}~\textup{(iii)}. Thus $Y^\sigma=Y_0^\sigma$, which is a contradiction. Since either orientation can be used, each side of a nondegenerate arc has a germ at some interior point.

If $D_p^\sigma(w)\subseteq Y_N^\sigma$, take a nondegenerate subarc of $(p,w)$ avoiding the branch points of $Y_N^\sigma$. A germ on one side at an interior point of this subarc lies in $D_p^\sigma(w)$, so that point would be a branch point of $Y_N^\sigma$. Thus some $w'\in D_p^\sigma(w)$ lies outside $Y_N^\sigma$, and by closedness so does a cut point in $(p,w')$.

Fix a side of a nondegenerate oriented arc $I$. Take a germ $D_{p_0}^\sigma(w_0)$ on that side at $p_0\in\inte_\circ I$ and start with $p_1\in(p_0,w_0)$. For $j\ge1$, take $N_j\ge j$ with $p_0,\,p_j\in Y_{N_j}^\sigma$. Since $p_j$ is a cut point, a germ at $p_j$ distinct from $D_{p_j}^\sigma(p_0)$ contains a cut point $p_{j+1}\notin Y_{N_j}^\sigma$ by the preceding argument. Then $p_j\in(p_0,p_{j+1})$ and $R=\bigcup_{j\in\N}[p_0,p_j]$ lies in $\{p_0\}\cup D_{p_0}^\sigma(w_0)$. The subtree $Y_{N_j}^\sigma$ contains $p_0$ and omits $p_{j+1}$, so $R\cap Y_{N_j}^\sigma\subseteq[p_0,p_{j+1}]$. Since every compact set lies in some $Y_{N_j}^\sigma$, the parametrization in the proof of Proposition~\ref{prop:rooted-representatives} makes $R$ a proper ray.

\smallskip
\textup{(ii)} If $Y_1^\sigma=Y_0^\sigma$, then~\eqref{eq:pi-exact-pullback} gives $Y_{n+1}^\sigma=f_\sigma^{-n}(Y_1^\sigma)=f_\sigma^{-n}(Y_0^\sigma)=Y_n^\sigma$ for $n\in\N_0$. Hence $Y^\sigma=Y_0^\sigma$. The converse follows from $Y_0^\sigma\subseteq Y_1^\sigma\subseteq Y^\sigma$.

\smallskip
\textup{(iii)} If $Y^\sigma=Y_0^\sigma$, each $T_X^\sigma$ is a nondegenerate subarc of a compact interval by~\textup{(i)}. Since $\Gamma_m^\sigma$ is a tree, two copies meet in at most one point by~\eqref{eq:ch4-stage-copy-intersection}. They meet only at common endpoints, with at most two copies at each such point. Hence each retained block has degree two. Distinct blocks incident to $X$ have distinct attachments by~\eqref{eq:ch4-incident-block-uniqueness} and Proposition~\ref{c:initial}~\textup{(ii)}, so at most two retained blocks are incident to $c(X)$.

If $Y^\sigma$ is hairy, a triod lies in some $Y_N^\sigma$ by Proposition~\ref{prop:pi-weak-tree}~\textup{(i)}. As in the choice of road nodes in the proof of Lemma~\ref{lem:ch5-marked-turns}, take an image $\beta_m^\sigma(b_i)$ of an interior node in each open arm at a common level $m>N$. Each $b_i$ is retained because the adjacent atoms have distinct tile vertices. If the induced graph were a path, the arcs joining successive block attachments in its tile copies would form an arc through all three points by~\eqref{eq:ch4-stage-copy-intersection}, contrary to their choice in distinct arms.

\smallskip
\textup{(iv)} The second assertion follows from~\textup{(iii)} with $m=1$. Suppose neither tree is hairy. Then every induced connection graph is a path by~\textup{(iii)}. Since only degree-one blocks are deleted, the circuit in the proof of Lemma~\ref{lem:ch5-marked-turns} visits each block $b$ at $\card\cA_m^\sigma(b)=\deg_{\Gamma_m^\sigma}(b)\le2$ angles. Its visits to the retained blocks follow the induced path once in each direction. Exchanging the two angles at each retained block therefore reverses their cyclic order.

If $\card\cA_m^\sigma(b)=2$, then $\cA_{m+1}^\sigma(\operatorname{Lift}_m^\sigma b)=\cA_m^\sigma(b)$ by~\eqref{eq:ch4-block-persistence} and the bound of two angles per class. Thus exchanging the two angles in each such class defines an involution $j_\sigma$ on their union $D_\sigma$. Any three angles in $D_\sigma$ lie in two-element classes at a common level, where the exchanges reverse their cyclic order. If $\mu(t)$ belongs to a level-$m$ pair, then $t\in\mathbf A^{m+1}$ and its class maps onto that pair by~\eqref{eq:finite-angle-preimages} and~\eqref{eq:finite-class-image}. The class of $t$ has at most two angles, so the map is bijective. Hence $j_\sigma\circ\mu=\mu\circ j_\sigma$ on $\mu^{-1}(D_\sigma)\subseteq D_\sigma$. The induced $1$-connection graph joins $d\ge2$ tiles, so $D_\sigma\ne\emptyset$. For $a\in D_\sigma$, the dense set $\bigcup_{k\in\N}\mu^{-k}(a)$ lies in $\mu^{-1}(D_\sigma)$. Thus $j_\sigma$ extends by order completion to an orientation-reversing homeomorphism of $S^1$ and commutes with $\mu$ by continuity. By Definition~\ref{oversetwsim} and density of $D_\sigma$, the closed relation of color $\sigma$ is the diagonal together with the graph of $j_\sigma$.

Since $j_\sigma(0)$ is fixed by $\mu$, $j_\sigma$ maps $\mu^{-k}(0)$ onto $j_\sigma(0)+\mu^{-k}(0)$ in reverse order. Thus $j_\sigma(t)=j_\sigma(0)-t$ on these sets for $k\ge1$, and hence on $S^1$ by density. The rotation $j_+\circ j_-$ has order dividing $d-1$, so $G=\langle j_+,j_-\rangle$ is finite. Its orbit relation $\overset{\white}{\sim}\vee\overset{\black}{\sim}$ is closed, so $S^1/G\cong S^2$ by Lemma~\ref{closureequi}~\textup{(i)}. Quotienting first by the rotations and then by the induced reflection makes $S^1/G$ a compact interval, which is a contradiction.
\end{proof}

\begin{lemma}\label{genuine_hairy}
Fix $\sigma\in\{+,\,-\}$ and suppose $Y^\sigma$ is hairy with the ribbon structure of Proposition~\ref{Yribbon}. Use the representatives $R_\omega$ rooted at $o^\sigma\in Y_0^\sigma$ from Proposition~\ref{prop:rooted-representatives}. For each side of a nondegenerate oriented arc $I\subseteq Y^\sigma$, there are $\omega\in\mathcal E_g(Y^\sigma)$ and $x\in\inte_\circ I\cap R_\omega$ such that $R_\omega\smallsetminus[o^\sigma,x]$ lies in a germ at $x$ on that side.
\end{lemma}

\begin{proof}
Put $p\=r_I^\sigma(o^\sigma)$ and take a nondegenerate closed subarc of $\inte_\circ I\smallsetminus\{p\}$ oriented as $I$. By hairiness, take a proper ray $R$ on the chosen side from an interior point $x$ of this subarc. The germ $D_x^\sigma(p)$ contains $[o^\sigma,\,x)$ and meets $I$, so it differs from the germ containing $R\smallsetminus\{x\}$. These germs are open, so pasting $[o^\sigma,x]$ and $R$ at $x$ gives a ray from $o^\sigma$.

For each interior subtree $K$, let $C_K$ be the component of $Y^\sigma\smallsetminus K$ containing a tail of $R$. For $K\subseteq L$, the component of $Y^\sigma\smallsetminus K$ containing $C_L$ also contains a tail of $R$. By the definition of $C_K$, this component is $C_K$. Hence $C_L\subseteq C_K$. Thus $\omega\=(C_K)_{K\in\mathcal K(Y^\sigma)}$ is an end. The joined ray represents $\omega$ and consists of cut points, namely its root and the interior points of its subarcs. Hence $R_\omega=[o^\sigma,x]\cup R$ by Proposition~\ref{prop:rooted-representatives}. Since $A_\omega$ contains the proper ray $R$, $\omega$ is genuine. The tail $R_\omega\smallsetminus[o^\sigma,x]=R\smallsetminus\{x\}$ lies in $D_x^\sigma(y)$.
\end{proof}

For a connection of $1$-tiles as in Lemma~\ref{degree3}, Lemma~\ref{genuine_hairy} applies to both trees by Proposition~\ref{lem:hairy}~\textup{(iv)}.

\subsection{Circular order on the end spaces}\label{sub_Circular_order_on_the_end_spaces}
In this subsection, we define the circular order on the end space of $Y^\sigma$ with respect to the ribbon structure

\begin{definition}[Circular order on ends]\label{def:end-circular-order}
    For distinct ends $\omega_1,\,\omega_2,\,\omega_3\in\mathcal E(Y^\sigma)$, suppose $\omega_i\=\{C_K^{\omega_i}\}$ for $i\in\{1,\,2,\,3\}$. Choose an interior subtree $K$ such that $C_K^{\omega_1},\,C_K^{\omega_2},\,C_K^{\omega_3}$ are pairwise distinct. Take points $x_i\in C_K^{\omega_i}$. Put $m\=\operatorname{med}_{Y^\sigma}(x_1,x_2,x_3)$ and let $G_i$ be the germ determined by $[m,x_i]_{Y^\sigma}$. Define $\langle\omega_1,\omega_2,\omega_3\rangle_\sigma\=\langle G_1,G_2,G_3\rangle_{\sigma}$.
\end{definition}

We now show that the above definition is independent of the choice of $K$ and $x_i\in C_K^{\omega_i}$.

\begin{prop}
    $\langle\omega_1,\omega_2,\omega_3\rangle_\sigma$ defined above is independent of the choice of $K$ and $x_i\in C_K^{\omega_i}$. Moreover, if $\card\mathcal E(Y^\sigma)\ge3$, these values define a circular order on $\mathcal E(Y^\sigma)$.
\end{prop}

\begin{proof}
    By the construction of $m\=\operatorname{med}_{Y^\sigma}(x_1,x_2,x_3)$ as above, we know that for $i\in\{1,\,2,\,3\}$, each $\omega_i$ can be presented by a ray $R_{\omega_i}$ rooted at $m$ as in Proposition~\ref{prop:rooted-representatives}. Moreover, we have $R_{\omega_1}\cap R_{\omega_2}\cap R_{\omega_3}=\{m\}$. Then the interior subtree $K$ with $C_K^{\omega_1},\,C_K^{\omega_2},\,C_K^{\omega_3}$ pairwise distinct contains $m$ and $x_i\in R_{\omega_i}$. It follows that the germ $G_i$ determined by $[m,x_i]_{Y^\sigma}$ is independent of the choice of $K$ and $x_i\in C_K^{\omega_i}$. Hence, $\langle\omega_1,\omega_2,\omega_3\rangle_\sigma$ is also independent of the choice of $K$ and $x_i\in C_K^{\omega_i}$.
    
    Since $\langle G_1,G_2,G_3 \rangle = {\rm sgn}(\tau)\langle \tau(G_1),\tau(G_2),\tau(G_3)\rangle$ for any permutation $\tau$ of $\{G_1,\,G_2,\,G_3\}$, we also have $\langle \omega_1,\omega_2,\omega_3 \rangle = {\rm sgn}(\tau)\langle \tau(\omega_1),\tau(\omega_2),\tau(\omega_3)\rangle$. 
    
    Suppose $\omega'\in\mathcal E(Y^\sigma)\smallsetminus\{\omega_1,\,\omega_2,\,\omega_3\}$ with $\langle\omega_1,\omega_2,\omega_3\rangle_\sigma=1$ and $\langle\omega_1,\omega_3,\omega'\rangle_\sigma=1$. We show that $\langle\omega_1,\omega_2,\omega'\rangle_\sigma=\langle\omega_2,\omega_3,\omega'\rangle_\sigma=1$. Suppose that $K$ is an interior subtree such that $C_K^{\omega_1},\,C_K^{\omega_2},\,C_K^{\omega_3},\,C_K^{\omega'}$ are pairwise distinct. We also choose $x'\in C_K^{\omega'}$ suppose that $m'\=\operatorname{med}_{Y^\sigma}(x_1,x_3,x')$. If $m'=m$, then we are done since $\langle\,\cdot\,,\,\cdot\,,\,\cdot\,\rangle_\sigma$ is a circular order so satisfies Definition~\ref{def_circularorder}~(ii). If $m'\neq m$, then $m'\in[x_1,m]_{Y^\sigma}$ or $m'\in[x_1,m]_{Y^\sigma}$. When $m'\in[x_1,m]_{Y^\sigma}$, we have $\operatorname{med}_{Y^\sigma}(x_1,x_2,x')=\operatorname{med}_{Y^\sigma}(x_1,m,x')=m'$, so $[m',x_3]_{Y^\sigma}$ and $[m',x_2]_{Y^\sigma}$ determine the same germ at $m'$ with $[m',m]_{Y^\sigma}$. It follows that $\langle\omega_1,\omega_2,\omega'\rangle_\sigma=\langle\omega_1,\omega_3,\omega'\rangle_\sigma=1$. When $m'\in[x_1,m]_{Y^\sigma}$, we also have $\operatorname{med}_{Y^\sigma}(x_2,x_3,x')=\operatorname{med}_{Y^\sigma}(x_2,x_3,m')=m$. It follows that $\langle\omega_2,\omega_3,\omega'\rangle_\sigma=\langle\omega_2,\omega_3,\omega_1\rangle_\sigma=1$. By Definition~\ref{def_circularorder}, we know that $\langle\,\cdot\,,\,\cdot\,,\,\cdot\,\rangle_\sigma$ on $\cE(Y^\sigma)$ is a circular order.
\end{proof}

Now we show that the order of three ends is determined if they are distinguished in a finite state $Y^\sigma_n$ for some $n\in\N$.

\begin{lemma}\label{lem:first-separation-order}
    Let $v\in Y_n^\sigma$ be a cut point of $Y^\sigma$. If three ends $\omega_i$ determine distinct nonempty germs $C_{\{v\}}^{\omega_i}\cap Y_n^\sigma$ at $v$ in $Y_n^\sigma$, the ribbon order of these germs is $\langle\omega_1,\omega_2,\omega_3\rangle_\sigma$.
\end{lemma}
\begin{proof}
    Take $x_i\in C_{\{v\}}^{\omega_i}\cap Y_n^\sigma$. Their median is $v$. Apply Definition~\ref{def:end-circular-order} with $K=\{v\}$. The compatibility of the finite-stage ribbon structures (see Proposition~\ref{Yribbon}~(i)) shows that this order is unchanged at every higher stage.
\end{proof}

\subsection{The end circle} \label{theendspace}

In what follows, we show that when $Y^+$ and $Y^-$ are hairy, the circular order completion of $\mathcal E(Y^+)$ and $\mathcal E(Y^-)$ are topological circles. Then we show that there is an orientation-reversing homeomorphism between the two circles.

Let $S^1(Y^\sigma)\=\overline{\mathcal E(Y^\sigma)}$ denote the circular order completion of $\mathcal E(Y^\sigma)$ (see Proposition~\ref{completionofcircularorder}). The main result of this subsection is the following.

\begin{theorem}
\label{thm:end-circle}
Fix $\sigma\in\{+,\,-\}$ and suppose $Y^\sigma$ is hairy. The circular order completion $S^1(Y^\sigma)$ is canonically an oriented topological circle. More precisely, the finite ribbon trees determine an inverse-limit circle $B_\infty^\sigma=\varprojlim(B_n^\sigma,\chi_n^\sigma)$ (here $B_n^\sigma$ and $\chi_n^\sigma$ will be defined later) and an order-preserving embedding $q^\sigma\:\mathcal E(Y^\sigma)\to B_\infty^\sigma$ whose restriction to the genuine ends has dense image. Consequently $q^\sigma$ extends uniquely to an orientation-preserving homeomorphism $\rho^\sigma\:S^1(Y^\sigma)\to B_\infty^\sigma$.
\end{theorem}

The proof of the theorem is organized as follows. We first construct boundary circles of finite trees and the collapse maps between them. We then show that their inverse limit is a topological circle and embed the end space into this circle. Finally, using the hairiness of $Y^\sigma$, we show that images of the genuine ends are dense in the circle.

Let $T$ be a finite tree with a ribbon structure. Since $T$ is finite, we can decompose $T$ into finitely many oriented edges $\{a_i\}_{i=1}^n$, and points in the interior of an edge have valence two in $T$. For a directed edge $a\=[0_{a},1_{a}]_T\in\{a_i\}_{i=1}^n$, let $\oa$ be its reverse and $0_{a},\,1_{a}$ be its initial and terminal vertices. At $0_{a}\in T$, let $\sigma_{0_{a}}(a)$ be the successor of $a$ on the directed edges issuing from $0_a$ with respect to the ribbon structure. More precisely, if $0_a$ has valence two, $\sigma_{0_a}(a)$ is the other edge issuing from $0_a$. If $0_a$ has valence one, $\sigma_{0_a}(a)$ is $a$. If $0_a$ has valence larger than two, $\sigma_{0_a}(a)$ is the edge such that for each edge $b$ different from $a$ and $\sigma_{0_a}(a)$, we have $\langle b,a,\sigma_{0_a}(a)\rangle_\sigma=1$. Set $\beta_T(a)\=\sigma_{1_a}(\oa)$. Let $I_a$ be an interval with the same orientation as $a$ and $I_{\oa}$ be an interval with the same orientation as $\oa$. Glue the terminal point of $I_a$ to the initial point of $I_{\beta_T(a)}$ for all $a\in \{a_i\}_{i=1}^n$. Then the resulting space $B(T)$, with natural projection $\tau_T\: B(T)\to T$, is a circle. Put $B_n^\sigma\=B(Y_n^\sigma)$.

Suppose $T'$ is another finite tree with a ribbon structure, and $\iota \: T \hookrightarrow T'$ is a ribbon-preserving inclusion. Subdivide $T'$ into $\{a_i'\}_{i=1}^{n'}$ as $T$. Moreover, we assume that $\{\iota(a_i)\}_{i=1}^n\subseteq\{a_i'\}_{i=1}^{n'}$.

\begin{lemma}\label{chidef}
    There exists a canonical orientation-preserving monotone map $\chi_{T', T}\: B(T')\to B(T)$ that maps $I_{\iota(a_i)}$ to $I_{a_i}$ for $i\in\{1,\,2,\,\dots,\,n\}$ and maps each component of $B(T')\smallsetminus \bigcup_{i=1}^nI_{\iota(a_i)}\cup \bigcup_{i=1}^nI_{\iota(\overline{a_i})}$ to a point.
\end{lemma}

\begin{proof}
    Every component $C$ of $T'\smallsetminus T$ attaches at one point, since two attachments would create two tree arcs between them. Moreover, the intervals $I_a\subseteq B(T')$ with $a\in C$ compose a whole interval $I_C$ of $B(T')$. Suppose $I$ is the interval between two intervals $I_{\iota(a_i)}$ and $I_{\iota(a_j)}$ for $i\neq j\in\{1,\,2,\,\dots,\,n\}$. If there is no other $I_{\iota(a)}\subseteq I$, then $I$ is a union of several $I_C$ where $C$ is a component of $T'\smallsetminus T$. Let $\chi_{T', T}(I)$ be the initial point of $I_{a_j}\subseteq B(T)$, that is also the terminal point of $I_{a_i}$. When restricted to $I_{\iota(a_i)}$, let $\chi_{T', T}$ be a homeomorphism between $I_{\iota(a_i)}$ to $I_{a_i}$. Then $\chi_{T', T}$ is the map we want.
\end{proof}

We call the map $\chi_n^\sigma\=\chi_{Y^\sigma_{n+1}, Y^\sigma_n}\:B_{n+1}^\sigma\to B_n^\sigma$ the \emph{stage collapse}, where $\chi_{Y^\sigma_{n+1}, Y^\sigma_n}$ is given by above lemma. Thus, the finite ribbon trees $Y^\sigma_n$ determine an inverse system of oriented circles.

Define $B_\infty^\sigma\=\varprojlim(B_n^\sigma,\chi_n^\sigma)$, with coordinate projections $p_n^\sigma$. For distinct $x,\,y,\,z\in B_\infty^\sigma $, define their cyclic order at each coordinate where the three images are distinct. Monotonicity and orientation preservation of the bonding maps $\chi_n^\sigma$ make this independent of the coordinate.

\begin{prop}\label{prop:inverse-limit-boundary-circle}
    The inverse-limit topology equals the circular-order topology on $B_\infty^\sigma$, and $B_\infty^\sigma$ is an oriented topological circle. Each $p_n^\sigma$ is an orientation-preserving monotone map.
\end{prop}

\begin{proof}
    Each $p_n^\sigma$ is surjective, and its fibers are inverse limits of connected fibers of finite monotone compositions, hence connected. For $a,\,b\in B_\infty^\sigma $, let $(a,b)\=\{x\in B_\infty^\sigma : a,\,x,\,b\text{ occur in the}\\ \text{ positive cyclic order}\}$ be an interval in $B_\infty^\sigma$. If $x\in(a,b)$, a coordinate where $a,\,x,\,b$ are distinct and a small arc around $p_N^\sigma(x)$ gives an inverse-limit neighborhood inside $(a,b)$; conversely, inverse images of small proper arcs under a monotone coordinate projection are intervals with respect to the circular order. Thus, the topologies agree.

   Since $B_\infty^\sigma$ is compact and metrizable, it is separable. Since the inverse-limit topology agrees with the circular-order topology, every closed interval is closed in $B_\infty^\sigma$. Hence, every nested sequence $I_1\supseteq I_2\supseteq\cdots$ of nonempty closed intervals has the finite-intersection property, so compactness gives $\bigcap_{j=1}^{+\infty}I_j\neq\emptyset$. Thus, $B_\infty^\sigma$ is order complete. We claim that it has no gaps. Given two distinct points $a,\,b$, choose a separating coordinate $p_N^\sigma(x)$. There is a point strictly between their images under $p_N^\sigma(x)$. Then the lift of this point is between $a$ and $b$. By \cite[Proposition~5.6]{Frankel13}, a complete separable circular order without gaps is order-isomorphic to a closed subset of $S^1$. Then the absence of gaps implies that $B_\infty^\sigma$ is homeomorphic to $S^1$. Since $p_n^\sigma$ is surjective, monotone, and preserves every triple with distinct images, it is orientation-preserving.
\end{proof}

We have therefore obtained a concrete oriented circle from the finite ribbon trees alone. It remains to identify this circle with the circular order completion of the end space.

Fix $o^\sigma\in Y^\sigma_0$. Let $\omega=\{C_K^\omega\}_{K\in \mathcal K_{ o^\sigma}(Y^\sigma)} \in \mathcal E(Y^\sigma)$ be an end. Let $R_\omega$ be the rooted ray representing $\omega$ as in Proposition~\ref{prop:rooted-representatives} and $A_\omega\=\overline {R_\omega}^{Y^\sigma}$. For integers $m>n$, write $\chi_{n,m}^\sigma\=\chi_n^\sigma\circ\chi_{n+1}\circ\cdots\circ\chi_{m-1}^\sigma$. If $A_\omega\not\subseteq Y_n^\sigma$, its tail lies in one component of $Y^\sigma\smallsetminus Y_n^\sigma$. Then $Y_m^\sigma$ meets that tail, and the corresponding component of $Y_m^\sigma\smallsetminus Y_n^\sigma$ determines an interval in $B_m^\sigma$ as in Lemma~\ref{chidef}, which is collapsed by $\chi_{n,m}^\sigma$ to a point. We denote this point by $q_n^\sigma(\omega)$. The point $q_n^\sigma(\omega)$ is independent of the integer $m$ since $(\chi_{n,m}^\sigma)^{-1}(q_n^\sigma(\omega))\subseteq (\chi_{n,m+1}^\sigma)^{-1}(q_n^\sigma(\omega))$ are both collapsed into the same point by Lemma~\ref{chidef}.

If $A_\omega\subseteq Y_n^\sigma$, then $\omega$ is a hard end with terminal point $h_\omega$. Since $h_\omega$ is a hard endpoint, it has a unique preimage under $\tau_{Y_m^\sigma} \: B_m^\sigma\to Y_m^\sigma$ for all integers $m\geq n$. Define $q_n^\sigma(\omega)$ as the preimage of $h_\omega$ under $\tau_{Y_n^\sigma}$ in $B_n^\sigma$. In both cases $\chi_n^\sigma(q_{n+1}^\sigma(\omega))=q_n^\sigma(\omega)$ and we set $q^\sigma(\omega)\in B_\infty^\sigma$ determined by $\{q_n^\sigma(\omega)\}_{n\in\N_0}$.

\begin{prop}\label{prop:end-embedding}
    The map $q^\sigma:\mathcal E(Y^\sigma)\to B_\infty^\sigma$ is injective and preserves the cyclic order of every triple of distinct ends.
\end{prop}

\begin{proof}
    First, we show that $q^\sigma$ is injective. For two distinct ends $\omega_1=\bigl\{C_K^{\omega_1}\bigr\}_{K\in \mathcal K_{ o^\sigma}(Y^\sigma)}$ and $\omega_2=\bigl\{C_K^{\omega_2}\bigr\}_{K\in \mathcal K_{ o^\sigma}(Y^\sigma)}$, take an interior subtree $K$ containing $o^\sigma$ separating the two ends. Take a stage $Y_n^\sigma$ containing $K$. Let $R_{\omega_i}$ be the rooted ray representing $\omega_i$ as in Proposition~\ref{prop:rooted-representatives} and $A_{\omega_i}\=\overline {R_{\omega_i}}^{Y^\sigma}$ for $i\in\{1,\,2\}$. The points $q_n^\sigma(\omega_1)$ and $q_n^\sigma(\omega_2)$ correspond to the last point of $A_{\omega_1}\cap Y_n^\sigma$ and $A_{\omega_2}\cap Y_n^\sigma$, which are different. Thus $q^\sigma$ is injective. 
    
    Next, we show that $q^\sigma$ preserves the cyclic order of every triple of distinct ends. For three distinct ends $\omega_i=\bigl\{C_K^{\omega_i}\bigr\}_{K\in \mathcal K_{ o^\sigma}(Y^\sigma)}$ and $i\in\{1,\,2,\,3\}$, take an interior subtree $K$ containing $o^\sigma$ separating the two ends. Take a stage $Y_n^\sigma$ containing $K$. Let $R_{\omega_i}$ be the rooted ray and $x_i$ be the last point of $A_{\omega_i}\cap Y_n^\sigma$. We assume that $x_i$ are pairwise distinct. Let $m\=\operatorname{med}_{Y^\sigma}(x_1,x_2,x_3)\in K$. We assume that $m$ is different from $x_i$ and the three germs at $m$ toward $x_i$ are distinct. Note that $[m,x_i]_{Y_n^\sigma}$ also determines an interval $I_i\subseteq B_n^\sigma$. Then $q^\sigma_n(\omega_i)$ is a point in $(\tau_{Y_n^\sigma})^{-1}(x_i)\in I_i$. Lemma~\ref{lem:first-separation-order} identifies the end order of $w_i$ in $Y_n^\sigma$ with the ribbon order of $[m,x_i]_{Y_n^\sigma}$ at $m$, and also identifies the order of $I_i$ on $B_n^\sigma$. This order also identifies the order of $q^\sigma_n(\omega_i)$ on $B_n^\sigma$.
\end{proof}

\begin{lemma}\label{lem:genuine-end-density}
    Suppose $Y^\sigma$ is hairy. The set $q^\sigma(\mathcal E_g(Y^\sigma))$ is dense in $B_\infty^\sigma$.
\end{lemma}

\begin{proof}
    Let $U\subseteq B_\infty^\sigma$ be nonempty and open. Choose $N$ and a proper open arc $V\subseteq B_N^\sigma$ with $(p_N^\sigma)^{-1}(V)\subseteq U$, and $V'\subseteq V$ such that $V'\in I_a$ where $a\in Y_N^\sigma$ is an oriented interval in the subdivision of $Y_N^\sigma$ as in the beginning of this subsection. By Lemma~\ref{genuine_hairy}, there exists a genuine end $\omega$ branching off from both sides of $a$. By definition, we have $q_N^\sigma(\omega)\in V$. Thus $q^\sigma(\omega)\in U$.
\end{proof}

\begin{proof}[Proof of Theorem~\ref{thm:end-circle}]
    Proposition~\ref{prop:inverse-limit-boundary-circle} identifies $B_\infty^\sigma$ with an oriented topological circle, while Proposition~\ref{prop:end-embedding} gives an order-preserving embedding $q^\sigma$ of the end space. By Lemma~\ref{lem:genuine-end-density}, the images of genuine ends are dense. Since $B_\infty^\sigma$ is order complete, the universal property of the circular order completion gives a unique orientation-preserving homeomorphism $\rho^\sigma\:S^1(Y^\sigma)\to B_\infty^\sigma$ extending $q^\sigma$ (see Proposition~\ref{completionofcircularorder}).
\end{proof}

Now, we show that a component in $Y^+\smallsetminus K$ for $K\in \mathcal K(Y^+)$ is associated to an interval in $S^1(Y^+)=\overline{\mathcal E(Y^+)}$. The statement is also true for $Y^-$.
\begin{lemma} \label{component_to_interval}
    Let $K\in \mathcal K(Y^+)$ be an interior subtree of $Y^+$ and $C$ be a connected component of $Y^+\smallsetminus K$. Define $I_C\=\{w=\{C_{L}^w\}_{L\in \mathcal K(Y^+)}\in \mathcal{E}(Y^+):C_{K}^w=C\}$ and denoted by $\overline {I_C}$ the closure of $I_C$ in $\overline{\mathcal E(Y^+)}$. Then $\overline{I_C}$ is a closed circular interval in $\overline{\mathcal E(Y^+)}$. Moreover, $\overline{I_C}\cap \mathcal E(Y^+)=I_C$. In particular, if $w\in \mathcal E(Y^+)$ satisfies $C_{K}^w\neq C$, then $w\notin \overline {I_C}$.
\end{lemma}
\begin{proof}
    Write $S\=\mathcal E(Y^+)$ and $\overline S\=\overline{\mathcal E(Y^+)}$. Fix a point $o^+\in K$. Then by Proposition~\ref{prop:rooted-representatives}, each end is identified with a ray rooted at $o^+$ for $\inte Y^+$. By Lemma~\ref{lem:end-in-component}, $I_C$ is nonempty. Since $Y^+$ is hairy and $S\smallsetminus I_C\neq \emptyset$, we can choose $d\in S\smallsetminus I_C$. Take a linear order $<_d$ in $\overline S\smallsetminus \{d\}$ as follows: $x<_d y$ if and only if $\langle d,x,y\rangle =+1$. Then $\oS\smallsetminus\{d\}$ is homotopic to the open interval $(0,1)$.

    We show that $I_C$ is convex with respect to $<_d$. Take $u\neq v\in I_C$ and let $x_0\in S$ satisfy $u<_d x_0<_dv$. Then we have that $\langle d,u,x_0\rangle_+=\langle d,x_0,v\rangle_+=1$. It follows that $\langle u,x_0,v\rangle_+=1$ and $\langle d,u,v\rangle_+=1$. We prove $x_0\in I_C$ by contradiction. Assume $x_0\notin I_C$ and take an arbitrary point $x\in S\smallsetminus I_C$. Consider the rooted representatives $R_u$, $R_v$ and $R_x$ of $u$, $v$ and $x$ respectively as in Proposition~\ref{prop:rooted-representatives}. Suppose $\oC\cap K=\{a\}$. Since $u,\,v\in I_C$, we have that $R_u\cap C\neq \emptyset$ and $R_v\cap C\neq \emptyset$. It follows that $[o^+,a]_{Y^+}\subseteq R_u\cap R_v$ and the rays $R_u$, $R_v$ branch at a point $w\in C$. Since $x\notin I_C$, we have $R_x\cap C=\emptyset$, where $\overline{C_K^x}\cap K=\{b\}$. This implies that $\operatorname{med}(u,x,v)=w\in C$. By Definition~\ref{def:end-circular-order}, we have that $\langle u,x,v\rangle_+=\langle a_u,a,a_v\rangle_+$, where $a_u\in R_u\smallsetminus [o^+,a]_{Y^+}$ and $a_v\in R_v\smallsetminus [o^+,a]_{Y^+}$. Since $d,\,x_0\notin I_C$, we have that $\langle u,x_0,v\rangle_+=\langle u,d,v\rangle_+=\langle a_u,a,a_v\rangle_+$, which contradicts that $\langle u,x_0,v\rangle_+=\langle d,u,v\rangle_+=1$. Therefore, $x_0\in I_C$.
    
    Then we prove that $I_C$ is open in $S$. Fix an arbitrary $w\in I_C$ and consider $R_{w}$. We can choose the two rays whose tails are contained in $C$ and branch from both sides of an oriented subarc of $R_{w}\cap C$ by the hairy property. These two rays determine two ends $w_-,\,w_+\in I_C$. Since the two rays branch from both sides of an oriented subarc of $R_{w}\cap C$, we have $w_-<_d w<_d w_+$. Let $(w_-,w_+)_S\=\{a\in S:w_-<_d a<_d w_+\}$. Since $I_C$ is convex with respect to $<_d$, we have that $(w_-,w_+)_S\subseteq I_C$. Then $(w_-,w_+)_S$ is an open interval neighborhood of $w$ in $I_C$. Thus, $I_C$ is open.

    Next we show that $\overline{I_C}$ is a closed interval. By the hairy property, there exists $d_-,\,d_+\in S$ with $d\in (d_-,d_+)_S$ and $(d_-,d_+)_S\cap I_C=\emptyset$. Therefore, $I_C$ is bounded away from $d$. Hence, we can define $\alpha\=\inf_{<_d} I_C$ and $\beta\=\sup_{<_d} I_C$, i.e., $\alpha$ is the infimum of $I_C$ and $\beta$ is the supremum of $I_C$ with respect to $<_d$. By the definition of order completions, we induce $\overline{I_C}=[\alpha,\beta]_{<_d}\=\bigl\{x\in \oS:\alpha<_d x <_d\beta\bigr\}\cup\{\alpha,\,\beta\}$, which is a closed circular interval in $\overline S$.

    It remains to show that $\overline {I_C}\cap S=I_C$. Take $w\in S\smallsetminus I_C$. Then $w\in I_D$ for some component $D\neq C$ of $Y^+\smallsetminus K$. Since $I_D$ is open, there exists an open interval neighborhood $U$ of $w$ in $S$ such that $w\in U\subseteq I_D$. Then $I_D\cap I_C=\emptyset$ and $w\notin \overline{I_C}$. Thus, $\overline {I_C}\cap S=I_C$.
\end{proof}

Now we recall the notion of ideal gap from \cite[Subsection~3.2]{CLWheels}. Note that $\cE(Y^\sigma)=\cE(\inte Y^\sigma)$ by Definition~\ref{def_endspace}, so we may consider $\inte Y^\sigma$ when necessary. The points of $ \overline{\cE(\inte Y^\sigma)}\smallsetminus \mathcal E(\inte Y^{\sigma})$ are represented by geometric objects called \emph{ideal gaps} as in Definition~\ref{idealgap}.

If $e^L,\,e^R \in \mathcal E(\inte Y^+)$, as above, we define $\bigl[e^L,e^R\bigr] \subseteq \mathcal E(\inte Y^+)$ to be the set
\begin{equation*}
\bigl[e^L,e^R\bigr] \= \bigl\{e^L,\,e^R\bigr\}
\cup \bigl\{e \in \mathcal E(Y^+): \bigl\langle e^L,e,e^R\bigr\rangle_+=1\bigr\}.
\end{equation*} 

\begin{lemma}[{Calegari \& Loukidou \cite[Lemma~2.15]{CL24}}]\label{nestedsequence}
Suppose there is an infinite nested sequence
\begin{equation*}
\cdots \subseteq \bigl[e_{i+1}^L,e_{i+1}^R\bigr]
\subseteq \bigl[e_i^L,e_i^R\bigr] \subseteq \cdots
\end{equation*}
whose intersection is empty. For each $i\in\N$, let $l_i \subseteq \inte Y^+$ denote the convex hull of $e_i^L$ and $e_i^R$ (this is a properly embedded line), oriented to run from $e_i^L$ to $e_i^R$. Then there is some $i\in\N$ so that the intersection $X \= \bigcap_{j = i}^{+\infty} l_j$ is either a point or a compact, nonempty interval.
\end{lemma}

We are now in a position to define ideal gaps.

\begin{definition}[{Calegari \& Loukidou \cite[Definition~2.16]{CL24}}]\label{idealgap}
A sequence of proper lines $l_i$ in $\inte Y^+$ which are the convex hulls of a sequence of nested intervals $\bigl[e_i^L,e_i^R\bigr]$ in $\mathcal E(\inte Y^+)$ is said to be \emph{end nested}.

A point $p \in \inte Y^+$ together with an end nested sequence of proper lines $l_i$ associated to a nested sequence of closed intervals in $\mathcal E(\inte Y^+)$ with empty intersection for which $p \in \bigcap_{i\in\N} l_i$ (as in the statement of Lemma~\ref{nestedsequence}) is called an \emph{ideal gap}.
Two ideal gaps are called \emph{equivalent} if they contain cofinal end nested subsequences of some end nested sequence.
\end{definition}

A point $p$ as above is said to be contained in the \emph{support} of the ideal gap (and by abuse of notation, the support of an equivalence class of the ideal gap). 

\begin{lemma}[{Calegari \& Loukidou \cite[Lemma~3.8]{CLWheels}}]\label{lemma_ordercompletion}
The union of $\mathcal E(\inte Y^\sigma)$ together with the equivalence classes of ideal gaps of $\inte Y^\sigma$ admits a natural circular ordering which is equal to the order completion $\overline{\mathcal E(\inte Y^\sigma)}$.
\end{lemma}

Since $Y^\sigma$ is hairy, $\inte Y^\sigma$ is also hairy. Then we recall the following lemma showing that the support of an equivalence class of ideal gap necessarily consists of a single point. 

\begin{lemma}[{Calegari \& Loukidou \cite[Lemma~3.11]{CLWheels}}]
Suppose $\inte Y^\sigma$ is a hairy zipper. Then the support of each equivalence class of ideal gap is a single point.
\end{lemma}

Now we show that there is an orientation-reversing homeomorphism from $S^1(Y^+)$ to $S^1(Y^-)$.

\begin{prop}\label{prop_reversesorientationheom}
    Let $S^1$ be the unit circle, $\gamma^n$ be the e approximation in Subsection~\ref{Pseudo_isotopies}, and $\mathbf A^n\=(\gamma^n)^{-1}(\V^n)$. For $\sigma\in\{+,\,-\}$, if $Y^\sigma$ is hairy, then there exists a homeomorphism $\mu_\sigma\:S^1\to S^1(Y^\sigma)$ such that $\mu_-\circ\mu_+^{-1}\:S^1(Y^+)\to S^1(Y^-)$ reverses orientation when $Y^+$ and $Y^-$ are both hairy.
\end{prop}

\begin{proof}
    We give the proof for a hairy $Y^+$. The proof for a hairy $Y^-$ interchanges the colors. The references below to $S^1(Y^-)$ and $\mu_-$ assume that $Y^-$ is also hairy.
    
    We give some notations used later. Let $t\in \mathbf A^n\subseteq S^1$ be an $n$-angle and denote $v\=\gamma^n(t)\in\V^n$. Suppose $a^n_-(t)$ and $a^n_+(t)$ are $n$-arcs such that $a^n_-(t)$ has terminal point $t$ and $a_+^n(t)$ has initial point $t$. By~\eqref{defZ0+} and~\eqref{defZpm}, we have that $v\in Z^+\cap Z^-$. Let $b_n^+\in\sB_n^+$ (resp.\ $b_n^-\in\sB_n^-$) be the vertex in white (resp.\ black) $n$-connection graph with a block at $v$ containing the succeeding edges $\gamma^n(a^n_-(t))$ and $\gamma^n(a^n_+(t))$. Put $u^+\=\beta^+_n(b_n^+)\in Y_n^+$ and $u^-\=\beta^-_n(b_n^-)\in Y_n^-$, where $\beta_n^\sigma\:\sB_n^\sigma\to Y_n^\sigma$ is the block attachment of Definition~\ref{def:ch4-tile-copy-data}.
    
    Now we define $\mu_+(t)$. If $u^+$ is a hard end of $Y^+$, let $\widehat u^+$ denote this hard end in $S^1(Y^+)$. Then we let $\mu_+(t)\=\widehat u^+$. If $u^+$ is a cut point of $Y^+$, let $I\=[0,1]$ and let $Y_{u^+}^+\=Y^+\sqcup I/(0\sim u^+)$ be a new topological $\R$-tree. Here we identify $u^+\in Y^+$ with $0\in I$. It follows that $I\smallsetminus\{0\}$ is a germ at $u^+$. Then we show that we can still define a ribbon structure at $u^+$, still denoted by $\langle\,\cdot\,,\,\cdot\,,\,\cdot\,\rangle_+$. Let $C_1,\,C_2,\,C_3$ be three germs at $u^+$. If each germ $C_i$ for $i\in\{1,\,2,\,3\}$ is different from $I\smallsetminus\{0\}$, let their order be the same as their order in $Y^+$. If there is $C_i=I\smallsetminus\{0\}$ for some $i\in\{1,\,2,\,3\}$, we assume that $C_1=I\smallsetminus\{0\}$ without loss of generality. Let $m\geq n\in\N$ be an integer such that $C_2'\=C_2\cap Y_m^+\neq\emptyset$ and $C_3'\=C_3\cap Y_m^+\neq\emptyset$. Since $\V^n\subseteq\V^m$ for $m\geq n\in\N$, the $n$-vertex $v$ is also an $m$-vertex. We also define $m$-arcs $a^m_-(t)$ and $a^m_+(t)$ as above such that $a^m_-(t)$ has terminal point $t$ and $a_+^m(t)$ has initial point $t$. Let $b_m^+\in\sB_m^+$ (resp.\ $b_m^-\in\sB_m^-$) be the vertex in white (resp.\ black) $m$-connection graph with a block at $v$ containing the succeeding edges $\gamma^m(a_-^m(t))$ and $\gamma^m(a_+^m(t))$. By Proposition~\ref{prop:compact-finite-trees}, we have that on $C_2$ (resp.\ $C_3$), there is a vertex $v_2$ (resp.\ $v_3$) such that $v_2=\beta(b_2)$ (resp.\ $v_3=\beta(b_3)$) for some $b_2\in\sB^+_m$ (resp.\ $b_2\in\sB^+_m$). Let $X^m_1$ be an black $m$-tile incident to the block of $b_m^-$. If $c\=\operatorname{med}_{\Gamma_m^+}(b_m^+,b_2,b_3)\neq b_m^+$, we define the order of $C_1,\,C_2,\,C_3$ by the order of $[c,b_m^+],\,[c,b_2],\,[c,b_3]$ around $c$ in $\Gamma_m^+$. If $c\=\operatorname{med}_{\Gamma_m^+}(b_m^+,b_2,b_3)= b_m^+$, we consider the embedding $\nu_m\:\Gamma_m^+\to S^2$ as in Proposition~\ref{Yribbon}. Let $w\in S^2$ denote a point in the component corresponding to $b_m^-$ of the geometric representation of the connection $\D(\pi_\w(v)\cup\pi_\b(v))$. Then we connected $w$ and $\nu_m(b_m^+)$ by an arc $l_w$ on $S^2$. Then we define the order of $C_1,\,C_2,\,C_3$ by the order of $[\nu_m(b_m^+),w],\,[\nu_m(b_m^+),\nu_m(b_2)],\,[\nu_m(b_m^+),\nu_m(b_3)]$ around $\nu_m(c)$ in $\nu_m(\Gamma_m^+)\cup l_w$. As the proof of the statement, every $\iota_n^\sigma$ is ribbon-preserving in Proposition~\ref{Yribbon}, we know that $\langle C_1,C_2,C_3\rangle_+$ is independent of $m\in\N$. It follows that we can define an order of two ends $w_1$, $w_2$ in $\cE(Y^+)$ and the hard end $C_1$ standing for $I\subseteq Y_{u^+}^+$. We can also define a complete order of every two ends $w_1$, $w_2$ in $\cE(Y^+)$. If $\langle C_1,w_1,w_2\rangle_+=1$, let $w_1<w_2$, or let $w_1>w_2$. This determines a gap of $\cE(Y^+)$ as well as a point in $S^1(Y^+)$, denoted by $\widehat u^+$. Let $\mu_+(t)\=\widehat u^+$. Similarly, we can define $\mu_-(t)$ for each $t\in \mathbf A^n\subseteq S^1$.

    Moreover, $\widehat u^+$ is an ideal gap with support $u^+$. For an integer $m>n$, since $Y_m^+$ is embedded in $Y^+$, we also have $\inte Y_m^+$ is embedded in $\inte Y^+$. It follows that each component $C_m$ of $\inte Y_m^+\smallsetminus\{u^+\}$ is contained in a component $C_m'$ of $\inte Y^+\smallsetminus\{u^+\}$. Let $I_m$ be an interval in $\overline{\cE(\inte Y^\sigma)}$ corresponding to $C_m'$ defined in Lemma~\ref{component_to_interval}. Note that $\widehat u^+\notin I_m$ for each $C_m$. Note that these finitely many intervals are pairwise disjoint and have a circular order. Then we can choose the two intervals $I_m^L$, $I_m^R$ immediately preceding and following $\widehat u^+$ in the circular order. Then we can choose ends $e_m^L\in I_m^L$, $e_m^R\in I_m^R$ recursively, such that $\bigl[e_{m+1}^L,e_{m+1}^R\bigr] \subseteq \bigl[e_m^L,e_m^R\bigr]$. Let $l_m \subseteq \inte Y^+$ denote the convex hull of $e_m^L$ and $e_m^R$. It follows that $u^+\in l_m$ and $\widehat u^+\in \bigl[e^L,e^R\bigr]$. Then we obtain a sequence $\cdots \subseteq \bigl[e_{m+1}^L,e_{m+1}^R\bigr] \subseteq \bigl[e_m^L,e_m^R\bigr] \subseteq \cdots$, and $u^+\in\bigcap_{m=n+1}^{+\infty} l_m$. Hence, $\widehat u^+$ is an ideal gap with support $u^+$. 

    We now show that $\mu_+$ (resp.\ $\mu_-$) preserves the orientation of $\mathbf A^n$, where $S^1$ is endowed with the orientation corresponding, via $\gamma^1$, to the orientation of $\gamma^1\bigl(S^1\bigr)$ induced by the orientation of $S^2$ (resp.\ its opposite), and $S^1(Y^+)$ (resp.\ $S^1(Y^-)$) is endowed with the orientation specified in Subsection~\ref{theendspace}. Suppose $t_1,\,t_2,\,t_3\in\mathbf A^{\infty}$ are three disjoint angles and $v_i\=\gamma(t_i)$. Since $\V^n\subseteq\V^m$ for $m\geq n\in\N$, we can assume that $t_1,\,t_2,\,t_3\in\mathbf A^{n}$ for a certain $n\in\N$. We also assume that $u_i^+\=\beta^+_n(b_{n,i}^+)\in Y_n^+$ (resp.\ $u_i^-\=\beta^-_n(b_{n,i}^-)\in Y_n^-$) for some $b_{n,i}^+\in\sB_n^+$ (resp.\ $b_{n,i}^-\in\sB_n^-$) and $i\in\{1,\,2,\,3\}$. We consider the embedding $\nu_n\:\Gamma_n^+\to S^2$ as in Proposition~\ref{Yribbon}. Let $w_i\in S^2$ denote a point in the component corresponding to $b_{n,i}^-$ of the geometric representation of the connection $\D(\pi_\w(v_i)\cup\pi_\b(v_i))$. Then we connected $w_i$ and $\nu_n(b_{n,i}^+)$ by an arc $l_{w_i}$ on $S^2$. Let $u_0\=\operatorname{med}(w_1,w_2,w_3)$. Then $\langle\mu_+(t_1),\mu_+(t_2),\mu_+(t_3)\rangle_+$ is equal to the order of $[u_0,w_1],\,[u_0,w_2],\,[u_0,w_3]$ around $u_0$ in $\nu_n(\Gamma_n^+)\cup l_{w_1}\cup l_{w_2}\cup l_{w_3}$. Note that the order of germs at $u_0$ is identified with the orientation of $S^2$, as well as with the orientation of $\gamma^n\bigl(S^1\bigr)$ induced by the orientation of $S^2$. Since the orientation of $S^1$ corresponding, via $\gamma^1$, to the orientation of $\gamma^1\bigl(S^1\bigr)$ is equal to the orientation of $S^1$ corresponding, via $\gamma^n$, to the orientation of $\gamma^n\bigl(S^1\bigr)$, $\mu_+$ preserves the orientation of $\mathbf A^n$. Similarly, we can prove the statement for $\mu_-$.
    
    Now we show that $\mu_+(\mathbf A^{\infty})$ (resp.\ $\mu_-(\mathbf A^{\infty})$) is dense in $S^1(Y^+)$ (resp.\ $S^1(Y^-)$). Suppose $e\= \{C_K^e\}_{K\in\cK_{o^+}(Y^+)}$ is an end in $\cE(Y^+)$. For each closed interval $U\ni e$ in $S^1(Y^+)$, there is a finite subtree $K_U\ni o^+$ of $Y^+$ such that $I_{C_{K_U}^e}\=\bigl\{a\in\cE(Y^+):C_{K_U}^a= C_{K_U}^e\bigr\}\subseteq U\cap \cE(Y^+)$. Then we choose a vertex $u_0=\beta_m(b)$ for $b\in\sB^+_m$ in $C_{K_U}^e$, with $v_0\=\pi^+(u_0)$ is an $m$-vertex. For $t\in \cA_m^+(b)$, since $\mu_+(t)$ is an ideal gap with support $u_0$, and there are ends in $I_{C_{K_U}^e}$ converge to $\mu_+(t)$ by the above construction. Then we have $\mu_+(t)\in \overline{I_{C_{K_U}^e}}\subseteq U$. Then $\mu_+(\mathbf A^{\infty})$ is dense in $S^1(Y^+)$. Similarly, we can show that $\mu_-(\mathbf A^{\infty})$ is dense in $S^1(Y^-)$.
    
    Since $\mu_+$ has been defined on a dense subset $\mathbf A^{\infty}\subseteq S^1$, $\mu_+(\mathbf A^{\infty})$ is dense in $S^1(Y^+)$, and $\mu_+$ preserves the orientation of $\mathbf A^{\infty}$, we conclude that $\mu_+$ is a homeomorphism. Similarly, $\mu_-$ is also a homeomorphism. By our discussion of the orientation above, it is easy to see that $\mu_-\circ\mu_+^{-1}\: S^1(Y^+)\to S^1(Y^-) $ is an orientation-reversing homeomorphism.

\end{proof}

\subsection{Dynamics on the end circle}\label{sub_Dynamics_on_the_end_circle}

For an integer $m\ge1$, $\sigma\in\{+,\,-\}$, and $X\in\X_\sigma^m$, recall $G^\sigma(X)$ from Definition~\ref{def:ch4-sector-data} and write $D_m(X)\=G^\sigma(X)\cap S^1$.

\begin{lemma}\label{lem:ch6-finite-positions-ref}
Fix $\sigma\in\{+,\,-\}$ and let $f_\sigma\:Y^\sigma\to Y^\sigma$ be the map of Definition~\ref{def:ch4-tree-maps}. For each integer $n\ge0$, let $r_n^\sigma\:Y^\sigma\to Y_n^\sigma$ be the first-point retraction of Definition~\ref{def:ch4-tree-maps}. There is a unique continuous map $x_n^\sigma\:S^1\to Y_n^\sigma$ such that, for each integer $m\ge1$, $b\in\sB_m^\sigma$, and $t\in\cA_m^\sigma(b)$, we have 
\begin{equation}
\label{eq:ch6-marked-values-ref}
x_n^\sigma(t)=r_n^\sigma(\beta_m^\sigma(b)),
\end{equation}
where $\beta_m^\sigma\:\sB_m^\sigma\to Y_m^\sigma$ is the block attachment of Definition~\ref{def:ch4-tile-copy-data}.
Moreover, for each integer $m\ge n$, $r_n^\sigma\circ x_m^\sigma=x_n^\sigma$. For each integer $m\ge1$ and $X\in\X_\sigma^m$, $x_m^\sigma(D_m(X))\subseteq T_X^\sigma$. Finally, for $n\in\N_0$ and $t\in S^1$, we have
\begin{equation}\label{eq:ch5-position-dynamics}
x_n^\sigma(\mu(t))=f_\sigma(x_{n+1}^\sigma(t)).
\end{equation}
\end{lemma}

\begin{proof}
For each integer $n\ge0$, put $x_n^\sigma\=r_n^\sigma\circ x^\sigma$ on $\mathbf A^\infty$, where $x^\sigma$ is defined in Subsection~\ref{subsec:hairiness}. Then~\eqref{eq:ch6-marked-values-ref} holds, and $r_n^\sigma\circ x_m^\sigma=x_n^\sigma$ for each integer $m\ge n$ by Proposition~\ref{prop:pi-weak-tree}~\textup{(ii)}.

For each integer $m\ge1$ and $X\in\X_\sigma^m$, equations~\eqref{eq:ch4-copy-stability},~\eqref{eq:ch4-refined-copy-intersection}, and~\eqref{eq:ch4-refined-mark} give $G_X^\sigma\cap Y_m^\sigma=T_X^\sigma$. Since $G_X^\sigma$ is a subtree, $r_m^\sigma(G_X^\sigma)=T_X^\sigma$. For each integer $k>m$ and $V\in\X_\sigma^k$ with $A_{m,k}(V)=X$, iteration gives $r_m^\sigma(T_V^\sigma)\subseteq T_X^\sigma$.

For $t\in D_m(X)\cap\mathbf A^\infty$, choose an integer $k>m$ and a $k$-arc $J\subseteq D_m(X)$ with endpoint $t$. Such an arc exists by~\eqref{eq:finite-gap-trace}. Let $V$ be the color-$\sigma$ $k$-tile incident to $\gamma^k(J)$. Then $J\subseteq D_k(V)\subseteq D_m(A_{m,k}(V))$ by~\eqref{eq:finite-gap-trace} and~\eqref{eq:ch4-gap-transport}. The interior of $J$ lies in $S^1\smallsetminus\mathbf A^m$, where each point belongs to a unique $m$-gap trace. Hence $A_{m,k}(V)=X$. Lemma~\ref{lem:finite-gap-incidence} and Proposition~\ref{prop:compact-finite-trees}~\textup{(i)} give $x^\sigma(t)\in T_V^\sigma$. Hence, $x_m^\sigma(t)\in T_X^\sigma$ and $x_n^\sigma(t)\in r_n^\sigma(T_X^\sigma)$ for $m\ge n$.

Take the closure of $\{(s,x_n^\sigma(s)):s\in\mathbf A^\infty\}$ in $S^1\times Y_n^\sigma$. Let $(t,a)$ and $(t,b)$ belong to this closure. Fix an integer $m>n$. If $t\in\mathbf A^m$, choose $X,\,Y\in\X_\sigma^m$ so that the two $m$-arcs with endpoint $t$ lie in $D_m(X)$ and $D_m(Y)$, respectively, by~\eqref{eq:finite-gap-trace}. If $t\notin\mathbf A^m$, choose $Y=X$ so that $D_m(X)$ contains the unique $m$-arc containing $t$. The chosen arcs cover a neighborhood of $t$. The preceding inclusion and closedness of the projected subtrees give $a,\,b\in r_n^\sigma(T_X^\sigma)\cup r_n^\sigma(T_Y^\sigma)$. If $X\ne Y$, then $t\in\mathbf A^m\cap D_m(X)\cap D_m(Y)$, so both projected subtrees contain $x_n^\sigma(t)$. In either case, $r_n^\sigma(T_X^\sigma)\cup r_n^\sigma(T_Y^\sigma)$ is connected and contains the arc $[a,b]$ in the tree.

If $a\ne b$, some road image contains a nondegenerate subarc $\eta_n^\sigma(I)\subseteq[a,b]$ with $I\subseteq I_\lambda^\sigma$. The $m$-atoms of $I_\lambda^\sigma$ have images in distinct copies by~\eqref{eq:compact-stage-stability},~\eqref{eq:ch4-stage-model-map}, and Lemma~\ref{prop:ch4-road-data}~\textup{(i)}. For each $W\in\{X,\,Y\}$, the set $r_n^\sigma(T_W^\sigma)$ is a point or equals $T_W^\sigma\cap Y_n^\sigma$. By~\eqref{eq:ch4-stage-copy-intersection}, the image of each $m$-atom labelled by a tile other than $W$ meets $T_W^\sigma$ in at most one point. Thus the preimage of $r_n^\sigma(T_W^\sigma)$ in $I_\lambda^\sigma$ is contained in at most one $m$-atom and finitely many points. Since $\eta_n^\sigma(I)\subseteq[a,b]\subseteq r_n^\sigma(T_X^\sigma)\cup r_n^\sigma(T_Y^\sigma)$, we obtain
\begin{equation*}
 |I|\le 2\max_{J\in\mathcal P_m^{\lambda,\sigma}}|J|\to 0
\end{equation*}
by Lemma~\ref{prop:ch4-road-data}~\textup{(ii)}, which is a contradiction.

The projection $(t,y)\mapsto t$ is injective on $\overline{\{(s,x_n^\sigma(s)):s\in\mathbf A^\infty\}}$, and its compact image contains the dense set $\mathbf A^\infty$. Hence the projection is a homeomorphism onto $S^1$, and the second coordinate of its inverse gives the unique continuous extension of $x_n^\sigma$. The equality $r_n^\sigma\circ x_m^\sigma=x_n^\sigma$ extends by continuity. Density on each $m$-arc and closedness of $T_X^\sigma$ give $x_m^\sigma(D_m(X))\subseteq T_X^\sigma$.

For $v\in Y^\sigma$, put $a\=r_{n+1}^\sigma(v)$. Since $[v,a]\cap Y_{n+1}^\sigma=\{a\}$,~\eqref{eq:pi-exact-pullback} gives $f_\sigma([v,a])\cap Y_n^\sigma=\{f_\sigma(a)\}$. By Proposition~\ref{prop:pi-weak-tree}~\textup{(iii)}, $[f_\sigma(v),f_\sigma(a)]\subseteq f_\sigma([v,a])$. Thus, $r_n^\sigma(f_\sigma(v))=f_\sigma(a)$, and
\begin{equation*}
r_n^\sigma\circ f_\sigma=f_\sigma\circ r_{n+1}^\sigma.
\end{equation*}

Recall the sets $\mathbf A^m$ from~\eqref{eq:finite-angle-preimages}. We have that $\mathbf A^\infty$ is dense in $S^1$ and $\mu(\mathbf A^\infty)\subseteq\mathbf A^\infty$. For $t\in\mathbf A^\infty$, by~\eqref{eq:ch5-angle-point-dynamics} we have
$x_n^\sigma(\mu(t))=(r_n^\sigma\circ f_\sigma)(x^\sigma(t))=(f_\sigma\circ r_{n+1}^\sigma)(x^\sigma(t))=f_\sigma(x_{n+1}^\sigma(t))$.
By continuity,~\eqref{eq:ch5-position-dynamics} holds on $\overline{\mathbf A^\infty}=S^1$.
\end{proof}

For a hairy $Y^\sigma$, let $\mu_\sigma\:S^1\to S^1(Y^\sigma)$ be the angle parametrization of Proposition~\ref{prop_reversesorientationheom}.

\begin{prop}\label{prop:ch5-position-dynamics}
Let $f\:S^2\to S^2$ be an expanding Thurston map of degree $d$. Fix an $f$-invariant Jordan curve $\cC\supseteq\post f$ and a connection of $1$-tiles satisfying \textup{(C1)}--\textup{(C2)} of Definition~\ref{Propertiesofconnections}. Let $Y^\sigma$ be the associated ribbon tree for a fixed $\sigma\in\{+,\,-\}$. Suppose that $Y^\sigma$ is hairy. Let $f_\sigma\:Y^\sigma\to Y^\sigma$ be the map of Definition~\ref{def:ch4-tree-maps}.

For $n\in\N_0$, let $x_n^\sigma\:S^1\to Y_n^\sigma$ be the map given by Lemma~\ref{lem:ch6-finite-positions-ref} and recall $r_n^\sigma\=r_{Y_n^\sigma}^\sigma\:Y^\sigma\to Y_n^\sigma$ from Definition~\ref{def:ch4-tree-maps}. Identify $S^1$ with $\R/\Z$ and write $\mu(t)\=d\,t\bmod 1$.

For each $n\in \N_0$, we have
\begin{equation*}
x_n^\sigma=\tau_{Y_n^\sigma}\circ p_n^\sigma\circ\rho^\sigma\circ\mu_\sigma.
\end{equation*}
\end{prop}

\begin{proof}
For $b\in B_{n+1}^\sigma$, write $v\=\tau_{Y_{n+1}^\sigma}(b)$ and $a=\tau_{Y_n^\sigma}(\chi_n^\sigma(b))$. By the construction in Subsection~\ref{theendspace}, $a=v$ if $v\in Y_n^\sigma$. If $v$ lies in a component $C$ of $Y_{n+1}^\sigma\smallsetminus Y_n^\sigma$, then $a$ is the attachment point of $C$. Thus $[v,a]\cap Y_n^\sigma=\{a\}$ in both cases, so $a=r_n^\sigma(v)$. For each integer $n\ge0$, define $y_n^\sigma\:S^1\to Y_n^\sigma$ by $y_n^\sigma\=\tau_{Y_n^\sigma}\circ p_n^\sigma\circ\rho^\sigma\circ\mu_\sigma$. The identity $a=r_n^\sigma(v)$ gives $\tau_{Y_n^\sigma}\circ\chi_n^\sigma=r_n^\sigma\circ\tau_{Y_{n+1}^\sigma}$. Since $p_n^\sigma=\chi_n^\sigma\circ p_{n+1}^\sigma$, it follows that $y_n^\sigma=r_n^\sigma\circ y_{n+1}^\sigma$. By Proposition~\ref{prop:pi-weak-tree}~\textup{(ii)}, $r_n^\sigma$ fixes $Y_n^\sigma$ and $r_n^\sigma\circ r_j^\sigma=r_n^\sigma$ for each integer $j\ge n$. Hence $r_n^\sigma\circ y_m^\sigma=y_n^\sigma$ for integers $m\ge n\ge0$.

For $t\in\mathbf A^\infty$, take an integer $N\ge1$ with $t\in\mathbf A^N$. By the definition of $x^\sigma$ in Subsection~\ref{subsec:hairiness}, the point $u=x^\sigma(t)$ lies in $Y_N^\sigma$. For each integer $m\ge N$, write $b_m\=p_m^\sigma(\rho^\sigma(\mu_\sigma(t)))$. If $u$ is a hard endpoint, write $\omega\=\mu_\sigma(t)$ for the hard end at $u$. Then $b_m=q_m^\sigma(\omega)$ because $\rho^\sigma$ extends $q^\sigma$. Since $A_\omega=[o^\sigma,u]\subseteq Y_m^\sigma$, the definition of $q_m^\sigma$ in Subsection~\ref{theendspace} yields $\tau_{Y_m^\sigma}(b_m)=u$.

Suppose $u$ is a cut point and $\tau_{Y_m^\sigma}(b_m)\ne u$. Take an open arc $J\subseteq B_m^\sigma$ containing $b_m$ with $u\notin\tau_{Y_m^\sigma}(J)$. By Lemma~\ref{lem:genuine-end-density} and the surjectivity of $p_m^\sigma$, take genuine ends $\omega_1,\omega_2$ with $s_i=q_m^\sigma(\omega_i)$ in different components of $J\smallsetminus\{b_m\}$. The tails $R_{\omega_i}\smallsetminus Y_m^\sigma$ attach at $\tau_{Y_m^\sigma}(s_i)$, so both lie in the component of $Y^\sigma\smallsetminus\{u\}$ containing $\tau_{Y_m^\sigma}(J)$. For large $k>m$, let $a_i\=q_k^\sigma(\omega_i)$ for $i=1,\,2$ and take $a_0\in B_k^\sigma$ with $\tau_{Y_k^\sigma}(a_0)=u$. The median of the three projections lies in that component. At this median, Subsection~\ref{theendspace} uses the direction toward $u$ to determine the cyclic order of $(\omega_1,\mu_\sigma(t),\omega_2)$. This is also the cyclic order of $(a_1,a_0,a_2)$. The maps $p_m^\sigma\circ\rho^\sigma$ and $\chi_{m,k}^\sigma$ send these triples to $(s_1,b_m,s_2)$ and $(s_1,\chi_{m,k}^\sigma(a_0),s_2)$, respectively. Since $\chi_{m,k}^\sigma(a_0)$ projects to $u$, it lies outside $J$. The image triples have opposite cyclic orders, contradicting order preservation. Thus, $\tau_{Y_m^\sigma}(b_m)=u$.

In either case, $y_m^\sigma(t)=u$ for each integer $m\ge N$. For each integer $n\ge0$, choose an integer $m\ge\max\{n,\,N\}$. Then $y_n^\sigma(t)=r_n^\sigma(u)=x_n^\sigma(t)$, where the last equality follows from~\eqref{eq:ch6-marked-values-ref}. Since $x_n^\sigma$ and $y_n^\sigma$ are continuous and $\mathbf A^\infty$ is dense in $S^1$, we have $y_n^\sigma=x_n^\sigma$ on $S^1$.
\end{proof}

\begin{theorem}\label{prop_z^d}
Let $f\:S^2\to S^2$ be an expanding Thurston map of degree $d$. Suppose $f$ admits an $f$-invariant Jordan curve $\cC\supseteq\post f$ and a connection of $1$-tiles satisfying \textup{(C1)}--\textup{(C2)} of Definition~\ref{Propertiesofconnections}. For a fixed $\sigma\in\{+,\,-\}$, suppose that the associated ribbon tree $Y^\sigma$ is hairy. Let $f_\sigma\:Y^\sigma\to Y^\sigma$ be the map of Definition~\ref{def:ch4-tree-maps}. Recall $\mu_\sigma\:S^1\to S^1(Y^\sigma)$ from Proposition~\ref{prop_reversesorientationheom}. Write $\eta_\sigma\=\mu_\sigma^{-1}\:S^1(Y^\sigma)\to S^1$ for the angle coordinate, and let $\mu(t)\=dt\bmod 1$ on $S^1\=\R/\Z$.

Then, for every $\omega\in\mathcal E_g(Y^\sigma)$ there is a unique $\omega'\in\mathcal E_g(Y^\sigma)$ such that, for every $n\in\N_0$ and every interior subtree $K\subseteq Y_n^\sigma$, $f_\sigma(R_\omega\smallsetminus Y_{n+1}^\sigma)\subseteq C_K^{\omega'}$.

Moreover, the map $\omega\mapsto\omega'$ has a unique continuous extension $\overline f_\sigma\:S^1(Y^\sigma)\to S^1(Y^\sigma)$. The extension satisfies $\eta_\sigma\circ\overline f_\sigma=\mu\circ\eta_\sigma$ and is a topologically exact covering of degree $d$. There is an orientation-preserving homeomorphism $\theta_\sigma\:S^1(Y^\sigma)\to S^1$ conjugating $\overline f_\sigma$ to $z\mapsto z^d$.
\end{theorem}

\begin{proof}
Fix $\omega\in\mathcal E_g(Y^\sigma)$ and put $t\=\eta_\sigma(\omega)$. The arcs $R_\omega\cap Y_n^\sigma=[o^\sigma,x_n^\sigma(t)]$ exhaust the proper ray $R_\omega$. Thus, for each fixed $j$, $x_{n+1}^\sigma(t)\notin Y_{j+1}^\sigma$ for all sufficiently large integer $n$. By~\eqref{eq:ch5-position-dynamics}, put $a_n\=x_n^\sigma(\mu(t))=f_\sigma(x_{n+1}^\sigma(t))$. By~\eqref{eq:pi-exact-pullback}, $a_n\notin Y_j^\sigma$ for all sufficiently large $n$.

Lemma~\ref{lem:ch6-finite-positions-ref} gives $r_n^\sigma(a_m)=a_n$ for each integer $m\ge n$, so the arcs $[o^\sigma,a_n]$ are nested. Parametrize their union $R$ by successive subarcs as in the proof of Proposition~\ref{prop:rooted-representatives}. Once $a_n\notin Y_j^\sigma$, we have $R\cap Y_j^\sigma\subseteq[o^\sigma,a_n]$. Since every compact subset of $Y^\sigma$ lies in some $Y_j^\sigma$ by Proposition~\ref{prop:pi-weak-tree}~\textup{(i)}, $R$ is a proper ray. Let $\omega'$ be the genuine end it represents.

Fix $n$ and an interior subtree $K\subseteq Y_n^\sigma$. By~\eqref{eq:pi-exact-pullback}, $f_\sigma(R_\omega\smallsetminus Y_{n+1}^\sigma)$ is a connected subset of $Y^\sigma\smallsetminus K$. It contains $a_m\in C_K^{\omega'}$ for all sufficiently large $m$, so it lies in $C_K^{\omega'}$. Note that every interior subtree $K$ lies in some $Y_n^\sigma$. For such an integer $n$, the nonempty image $f_\sigma(R_\omega\smallsetminus Y_{n+1}^\sigma)$ determines $C_K^{\omega'}$ uniquely. Hence, $\omega'$ is unique.

Put $b\=\rho^\sigma\bigl(\eta_\sigma^{-1}(\mu(t))\bigr)$. For each integer $n$ and all sufficiently large integers $m>n$, we have $\tau_{Y_m^\sigma}(p_m^\sigma(b))=a_m\in R\smallsetminus Y_n^\sigma$. By the definition of $q_n^\sigma$, $p_n^\sigma(b)=\chi_{n,m}^\sigma(p_m^\sigma(b))=q_n^\sigma(\omega')$. Thus, $b=q^\sigma(\omega')=\rho^\sigma(\omega')$, and $\eta_\sigma(\omega')=\mu(t)$.

Consequently, $\overline f_\sigma=\eta_\sigma^{-1}\circ\mu\circ\eta_\sigma$ is a continuous extension of $\omega\mapsto\omega'$. It is unique since $\mathcal E_g(Y^\sigma)$ is dense in $S^1(Y^\sigma)$ by Lemma~\ref{lem:genuine-end-density}.

Let $J\subseteq S^1$ be an open arc with a lift $I\subseteq\R$ of length $0<\ell<1$. The preimage $\mu^{-1}(J)$ consists of the $d$ disjoint arcs $(I+j)/d\bmod 1$ for $0\le j<d$. On each of these arcs, $\mu$ is an orientation-preserving homeomorphism onto $J$. Also, $\mu^k(J)=S^1$ once $d^k\times \ell>1$. Thus $\mu$ is a topologically exact covering of degree $d$. By conjugacy with $\mu$, $\overline f_\sigma$ is a topologically exact covering of degree $d$. Thus, for every $z\in S^1$, $\bigl(\theta_\sigma\circ\overline f_\sigma\circ\theta_\sigma^{-1}\bigr)(z)=z^d$.
\end{proof}

Now, we are ready to prove Theorem~\ref{B}.

\begin{proof}[\textbf{Proof of Theorem~\ref{B}}]
By Proposition~\ref{lem:hairy}~\textup{(i), (iv)}, we have~\textup{(i)}.

Fix $\sigma\in\{+,\,-\}$ and suppose that $Y^\sigma$ is hairy. By Theorem~\ref{thm:end-circle}, $S_\sigma=S^1(Y^\sigma)$ is an oriented topological circle and the genuine ends of $Y^\sigma$ are dense in it. The homeomorphism $\mu_\sigma\:S^1\to S_\sigma$ of Subsection~\ref{theendspace} is orientation preserving for $\sigma=+$ and orientation reversing for $\sigma=-$. Its inverse $\eta_\sigma=\mu_\sigma^{-1}$ has the same orientation property. By Theorem~\ref{prop_z^d}, the map induced by $f_\sigma$ on genuine ends extends uniquely to a topologically exact covering $\overline f_\sigma=\eta_\sigma^{-1}\circ\mu_d\circ\eta_\sigma$ of degree $d$. Thus $\eta_\sigma\circ\overline f_\sigma=\mu_d\circ\eta_\sigma$, proving~\textup{(ii)}.

If both trees are hairy, $h=\eta_-^{-1}\circ\eta_+$ pairs points with the same angle and reverses orientation. Moreover, we have $h\circ\overline f_+=\eta_-^{-1}\circ\mu_d\circ\eta_+=\overline f_-\circ h.$
This proves~\textup{(iii)}.
\end{proof}

\section{The end circle and laminations}
\label{The_end_circle_and_laminations}

In this section, we use the maps $x_n^\sigma$ of Lemma~\ref{lem:ch6-finite-positions-ref} to define the equivalence relations $R_+$ and $R_-$ on $S^1$. We identify them with Meyer's closed white and black relations, show that both are laminar, and recover the invariant Peano curve $\gamma$, proving Theorem~\ref{C}.

\begin{definition}
\label{def_equivalencerelationsassociatedtoY-ref}
Fix $\sigma\in\{+,\,-\}$. For each integer $n\ge0$, let $x_n^\sigma\:S^1\to Y_n^\sigma$ be the continuous map of Lemma~\ref{lem:ch6-finite-positions-ref}, and write
\begin{equation*}
    \ker x_n^\sigma\=\bigl\{(s,t)\in S^1\times S^1: x_n^\sigma(s)=x_n^\sigma(t)\bigr\}.
\end{equation*}
Define $R_\sigma\=\bigcap_{n=0}^{+\infty}\ker x_n^\sigma$.
\end{definition}

\begin{prop}\label{relationsameprop}
Let $f\:S^2\to S^2$ be an expanding Thurston map with an invariant Jordan curve $\cC\supseteq\post f$. For a connection of $1$-tiles satisfying \textup{(C1)}--\textup{(C2)} of Definition~\ref{Propertiesofconnections}, the relations $R_\sigma$ of Definition~\ref{def_equivalencerelationsassociatedtoY-ref} are Meyer's closed relations of the respective colors.
\end{prop}

\begin{proof}
Let $q_\sigma\:S^1\to S^1/{\sim_\sigma}$ be the quotient induced by Meyer's color-$\sigma$ relation. Choose a compatible metric on $S^1/{\sim_\sigma}$. By~\eqref{eq:ch6-marked-values-ref}, $R_\sigma$ contains all finite angle relations. Since $R_\sigma$ is closed, $\ker q_\sigma\subseteq R_\sigma$ by Definition~\ref{oversetwsim}.

For each integer $m\ge1$ and $X\in\X_\sigma^m$, list the $\card\post f$ $m$-arcs of $D_m(X)$ in their cyclic order along $\partial G^\sigma(X)$. For two successive arcs in this list, one endpoint of each belongs to the same level-$m$ angle class of color $\sigma$ \cite[proof of Lemma~5.9~\textup{(G2)} and Remark 5.6]{Meyer09}. Their images under $q_\sigma$ therefore intersect. The triangle inequality gives
\begin{equation*}
\varepsilon_m:=\max_{X\in\X_\sigma^m}\diam q_\sigma(D_m(X))
 \le (\card\post f)\max_{\alpha\in\mathbf A^m}
\diam q_\sigma(a_+^m(\alpha)).
\end{equation*}
By~\eqref{diamadn} and uniform continuity of $q_\sigma$, the right-hand side tends to zero as $m\to+\infty$.

For $(s,t)\in R_\sigma$ and each integer $m\ge1$, choose $X,Y\in\X_\sigma^m$ with $s\in D_m(X)$ and $t\in D_m(Y)$. Lemma~\ref{lem:ch6-finite-positions-ref} gives $x_m^\sigma(s)\in T_X^\sigma\cap T_Y^\sigma$. If $X\ne Y$, their common block $b$ from~\eqref{eq:ch4-stage-copy-intersection} satisfies $q_\sigma(\cA_m^\sigma(b))\subseteq q_\sigma(D_m(X))\cap q_\sigma(D_m(Y))$ by Lemma~\ref{lem:finite-gap-incidence}. Thus, $\operatorname{dist}(q_\sigma(s),q_\sigma(t))\le2\varepsilon_m$ for all $m$, also when $X=Y$. Hence $R_\sigma=\ker q_\sigma$.
\end{proof}

Now we recall the definition of the laminar equivalence relation.
\begin{definition}[{Calegari \& Loukidou \cite[Definition~1.16]{CLWheels}}]
An equivalence relation on $S^1$ is \emph{laminar} if the equivalence classes are closed as subsets of $S^1$, if the relation is closed in the space $\mathcal M$ of unordered distinct pairs, and if distinct equivalence classes are non-crossing.
\end{definition}

\begin{prop}\label{lemma_rker-ref}
For each $\sigma\in\{+,\,-\}$, the relation $R_\sigma$ of Definition~\ref{def_equivalencerelationsassociatedtoY-ref} is laminar.
\end{prop}

\begin{proof}
Fix $\sigma$ and let $\mathcal M$ be the space of unordered pairs of distinct points of $S^1$. Let $\Lambda$ be the closure in $\mathcal M$ of the geodesic sides of the convex hulls of the finite angle classes of color $\sigma$. For a two-point class, take the single geodesic joining its points.

For integers $1\le m\le k$, each level-$m$ hull lies in a level-$k$ hull by~\eqref{eq:ch4-block-persistence}. Its sides cannot cross a side of the containing hull. Distinct level-$k$ hulls are disjoint \cite[Lemma~5.7~\textup{($L^n2$)} and Remark 5.6]{Meyer09}, so no level-$m$ side crosses a level-$k$ side. If two geodesics in $\Lambda$ crossed, choose disjoint open arcs around their four alternating endpoints. By the definition of $\Lambda$, two finite-level sides have endpoints in the corresponding arcs and would also cross. Thus the closed set $\Lambda$ is a lamination.

Let $\operatorname{rel}(\Lambda)$ denote the smallest closed equivalence relation on $S^1$ that identifies the endpoints of each leaf of $\Lambda$ \cite[Definition~2.3]{CLWheels}. For each integer $m\ge1$ and $b\in\sB_m^\sigma$, the sides of the hyperbolic convex hull $L(\cA_m^\sigma(b))$ connect all angles in $\cA_m^\sigma(b)$. Thus identifying the endpoints of all these sides generates the finite angle relations of color $\sigma$. Since $\Lambda$ is the closure of these sides in $\mathcal M$, $\operatorname{rel}(\Lambda)$ is the smallest closed equivalence relation containing those finite angle relations. By Proposition~\ref{relationsameprop} and Definition~\ref{oversetwsim}, we have $R_\sigma=\operatorname{rel}(\Lambda)$. Hence $R_\sigma$ is laminar by \cite[the paragraph following Definition~2.3]{CLWheels}.
\end{proof}

Now, we can prove Theorem~\ref{C}.
\begin{proof}[\textbf{Proof of Theorem~\ref{C}}]
Proposition~\ref{relationsameprop} and Lemma~\ref{closureequi}~\textup{(i)} give $R\=\widehat{R_+\vee R_-}=\ker\gamma$. Thus $H([t])=\gamma(t)$ induces a homeomorphism from the compact quotient onto $S^2$.

Let $q\colon S^1\to S^1/R$ be the quotient map. By~\eqref{eq:ch5-position-dynamics}, for each $\sigma\in\{+,\,-\}$, if $(t,t')\in R_\sigma$, then $(\mu_d(t),\mu_d(t'))\in R_\sigma$. Hence, $\mu_d$ preserves $R_+\vee R_-$. Since $\mu_d$ is continuous and $R=\widehat{R_+\vee R_-}$, it preserves $R$. Thus, $F([t])=[\mu_d(t)]$ is well-defined on $S^1/R$ and is continuous by the quotient property. It satisfies $F\circ q=q\circ\mu_d$.

Let $f'\=H\circ F\circ H^{-1}\colon S^2\to S^2$. For each $t\in\mathbf A^\infty$, we have
\begin{equation*}
f'(\gamma(t))=H(F(q(t)))=H(q(\mu_d(t)))=\gamma(\mu_d(t))=f(\gamma(t)),
\end{equation*}
where the last equality follows from~\eqref{commutediaav}. Since $\gamma(\mathbf A^\infty)$ is dense in $S^2$, continuity gives $f'=f$, i.e., $H\circ F=f\circ H$.
\end{proof}

\section{Relations to the CaTherine wheels}\label{Connection_graphs_with_branching}
In this section, we establish a property of branches $Z^+$ and $Z^-$ of the P-zipper under the assumption that the connection of $1$-tiles satisfies the property in Lemma~\ref{degree3}. This property is an analogy of the construction of $Z^+$ and $Z^-$ in \cite[Section~1]{CLWheels}.

We now fix an expanding Thurston map $f$ and an invariant Jordan curve $\cC$ as in Lemma~\ref{degree3}. We also fix a connection of $1$-tiles satisfying the properties stated in Lemma~\ref{degree3}. By combining $l_{p,q}$ with the invariant Peano curve $\gamma$ from Theorem~\ref{Peanocurves}, we establish some properties of $Z^+$ and $Z^-$. The main result of this subsection is the following proposition.

\begin{prop}\label{boundaryzipper}
    For $p,\,q\in\post f$, let $F^{p,q}$ and $l_{p,q}\=F^{p,q}(I)$ be given as in Subsection~\ref{Construction_of_primary_zippers}. Let $G^{p,q}$ be the path for black tiles and set $g_{p,q}\=G^{p,q}(I)$. Suppose $E\in S^1$ is a $0$-arc with $\gamma^0(E)\=[p,q]\subseteq\cC$. Then we have that $\partial \gamma(E)=l_{p,q}\cup g_{p,q}$.
\end{prop}

To prove this proposition, we establish several lemmas. Recall the maps $g_{n,m}$ from Lemma~\ref{lem:ch4-sector-ancestors}. By~\eqref{diamgmsmall}, there exists a constant $C_2>0$ such that $\rho(x,g_{n,m}(x))<C_2\Lambda^{-m}$ for all integers $n\ge m\ge1$ and all $x\in S^2$.

\begin{lemma}\label{existleaftile}
    Let $E^0$ be a $0$-edge of $(f,\cC)$. Suppose the connection of $1$-tiles is given and satisfies the conditions in Definition~\ref{Propertiesofconnections}. Suppose the induced white (resp.\ black) $1$-connection graph contains a vertex of degree at least $3$. Then there exist $k\in\N$ and a connected subarc $U\subseteq E^0$ such that $g_{k-1,0}(U)$ coincides with the boundary of a white (resp.\ black) $k$-tile $X$, which is the union of all $k$-edges on its boundary.
\end{lemma}
\begin{proof}
    We only prove the statement for the white case, since the black case can be treated similarly. Let $E_1^0,\,E_2^0,\,\dots,\,E_k^0$ be the $0$-edges in their cyclic order on $\cC$, where $k\=\card\post f$. Let $d$ be the degree of $f$. We consider the leaf $c(X)$ (i.e., vertices of degree $1$) of the induced $n$-connection graph for $n\in\N$. The vertex $c(X)$ represents the $n$-tile $X$ and is incident to exactly one edge. 
    
    First, we show that the number of leaves in the induced $n$-connection graph tends to infinity as $n\to+\infty$. Note that by Lemma~\ref{degree3}, our induced white $1$-connection graph has a branched vertex with degree at least $3$. In other words, there is a white $1$-tile that connects to three other white $1$-tiles. By Lemma~\ref{lem:ch4-sector-copies}, for each white $n$-tile $X^n$ and white $n+1$-tiles $X^{n+1},\,Y^{n+1}$ with $c\bigl(X^{n+1}\bigr),\,c\bigl(Y^{n+1}\bigr)\in\Gamma_{n+1}^+[X^n]$, $X^{n+1},\,Y^{n+1}$ share a vertex $v$ are connected at $v$ if and only if $1$-tiles $f^n\bigl(X^{n+1}\bigr)$ and $f^n\bigl(Y^{n+1}\bigr)$ are connected at $f^n(v)$. Hence there is a white $(n+1)$-tile in $\Gamma_{n+1}^+[X^n]$ (here we mean that the corresponding tile vertex is contained in $\Gamma_{n+1}^+[X^n]$) that connects to three other white $1$-tiles in $\Gamma_{n+1}^+[X^n]$. It follows that the number of branched vertices in the induced white $(n+1)$-connection graph is at least $d^{n}$, which is the number of white $n$-tiles. It follows that the number of leaves in the induced $n$-connection graph tends to infinity as $n\to+\infty$ at the same time.

    Next we show that the statement holds for a certain $0$-edge. Suppose there is $N_1\in\N$ such that for an integer $m>N_1$ there are more than $2k+1$ leaves in the induced $m$-connection graph as shown above. Then we may choose $c(X^m)$ such that it is a leaf of the induced white $m$-connection graph and the $m$-tile $X^m$ does not contain one of the marked succeeding $m$-edges, because the number of the marked succeeding $m$-edges is at most $2k$. Then there is a $0$-edge $E_{i_0}^0$ such that all the $m$-edges on the boundary of $X^m$ are contained in $g_{m-1,0}\bigl(E_{i_0}^0\bigr)$. Suppose $v^m\in\partial_{o}X^m$ is the unique $m$-vertex at which $X^m$ is connected to other $m$-tiles. Since $c(X^m)$ is a leaf, every pair of adjacent $m$-edges on $\partial_{o}X^m$, except for the pair containing $v^m$, consists of succeeding edges. Thus, the union of all $m$-edges on $\partial_{o}X^m$ is the image under $g_{m-1,0}$ of a connected subarc of $E_{i_0}^0$.

    Finally, we show that the statement holds for all $0$-edges, including $E^0$. By Lemma~\ref{lem:cell mapping properties of Thurston map}, we have that $\delta_n\=\max\{ \diam_{\rho}(E) : E \in \E^n(f,\mathcal{C})\}\to 0$ as $n\to +\infty$. Here $\rho$ is a visual metric for $f$. Since $g_{n-1,0}\bigl(E^0_i\bigr)$ is a union of several $n$-edges, contains $p$ and $q$, and is connected (in the sense of standard topology), the number of $n$-edges in $g_{n-1,0}\bigl(E^0_i\bigr)$ is larger than $\rho(p,q)/\delta_n$. Then the number of $n$-edges contained in $g_{n-1,0}\bigl(E^0_i\bigr)$ tends to infinity as $n\to+\infty$. We label an $n$-edge $E^n$ by $j$ if $f^n(E^n)=E_j^0$. Note that the labels of $n$-edges contained in $g_{n-1,0}\bigl(E^0_i\bigr)$ are in their cyclic order on $g_{n-1,0}\bigl(E^0_i\bigr)$ since $f$ is orientation-preserving and the white $n$-tiles are always on one side of $g_{n-1,0}\bigl(E^0_i\bigr)$. Then there is $N_2\in\N$ such that for integer $n>N_2$ there is an $n$-edge $E^n$ with label $i_0$ contained in $g_{n-1,0}\bigl(E^0\bigr)$. Suppose $U^0\subseteq E^0$ is a subarc and $g_{n-1,0}\bigl(U^0\bigr)=E^n$. Since $(f^n\circ g_{n+m-1,0})\bigl(U^0\bigr)=\bigl(f^n\circ H^{n+m-1}_1\circ H_1^{n+m-2}\circ\dots\circ H_1^n\bigr)(E^n)=\bigl(H^{m-1}_1\circ H_1^{m-2}\circ\dots\circ H_1^0\circ f^n\bigr)(E^n)=g_{m-1,0}\bigl(E_{i_0}^0\bigr)$, there is a connected subset $U\subseteq E^0$ such that for $k'=n+m$, $g_{k'-1,0}(U)$ is the union of $k'$-edges on the boundary of a white $k'$-tile $X'
    \subseteq f^{-n}(X^m)$. 
\end{proof}

\begin{remark}
    Let $U'\subseteq S^1$ be the preimage of $U$ under $\gamma^0$. Then $U'$ is a subarc in $S^1$ with endpoints $a$ and $b$. Since $g_{k'-1,0}(U)$ is the union of $k'$-edges on the boundary of a white $k'$-tile $X$, we have that $\gamma^{k'}(a)=\gamma^{k'}(b)$. Hence $a\overset{k',\white}\sim b$.
\end{remark}

\begin{lemma}
    Let $E^0$ be a $0$-edge and $p\neq q\in\post f$ be the endpoints of $E^0$. Let $E_0$ be the~$0$-arc with $\gamma^0(E_0)=E^0$. Then $p\in S^2$ is on the boundary of $\gamma(E_0)$.
\end{lemma}
\begin{proof}
    First we show that in an arbitrary subarc $I_0\subseteq S^1$, there are distinct points in the same equivalence class under $\overset{\white}\sim$. By Lemma~\ref{existleaftile}, we know that between two adjacent postcritical points $p'$ and $q'$, there exist two distinct points that are in the same equivalence class under the equivalence relation $\overset{k,\white}\sim$. We denote the image of the $0$-arc $E_1$ under $\gamma^0$ as the $0$-edge $E_1^0=\bigl[p',q'\bigr]$. Now we consider an arbitrary subarc $I_0\subseteq S^1$. Recall~\eqref{diamadn} that there is a constant $C(\lesssim)>0$ such that $\diam_{S^1} a^m\lesssim d^{-m}$
    for each $m$-arc $a^m$. Here $d$ is the degree of $f$. Then when $m$ is large enough, there is an $m$-arc $a$ contained in $I$ such that $\mu^m(a)=E_1$. Then there will be two distinct points in $a\subseteq I_0$ in the same equivalence class under the equivalence relation $\overset{k+m,\white}\sim$. Recall from Subsection~\ref{wbequivalencerelations} that $\overset{n,\white}\sim\leq \overset{\white}\sim$ for all $n\in\N$. Then we have proved that, in an arbitrary subarc $I_0\subseteq S^1$, there are distinct points in the same equivalence class under $\overset{\white}\sim$. Similarly, we have that this is also true for $\overset{\black}\sim$. 

    Then we show that there is a point in $I_0$ that cannot be equivalent to points on $S^1\smallsetminus I_0$ under $\sim$. Suppose $p_1,\,q_1\in I_0$ are in the same equivalence class under $\overset{\white}\sim$. Denote the interval $[p_1,q_1]\in I_0$ by $I_1$. Choose $I_2\=[p_2,q_2]\subseteq I_1$ such that $p_2,\,q_2$ are in the same equivalence class under $\overset{\black}\sim$. We continue this procedure and obtain nested intervals $\{I_i\}_{i=0}^{+\infty}$ with $I_i=[p_i,q_i]$, satisfying that $p_{2j+1}\overset{\white}\sim q_{2j+1}$ and $p_{2j}\overset{\black}\sim q_{2j}$ for $j\in\N$. 

    Let $x_0\in S^2\smallsetminus\post  f$ be a point with the smallest number of preimages under $\gamma$, and enumerate $\gamma^{-1}(x_0)=\{t_1,\,t_2,\,\dots,\,t_M\}$ by Theorem~\ref{finitefiber}. Consider preimages, and set $\{s_1,\,s_2,\,\dots,\,s_{dM}\}\=\mu^{-1}(\{t_1,\,t_2,\,\dots,\,t_M\})$, and $\{w_1,\,w_2,\,\dots,\,w_d\}\=f^{-1}(x_0)$. Since $f\circ\gamma=\gamma\circ\mu$, it follows that $\gamma(\{s_1,\,s_2,\,\dots,\,s_{dM}\})=\{w_1,\,w_2,\,\dots,\,w_d\}$ and $\gamma^{-1}(\{w_1,\,w_2,\,\dots,\,w_d\})=\{s_1,\,s_2,\,\dots,\,s_{dM}\}$. By the minimality of $\card\gamma^{-1}(x_0)$, it follows that $\card\gamma^{-1}(w_j)\geq M$ for all $j\in\{1,\,2,\dots,\,d\}$. Thus $\card\gamma^{-1}(w_j)=M$ for all $j\in\{1,\,2,\,\dots,\,d\}$. The same argument yields that for all $w\in f^{-n}(x_0)$ there are $M$ points in $\mu^{-n}(\{t_1,\,t_2,\,\dots,\,t_M\})$ that are mapped by $\gamma$ to $w$. It follows that $V_M\=\bigl\{t:\card\gamma^{-1}(\gamma(t))=M\bigr\}$ is dense in $S^1$.

    Then we consider $I_{2M}$ and show that the point in $I_{2M}\cap V_M$ can not be equivalent to points on $S^1\smallsetminus I_0$ under $\sim$. We argue by contradiction and suppose that there are $t\in I_{2M}\cap V_M$ and $s\in S^1\smallsetminus I_0$ such that $t\sim s$. Suppose $\gamma^{-1}(\gamma(t))\=\{t_1,\,t_2,\,\dots,\,t_M\}$ and $t_1=t,\,t_M=s$. We choose disjoint open intervals $J_i\ni t_i$ for all $i\in\{1,\,2,\,\dots,\,M\}$ such that $p_j,\,q_j\notin J_i$ unless $p_j=t_i$ or $q_j=t_i$ for all $j\in\{1,\,\dots,\,2M\}$. Since $s\notin I_0$ we can also assume that $J_M\cap I_0=\emptyset$. Since $S^1\smallsetminus \bigl(\bigcup_{i=0}^{M}J_i \bigr)$ is closed and $\gamma$ is continuous, there is a neighborhood $O\ni\gamma(t)$ such that $O\cap\gamma\bigl(S^1\smallsetminus \bigl(\bigcup_{i=1}^{M}J_i \bigr)\bigr)=\emptyset$. Then $\gamma^{-1}(O)\subseteq \bigcup_{i=1}^{M}J_i$. Let $J_i(O)\=J_i\cap \gamma^{-1}(O)$.

    By Lemma~\ref{closureequi}, we have that $\sim$ is the closure of $\overset{\white}{\sim}\vee\overset{\black}{\sim}$. Let $j_0\in\{1,\,2,\,\dots,\,2M\}$ such that $I_{j_0}\smallsetminus I_{j_0+1}=(p_{j_0},p_{j_0+1})\cup(q_{j_0+1},q_{j_0})$ has no intersection with $\bigcup_{i=1}^M J_i$, then it has no intersection with $\gamma^{-1}(O)$. Let $L(O)\=I_{j_0}\cap \gamma^{-1}(O)$ and $R(O)\=\gamma^{-1}(O)\smallsetminus L(O)$. Suppose there are $s_1,\,\dots,\,s_N\in S^1$ such that $s'=s_1\overset{\white}\sim s_2\overset{\black}\sim s_3\cdots s_{N-2}\overset{\white}\sim s_{N-1}\overset{\black}\sim s_N=t'$ with $t'\in L(O)\subseteq I_{j_0}$ and $s'\in R(O)$. Then $\gamma(s_i)\in O$ for all $i\in\{1,\,2,\,\dots,\,N\}$. It follows that $s_i\in J_j(O)$ for some $j\in\{1,\,2,\,\dots,\,M\}$. Recall the equivalence relation $\overset{\white}{\sim}$ and $\overset{\b}{\sim}$ has non-crossing property (see Lemma~\ref{lemma_wbnoncrossing}). It follows that if $a\in I_{j_0+1}$ and $b\in S^1\smallsetminus I_{j_0}$, then $a$ and $b$ can not belong to the same equivalence class under $\overset{\white}{\sim}$ or $\overset{\b}{\sim}$. Then there exist $s_i\in I_{j_0}\smallsetminus I_{j_0+1}$ for some $i\in\{1,\,2,\,\dots,\,N\}$, which contradicts $s_i\in J_j(O)$. Hence if $a\in L(O)$ and $b\in R(O)$, then $a$ and $b$ can not belong to the same equivalence class under $\overset{\white}{\sim}\vee\overset{\black}{\sim}$. Set $E\=\bigl\{(s,u)\in S^1\times S^1:\gamma(s)=\gamma(u)\bigr\}$ and $E'\=E\smallsetminus((L(O)\times R(O))\cup(R(O)\times L(O)))$. The removed set is open in $\bigl(S^1\bigr)^2$, hence $E'$ is closed. Above $O$, each $\gamma$-fiber is divided into its $L(O)$-part and $R(O)$-part, while the fibers outside $O$ remain unchanged. Note that $(a,b)\notin ((L(O)\times R(O))\cup(R(O)\times L(O)))$ if $a$ and $b$ are in the same equivalence class under $\overset{\white}{\sim}\vee\overset{\black}{\sim}$. Thus, $E'$ is a closed equivalence relation bigger than $\overset{\white}{\sim}\vee\overset{\black}{\sim}$, contradicting the minimality of $E$.

    Now we come back to this lemma. We prove it by contradiction and assume that $p$ is in the interior of $\gamma(E_0)$. Then there is an open neighborhood $U$ of $p$ that is contained in $\gamma(E_0)$. Since $\gamma$ is continuous, $\gamma^{-1}(U)$ is an open set containing $\bigl(\gamma^0\bigr)^{-1}(p)$. Let $I\subseteq \gamma^{-1}(U)$ be an interval and $I\cap E_0=\emptyset$. By the above discussion, we have that there exists $t_0\in I$ that cannot be equivalent to points on $S^1\smallsetminus I$ under $\sim$, which means that $\gamma(t_0)\notin \gamma(E_0)$. This is a contradiction.
\end{proof}

Now we are ready to prove Proposition~\ref{boundaryzipper}.

\begin{proof}[Proof of Proposition~\ref{boundaryzipper}]
    For an integer $n\ge1$ and $t\in S_n^{p,q}$, equation~\eqref{eq:node-realization} gives that $F_n^{p,q}(t)=v_n$ is an $n$-vertex. Suppose $t_0\in S^1$ and $\gamma(t_0)=v_n$. By our construction, $F_n^{p,q}(I)\cup G_n^{p,q}(I)$ is a closed curve (maybe not simple) and $\gamma^n(E)$ is on one side of it. Since $v_n\in F_n^{p,q}(I)\cup G_n^{p,q}(I)$, we have that there exists an $n$-arc $e$ with $\gamma^n(e)$ incident to $v_n$ that is not contained in $\gamma(E)$. As discussed above, for an arbitrary subarc $I\subseteq e\smallsetminus\{t_0\}$, there is a point $t\in I$ that can not be equivalent to points on $E$ under $\sim$. It follows that $v_n\in\partial\gamma(E)$.

    By Lemma~\ref{prop:ch4-road-data}~\textup{(ii)} and Lemma~\ref{lem:ch3-road-limits}, for each $t\in[0,1]$ there is a sequence $\{v_n\}_{n=1}^{+\infty}$ such that $v_n\in F_n^{p,q}(I)$ and $v_n\to F^{p,q}(t)$ as $n\to+\infty$. Then we have that $F^{p,q}(t)\in\partial\gamma(E)$. Similarly, we have $G^{p,q}(I)\subseteq\partial\gamma(E)$. Thus, we have $l_{p,q}\cup g_{p,q}\subseteq\partial \gamma(E)$.

    For $v\in\partial\gamma(E)$, suppose $t_0\in E$ such that $\gamma(t_0)=v$. Consider a neighborhood $U\ni v$ on $S^2$. Since $v\in\partial\gamma(E)$ and $\gamma^n$ converge to $\gamma$ uniformly, there is $m$-vertex $v_m\in U$ such that the $m$-edges incident to it are contained in both $\gamma^m(E)$ and $\gamma^m\bigl(\overline{S^1\smallsetminus E}\bigr)$. Since $F_m^{p,q}(I)\cup G_m^{p,q}(I)$ is a closed curve (maybe not simple) and $\gamma^m(E)$ is on one side of it, we have that $v_m\in F_m^{p,q}(I)\cup G_m^{p,q}(I)$. We also have that $v_m\in l_{p,q}\cup g_{p,q}$ since it is an $m$-vertex. Thus, $v$ is in the closure of $l_{p,q}\cup g_{p,q}$. Since $l_{p,q}\cup g_{p,q}$ is already closed, we have $v\in l_{p,q}\cup g_{p,q}$. Thus, $\partial\gamma(E)\subseteq l_{p,q}\cup g_{p,q}$. Therefore, we have $\partial \gamma(E)=l_{p,q}\cup g_{p,q}$.
\end{proof}
\begin{remark}
    Since $\gamma(E)$ is closed, $\partial \gamma(E)$ is nowhere dense on $S^2$. Thus, $Z_0^+$ and $Z_0^-$ are nowhere dense.
\end{remark}

As noted in \cite[Subsection~7.1.4]{CLWheels}, the laminations associated with expanding Thurston maps generally have perfect fits. Thus, in general, Meyer's invariant Peano curves are not CaTherine wheels, and the analogy established above does not extend to other properties of Catherine wheels.

\end{document}